# MORSE THEORY AND HIGHER TORSION INVARIANTS I

Sebastian Goette

Abstract. We compare the higher analytic torsion $\mathcal{T}$ of Bismut and Lott of a fibre bundle $p\colon M \to B$ equipped with a flat vector bundle $F \to M$ and a fibre-wise Morse function $h$ on $M$ with a higher torsion $T$ that is constructed in terms of a families Thom-Smale complex associated to $h$ and $F$, thereby extending previous joint work with Bismut. Under additional conditions on $F$, the torsion $T$ is related to Igusa's higher Franz-Reidemeister torsion. As an application, we use the higher analytic torsion to detect infinite families of smooth bundles $p_i\colon M \to B$ with diffeomorphic fibres that are homeomorphic but not diffeomorphic as bundles.

This is the first of two papers devoted to a comparison of Igusa's and Klein's higher Franz-Reidemeister torsion $\tau$ with Bismut's and Lott's higher analytic torsion $\mathcal{T}$. In this paper, we evaluate the Bismut-Lott torsion form $\mathcal{T}$ for families that carry a fibre-wise Morse function, thus extending earlier work with Bismut ([BG1]). In the second part of the series, we relate $\tau$ and $\mathcal{T}$ for families with fibre-wise Morse functions in those situations where both invariants are defined. Let us recall the development of higher torsion invariants.

Franz and Reidemeister constructed a numerical invariant $\tau_{\mathrm{FR}}$ of chain complexes in [R], [F], and used it to detect homeomorphism types of homotopy-equivalent Lense-spaces. The higher Franz-Reidemeister torsion $\tau$ is an extension of $\tau_{\mathrm{FR}}$ to families of manifolds $p\colon M \to B$, which takes values in the cohomology of $B$. It was first constructed by John Klein in [K] using a variation of Waldhausen's $A$-theory. Other describtions of $\tau$ were later given by Igusa and Klein, see [I2] for more detailled references. In [IK], they computed $\tau$ in the special case of circle bundles.

We will constantly refer to the construction of $\tau$ in [I2]. Given a family $p\colon M \to B$ of smooth manifolds, one first finds a function $h\colon M \to \mathbb{R}$ that has only Morse-type and cubical singularities along each fibre of $p$, such that the unstable manifolds of the fibre-wise singularities are trivialised in a compatible way. By [I1], such a "framed" function always exists if $\dim M > 2\dim B$. If $\dim M \leq 2\dim B$, one may replace $M$ by $M \times \mathbb{R}P^{2N}$ for some sufficiently large $N$. Let $F \to M$ be a unitarily flat complex vector bundle that is fibre-wise acyclic, then $F$ and $h$ give rise to a functor from the category of generic small simplices on $B$ to the simplicial Whitehead category $Wh^h(M(\mathbb{C}), U)$ associated to the infinite matrix ring $M(\mathbb{C}) = \lim_{\to} M_n(\mathbb{C})$ and the infinite unitary group $U = \lim_{\to} U(n)$. Now, the $k$-th higher torsion $\tau_k(M/B; F)$ is defined as the pull-back of a certain cohomology class $D_{2k} \in H^{2k}\bigl(Wh^h(M(\mathbb{C}), U), \mathbb{R}\bigr)$.

On the other hand, Ray and Singer defined an analytic torsion $\mathcal{T}_{\mathrm{RS}}$ of unitarily flat complex vector bundles on compact manifolds in [RS] and conjectured that $\mathcal{T}_{\mathrm{RS}} = \tau_{\mathrm{FR}}$. This conjecture was established independendly by Cheeger ([C]) and Müller ([M1]), and further generalised by Müller ([M2]) to unimodular flat bundles, and by Bismut and Zhang ([BZ1]) to arbitrary flat complex vector bundles.

2000 *Mathematics Subject Classification.* Primary 58J52; Secondary 57R22, 57Q10.

Research at MSRI is supported in part by NSF-grant DMS-9701755

In [BL], Bismut and Lott constructed characteristic differential forms $\mathrm{ch}^{\mathrm{o}}(\nabla^F, g^F)$ of flat vector bundles $(F, \nabla^F)$ with arbitrary metrics $g^F$, which provide obstructions towards finding a parallel metric. In fact, the forms defined in [BL] have different coefficients, but we prefer to state the results of [BL] in the so-called "Chern normalisation" that was introduced in [BG1], for reasons explained below. Let $p\colon M \to B$ be a family of compact manifolds with a horizontal distribution $T^HM \subset TM$ and a metric $g^{TX}$ on the vertical tangent bundle $TX$, and let $(F, \nabla^F) \to M$ be a flat vector bundle with metric $g^F$. Then Bismut and Lott defined a higher analytic torsion form $\mathcal{T}(T^HM, g^{TX}, \nabla^F, g^F)$ on $B$ such that

$$(*) \qquad d^B\mathcal{T}\bigl(T^HM, g^{TX}, \nabla^F, g^F\bigr) = \int_{M/B} e\bigl(\nabla^{TX}\bigr)\ \mathrm{ch}^{\mathrm{o}}\bigl(\nabla^F, g^F\bigr) - \mathrm{ch}^{\mathrm{o}}\bigl(\nabla^H, g^H\bigr)\ ,$$

where $e(\nabla^{TX})$ is the Euler form of the vertical tangent bundle with its natural connection $\nabla^{TX}$, and $H = H^*(M/B; F) \to B$ is the bundle of the fibre-wise cohomology of $F$, equipped with its natural connection $\nabla^H$ and the $L_2$-metric $g^H$ on harmonic forms. They also showed that the component of $\mathcal{T}(T^HM, g^{TX}, \nabla^F, g^F)$ in degree 0 equals the Ray-Singer analytic torsion of the fibres.

Formula $(*)$ above in particular implies a "Grothendieck-Riemann-Roch" theorem for flat vector bundles at a cohomological level. In [DWW], Dwyer, Weiss and Williams proved such a theorem at the level of algebraic $K$-theory. They also suggested another higher torsion invariant. We do not know wether their higher torsion is related to one of the above.

In [L], Lott defined a secondary family index using the analytic torsion form, which is a real analogue of the push-forward in Arakelov geometry. Ma has studied the behaviour of $\mathcal{T}$ under iterated fibrations in [Ma]. As a consequence, Bunke showed in [Bu3] that Lott's secondary index is indeed functorial. One can use Theorem 0.1 below to evaluate Lott's index in certain cases.

If the family $p\colon M \to B$ has compact structure group $G$ and if $F$ is induced from a $G$-equivariant flat vector bundle on the typical fibre $X$, then the higher analytic torsion of Bismut and Lott can be treated with equivariant methods, see [BL], Chapter 4. In [Bu1] and [Bu2], Bunke computed both the equivariant torsion and the higher analytic torsion for compact structure groups in degree $> 0$ for fibre-wise acyclic, unitarily flat vector bundles from the equivariant Euler characteristic of $X$. This was used in [Bu4] to construct cohomology classes of $\mathcal{D}\mathit{iff}(S^{2n-1})$ both with the smooth and with the discrete topology. Bunke's results also imply a relation between higher analytic torsion for families with compact structure group and equivariant torsion in the special case above, which was proved in general in [BG2]. Since this result depends on the Chern-normalisation of the analytic torsion forms, while all other results discussed here hold for any possible choice of coefficients, we have chosen to work with the particular normalisation $\mathcal{T}$.

In [BG1], we calculated $\mathcal{T}(T^HM, g^{TX}, \nabla^F, g^F)$ for families $p\colon M \to B$ with a function $h\colon M \to \mathbb{R}$ that is fibre-wise Morse, such that there exists a fibre-wise gradient field for $h$ that satisfies Smale's transversality condition on every fibre of $p$. In this setting, the fibre-wise Thom-Smale complexes form a flat $\mathbb{Z}$-graded vector bundle $(V, \nabla^V)$ over $B$, the metric $g^F$ defines a metric $g^V$ on $V$, and the coboundary map $v$ becomes a parallel section of $\operatorname{End} V$. In particular, there is a flat superconnection $A' = \nabla^V + v$ of total degree 1, and we showed that $\mathcal{T}(T^HM, g^{TX}, \nabla^F, g^F)$ is related to the torsion form $T(A', g^V)$ defined in [BL], section 2. If $H$ still denotes the fibre-wise cohomology and $g^H_V$ denotes the metric induced by identifying $H$ with a subbundle of $V$ using finite-dimensional Hodge theory, this torsion form satisfies

$$d^BT\bigl(A', g^V\bigr) = \mathrm{ch}^{\mathrm{o}}\bigl(\nabla^V, g^V\bigr) - \mathrm{ch}^{\mathrm{o}}\bigl(\nabla^H, g^H_V\bigr)\ .$$

When [BG1] was written, Igusa's framing principle of [I2], chapter 6, was not yet available. Thus in order to compare our result with the theory of Igusa and Klein, we could only regard unitarily flat bundles $F \to M$ such that the fibre-wise cohomology forms a trivial flat bundle over $B$, and only framed Morse functions $h$, i.e., Morse functions whose fibre-wise unstable manifolds are trivialised. Moreover, we still needed the existence of a fibre-wise Smale gradient field, in which case there always is a deformation of $h$ to a self-indexing Morse function on $M \to B$. In this very special case, both the higher Franz-Reidemeister torsion of Igusa and Klein and the higher analytic torsion form vanish.

In the present paper, we generalise the results of [BG1] to families of compact manifolds $M \to B$ that carry a fibre-wise Morse function $h: M \to \mathbb{R}$, but we do not require the existence of a fibre-wise Smale gradient field. Let $\hat{p}: C \to B$ denote the covering of $B$ by the fibre-wise critical points of $h$, and let $o(T^uX) \to C$ denote the orientation bundle of the unstable manifolds along $C$. We still define a vector bundle

$$V = \hat{\pi}_*\big(F|_C \otimes o(T^uX)\big)$$

of fibre-wise Thom-Smale complexes. A generic fibre-wise gradient field will satisfy the Smale condition everywhere outside a stratified subset of $B$ of codimension 1. Whenever one passes this subset, the coboundary map $v$ on the bundle $V$ will change by an elementary isomorphism. This elementary isomorphism will change whenever we pass the stratum of codimension 2, and the composition of the elementary isomorphisms that we encounter when we go around the stratum of codimension 2 will be homotopic to the identity, but not necessarily equal to it. If we fix homotopies on the stratum of codimension 2, then these will change when we pass the stratum of codimension 3, and the process continues. The algebraic structure arising from these coboundary maps, elementary isomophisms, homotopies etc. is called a "system of higher homotopies" in [I2]. We encounter such structures in Section 1, where we interpret them as "simplicial superconnections", since they behave similar as classical superconnections, except that the exterior differential $d$ gets replaced by the simplicial coboundary operator $\delta$.

In order to define an analytic torsion form for a family $V$ of Thom-Smale complexes where the fibre-wise coboundary map $v$ jumps as described above, we replace the superconnection $\nabla^V + v$ of [BG1] by a flat superconnection $A'$ of total degree 1 on $V$ that may now contain terms of arbitrarily high degree with respect to $\Omega^*(B)$. This unfortunately implies that the analytic torsion form $T(A', g^V)$ of Bismut-Lott is no longer defined. However, the Morse function $h$ induces an endomorphism $h^V$ of $V$ by multiplication with $h|_C$, and we can construct an analytic torsion form $T(A', g^V, h^V)$ that depends explicitly on $h^V$, such that still

$$d^B T\big(A', g^V, h^V\big) = \mathrm{ch}^{\mathrm{o}}\big(\nabla^V, g^V\big) - \mathrm{ch}^{\mathrm{o}}\big(\nabla^H, g^H_V\big)\ .$$

In contrast to the situation in [BG1], the form $T(A', g^V, h^V)$ can be non-trivial even if the bundle $F$ is fibre-wise acyclic and $g^F$ is parallel. For example, it will follow from [I2] and our results in [G] that $T(A', g^V, h^V)$ restricts to the Borel regulator classes on the Volodin $K$-theory space.

In order to state the main theorem of this paper, let $\widetilde{\mathrm{ch}}^{\mathrm{o}}(H, g^H_V, g^H_{L_2}) \in \Omega^*(B)/d^B\Omega^*(B)$ be the natural Chern-Simons class satisfying

$$d^B \widetilde{\mathrm{ch}}^{\mathrm{o}}\big(H, g^H_V, g^H_{L_2}\big) = \mathrm{ch}^{\mathrm{o}}\big(H, g^H_{L_2}\big) - \mathrm{ch}^{\mathrm{o}}\big(H, g^H_V\big)\ .$$

Let $\psi(\nabla^{TX}, g^{TX})$ denote the Mathai-Quillen current on the total space of $TX \to M$, such that

$$d^M\Big(\big(\nabla^{TX}h\big)^*\psi\big(\nabla^{TX}, g^{TX}\big)\Big) = e\big(TX, \nabla^{TX}\big) - (-1)^{\mathrm{ind}_h}\,\delta_C\ ,$$

where $\delta_C$ denotes the current of integration over the critical set $C$. Finally, let ${}^0J$ denote the additive genus defined in [BG1], cf. (3.6) below.

**0.1. Theorem.** *Modulo exact forms on $B$,*

$$\begin{aligned}\mathcal{T}\big(T^HM,g^{TX},\nabla^F,g^F\big)-T\big(A',g^V,h^V\big)-\widetilde{\mathrm{ch}}^{\mathrm{o}}\big(\nabla^H,g^H_{L_2},g^H_V\big)\\ =-\int_{M/B}\mathrm{ch}^{\mathrm{o}}\big(\nabla^F,g^F\big)\,\big(\nabla^{TX}h\big)^*\psi\big(\nabla^{TX},g^{TX}\big)+\hat p_*\big((-1)^{\mathrm{ind}_h}\,{}^0J(T^sX-T^uX)\big)\,\mathrm{rk}\,F\;.\end{aligned}$$

The formulation of Theorem 0.1 precisely equals Theorem 0.1 in [BG1] in the non-equivariant case, only $T(A',g^V)$ has been replaced by $T(A',g^V,h^V)$. In fact, a large portion of the proof of Theorem 0.1 is identical to the corresponding parts of [BG1] and has therefore been omitted. However, the computation of the limit of the analytic torsion under the Witten deformation in [BG1] relies on Smale's transversality condition; we give a replacement for this part of the proof in Section 4.

As an application of Theorem 0.1, we use the higher analytic torsion to detect homeomorphic, but not diffeomorphic fibre bundles. More precisely, we say that two maps $p$, $p'\colon M\to B$ are *homeomorphic (diffeomorphic) as maps* iff there exist homeomorphisms (diffeomorphisms) $\Phi$ of $M$ and $\Psi$ of $B$ such that the diagramme

$$\begin{array}{ccc} M & \xrightarrow[\sim]{\Phi} & M\\ {\scriptstyle p}\big\downarrow & & \big\downarrow{\scriptstyle p'}\\ B & \xrightarrow[\sim]{\Psi} & B\end{array}$$

commutes. If $p$ and $p'$ are continuous (smooth) fibre bundles, we say that $p$ and $p'$ are *homeomorphic (diffeomorphic) as fibre bundles* iff they are equivalent as continuous (smooth) maps and we can take $\Psi=\mathrm{id}_B$ in the diagramme above. The following result generalises Igusa's Corollary 6.5.8 in [I2] for smooth bundles that are homeomorphic to trivial bundles.

**0.2. Theorem.** *For each integer $k>0$ there exists $l_0>0$ such that the following holds. Let $p\colon M\to B$ be a smooth fibre bundle with compact, connected, oriented base $B$ of dimension $4k$ and with compact, connected fibres of dimension $2l-1$, where $l\ge l_0$. Then there exist smooth fibre bundles $p_j\colon M\to B$ for all $j\in Z$ with $p_0=p$, which are all homeomorphic to $p$ as fibre bundles and have diffeomorphic fibres, but which are pairwise not diffeomorphic as fibre bundles.*

*Moreover, if the fibre-wise cohomology bundle $H^*(M/B;\mathbb{C})$ admits a parallel metric with respect to the Gauß-Manin connection, then the $p_j$ are pairwise not diffeomorphic as maps.*

We give a more general "quantitative" version in Theorem 5.29 below. The bundles $p_i$ are constructed by a gluing construction involving homeomorphically trivial bundles with base $S^{4k}$ and fibre $S^{2l-1}$. The construction of these bundles goes back to Hatcher and is discussed in Section 5. Note that the second statement of Theorem 0.2 can be proved using an extension of Igusa's higher Franz-Reidemeister torsion that we will introduce in [G]. Our proof of the first statement however relies on the fact that the higher analytic torsion forms can be defined even if the cohomology bundle does not admit a parallel metric.

In the second paper [G] of this series, we will show that the construction of the bundle $(V,A',h^V)$ over $B$ gives rise to a map $\xi_{p,h,F}$ from $B$ to Igusa's Whitehead space $Wh(M_r(\mathbb{C}),GL(r,\mathbb{C}))$, which is homotopic to a similar map constructed in Chapter 4 of [I2]. Conversely, each such functor $\xi_{p,h,F}$ up to homotopy determines a $\mathbb{Z}$-graded flat vector bundle $V\to B$ with a flat superconnection $A'=\nabla^V+a'$ and an endomorphism $h^V$. Thus, $\xi_{p,h,F}$ is the analogue of a classifying map for flat vector bundles with the extra structure described above. The bundle $V$ admits a parallel metric $g^V$ that is adapted to $h^V$ and $A'$ and the fibre-wise cohomology $H=H^*(V,a'_0)\to B$ vanishes iff the image of $\xi_{p,h,F}$ lies in the subcategory $Wh^h(M_r(\mathbb{C}),U(r)) \subset Wh(M_r(\mathbb{C}),GL(r,\mathbb{C}))$. In this case, the torsion form $T(A',g^V,h^V)$ becomes a cohomology class on $B$, which is the pull-back of a cohomology class $T \in H^*\big(Wh^h(M_r(\mathbb{C}),U(r))\big)$. This class $T$ will be shown to equal Igusa's torsion class $\tau$.

We also construct a Whitehead space $Wh^u\big(M_r(\mathbb{C}),U(r))\big)$ that classifies complexes $(V,A')$ as above where both $F$ and $H$ carry parallel metrics. Then $T=\tau$ extends to a cohomology class of $Wh^u\big(M_r(\mathbb{C}),U(r))\big)$. Igusa's general construction gives a map $\xi\colon B \to Wh^u(M_r(\mathbb{C}),U(r))$, and the pull-back $\xi^*T$ yields a new class $\tau(M/B;F) \in H^*(B)$ for all families $M \to B$ of smooth manifolds with a unitarily flat vector bundles $(F,\nabla^F,g^F) \to M$ with unitarily flat fibre-wise cohomology $(H,\nabla^H,g^H)$—the existence of a fibre-wise Morse function is no longer required, and compactness of the fibres can be relaxed as in [I2].

In the case that the fibres of $p\colon M \to B$ are compact and the bundles $F \to M$ and $H \to F$ admit parallel metrics, the form

$$\mathcal{T}(M/B;F) = \Big(\mathcal{T}\big(T^HM,g^{TX},\nabla^F,g^F\big) + \widetilde{\mathrm{ch}}^{\mathrm{o}}\big(\nabla^H,g^H,g^H_{L_2}\big)\Big)^{[\geq 2]}$$

of Definition 2.88 is closed, and its cohomology classes in $H^*(B)$ is independent of $T^HM$, $g^{TX}$, $g^F$, and $g^H$. Using Igusa's framing principle and Theorem 0.1, we will prove the following result.

**0.3. Theorem** ([G]). *Let $M \to B$ be a family of compact manifolds, and let $(F,\nabla^F) \to M$ be a flat complex vector bundle with a parallel Hermitian metric, such that $(H,\nabla^H)$ also admits a parallel metric. Assume that there exists a function $h\colon M \to \mathbb{R}$ that is fibre-wise Morse, then*

$$\mathcal{T}(M/B;F) = \tau(M/B;F)^{[\geq 2]} + \int_{M/B} e(TX)\ {}^0J(TX)\ \mathrm{rk}\,F \quad \in H^*(B)\,.$$

Note that the fibre-wise Morse function does not appear in the conclusion. We conjecture that Theorem 0.3 holds without the assumption that such a function exists. In fact, by [I1], after taking the fibre-wise product with $\mathbb{R}P^{2N}$ for $N$ sufficiently large, there exists a framed function $h\colon M \to \mathbb{R}$ whose fibre-wise singularities are either of Morse type or cubical. In order to treat the general case, it thus remains to incorporate cubical singularities into Theorems 0.1 and 0.3.

The present paper is organised as follows. In Section 1, we construct a CW-structure on $M$ over a fixed simplicial structure on $B$, such that the cells of $M$ over each simplex of $B$ give a decomposition of the fibres that is adapted to $h$. We then describe how to turn the dicrete coboundary operator of the CW-complex into a smooth flat superconnection $A'$ on the bundle $V$. We also construct a de-Rham or "integration" map $I$ that gives a quasi-isomorphism from the total de Rham complex $(\Omega^*(M;F),\nabla^F)$ to $(\Omega^*(B;V);A')$ and has certain geometric properties that are needed in Section 4.

In Section 2, we recall the definition of the characteristic form $\mathrm{ch}(\nabla^F,g^F)$ and of the analytic torsion classes $T(A',g^V)$ and $\mathcal{T}(T^HM,g^{TX},\nabla^F,g^F) \in \Omega^*(B)/d^B\Omega^*(B)$ from [BL] and [BG1]. We describe the behaviour of the superconnection $A'$ and the endomorphism $h^V$ on $V$, and use it to define the torsion class $T(A',g^V,h^V) \in \Omega^*(B)/d^B\Omega^*(B)$. This definition already involves a finite-dimensional Witten-deformation on $V$, which is similar to the Witten deformation used in the proof of Theorem 0.1. We also define the higher torsion classes $T(M/B;F,h)$ and $\mathcal{T}(M/B;F)$ if $F$ and $H$ are unitarily flat.

In Section 3, we show that Theorem 0.1 is compatible with several known results on higher analytic torsion. We apply our result to study Lott's $\overline{K}^0_R$-pushdown. We then recall the outline of the proof of Theorem 0.1 in [BG1], stating all the intermediate results needed to prove it, and noting that all but one of these intermediate steps hold unchanged in our situation.

In Section 4, we study the behaviour of the analytic torsion form $\mathcal{T}(T^HM, g^{TX}, \nabla^F, g^F)$ under Witten deformation, without using the Smale property. Here, we use all the properties of the superconnection $A'$ on $V$, the map $I\colon \Omega^*(M;F) \to \Omega^*(B;V)$, and the higher torsion class $T(A', g^V, h^V)$ that were established in Sections 1 and 2.

Finally, in Section 5, we use Hatcher's example of a non-linear disc bundle to construct a certain sphere bundle $p\colon M \to B$. Using the higher torsion invariants, we see that $p\colon M \to B$ is homeomorphic but not diffeomorphic to a trivial sphere bundle. Thus we recover a result by Bökstedt, cf. [I2], Chapter 6. We then use Hatcher's example prove Theorem 0.2 above. We also discuss Lott's conjecture that the $\overline{K}^0_{\mathbb{R}}$-pushdown is related to a Becker Gottlieb transfer.

We have not attempted to keep this paper self-contained, because we rely heavily on the methods and results of [BG1]. However, at the beginnings of Sections 1 and 2, we recall all necessary constructions from [BL] and [BG1], so that at least the part of this paper that constructs $T(A', g^V, h^V)$ is reasonably self-contained.

We are grateful to J.-M. Bismut, U. Bunke, and P. Teichner for fruitful discussions and their interest in this project. We are in great debt to K. Igusa for giving us some topological background of his work [I2] on higher FR-torsion, and for explaining Hatcher's example, cf. Section 5. We also thank the Mathematical Sciences and Research Institute at Berkeley and Igusa's family for their hospitality.

## 1. Generalised Thom-Smale complexes

Let $M \to B$ be a smooth bundle with compact fibres, and let $h\colon M \to \mathbb{R}$ be a Morse function on each fibre. We show that one can patch together locally defined Thom-Smale complexes for $h$ to obtain a finite-dimensional vector bundle $V \to B$ with a flat superconnection $A'$, and a cochain map $I\colon (\Omega^*(M;F), d^M) \to (\Omega^*(B;V), A')$ Together with the action of the Morse function on $V$, we obtain an example of a family Thom-Smale complex as defined in Definition 2.35.

We start by reviewing in Section 1.a the construction of a Thom-Smale complex $(\Omega^*(B;V), A')$ associated to a family of fibre-wise Morse-Smale gradient fields, cf. [BG1], Chapter 5. In Section 1.b, we state some simplifying assumptions on the local geometry of $M$ and $F$ near the set of fibre-wise critical points. Dropping the Smale transversality condition, we construct in Section 1.c a particular CW-complex that is adapted to a simplicial complex $S$ on the base space $B$. The coboundary operator of the associated cochain complex naturally leads to the notion of a flat "simplicial" superconnection $\delta + a$, see Proposition 1.31. This is nothing but a "system of higher homotopies" as defined in chapter 3 of [I2]. We also construct a quasi-isomorphism relating $\nabla^F$ to $\delta + a$.

In the next step, we lift the simplicial superconnection to a neighbourhood of the diagonal in $B \times B$ in Section 1.d. In Section 1.e, we extend it to a flat "mixed" simplicial and differential superconnection. The construction in Proposition 1.47 is reminiscent of the spectral sequence proof of de Rham's theorem that simiplicial cohomology and de Rham cohomology are isomorphic. Nevertheless, the situation here is more complicated since flat superconnections obey a quadratic equation, while the equation for cohomology is linear. Note also that we have to construct an explicit quasi-isomorphism from $\Omega^*(M;F)$ to $\Omega^*(B;F)$ with certain geometric properties that will be crucial in Section 4; this is done in Proposition 1.54. Finally in Section 1.f, we pull the purely differential part of the mixed superconnection back to $B$ by a partition of unity, to obtain a flat smooth superconnection on the Thom-Smale bundle $V \to B$ in Definition 1.60. The main properties of our construction are summarized in Theorem 1.61.

**1.a. Local Thom-Smale complexes.** Let $p\colon M \to B$ be a fibre-bundle with compact fibre $X$, and let $F \to M$ be a flat vector bundle, giving rise to a de Rham complex $\Omega^*(M;F)$. Let $h\colon M \to \mathbb{R}$ be a family of Morse functions on the fibres of $p$. and let $C \subset M$ be the fibre-wise critical set of $h$. Then $C$ forms a finite covering $\widehat{p} = p|_C\colon C \to B$. If $C^j$ denotes the critical set of index $j$, then $\widehat{p}^j = p|_{C^j}\colon C^j \to B$ is also a covering, and $\widehat{p}$ is the disjoint union of the coverings $\widehat{p}_0$, $\widehat{p}_1$, ...

Let $TX \to M$ denote the vertical tangent bundle, and let $T^sX \oplus T^uX \to C$ be any decomposition of $TX|_C$ such that the vertical second derivative $d^2h\colon (TX)^2 \to \mathbb{R}$ of $h$ becomes positive (negative) definite on $T^sX$ (on $T^uX$), and such that $d^2h$ vanishes on $T^sX \times T^uX$. Then the orientation bundle $o(T^uX) = \Lambda^{\max} T^uX$ is independent of the choice of such splitting.

**1.1. Definition.** Define vector bundles

$$\widehat{V}^j = F|_{C^j} \otimes o(T^uX) \longrightarrow C^j\ , \qquad \widehat{V} = \bigcup_{j\geq 0} \widehat{V}^j = F|_C \otimes o(T^uX) \longrightarrow C\ ,$$

$$V^j = \widehat{p}^j_* \widehat{V}^j \longrightarrow B\ , \qquad V = \bigoplus_{j\geq 0} V^j = \widehat{p}_* \widehat{V} \longrightarrow B\ ,$$

$$\text{and} \qquad \operatorname{End}^k V = \bigoplus_{j\geq 0} \operatorname{Hom}\bigl(V^j, V^{j+k}\bigr) \subset \operatorname{End} V\ .$$

Let $\nabla^{\widehat{V}}$, $\nabla^V$ and $\nabla^{\operatorname{End} V}$ be the flat connections on $\widehat{V}$ and $V$ respectively that are induced by $\nabla^F$.

**1.2. Definition.** Let $N^{\widehat{V}} \in \Gamma(\operatorname{End} \widehat{V})$ denote multiplication by $j$ on $\widehat{V}^j$, and let $N^V \in \Gamma(\operatorname{End}^0 V)$ be its pushdown to $B$. Multiplication by $h|_C$ on $\widehat{V}$ will be denoted by $h^{\widehat{V}}$, and its pushdown to $B$ by $h^V \in \Gamma(\operatorname{End}^0 V)$. Then both $N^V$ and $h^V$ extend to $\Omega^*(B)$-linear endomorphisms of $\Omega^*(B, V)$.

An endomorphism $e \in \Omega^*(B, \operatorname{End} V)$ will be called *strictly $h^V$-increasing* if $e$ maps the $\lambda$-eigenspace of $h^V$ to the sum of eigenspaces of eigenvalues strictly greater than $\lambda$. In other words, such an operator maps subspaces $\widehat{V}_c \subset V_b$ for $c \in \widehat{p}^{-1}b$ to the sum of subspaces $\widehat{V}_{c'}$ with $c' \in \widehat{p}^{-1}b$ and $h(c') > h(c)$.

**1.3. Definition.** For $b \in B$, let $\operatorname{End}_+ V_b = \bigoplus_{k\in\mathbb{Z}} \operatorname{End}^k_+ V_b$ be the space of strictly $h^V$-increasing endomorphisms of $V_b$. Let

$$\Gamma(\operatorname{End}_+ V) = \bigl\{\, e \in \Gamma(\operatorname{End} V) \bigm| e_b \in \operatorname{End}_+ V_b \text{ for all } b \in B \,\bigr\}\ ,$$

$$\text{and} \qquad \Omega^*(B; \operatorname{End}_+ V) = \Omega^*(B) \otimes_{C^\infty(B)} \Gamma(\operatorname{End}_+ V)\ .$$

Note however that the spaces $\operatorname{End}_+ V_b$ may not form vector bundles over $B$ because $h$ is not locally constant on $C$.

**1.4. Definition.** Let $F^0(M; F) \supset F^1(M; F) \supset \dots$ and $F^0(B; V) \supset F^1(B; V) \supset \dots$ be the filtrations of $\Omega^*(M; F)$ and $\Omega^*(B, V)$ respectively by horizontal degree, given by

$$F^p(M; F) = \bigoplus_k F^{p,k}(M; F)\ ,$$

$$\text{with} \qquad F^{p,k}(M; F) = \bigl\{\, \alpha \in \Omega^k(M; F) \bigm| \iota_{X_0} \dots \iota_{X_{k-p}} \alpha = 0 \text{ for all } X_0, \dots, X_{k-p} \in \Gamma(TX) \,\bigr\}$$

$$\text{and} \qquad F^p(B; V) = \bigoplus_{q\geq p} \Omega^q(B; V)\ .$$

We can choose an open cover $\mathcal{U} = \{U_i \mid i \in \mathcal{I}\}$, such that we have vertical metrics $g_i^{TX}$ on $TX|_{p^{-1}(U_i)}$ over $U_i$, for which $h$ becomes a fibre-wise Morse-Smale function. Let $M_i^u \to \widehat{p}^{-1}(U_i) \subset C$ with $M_i^u|_c \subset M_{\widehat{p}(c)}$ denote the bundle of unstable cells of the fibre-wise negative gradient flow with respect to the metric $g_i^{TX}$.

**1.5. Definition.** Set

$$\hat I_i\colon\ \Omega^*(M;F)\ \longrightarrow\ \Omega^*\bigl(\widehat p^{-1}(U_i);\widehat V|_{\widehat p^{-1}(U_i)}\bigr)\ ,\qquad \alpha\ \longmapsto\ \int_{M_i^u/C}\alpha\ ,$$

$$\text{and}\qquad I_i\colon\ \Omega^*(M;F)\ \longrightarrow\ \Omega^*\bigl(U_i;V|_{U_i}\bigr)\ ,\qquad \alpha\ \longmapsto\ \int_{M_i^u/B}\alpha\ .$$

Note that $\hat I_i$ and $I_i$ are $\Omega^*(B)$-linear in the sense that

$$\hat I_i(p^*\alpha\wedge\beta)=\widehat p^*\alpha\wedge\hat I_i(\beta)\qquad\text{and}\qquad I_i(p^*\alpha\wedge\beta)=\alpha\wedge I_i(\beta)\ .$$

By [BG1], Theorem 5.5, these maps are $C^\infty$ under the simplifying assumptions of Section 1.b.

Summarizing the results of [BG1], Sections 5.4 and 5.5, in particular Theorem 5.8, we obtain the following characterisation of $I_i$.

**1.6. Theorem.** *There exists a unique coboundary operator $\partial_i\in\Gamma(\operatorname{End}^1_+V|_{U_i})$ such that the superconnection $A_i'=\partial_i+\nabla^V$ on $V$ is flat and*

$$I_i\Bigl(\Omega^*\bigl(p^{-1}(U_i),F|_{p^{-1}(U_i)}\bigr),\nabla^F\Bigr)\ \longrightarrow\ \bigl(\Omega^*(U_i,V|_{U_i}),A_i'\bigr)$$

*is a quasi-isomorphism of complexes respecting the filtrations by horizontal degree. If we denote the spectral sequences associated to the filtered complexes above by $(E_k^*(M;F),d_k)$ and $(E_k^*(U_i;V),d_k)$ respectively, then $I_i$ induces isomorphisms*

$$I_{i*}\colon\ \bigl(E_k^*(M;F),d_k\bigr)\xrightarrow{\ \sim\ }\bigl(E_k^*(U_i;V),d_k\bigr)\qquad\text{for } k\ge 1.$$

We will be mainly interested in the $E_0$- and $E_1$-terms of these spectral sequences. Let $T^HM\subset TM$ be a horizontal subbundle complementary to $TX$, and let $\Omega\in\Omega^2(B,p_*TX)$ be the associated fibre bundle curvature. Then the operators $d^M$ on $\Omega^*(M)$ and $\nabla^F$ on $\Omega^*(M;F)$ split naturally as

$$d^M=d^X+\nabla^{\Omega^*(X)}-\iota_\Omega\qquad\text{and}\qquad\nabla^F=\nabla^{F|_X}+\nabla^{\Omega^*(X;F)}-\iota_\Omega$$

with respect to the splitting $\Omega^*(M)=\Omega^*(B;\Omega^*(X))$ and $\Omega^*(M;F)=\Omega^*(B;\Omega^*X;F)$ defined by $T^HM$. The operators $d^X$ and $\nabla^{F|_X}$ induce operators on the quotients

$$F^p(M)/F^{p+1}(M)\cong\Omega^p\bigl(B;\Omega^*(X)\bigr)\qquad\text{and}\qquad F^p(M;F)/F^{p+1}(M;F)\cong\Omega^p\bigl(B;\Omega^*(X;F)\bigr)$$

respectively, which are independent of $T^HM$. We have

$$\bigl(E_0^p(M;F),d_0\bigr)\cong\Bigl(\Omega^p\bigl(B;\Omega^*(X;F)\bigr),\nabla^{F|_X}\Bigr)\ ,\qquad\text{and}\qquad\bigl(E_0^p(U_i;V),d_0\bigr)\cong\bigl(\Omega^p(U_i;V),\partial_i\bigr)\ . \tag{1.7}$$

**1.8. Definition.** Let $H=H^*_{\mathrm{dR}}(M/B;F)$ denote the bundle of the fibre-wise cohomology with its natural flat Gauß-Manin connection $\nabla^H$. Similarly, let $H_i=H^*(V;\partial_i)$, and let $\nabla^{H_i}$ be the connection induced by $\nabla^V$.

By flatness of $A_i'$, the operator $\partial_i$ is parallel. If follows that $\nabla^{H_i}$ is well-defined and flat. As in [BG1], we have

$$\bigl(E_1^*(M;F),d_1\bigr)\cong\bigl(\Omega^*(B;H),\nabla^H\bigr)\ ,\qquad\text{and}\qquad\bigl(E_1^*(U_i;V),d_1\bigr)\cong\bigl(\Omega^*(U_i;H_i),\nabla^{H_i}\bigr)\ .$$

By Theorem 1.6, the integration map gives an isomorphism of these $E_1$-terms. In particular, there is an isomorphism of flat vector bundles

$$I_{i*}\colon\ \bigl(H|_{U_i},\nabla^H\bigr)\ \longrightarrow\ \bigl(H_i,\nabla^{H_i}\bigr)\ .$$

**1.b. Some simplifying assumptions.** In this section, we construct vertical metrics $g^{TX}$, horizontal subbundles $T^H M$, and metrics $g^F$ on $F$, satisfying certain properties that we will use throughout this section and Section 4.

Let $p\colon M \to B$ be a smooth bundle with compact fibres, and let $h\colon M \to \mathbb{R}$ be a smooth function that is Morse on each fibre of $p$. Let $\widehat{p} = p|_C\colon C \to B$ be the finite covering of $B$ by the fibre-wise critical points of $h$. Let $(F, \nabla^F) \to M$ be a flat vector bundle.

For a Euclidean vector bundle $\pi\colon (E, g^E) \to M$, define

$$D_r E = \left\{ e \in E \;\middle|\; |e| < r \right\} .$$

A Euclidean connection $\nabla^E$ on $E$ with curvature $F^E$ naturally defines a horizontal subbundle $T^H E$ of the tangent bundle $TE$ to its total space, such that the fibre bundle curvature $\Omega^E\colon \pi^* TM \to E$ takes the form

$$\Omega^E(v, w)_e = -[\bar{v}, \bar{w}]_e^\perp = F^E_{v,w} e \; , \tag{1.9}$$

where $\bar{v}$, $\bar{w} \in \Gamma(T^H E)$ are horizonttal lifts of vector fields $v$, $w$ on $M$, and $e \in E$ is a point on the total space.

If $g^{TX}$ is a vertical metric on $TX \to M$, let $T^s X \subset TX|_C$ ($T^u X \subset TX|_C$) denote the sum of all positive (negative) eigenspaces of the second derivative $\nabla\nabla h \in \Gamma(\operatorname{End} TX|_C)$, so the bundle $\pi\colon TX|_C \to C$ splits as $TX|_C = T^s X \oplus T^u X$. Let $U_r(C)$ denote the fibre-wise neighbourhood of $C$ of radius $r$, and let $\exp^X\colon TX \to M$ be the fibre-wise Riemannian exponential map.

**1.10. Proposition.** *There exists a horizontal subbbundle $T^H M \subset TM$, a vertical metric $g^{TX}$ on $TX$, a metric $g^F$ on $F$, and a Euclidean connection $\nabla^{TX|_C}$ respecting the splitting $TX|_C = T^s X \oplus T^u X$, with the following property. For each compact subset $K \subset B$ there exists an $\varepsilon > 0$, such that*

1. *over $K$, the fibre-wise exponential map $\exp^X|_{D_{2\varepsilon} TX|_C}$ induces an isometry onto its image $U_{2\varepsilon}(C)$, and its differential maps $T^H(TX)$ to $T^H M|_{U_{2\varepsilon}(C)}$;*
2. *the Morse function $h|_{U_{2\varepsilon}(C)}$ is for $c \in C|_K$ and $(z^s, z^u) \in D_{2\varepsilon} TX_c$ given by*
$$h\bigl(\exp_c^X(z^s, z^u)\bigr) = h(c) + \frac{|z^s|^2 - |z^u|^2}{2} \; ;$$
3. *over $K$, the flat bundle $(F, \nabla^F, g^F)|_{U_{2\varepsilon}(C)}$ is isometric to $\pi^*(F|_C, \nabla^F|_C, g^F|_C)$.*

*Proof.* This follows easily because the classical Morse lemma holds for families. □

If conditions (1) – (3) above hold for $T^H M$, $g^{TX}$ and $g^F$, we will say that $T^H M$, $g^{TX}$ and $g^F$ satify the simplifying assumptions of Section 1.b. We assume throughout this section that this is the case. For simplicity, we may assume (at least for the purpose of estimates) that $B$ is compact, so that we may fix an $\varepsilon > 0$ as in Proposition 1.10 for all of $B$.

In particular, the closures of the $2\varepsilon$-tubular neighbourhoods of different leaves of the critical covering $\widehat{p}\colon C \to B$ do not intersect, and $h$ differs by at most $2\varepsilon^2$ from $h|_C$ along the fibres of these tubular neighbourhoods. For each flow line $\gamma$ of the negative gradient flow, this implies that $h(c_-) - h(c_+) > 4\varepsilon^2$, where $c_\pm = \lim_{t \to \pm\infty} \gamma(t) \in C$.

**1.c. A CW complex for a fibre-wise Morse function.** In Section 1.a, we recalled from [BG1] how to associate a family of CW complexes to a fibre-wise Morse function $h$ in the special case that there exists a fibre-wise Morse-Smale gradient field for $h$. In this subsection, we similarly construct a CW complexes on $M$ over a simplicial complex $S$ on $B$ for general fibre-wise Morse functions.

With respect to these complexes, the map $p\colon M \to B$ will be cellular, and the preimage of any simplex $|\sigma|$ in $M$ will consist of one cell $M_\alpha(\sigma)$ of dimension $k + \operatorname{ind} h$ for each leaf $C_\alpha(\sigma)$ of the critical covering $\widehat{p}\colon C \to B$.

Let the *standard $k$-simplex* be a fixed closed affine $k$-simplex $\Delta_k \subset \mathbb{R}^k$, whose corners will be enumerated $\Delta_k(0), \dots, \Delta_k(k)$. Then for each finite sequence of integers $0 \le i_0 \le \dots \le i_l \le k$ there is a unique affine map $f_{i_0 \dots i_l}\colon \Delta_l \to \Delta_k$ mapping $\Delta_l(j)$ to $\Delta_k(i_j)$. We call $f_{i_0 \dots i_l}$ a *face* of $\Delta_k$ if $i_0 < \dots < i_l$, and a *degeneracy* otherwise.

Let $x_0, \dots, x_k$ are barycentric coordinates on $\Delta_k$, i.e., $x_j$ is an affine function on $\Delta_k \subset \mathbb{R}^k$ that takes the value $\delta_{ij}$ at $\Delta_k(i)$, and $X_0 + \dots + x_k = 1$. Let $w_{ij}$ denote the constant vector field on $\mathbb{R}^k$ with $w_{ij} = \Delta_k(j) - \Delta_k(i) \in \mathbb{R}^k$, so

$$w_{ij}(x_l) = \begin{cases} -1 & \text{for } l = i, \\ 1 & \text{for } l = j, \text{ and} \\ 0 & \text{otherwise.} \end{cases} \tag{1.11}$$

We recall that $\Delta_k$ is oriented by fixing the oriented base $w_{0,1}, \dots, w_{0,k}$, and that $\partial\Delta_k$ is oriented such that prefixing an oriented base of $\partial\Delta_k$ with an outward normal vector gives an oriented base of $\Delta_k$. Then we have

$$\partial\Delta_k = \sum_{j=0}^{k} (-1)^j f_{0 \dots \widehat{j} \dots k} \Delta_{k-1} .$$

**1.12. Definition.** A *smoothly embedded $k$-simplex* in $B$ is a map $\sigma\colon \Delta_k \to B$ that extends to a smooth embedding of some small neighbourhood of $\Delta_k$ in $\mathbb{R}^k$ into $B$. Let $|\sigma|$ be the image of $\sigma$, and write $\sigma_j = \sigma(\Delta_k(j)) \in B$ for $0 \le j \le k$.

A *smooth ordered simplicial complex* $S$ on $B$ consists of sets $S_k$ of smooth $k$-simplices in $B$ for each $k \ge 0$, and of *face maps* $F_{i_0,\dots,i_l}\colon S_k \to S_l$ for each finite sequence of integers $0 \le i_0 < \dots < i_l \le k$, such that

$$F_{i_0,\dots,i_l}(\sigma) = \sigma \circ f_{i_0,\dots,i_l}\colon\ \Delta_l \longrightarrow B ,$$

and $B$ is the disjoint union of the interiors of the images of all simplices in $S$. We say that $\sigma' \in S_l$ is a *face* of $\sigma \in S_k$, if $\sigma' = F_{i_0,\dots,i_l}(\sigma)$ for some face map $F_{i_0,\dots,i_l}$. A *simplicial subset* of $B$ is a closed subset that is a union of images of simplices of $S$.

We will now fix a smooth ordered simplicial complex $S$ on $B$. Let us assume that $S$ is so fine that each simplex in $S$ is uniquely determined by the set of its vertices. We will write shortly $\sigma_{i_0,\dots,i_l}$ for the face $F_{i_0,\dots,i_l}(\sigma)$, thereby identifying vertices $\sigma_j$ with their images in $B$. Note that we will not need degerated simplices or degeneracy maps for most of the following constructions.

*1.13. Remark.* By the definition of a smoothly embedded $k$-simplex $\sigma$ in $B$, we can pull the bundle $M \to B$ back to a neighbourhood $U(\Delta_k)$ of $\Delta_k$ in $\mathbb{R}^k$. In particular, we can say a vector field on $p^{-1}(|\sigma|) \subset M$ is smooth if its pullback extends smoothly to the pullback of $M$ to $U(\Delta_k)$. If two smooth submanifolds of $p^{-1}(|\sigma|)$ intersect, we can say that they are transversal if their pullbacks extend to smooth submanifolds over $U(\Delta_k)$ that intersect transversally. This makes sense even if the original intersection took place in the boundary of $p^{-1}(|\sigma|)$.

For each $k$-simplex $\sigma \in S_k$, let $P(\sigma) = C|_{|\sigma|}/\,|\sigma|$ be the set of leaves of the covering $\widehat{p}\colon C \to B$ restricted to $|\sigma|$. For $\alpha \in P(\sigma)$, let $C_\alpha(\sigma)$ be the corresponding leaf of $C|_{|\sigma|}$. Let $h_\alpha$ be restriction of $h$ to $C_\alpha(\sigma)$, regarded as a function on $|\sigma|$, and let $\operatorname{ind}_\alpha$ denote the fibre-wise Morse index of $h$ at $C_\alpha(\sigma)$. We choose the simplical structure on $B$ fine enough such that

$$|h_\alpha(q') - h_\alpha(q)| < \frac{\varepsilon^2}{2} \qquad \text{for all } \alpha \in P(\sigma) \text{ and all } q,\, q' \in |\sigma|. \tag{1.14}$$

Let us define a relation $\prec_\sigma$ on $P(\sigma)$ by setting

$$\alpha \prec_\sigma \beta \qquad \text{iff} \qquad \big(h_\beta - h_\alpha\big)(q) > 2\varepsilon^2 \qquad \text{for some } q \in |\sigma|. \tag{1.15}$$

*1.16. Remark.*

(1) Assumption (1.14) guarantees that $\prec_\sigma$ is a partial ordering on $P(\sigma)$.
(2) If a $\sigma'$ is a face of $\sigma$, then $\prec_\sigma$ refines on $\prec_{\sigma'}$.
(3) Because of the simplifying assumptions of Section 1.b, the unstable cell $M^u(C_\beta(\sigma))$ can only intersect the closed $2\varepsilon$-tubular neighbourhood of $C_\alpha(\sigma)$ if $\alpha \prec_\sigma \beta$. In particular, there is no negative gradient flow line with respect to $g_0^{TX}$ from $C_\beta(\sigma)$ to $C_\alpha(\sigma)$ over some $q \in |\sigma|$ unless $\alpha \prec_\sigma \beta$.

Using barycentric coordinates and the vector fields $w_{ij}$ of (1.11), we construct compatible vector fields on all simplices of $S$.

**1.17. Proposition.** *The vector fields*

$$W_k = \sum_{0 \le i < j \le k} x_i\, x_j\, w_{ij}$$

*on the standard simplices $\Delta_k$ have the following properties:*

(1) *whenever $0 \le i_0 \le \cdots \le i_l \le k$, we have*

$$W_k \circ f_{i_0,\dots,i_k} = df_{i_0,\dots,i_k} \circ W_l\ ,$$

*in particular, $W_k$ is tangential to all faces of $\Delta_k$;*
(2) *$W_k$ vanishes only at the vertices $\Delta_k(j)$ of $\Delta_k$, where it has non-degenerate zeros, and the stable and unstable tangential subspaces are spanned by $w_{0,j}$, … , $w_{j-1,j}$ and by $w_{j,j+1}$, … , $w_{j,k}$ respectively; and*
(3) *the stable and unstable cells of $W_k$ at $\Delta_k(j)$ are given by*

$$M^s\big(\Delta_k(j)\big) = |f_{0,\dots,j}| \setminus |f_{0,\dots,j-1}| \qquad \text{and} \qquad M^u\big(\Delta_k(j)\big) = |f_{j,\dots,k}| \setminus |f_{j+1,\dots,k}|\ .$$

By [S], the vector field $W_k$ can be obtained as a gradient field of a Morse function on $\Delta_k$ with respect to some Riemannian metric.

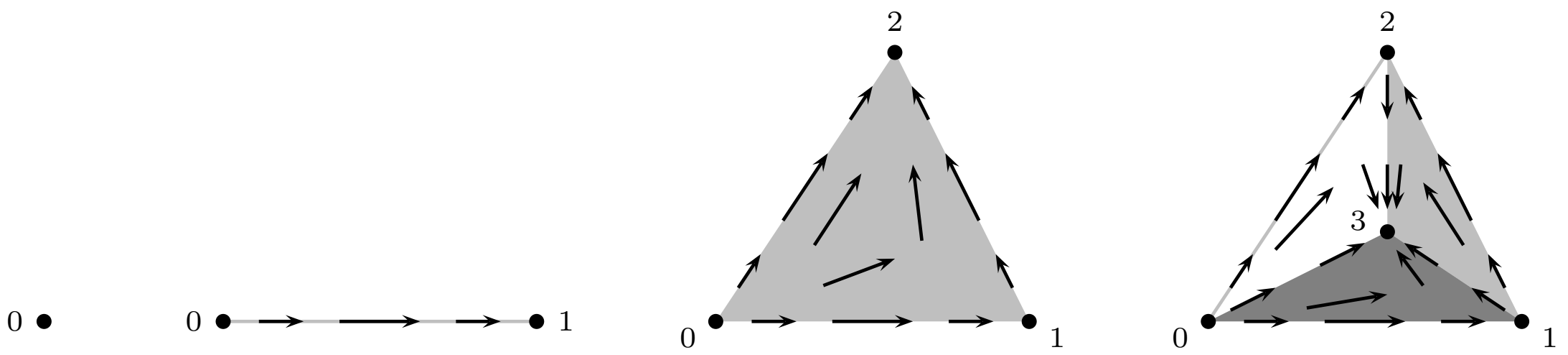

Figure 1.18. The vector fields $W_k$ on small standard simplices.

*Proof.* Because $w_{ij}$ is tangent to all faces of $\Delta_k$ containing $i$ and $j$, and because $x_ix_j$ vanishes on all other faces of $\Delta_k$, the vector field $W_k$ is tangential to all faces of $\Delta$. Now (1) follows by an elementary calculation. Also, it is clear that $W_k$ has non-degenerate zeros at the vertices of $\Delta_k$ of the type described in (2).

To prove that $W_k$ has no other zeros and to determine the stable and unstable cells, regard the affine function $h_k = \sum_{j=0}^k j\,x_j$. We have

$$W_k(h_k) = \sum_{0 \le i < j \le k} (j-i)\, x_i\, x_j \ge 0\;, \tag{1.19}$$

with "=" only at the vertices of $\Delta_k$, which completes the proof of (2).

By (1.19), each trajectory $\gamma(t)$ of $W_k$ converges for $t \to \pm\infty$ to a zero of $W_k$. By (1), each face of $\Delta_k$ is invariant under the flow of $W_k$. Since the dynamics of $W_k$ near the vertices are given by (2), (3) now follows easily. □

Property (1) guarantees the existence of a piece-wise smooth vector field $W_B$ on $B$ such that

$$W_B \circ \sigma = D\sigma \circ W_k$$

for all $k$ and all $\sigma \in S_k$. Note that $W_B$ looks like the negative gradient vector field of a Morse function when restricted to the image of a single simplex.

Let $\overline{W}_B$ be the horizontal lift of $W_B$ to $M$, so $\overline{W}_B$ is now a piece-wise smooth horizontal vector field on $M$. We fix a smooth cut-off function $\varphi_M\colon M \to [0,1]$ with support in $U_\varepsilon(C)$ and $\varphi_M|_{U_{\varepsilon/2}(C)} = 1$, and we define vector field

$$W = -\nabla^{TX} h + \varphi_M\, \overline{W}_B\;. \tag{1.20}$$

Then $W$ is smooth on $p^{-1}(|\sigma|)$ in the sense of Remark 1.13 for each simplex $\sigma \in S$. We denote by $M^s_W(x,\sigma)$, $M^u_W(x,\sigma)$ the stable and unstable cells of the flow of $W$ restricted to $p^{-1}(|\sigma|)$ at a critical point $x$ of $W$, while $M^s$, $M^u$ still denotes the fibre-wise stable and unstable cells of $-\nabla^{TX} h$. For $c \in \mathbb{R}$, we define the sets

$$\begin{aligned} h^{-1}(h_\alpha + c) &= \bigl\{\, x \in p^{-1}(|\sigma|) \bigm| h(x) = h_\alpha\bigl(p(x)\bigr) + c \,\bigr\}\;, \\ h^{-1}(h_\alpha + c, \infty) &= \bigl\{\, x \in p^{-1}(|\sigma|) \bigm| h(x) > h_\alpha\bigl(p(x)\bigr) + c \,\bigr\}\;, \\ \text{and} \qquad h^{-1}(-\infty, h_\alpha + c) &= \bigl\{\, x \in p^{-1}(|\sigma|) \bigm| h(x) < h_\alpha\bigl(p(x)\bigr) + c \,\bigr\} \subset p^{-1}(|\sigma|)\;. \end{aligned}$$

**1.21. Proposition.** *The zeros of $W$ are precisely the points $\widehat{p}^{-1}(|\sigma|) \in C$ for $\sigma \in S_0$. For each $k$ and each $k$-simplex $\sigma \in S_k$, the flow of $W$ preserves both $p^{-1}(|\sigma|) \subset M$ and $C(|\sigma|)$. If $W$ vanishes at $x \in p^{-1}(|\sigma|)$, then there is a coordinate chart $\psi\colon U \to V$ of $M$ around $x$ that maps $p^{-1}(|\sigma|) \cap U$ diffeomorphically onto $([0,\infty)^k \times \mathbb{R}^n) \cap V$, such that $W$ is in this coordinates given as*

$$(D\psi \circ W)|_{p^{-1}(|\sigma|)\cap U} = \left( \pm \psi^1 \frac{\partial}{\partial \psi^1} \pm \psi^2 \frac{\partial}{\partial \psi^2} \right) \circ W|_{p^{-1}(|\sigma|)\cap U}\;. \tag{1}$$

*For $\alpha \in P(\sigma)$, let $x \in C_\alpha(|\sigma|)$ project to the $j$-th vertex $\sigma_j$ of $|\sigma|$. Then we have*

$$\begin{aligned} M^s_W(x,\sigma) \cap h^{-1}\left(-\infty, h_\alpha + \frac{3}{2}\,\varepsilon^2\right) &= D_{\sqrt{3}\,\varepsilon} T^s X|_{C_\alpha(|\sigma_{0,\dots,j}| \setminus |\sigma_{0,\dots,j-1}|)}\;, \\ \textit{and} \qquad M^u_W(x,\sigma) \cap h^{-1}\left(h_\alpha - \frac{3}{2}\,\varepsilon^2, \infty\right) &= D_{\sqrt{3}\,\varepsilon} T^u X|_{C_\alpha(|\sigma_{j,\dots,k}| \setminus |\sigma_{j+1,\dots,k}|)}\;. \end{aligned} \tag{2}$$

*Moreover, if $y$ is another critical point of $W$ and $M^u_W(x,\sigma) \cap M^s_W(y,\sigma) \neq \emptyset$, then either $\beta = \alpha$ or $\beta \prec_\sigma \alpha$.*

*Proof.* The first two claims of the proposition follow directly from our construction. Let $x$ be a zero of $W$ with $p(x)=\sigma_j$. We can choose a coodinate chart of $B$ that restricts to the barycentric coordinates of $|\sigma|$, except the $j$th. Composing with $\hat p$, we obtain a chart $\psi_0$ of $C$ around $x$. Finally, we obtain the desired chart by trivialising $T^sX$ and $T^uX\to C$ by parallel translation with respect to the connection $\nabla^{TX|_C}$ of Proposition 1.10, along paths starting from $x$ that are radial in the chart $\psi_0$.

Let $\sigma\in S_k$ be a $k$-simplex, and let $h_k$ be the function constructed in the proof of Proposition 1.17, which we now regard as a function on $|\sigma|\subset B$. Then for some sufficiently large $C>0$, we clearly have

$$W(h-C\,h_k)\le 0\;,$$

with "=" only at the zeros of $W$. This shows that each non-constant flow line $\gamma(t)$ converges to zeros of $W$ for $t\to\pm\infty$.

We will restrict our attention to the unstable cells of $W$ and $-\nabla^{TX}h$, the stable case being completely analogous. By our simplifying assumptions and because $\varphi_M\,\overline{W}_B$ is horizontal, $W$ is tangential to the subset $D_{2\varepsilon}T^uX$ of the $2\varepsilon$-tubular neighbourhood $U_{2\varepsilon}(C_\alpha(|\sigma|))$, and restricted to $U_{2\varepsilon}(C_\alpha(|\sigma|))$, the unstable cells of $x$ is

$$D_{2\varepsilon}T^uX|_{C_\alpha(|\sigma_{0,\dots,j}|\setminus|\sigma_{j+1,\dots,k}|)}\;.$$

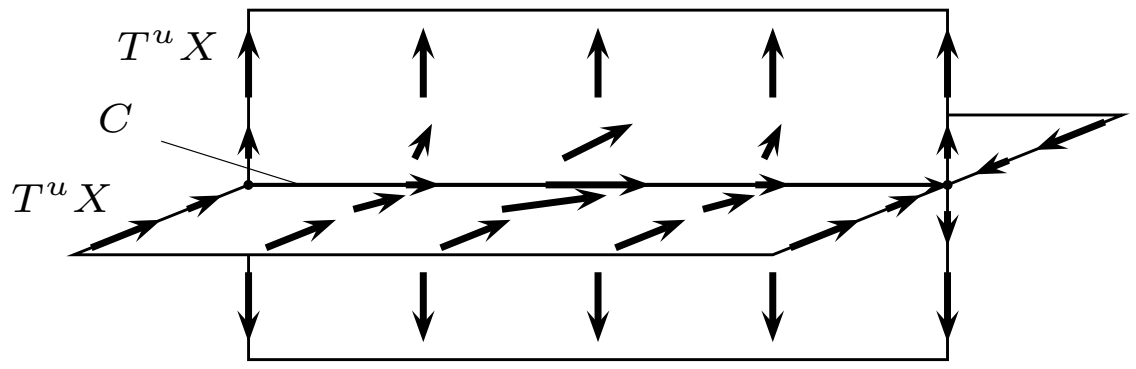

Figure 1.22. The vector field $W$ near the fibre-wise critical set of $h$.

Now regard a typical flow-line $\gamma(t)$ of $W$ with $\lim_{t\to-\infty}\gamma(t)=x$. Unless $\gamma(t)\in C_\alpha(|\sigma|)$ for one $t\in\mathbb{R}$ (and thus for all $t$), the flow-line $\gamma$ will leave $U_{2\varepsilon}(C_\alpha(|\sigma|))$ at a certain times $t_0$ with

$$h\big(\gamma(t_0)\big)<h_\alpha\big(p(\gamma(t_0))\big)-2\varepsilon^2\;.$$

Because $\operatorname{supp}(\varphi_M)\subset U_\varepsilon(C)$, the trajectory $\gamma$ will now follow the fibre-wise flow of $-\nabla^{TX}h$ until it enters $U_\varepsilon(C_\beta(|\sigma|))$ for some $\beta$ Let $t_1>t_0$ be the time where $\gamma(t)$ enters the larger neighbourhood $U_{2\varepsilon}(C_\beta(|\sigma|))$, and let $\gamma(t_1)=(z^s,z^u)\in\partial D_{2\varepsilon}TX|_{C_\beta(|\sigma|)}$. Since $\gamma(t)$ hits $U_\varepsilon(C_\beta(|\sigma|))$, clearly

$$|z^u|<\varepsilon\qquad\text{and}\qquad|z^s|=\big(4\varepsilon^2-|z^u|^2\big)^{\frac12}>\sqrt3\,\varepsilon\;,$$

so by Proposition 1.10 (2),

$$h_\alpha\big(p(\gamma(t_0))\big)-\varepsilon^2>h\big(\gamma(t_0)\big)>h\big(\gamma(t_1)\big)>h_\beta\big(p(\gamma(t_1))\big)+\varepsilon^2\;,$$

where of course $p(\gamma(t_0))=p(\gamma(t_1))\in|\sigma|$. In particular, we have $\beta\prec_\sigma\alpha$ by (1.15).

Repeating this argument starting from the time where $\gamma(t)$ leaves $U_{2\varepsilon}(C_\beta(|\sigma|))$, we check that $\gamma(t)$ passes the $\varepsilon$-tubular neighbourhoods of a finite number of critical levels that form a decreasing sequence with respect to the partial ordering $\prec_\sigma$. The vertical movement inside the larger $2\varepsilon$-tubular neighbourhood of any of these fixpoints will decrement $h(\gamma(t))$ by at least $2\varepsilon^2$, while the horizontal movement introduced by $\varphi_M\,W_B$ can increment $h(\gamma(t))$ by at most $\varepsilon^2/2$ by (1.14). This shows in particular that $h(\gamma(t))<h_\alpha(p(\gamma(t)))$ for $t>t_0$, and as claimed,

$$M^u_W(x,\sigma)\cap h^{-1}\left(h_\alpha-\frac32\,\varepsilon^2,\infty\right)=D_{\sqrt3\,\varepsilon}T^uX|_{C_\alpha(|\sigma_{0,\dots,j}|\setminus|\sigma_{j+1,\dots,k}|)}\;.\qquad\square$$

In order to describe the closures of the stable and unstable cells geometrically, we need the notion of submanifolds with conical singularities. We recall the recursive definition in [La].

**1.23. Definition.** A stratified subset $\Sigma = (\Sigma_0, \dots, \Sigma_k)$ is a *submanifold with conical singularities* of dimension $k$, if

(1) for any $0 \le i \le k$, the set $\Sigma_i \setminus \Sigma_{i+1}$ is a smooth submanifold of $M$,
(2) around each point $x \in \Sigma_i \setminus \Sigma_{i+1}$, there is a $C^\infty$-coordinate chart $U_x \to \mathbb{R}^m$ of $M$ that maps $U_x \cap (\Sigma_0, \dots, \Sigma_i)$ to
$$\mathbb{R}^{k-i} \times \big(T_0, \dots, T_i\big) ,$$
where $T = (T_0, \dots, T_i,) \subset \mathbb{R}^{m-k+i}$ is a submanifold with conical singularities, and
(3) around each $x \in \Sigma_k$, there is a $C^1$-coordinate chart $U'_x \to \mathbb{R}^m$ of $M$ that maps $\Sigma \cap U'_x$ to the cone of a $(k-1)$-dimensional submanifold of $S^{m-1}$ with conical singularities.

Two submanifolds $\Sigma$ and $T \subset M$ with conical singularities are said to intersect *transversally* iff $\Sigma_i \setminus \Sigma_{i+1}$ intersects $T_j \setminus T_{j+1}$ transversally for all $i$, $j$.

*1.24. Remark.* Let $\Sigma = (\Sigma_0, \dots, \Sigma_k)$ be a submanifold with conical singularities of $M$, and let $T \subset M$ be a submanifold. Suppose that $T$ intersects a stratum $\Sigma_j \setminus \Sigma_{j+1}$ transversally. Then there is a neighbourhood $U$ of $T \cap (\Sigma_j \setminus \Sigma_{j+1})$ in $M$ such that $T \cap U$ intersects $\Sigma_i \cap U$ transversally in $U$ for all $i \le j$.

Note that $|\sigma| \subset B$ and $p^{-1}(|\sigma|) \subset M$ are submanifolds with conical singularities. We now want to take the unstable cells of $W$ to form a CW-complex on $M$ whose cells are submanifolds with conical singularities. By the following generalisation of a result of Laudenbach, it suffices to deform $W$ in such a way that a suitable version of the Smale transversality condition along the lines of Remark 1.13 holds. For a fixed simplex $\sigma \in S$, let $h' \colon |\sigma| \to \mathbb{R}$ be a smooth function such that $h'(q)$ is a regular value of $h|_{p^{-1}(q)}$ for all $q \in |\sigma|$.

**1.25. Proposition.** *Let $\sigma \in S$ be a fixed simplex. If the stable cells $M^u_W(x, \sigma)$ intersect the unstable cells $M^s_W(y, \sigma)$ of $W$ transversally in $p^{-1}(|\sigma|) \subset M$ in the sense of Remark 1.13, then $\overline{M^u_W(x, \sigma)}$ and $\overline{M^s_W(y, \sigma)}$ are submanifolds of $M$ with conical singularities.*

*If for a fixed zero $x$ of $W$ over $|\sigma|$, all $M^u_W(x', \sigma)$ with $x' \in \overline{M^u_W(x, \sigma)} \cap h^{-1}(h', \infty)$ meet all $M^s_W(y, \sigma)$ with $y \in h^{-1}(h', \infty)$ transversally in the sense of Remark 1.13, then $\overline{M}^u_W(x, \sigma) \cap h^{-1}(h')$ is a submanifold of $h^{-1}(h')$ with conical singularities. The strata of $\overline{M}^u_W(x, \sigma) \cap h^{-1}(h')$ are of the form $M^u_W(x', \sigma') \cap h^{-1}(h')$, where $\sigma'$ is a face of $\sigma$ and $x' \in \overline{M^u_W(x, \sigma)} \cap p^{-1}(|\sigma'|) \cap h^{-1}(h', \infty)$.*

*Proof.* First, assume that for one of the unstable cells $M^u_W(x, \sigma)$, all $M^u_W(x', \sigma)$ with $x' \in \overline{M^u_W(x, \sigma)}$ meet all $M^s_W(y, \sigma)$ transversally in the sense of Remark 1.13, for a fixed $\sigma \in S_k$ and a zero $x$ of $W$ over $p^{-1}(|\sigma|)$. In particular, we can pull $M$ and $W$ back to $\Delta_k$ by $\sigma$, and we can extend the pullbacks $\bar{p} \colon \overline{M} \to \Delta_k$ and $\overline{W}$ smoothly over a neighbourhood $U$ of $\Delta_k$ in $\mathbb{R}^k$, preserving all transversal intersections. Because $\overline{W}$ takes the form of Proposition 1.21 (1) near zeros of $\overline{W}$, Laudenbach's arguments in [La] imply that the closure $\overline{\overline{M}^u_{\overline{W}}(\bar{x}, U)}$ of the unstable cell at the pullback $\bar{x}$ of $x$ is a submanifold with conical singularities whose strata are unions of other unstable cells.

Let $H$ be an affine hyperplane that contains a codimension 1 face of $\Delta_k$. Me may assume that $\bar{p}^{-1}(H)$ is preserved by the flow. Assume first that $\bar{p}(\bar{x}) \notin H$, then $\overline{\overline{M}^u_{\overline{W}}(\bar{x}, U)}$ may hit $\bar{p}^{-1}(H)$, but it cannot extend across $\bar{p}^{-1}(H)$. On the other hand, assume that $\bar{p}(\bar{x}) \in H$. Then it is easy to see that either $\overline{M}^u_{\overline{W}}(\bar{x}, U)$ is contained in $\bar{p}^{-1}(H)$, or $\overline{M}^u_{\overline{W}}(\bar{x}, U)$ intersects $\bar{p}^{-1}(H)$ transversally. Similar statement hold for all zeros $\bar{x}' \in \overline{\overline{M}^u_{\overline{W}}(\bar{x}, U)}$. This implies that $\overline{\overline{M}^u_{\overline{W}}(\bar{x}, U)} \cap \bar{p}^{-1}(\Delta_k)$ is again a submanifold with conical singularities, which has some extra strata given as intersections of unstable cells with $\bar{p}$-preimages of faces of $\Delta_k$.

Going back to $M$, we know now that the strata of $\overline{M^u_W(x, \sigma)}$ are unions of sets of the form

$$M^u_W(x', \sigma) \cap p^{-1}(|\sigma'|) = M^u_W(x', \sigma') , \tag{1.26}$$

where $\sigma'$ is a face of $\sigma$, and $x' \in \overline{M^u_W(x,\sigma)} \cap p^{-1}(|\sigma'|)$ is a zero of $W$. A similar statement holds if we reverse the roles of the stable and unstable cells.

Finally assume that $h'\colon |\sigma| \to \mathbb{R}$ is a smooth function such that $h'(q)$ is a regular value of $h|_{p^{-1}(q)}$ for all $q \in |\sigma|$. We also assume that for one of the unstable cells $M^u_W(x,\sigma)$, all $M^u_W(x',\sigma)$ with $x' \in \overline{M^u_W(x,\sigma)} \cap h^{-1}(h',\infty)$ meet all $M^s_W(y,\sigma)$ with $y \in h^{-1}(h',\infty)$ transversally in the sense of Remark 1.13. Then the same argument as above implies that $\overline{M}^u_W(x,\sigma) \cap h^{-1}(h')$ is a submanifold of $h^{-1}(h')$ with conical singularities, whose strata are of the form $M^u_W(x',\sigma') \cap h^{-1}(h')$, for $\sigma'$ a face of $\sigma$ and $x' \in \overline{M^u_W(x,\sigma)} \cap p^{-1}(|\sigma'|) \cap h^{-1}(h',\infty)$. □

We use Proposition 1.25 inductively to vertically deform $W$ such that Smale's transversality condition holds in $p^{-1}(|\sigma|)$ in the sense of Remark 1.13.

**1.27. Proposition.** *After a small vertical perturbation of $W$ outside $U_{2\varepsilon}(C)$ that is smooth on each $p^{-1}(|\sigma|)$, the closures of the stable and unstable cells $\overline{M^u_W(x,\sigma)}$, $\overline{M^s_W(y,\sigma)}$ are submanifolds of $M$ with conical singularities. The strata of $\overline{M^u_W(x,\sigma)}$ are unions of sets $M^u_W(x',\sigma')$, where $\sigma'$ is a face of $\sigma$, and $x'$ is a zero of $W$. If $W$ already satisfies the transversality condition of Remark 1.13 over each closed simplex of a simplicial subset $B_0 \subset B$, then $W$ can be deformed in such a way that $W|_{\pi^{-1}(B_0)}$ is not affected.*

*Proof.* By Proposition 1.25, it suffices to establish transversality in the sense above. We proceed by induction of the degree $k$ of the simplices of $S$, modifying $W$ only over those $k$-simplices where the the stable and unstable cells do not yet meet transversally. Then the relativity statement above is automatically satisfied. For $k = 0$, we can achieve transversality by the classical argument of Smale ([S]).

Thus assume that $\sigma \in S_k$ with $k \geq 1$, and that transversality in the sense of Remark 1.13 already holds over each face of $\partial\,|\sigma|$. Let $x$, $y$ be two zeros of $W$ on $p^{-1}(|\sigma|)$ with $\{x\} = C_\alpha(\sigma_i)$ and $\{y\} = C_\beta(\sigma_j)$. If $M^u_W(x,\sigma) \cap M^s_W(y,\sigma) \neq \emptyset$, then $0 \leq i \leq j \leq k$, and

$$p\big(M^u_W(x,\sigma) \cap M^s_W(y,\sigma)\big) \subset |\sigma_{i,\dots,j}| \setminus \big(|\sigma_{i+1,\dots,j}| \cup |\sigma_{i,\dots,j-1}|\big) .$$

Unless $i = 0$ and $j = k$, the simplex $\sigma_{i,\dots,j}$ is a proper face of $\sigma$, and we know by induction that the intersection of $M^u_W(x,\sigma)$ and $M^s_W(y,\sigma)$ is transversal over $|\sigma_{i,\dots,j}|$. Let $\gamma(t) \in M^u_W(x,\sigma) \cap M^s_W(y,\sigma)$ for some flow-line $\gamma(t)$ of $W$, then $\lim_{t\to-\infty} = x$ and $\lim_{t\to\infty} = y$. By definition of $W_k$ and $W$ and using Proposition 1.17, it is easy to check that

(1.28)
$$\begin{aligned} &T_{\gamma(t)} M^u_W(x,\sigma) \equiv T_{\gamma(t)}\big(p^{-1}\,|\sigma_{i,\dots,k}|\big) \qquad \operatorname{mod} T_{\gamma(t)}\big(p^{-1}\,|\sigma_{i,\dots,j}|\big) \\ \text{and}\qquad &T_{\gamma(t)} M^s_W(y,\sigma) \equiv T_{\gamma(t)}\big(p^{-1}\,|\sigma_{0,\dots,j}|\big) \qquad \operatorname{mod} T_{\gamma(t)}\big(p^{-1}\,|\sigma_{i,\dots,j}|\big) \end{aligned}$$

for all $t \in \mathbb{R}$. Together with transversality over $|\sigma_{i,\dots,j}|$, i.e.,

$$T_{\gamma(t)}\big(p^{-1}\,|\sigma_{i,\dots,j}|\big) \subset T_{\gamma(t)} M^u_W(x,\sigma) + T_{\gamma(t)} M^s_W(y,\sigma) ,$$

we obtain transversality in $p^{-1}(|\sigma|)$.

Now assume that $x \in C_\alpha(|\sigma|)$ and $y \in C_\beta(|\sigma|)$ with $p(x) = \sigma_0$ and $p(y) = \sigma_k$. By Proposition 1.21, we must have either $\beta = \alpha$ or $\beta \prec_\sigma \alpha$. In the first case, the intersection takes place along

$$C_\alpha\big(|\sigma| \setminus (|\sigma_{1,\dots,k}| \cup |\sigma_{0,\dots,k-1}|)\big) \tag{1.29}$$

and is clearly transversal by Proposition 1.21.

We are left with the case $\{x\} = C_\alpha(\sigma_0)$, $\{y\} = C_\beta(\sigma_k)$ with $\beta \prec_\sigma \alpha$. We start another induction over $\alpha$, where we complete $\prec_\sigma$ to a total ordering $<_\sigma$. Fix $\alpha$ and assume that $M^u_W(x',\sigma)$ and $M^s_W(y,\sigma)$ intersect transversally for all zeros $x$, $y$ of $W$ with $x \in C_{\alpha'}(|\sigma|)$ for some $\alpha' <_\sigma \alpha$, and with $y \in C(|\sigma|)$. This is in particular true for the smallest $\alpha$.

Let $S_{2\varepsilon}T^uX = \partial D_{2\varepsilon}T^uX$ be the sphere bundle in $T^uX$ of radius $2\varepsilon$. We may fix a small open neighbourhood $U_\alpha$ of $S_{2\varepsilon}T^uX|_{C_\alpha(|\sigma|)}$ in the hypersurface

$$h^{-1}\big(h_\alpha - 2\varepsilon^2\big) \cap p^{-1}(|\sigma|) = \big\{\, z \in p^{-1}(|\sigma|) \,\big|\, h(z) > h_\alpha\big(p(z)\big) - 2\varepsilon^2 \,\big\}\ ,$$

such that $U_\alpha \cap U_{2\varepsilon}((C'_\alpha(|\sigma|)) = \emptyset$ for all $\alpha' <_\sigma \alpha$. Then we have in particular $U_\alpha \cap M^u_W(x',\sigma) = \emptyset$ for all $x' \in C_{\alpha'}(|\sigma|)$ with $\alpha' <_\sigma \alpha$, so a small modification of $W$ along $U_\alpha$ over the interior $\mathring{\sigma}$ of $|\sigma|$ will not destroy any transversal intersection established in a previous step. Note that for $U_\alpha$ small enough, we get

$$M^u_W(x,\sigma) \cap U_\alpha = S_{2\varepsilon}T^uX|_{C_\alpha(|\sigma|)}\ .$$

Our goal is to find a diffeomorphism $F_\alpha$ of $U_\alpha$ that is an arbitrarily small perturbation of $\mathrm{id}_{U_\alpha}$, that equals $\mathrm{id}_{U_\alpha}$ outside some compact subset of $U_\alpha \cap p^{-1}(\mathring{\sigma})$, and that preserves fibres of $p$, such that the intersections of

$$F_\alpha\big(S_{2\varepsilon}T^uX|_{C_\alpha(|\sigma|)}\big)$$

with all $M^s_W(y,\sigma) \cap U_\alpha$ become transversal in $U_\alpha$. We will then cut $M$ along $U_\alpha \subset h^{-1}(h_\alpha - 2\varepsilon^2)$ and glue the two sides back using $F_\alpha$. Clearly, we can get the same effect by a suitable perturbation of $W$ "just below" $h^{-1}(h_\alpha - 2\varepsilon^2)$, so that $W|_{U_{2\varepsilon}(C)}$ is unaffected. Moreover, Proposition 1.21 still holds for the modified vector field $W$.

By induction over $\alpha$ and by Proposition 1.25, we know that $M^s_W(y,\sigma) \cap U_\alpha$ is a submanifold of $h^{-1}(h_\alpha - 2\varepsilon^2)$ with conical singularities for all zeros $y$ of $W$, whose strata are unions of sets of the form $M^s_W(y',\sigma') \cap U_\alpha$. Let $\sigma'$ be a proper face of $\sigma$, then we know by induction that all $M^s_W(y,\sigma') \cap U_\alpha$ intersect $S_{2\varepsilon}T^uX|_{C(|\sigma|')}$ transversally in $p^{-1}(|\sigma'|)$. As above, this implies that $M^s_W(y,\sigma') \cap U_\alpha$ intersects $S_{2\varepsilon}T^uX|_{C(|\sigma|)}$ transversally in $p^{-1}(|\sigma|)$.

By Remark 1.24, the intersection of $M^s_W(y,\sigma) \cap U_\alpha$ with $S_{2\varepsilon}T^uX|_{C(|\sigma|)}$ is transversal in a neighbourhood of $p^{-1}(\partial\,|\sigma|)$. There clearly is a diffeomorphism $F_\alpha$ as above such that $F_\alpha(S_{2\varepsilon}T^uX|_{C_\alpha(|\sigma|)})$ intersects all $M^s_W(y',\sigma)$ transversally, if we allow $F_\alpha$ to mix the fibres of $p$. However, because $p$ restricts to a submersion $S_{2\varepsilon}T^uX|_{C_\alpha(|\sigma|)} \to |\sigma|$, we can always take $F_\alpha$ such that it preserves the fibres of $p$. After the gluing construction described above, $M^u_W(x,\sigma)$ now intersects all $M^s(y,\sigma)$ transversally in a neighbourhood of $U_\alpha$. Moreover, the glueing construction is equivalent to a small vertical modification of $W$ over $\mathring{\sigma}$, and "just below" $h^{-1}(h_\alpha - 2\varepsilon^2)$. Since transversality is preserved by the flow of $W$, these intersections are everywhere transversal. This finishes the last step of both inductions. $\square$

For each $k$, each $k$-simplex $\sigma \in S_k$ and each $\alpha \in P(\sigma)$, let

$$M_\alpha(\sigma) = \overline{M^u_W(x,\sigma)} \tag{1.30}$$

denote the closure of the $W$-unstable manifold of the critical point $x \in C_\alpha(\sigma_0)$ over $|\sigma|$. Then $M_\alpha(\sigma)$ is a submanifold of $M$ with conical singularities whose interior is homeomorphic to a $k + \mathrm{ind}_\alpha$-dimensional disc. The following proposition shows that the $M_\alpha(\sigma)$ form a CW-structure on $M$ that behaves nicely with respect to the projection $p\colon M \to B$ and the simplicial complex $S$ on $B$. We write $M_\alpha(\sigma)$ for the generators of the corresponding chain complex, and use $\partial$ to denote its boundary operator. Let us for the moment assume that the bundle $T^uX \to C$ is oriented, so that all cells $M_\alpha(\sigma)$ carry a natural orientation that agrees with the orientation of $|\sigma| \times T^uX_{\sigma_0}$ at $\sigma_0$.

**1.31. Proposition.** *For $\alpha$, $\beta \in P(\sigma)$, there exist an integer $a_{\alpha\beta}(\sigma)$, such that*

(1) $$a_{\alpha\beta}(\sigma) = 0 \qquad \text{unless } \beta \prec_\sigma \alpha \text{ and } \operatorname{ind}_\alpha + k = \operatorname{ind}_\beta + 1,$$

(2) $$\partial M_\alpha(\sigma) = \sum_j (-1)^j M_\alpha(\sigma_{0\dots\widehat{j}\dots k}) + \sum_j \sum_\beta (-1)^{k(j-1)} a_{\alpha\beta}(\sigma_{0\dots j})\, M_\beta(\sigma_{j\dots k}) \; ,$$

(3) *and* $$0 = \sum_j (-1)^j a_{\alpha\beta}(\sigma_{0\dots\widehat{j}\dots k}) + \sum_j \sum_\gamma (-1)^{k(j-1)} a_{\alpha\gamma}(\sigma_{0\dots j})\, a_{\gamma\beta}(\sigma_{j\dots k}) \; .$$

Note that equation (1) makes sure that only cells of dimension $k + \operatorname{ind}_\alpha - 1$ enter the second sum in (2). The right hand side of (2) reads like the dual of a "simplicial" superconnection, consisting of the simplical boundary operator on $B$ and a sum of products with simplicial cochains valued in the endomorphisms of a bundle $\hat{p}_*(C \times \mathbb{Z})$ of free abelian groups. Then (3) simply says that this superconnection is flat.

*1.32. Remark.* In [I2], Definition 3.1.1, Igusa defines *system of higher homotopies* to be systems of endomorphism $e(v_0, \dots, v_k) \in \operatorname{End}(V)$ for $n \geq v_0 \geq \dots \geq v_k \geq 0$, satisfying a non-degeneracy axiom and axioms similar to Proposition 1.31 (1) and (3), where (3) above is replaced by

$$\sum_{j=0}^{k} (-1)^j \big( e(v_0, \dots, v_j)\, e(v_j, \dots, v_k) - e(v_0, \dots, \widehat{v_j}, \dots, v_k) \big) \; .$$

In particular, if we fix a simplex $\sigma \in S_n$ and set

$$e_\sigma(v_0, \dots, v_k) = \begin{cases} (-1)^{\frac{k(k-1)}{2}}\, a(\sigma_{v_0,\dots,v_k}) & \text{if } v_0 < \dots < v_k \text{, and} \\ 0 & \text{otherwise,} \end{cases}$$

we obtain a system of higher homotopies in the sense of Igusa.

If $T^u X \to C$ is orientable, let $G = \{1\}$ be the trivial group, otherwise take $G = \{\pm 1\}$. Then the structure group of $o(T^u X)$ is $G$, and the structure group of $\hat{p}_* o(T^u X) \to B$ is a subgroup of the group of $G$-monomial matrices, i.e., the group of invertible matrices having precisely one non-zero element in each row and each column, which is an element of $G$. Putting the $e_\sigma$ together for all simplices $\sigma \in S$, we get a functor $\xi$ from the category of simplices of $S$ to Igusa's Whitehead space $Wh(\mathbb{Z}, G)$ (which is constructed as a simplicial category), and $\xi$ is *$G$-monomial* in the sense of Chapter 3.3 in [I2].

*Proof* of Proposition 1.31. Let $\sigma \in S_k$, and let $x$ be the critical point in $C_\alpha$ over the vertex $\sigma_0$. Recall that by Proposition 1.27, the strata of $\overline{M^u_W(x,\sigma)}$ are unions of sets $M^u_W(x', \sigma')$, with $\sigma'$ a face of $\sigma$, and $x' \in \overline{M^u_W(x,\sigma)} \cap p^{-1}(|\sigma'|)$. We now have to identify the codimension-1-stratum of $\overline{M^u_W(x,\sigma)}$. It is a union of two types of unstable cells $M^u_W(x', \sigma')$. The first type has $\sigma' = \sigma$, so assume that $x' \in C_\beta$ lies over the vertex $\sigma_j$. Then the dimension of $M^u_W(x', \sigma) = M^u_W(x', \sigma_{j,\dots,k})$ is $\operatorname{ind}_\beta + k - j$. In particular, $M^u(x', \sigma)$ is of codimension 1 iff $\operatorname{ind}_\alpha + j = \operatorname{ind}_\beta + 1$ as in (1). This way we obtain all terms in the second sum in (2) for $\beta \prec_\sigma \alpha$, and the very first term $M_\alpha(\sigma_{1,\dots,k})$ in the first sum, cf. (1.29).

The remaining contributions to the codimension-1-stratum arise as intersections of unstable cells with $p^{-1}(\partial\,|\sigma|)$. The only such intersections of codimension 1 in $\overline{M^u_W(x,\sigma)}$ are of the form $\overline{M^u_W(x,\sigma')}$ with $\sigma' = \sigma_{0,\dots,\widehat{j},\dots,k}$ for $1 \leq j \leq k$, accounting for the rest of the first sum in (2). Note that $j = 0$ is impossible here, because for a transversal intersection, we had to assume that $p(x) = \sigma_0 \in |\sigma'|$ in the proof of Proposition 1.27.

The preceding arguments show that the boundary of $M_\alpha(\sigma)$ as a chain is given by (2), with constants $a_{\alpha\beta}(\sigma_{0\dots j})$ as in (1), which might, however, still depend on $\sigma$. To see that this is not the case, and to explain the sign factor $(-1)^{k(j-1)}$, we define $a_{\alpha\beta}(\sigma_{0\dots j})$ to be the coefficient of $M_\beta(\sigma_j)$ in the (chain) boundary of $M_\alpha(\sigma_{0\dots j})$. Suppose that we orient both $M_\alpha(\sigma)$ and $M_\beta(\sigma_{j\dots k})$ by placing an oriented base of $\sigma_{j\dots k}$ to the right of oriented bases of $M_\alpha(\sigma_{0\dots j})$ and $M_\beta(\sigma_j)$, respectively. By the arguments after (1.28) in the proof of Proposition 1.27, the coefficient of $M_\beta(\sigma_{j\dots k})$ in $\partial M_\alpha(\sigma)$ is still $a_{\alpha\beta}(\sigma_{0\dots j})$ with respect to these orientations. However, since we actually have to fit the oriented base of $\sigma_{j\dots k}$ to the left of the bases of $T^uX|_{C_\alpha}$ and $T^uX|_{C_\beta}$ according to the conventions fixed at the beginning of this subsection and before Proposition 1.31, these new orientations differ from the actual ones by factors of $(-1)^{(k-j)\operatorname{ind}_\alpha}$ and $(-1)^{(k-j)\operatorname{ind}_\beta}$, respectively. By (1), $\operatorname{ind}_\alpha - \operatorname{ind}_\beta = 1-j$, so we get the sign factor $(-1)^{(k-j)(j-1)} = (-1)^{k(j-1)}$ as claimed.

Finally, because $\partial^2 = 0$,

$$\begin{aligned}
0 &= \partial^2 M_\alpha(\sigma)\\
&= \sum_j (-1)^j \partial M_\alpha(\sigma_{0\dots\widehat{j}\dots k})\\
&\qquad + \sum_j \sum_\beta (-1)^{k(j-1)} a_{\alpha\beta}(\sigma_{0\dots j})\, \partial M_\beta(\sigma_{j\dots k})\\
&= \sum_{i<j} \sum_\beta (-1)^j\, (-1)^{(k-1)(i-1)} a_{\alpha\beta}(\sigma_{0\dots i})\, M_\beta(\sigma_{i\dots\widehat{j}\dots k})\\
&\qquad + \sum_{i>j} \sum_\beta (-1)^j\, (-1)^{(k-1)i} a_{\alpha\beta}(\sigma_{0\dots\widehat{j}\dots i})\, M_\beta(\sigma_{i\dots k})\\
&\qquad + \sum_{i\ge j} \sum_\beta (-1)^{k(j-1)}\, (-1)^{i-j} a_{\alpha\beta}(\sigma_{0\dots j})\, M_\beta(\sigma_{j\dots\widehat{i}\dots k})\\
&\qquad + \sum_{i\ge j} \sum_{\beta,\gamma} (-1)^{k(j-1)}\, (-1)^{(k-j)(i-j-1)} a_{\alpha\beta}(\sigma_{0\dots j})\, a_{\beta\gamma}(\sigma_{j\dots i})\, M_\gamma(\sigma_{i\dots k})\\
&= \sum_j \sum_\beta (-1)^{(k-1)j} \Bigg( \sum_{i\le j} (-1)^i a_{\alpha\beta}(\sigma_{0\dots\widehat{i}\dots j})\\
&\qquad\qquad \sum_{i\le j} \sum_\gamma (-1)^{j(i-1)} a_{\alpha\gamma}(\sigma_{0\dots i})\, a_{\gamma\beta}(\sigma_{i\dots j}) \Bigg) M_\beta(\sigma_{j\dots k})\;,
\end{aligned}$$

which implies (3). □

**1.d. Simplicial complexes and good coverings.** In the previous subsections, we constructed a CW complex on $M$ over a simplicial complex $S$ on $B$ and gave an adapted basis $M_\alpha(\sigma)$ for the associated chain complex. Each chain $M_\alpha(\sigma)$ lived over a simplex $\sigma \in S$ of $B$. For the following constructions however, we need families of cells and chains over subsets of $B$, such that we can integrate over the fibres of these cells and chains locally on $B$. We therefore modify the constructions of the previous subsections. We will need a sufficiently fine simplicial complex, in the following sense.

*1.33. Remark.* Let $\pi_1, \pi_2\colon B \times B \to B$ denote the two product projections, and let $\varepsilon > 0$ be defined as in Proposition 1.10. Then there exists a neighbourhood $U_{\Delta B}$ of the diagonal $\Delta B$ in $B \times B$ and a diffeomorphism $F\colon (\pi_2^*M)|_{U_{\Delta B}} \to (\pi_1^*M)|_{U_{\Delta B}}$ of fibre-bundles extending $\mathrm{id}_M\colon (\pi_2^*M)|_{\Delta B} \to (\pi_1^*M)|_{\Delta B}$, such that

(1) $|h \circ F_{q_1,q_2} - h| < \frac{\varepsilon^2}{2}$ for all $(q_1, q_2) \in U_{\Delta B}$, and

(2) $F_{q_1,q_2}(\exp_c^X(z^s,z^u)) = \exp^X_{F_{q_1,q_2}c}(F_{q_1,q_2,*c}z^s, F_{q_1,q_2,*c}z^u)$ for $(q_1,q_2)$ as above, $c \in \widehat{p}^{-1}(q_2)$, and $z^{s,u} \in T_c^{s,u}X$.

We may choose the smooth ordered simplicial complex $S$ so fine that

$$|\sigma| \times |\sigma| \subset U_{\Delta B} \qquad \text{for all } \sigma \in S. \tag{1.34}$$

For $\sigma \in S$, let $\mathring{\sigma}$ denote the interior of $|\sigma|$ as a submanifold of $M$, and set

$$U(\sigma) = \bigcup_{\substack{\sigma' \in S \\ \sigma \text{ is a face of } \sigma'}} \mathring{\sigma}' \tag{1.35}$$

for each $\sigma \in S$. We assume that $S$ is so fine that the closure $\overline{U(\sigma)}$ is contractible for all $\sigma \in S$. In particular,

$$\mathcal{U}(S) = \left\{ U(\sigma) \;\middle|\; \sigma \in S_0 \right\} \tag{1.36}$$

is a good open cover of $B$. Let us define

$$\overline{B}(\sigma) = \overline{U(\sigma)} \times |\sigma| \qquad \text{and} \qquad \overline{B} = \bigcup_{\sigma \in S} \overline{B}(\sigma) \;\subset\; B \times B\;, \tag{1.37}$$

then $\overline{B}(\sigma)$ is a manifold with corners for all $\sigma$. Note that by (1.34),

$$\overline{B} = \bigcup_{\sigma \in S} |\sigma| \times |\sigma| \;\subset\; U_{\Delta B}\;.$$

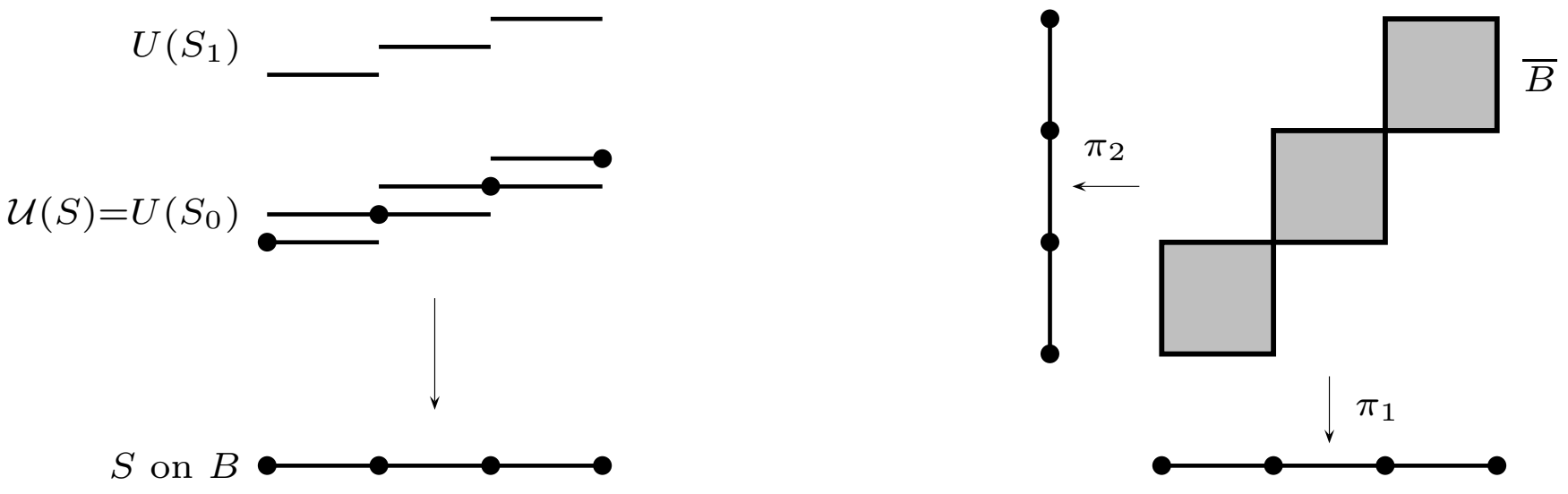

Figure 1.38. A one-dimensional simplical complex $S$, the cover $\mathcal{U}(S)$, and $\overline{B}$.

Consider the bundle $\overline{M} = (\pi_1^*M)|_{\overline{B}} \to \overline{B}$ with its fibre-wise Morse function $\overline{h} = \pi_1^*h$. Then

$$\overline{M}(\sigma) = (\pi_1^*M)|_{\overline{B}(\sigma)} \;\longrightarrow\; \overline{U(\sigma)} \tag{1.39}$$

is a bundle with fibre $|\sigma| \times X$. For each $\sigma$, each $\alpha \in P(\sigma)$, and each $q_1 \in \overline{U(\sigma)}$, we define a submanifold

$$\begin{aligned} \overline{M}_\alpha(\sigma)|_{q_1} &= F_{q_1,\cdot}\big(\{q_1\} \times M_\alpha(\sigma)\big) \\ &= \left\{ F_{q_1,q_2}x \;\middle|\; x \in M_\alpha(\sigma), \text{ and } q_2 = p(x) \right\} \;\subset\; \overline{M}(\sigma)|_{q_1} \end{aligned} \tag{1.40}$$

with conical singularities. We still assume that the bundle $T^uX|_C$ is oriented.

**1.41. Proposition.** *The submanifolds $\overline{M}_\alpha(\sigma)|_{q_1} \subset \overline{M}(\sigma)|_{q_1}$ with conical singularities form a bundle $\overline{M}_\alpha(\sigma) \to \overline{U(\sigma)}$, such that*

(1) *for each $\sigma \in S_k$, the restrictions of the $\overline{M}_\alpha(\sigma')$ to $\overline{U(\sigma)}$ for all $\alpha$ and all faces $\sigma'$ of $\sigma$ form a family of fibre-wise CW-complexes on $\overline{M}(\sigma) \to \overline{U(\sigma)}$;*

(2) *let $a_{\alpha\beta}(\sigma)$ be the coefficients of Proposition 1.31, then the boundary of $\overline{M}_\alpha(\sigma)$ as a fibre-wise chain is given by*

$$\partial \overline{M}_\alpha(\sigma) = \sum_j (-1)^j \overline{M}_\alpha(\sigma_{0\dots\hat{j}\dots k})|_{\overline{U(\sigma)}} + \sum_j \sum_\beta (-1)^{k(j-1)} a_{\alpha\beta}(\sigma_{0\dots j})\, \overline{M}_\beta(\sigma_{j\dots k})|_{\overline{U(\sigma)}} \,;$$

(3) *we have*

$$\overline{M}_\alpha(\sigma) \cap \left\{\, (q_1, q_2, x) \in \overline{M} \;\middle|\; h(x) > h_\alpha(q_1) - \varepsilon^2 \,\right\} = \left(\pi_1^* D_{\sqrt{2}\,\varepsilon} T^u X\right)|_{\overline{B}(\sigma)} \,.$$

*Proof.* From the construction (1.40), it follows that $\overline{M}_\alpha(\sigma) \to \overline{U(\sigma)}$ is a bundle with fibre diffeomorphic to the $M_\alpha(\sigma)$ of Proposition 1.31. Also, the $\overline{M}_\alpha(\sigma_{i_0\dots i_j})|_{\overline{U(\sigma)}}$ form a fibre-wise CW-complex on $\overline{M}(\sigma) \to \overline{U(\sigma)}$ with fibre diffeomorphic to the CW-complex on $p^{-1}(|\sigma|)$ formed by the $M_\alpha(\sigma_{i_0\dots i_j})$ for all $\alpha$ and all $0 \le i_0 < \dots < i_j \le k$, so we get (1). Moreover, from Proposition 1.31 (2), we get (2). Finally, (3) follows from Proposition 1.21 by the specifications in Remark 1.33 and (1.40). □

We want to state a dual version of Proposition 1.31 and Proposition 1.41. Let $F \to M$ be a flat vector bundle. We no longer assume that the bundle $T^u X \to C$ is oriented, and set $V_\alpha = o(T^u X) \otimes F|_{C_\alpha} \to B$ over any simplex $|\sigma|$. We define $V \to B$ as in Definition 1.1, so $V = \bigoplus_\alpha V_\alpha$.

Let $\overline{F} = (\pi_1^* F)|_{\overline{B}} \to \overline{M}$ be the flat bundle over $\overline{M}$ induced by $F$, and let $\overline{F}(\sigma) = \overline{F}|_{\overline{M}(\sigma)}$. We have natural pull-back maps

$$\pi_1^*(\sigma)\colon\ \Omega^*\bigl(p^{-1}(\overline{U(\sigma)}); F\bigr) \;\longrightarrow\; \Omega^*\bigl(\overline{M}(\sigma); \overline{F}\bigr) \,,$$

and we define integration maps by integrating over the fibres of $\overline{M}_\alpha(\sigma)/\overline{U(\sigma)}$,

(1.42)
$$\begin{aligned} I_\alpha(\sigma) &= \int_{\overline{M}_\alpha(\sigma)/\overline{U(\sigma)}} \circ\, \pi_1^*(\sigma)\colon\ \Omega^*\bigl(p^{-1}(\overline{U(\sigma)}); F\bigr) \;\longrightarrow\; \Omega^*\bigl(\overline{U(\sigma)}; V_\alpha\bigr) \\ \text{and}\qquad I(\sigma) &= \bigoplus_\alpha I_\alpha(\sigma)\colon\ \Omega^*\bigl(p^{-1}(\overline{U(\sigma)}); F\bigr) \;\longrightarrow\; \Omega^*\bigl(\overline{U(\sigma)}; V\bigr) \,. \end{aligned}$$

**1.43. Corollary.** *For all $k$ and all $\sigma \in S_k$, there exist $\nabla^V$-parallel endomorphisms $a(\sigma) \in \Gamma(\operatorname{End}^{1-k} V|_{\overline{U(\sigma)}})$ such that we have*

(1)
$$\begin{aligned} \nabla^V \circ I(\sigma) = (-1)^k\, I(\sigma) \circ \nabla^F &+ \sum_j (-1)^j\, I(\sigma_{0\dots\hat{j}\dots k})|_{\overline{U(\sigma)}} \\ &+ \sum_j (-1)^{k(j-1)}\, a(\sigma_{0\dots j}) \circ I(\sigma_{j\dots k})|_{\overline{U(\sigma)}} \,, \end{aligned}$$

(2)
$$\begin{aligned} \textit{and}\qquad 0 = &\sum_j (-1)^j\, a(\sigma_{0\dots\hat{j}\dots k})|_{\overline{U(\sigma)}} \\ &+ \sum_j (-1)^{k(j-1)}\, a(\sigma_{0\dots j}) \circ a(\sigma_{j\dots k})|_{\overline{U(\sigma)}} \,. \end{aligned}$$

*Suppose that* $v \in V_\beta|_{\overline{U(\sigma)}}$, *then*

$$a(\sigma)\, v \in \bigoplus_{\beta \prec_\sigma \alpha} V_\alpha \; . \tag{3}$$

*Suppose that* $\omega \in \Omega^*(\overline{U(\sigma)}; F)$ *has support in* $h^{-1}(h_\alpha - \varepsilon^2, \infty) \subset M|_{U_\sigma}$, *then*

$$I_\alpha(\sigma)\, \omega = \begin{cases} \int_{D_{2\varepsilon} T^u X|_{C_\alpha(\sigma)}/\overline{U(\sigma)}} \omega & \text{if } \sigma \in S_0, \text{ and} \\ 0 & \text{otherwise.} \end{cases} \tag{4}$$

Note that (4) implies in particular that if $\operatorname{supp} \omega \subset h^{-1}(h_\alpha - \varepsilon^2, \infty)$ and $\beta \prec_\sigma \alpha$, then $I_\beta(\sigma)\, \omega = 0$.

*1.44. Remark.* For each simplex $\sigma \in S$, the endomorphisms $a(\sigma')$ for all faces $\sigma'$ of $\sigma$ give rise to a system $e_\sigma$ of higher homotopies in the sense of Igusa, just as in Remark 1.32. Suppose that $R$ is an associative ring with unit, with a fixed complex representation $\rho\colon \mathrm{End}(\mathbb{C}^r)$ for some $r \geq 1$, and that $G$ is a subgroup of $\mathrm{GL}_s(R)$ for some $s \geq 1$, such that the flat bundle $F$ has structure group $\rho(G) \subset \mathrm{GL}_{rs}(\mathbb{C})$ (if $T^u X \to C$ is not orientable, we assume that $-1 \in G$). Then the coefficients $a_{\alpha\beta}(\sigma)\colon V_\beta \to V_\alpha$ of $a(\sigma)$ naturally take values in the matrices $M_s(R)$ over $R$, and the structure group of $V$ is a subgroup of the group of $G$-monomial matrices (matrices having exactly one non-zero entree in each row or column, which is an element of $G$). Then the $e_\sigma$ above form a $G$-monomial functor from the category of simplices of $S$ to the Whitehead space $Wh(R, G)$.

*Proof* of Corollary 1.43. Consider the families of cochain-complexes associated to the bundles of fibre-wise CW-complex on $\overline{M}(\sigma) \to \overline{U(\sigma)}$ constructed in Proposition 1.41, with values in $\overline{F}$. Since we now have to take the possible orientations of $T^u X$ into account, each cell $\overline{M}_\alpha(\sigma')$ contributes a copy $V_\alpha(\sigma')$ of $V_\alpha|_{\overline{U(\sigma)}} \to \overline{U(\sigma)}$ to the family of complexes over $\overline{U(\sigma)}$, where $\sigma'$ is a face of $\sigma$. For any fibre-wise cochain $s$, let $s_\alpha(\sigma')$ denote the component of $s$ in $V_\alpha(\sigma')$. By Proposition 1.31 and Proposition 1.41, the coboundary map of this complex takes the form

$$(\delta s)_\alpha(\sigma') = \sum_j (-1)^j s_\alpha(\sigma'_{0\dots\hat{j}\dots k}) + \sum_j \sum_\beta (-1)^{k(j-1)} a_{\alpha\beta}(\sigma'_{0\dots j})\, s_\beta(\sigma'_{j\dots k}) \tag{1.45}$$

if $\sigma' \in S_k$ is a face of $\sigma \in S$, where the $a_{\alpha\beta}(\sigma'_{0\dots j})\colon V_\beta \to V_\alpha$ are parallel homomorphisms satisfying Proposition 1.31 (1) and (3).

Let us define

$$a(\sigma') = \bigoplus_{\alpha,\beta} a_{\alpha\beta}(\sigma')\colon\; V|_{U(\sigma')} \longrightarrow V|_{U(\sigma')} \; ,$$

then (1) follows from (1.45) and Stokes' theorem, while (2) and (3) follow from Proposition 1.31 (1) and (3).

If $\operatorname{supp} \omega \subset h^{-1}(h_\alpha - \varepsilon^2, \infty)$, then by Proposition 1.41 (3) and (1.42),

$$I_\alpha(\sigma), \omega = \int_{(\pi_1^* D_{2\varepsilon} T^u X)|_{\overline{B}(\sigma)}/\overline{U(\sigma)}} \pi_1^*(\sigma)\, \omega \; ,$$

and $\pi_1^*(\sigma)\omega$ vanishes in the directions of $|\sigma|$, so this integral vanishes unless $\sigma \in S_0$. $\square$

**1.e. A "mixed" superconnection.** Starting from the endomorphisms $a(\sigma)$ of the previous subsection, we construct a family of compatible flat superconnections $\nabla^{\overline{V}} + a'(\sigma,\emptyset)$ on $V|_{|\sigma|}$, where $a'(\sigma,\emptyset) \in \Omega^*(|\sigma|\,,\mathrm{End}_+ V)$ is an $\mathrm{End}\,V$-valued differential form of total degree 1. We also construct a family of compatible "integration maps" $I'(\sigma,\emptyset)\colon\Omega^*(M;F) \to \Omega^*(|\sigma|\,;V)$ of total degree zero that satisfy $I'(\sigma,\emptyset)\circ\nabla^F = (\nabla^{\overline{V}} + a'(\sigma,\emptyset))\circ I'(\sigma,\emptyset)$. Using a partition of unity, we pull $\nabla^{\overline{V}} + a'(\emptyset)$ and $I'(\emptyset)$ back to a flat superconnection $A'$ on $V \to B$ of total degree one and an integration map $I\colon\Omega^*(M;F) \to \Omega^*(B;V)$ satisfying $I\circ\nabla^F = A'\circ I$.

If $\sigma' = \sigma_{i_0\dots i_j}$ is a face of $\sigma \in S_k$, set

$$|\sigma \setminus \sigma'| = \left|\sigma_{i_j,i_j+1,\dots,k}\right| \;, \tag{1.46}$$

so $|\sigma\setminus\sigma'|$ is the simplex spanned by all vertices of $\sigma$ that are "higher than or equal to" all vertices of $\sigma'$. For convenience, we introduce the *empty simplex* $\emptyset$ of formal dimension $-1$, which we regard as a face $\emptyset = \sigma_\emptyset$ of all $\sigma \in S$, so that we may write $|\sigma| = |\sigma\setminus\emptyset|$.

Let $a(\sigma) \in \Gamma(\mathrm{End}^{1-k} V|_{\overline{U(\sigma)}})$ be the parallel endomorphisms constructed in Corollary 1.43.

**1.47. Proposition.** *For each $0 \le k \le l+1$ and each (possibly empty) face $\sigma' \in S_{k-1}$ of $\sigma \in S_l$, there exist an $\mathrm{End}\,V$-valued differential form $a'(\sigma,\sigma')$ on $|\sigma\setminus\sigma'|$ of total degree $1-k$, such that*

1. *if $\sigma$ is a face of $\tau$, then*
$$a'(\tau,\sigma')|_{|\sigma\setminus\sigma'|} = a'(\sigma,\sigma')\;;$$
2. *the form $a'(\sigma,\sigma)$ on $|\sigma_l|$ vanishes;*
3. *if $\sigma \in S_l$, then*
$$a'(\sigma,\sigma_{0,\dots,l-1})|_{|\sigma_l|} = a(\sigma)|_{|\sigma_l|}\;;$$
4. *we have $0 = \left(\nabla^V + a'(\sigma,\emptyset)\right)^2$, and for $\sigma' \in S_k$ with $k \ge 1$,*
$$\begin{aligned}0 = (-1)^k \left[\nabla^V, a'(\sigma,\sigma')\right] &+ \sum_{j=0}^{k-1}(-1)^j\, a'(\sigma,\sigma'_{0,\dots,\hat{j},\dots,k-1}) + (-1)^k\, a(\sigma')\\ &+ \sum_{j=0}^{k-1}(-1)^{k(j-1)}\, a(\sigma,\sigma_{0\dots j})\, a(\sigma,\sigma_{j,\dots,k-1}) + a'(\sigma,\sigma')\, a'(\sigma,\emptyset)\;;\end{aligned}$$
5. *the component $a'_{\alpha\beta}(\sigma,\sigma')$ in $\Omega^*\big(|\sigma\setminus\sigma'|\,;\mathrm{Hom}(V^\beta,V^\alpha)\big)$ vanishes unless $\beta \prec_\sigma \alpha$.*

*1.48. Remark.* If $a'(\sigma,\sigma')$ was a function, Proposition 1.47 (1) would mean that $a'(\sigma,\sigma')$ is the restriction to $|\sigma\setminus\sigma'|$ of a piece-wise smooth function $a'(\sigma')$ on the union of all $|\sigma\setminus\sigma'|$ for all $\sigma$ having $\sigma'$ as a face. Note however that this union is not a differentiable manifold (not even a manifold with corners). Moreover, we are dealing with differential forms. Suppose that $\sigma$ is a codimension-1-face of two simplices $\tau$ and of $\tau'$. Then in local coordinates, those coefficients of $a'(\tau,\sigma')$ and $a'(\tau',\sigma')$ that contain the direction normal to $|\sigma|$ may jump at $|\sigma|$, so that $a'(\sigma')$ would not even be well-defined along $|\sigma|$. We will nevertheless write $a'(\sigma')$ instead of $a'(\sigma,\sigma')$ as long as no confusion about $\sigma$ arises.

If we introduce the notation

$$a(\sigma_{0,\dots,k-1,\,\cdot\,}) = a'(\sigma_{0,\dots,k-1}) \qquad\text{and}\qquad a(\sigma_{\,\cdot\,}) = a'(\emptyset)$$

then (4) above without the term involving $\nabla^{\overline{V}}$ looks precisely as Corollary 1.43 (2), where the index "$k$" has been replaced by "$\cdot$". The placement of "$\cdot$" to the right of $k-1$ in the index above is compatible with the fact that we may only evaluate $a'(\sigma_{0,\dots,k-1})$ on the set of points "behind" $\sigma'_{k-1}$ in $\sigma$. In fact, these two analogies with the coefficients $a(\sigma)$ of Corollary 1.43 were the motivation for our construction of $a'(\sigma')$.

While we could regard the $a(\sigma)$ of Corollary 1.43 as coefficients of a purely combinatorial superconnection, (or as a "system of higher homotopies" as in. [I2], chapter 3), we now have a mixed superconnection which is partly combinatorial and partly differential. Although the forms $a'(\sigma')$ do not pull back from endomorphisms defined on $B$, we can pull them back to $B$ using a section of $\pi_1\colon \overline{B} \to B$ constructed from a partition of unity, to obtain a differential superconnection on $V \to B$.

*Proof* of Proposition 1.47. We define the forms $a'(\sigma,\sigma')$ on $|\sigma \setminus \sigma'|$ for $\sigma \in S_l$ by induction over $l$. For all $\sigma \in S_0$, the endomorphisms $a'(\sigma,\sigma) = 0$ and $a'(\sigma,\emptyset) = a(\sigma)$ are determined by (2) and (3), and (4) holds on $\sigma$ by Corollary 1.43 (2).

Now fix $\sigma \in S_l$ for $l > 0$ and suppose that all $a(\tau,\sigma')$ have been constructed on all faces $\tau$ of $\sigma$, such that (1)–(5) are satisfied. If $\sigma' = \sigma_{i_0,\dots,i_{k-1}} \in S_{k-1}$ is a face of $\sigma$ with $i_{k-1} > k-1$, then the simplex $\tau = \sigma_{i_0,\dots,i_{k-1},i_{k-1}+1,\dots,l}$ is a proper face of $\sigma$, and $|\sigma \setminus \sigma'| = |\tau \setminus \sigma'|$. In particular, the form $a'(\sigma,\sigma') = a'(\tau,\sigma')$ is already determined by (1) and induction, and (4) holds for $a'(\sigma,\sigma')$.

Now fix $\sigma' = \sigma_{0,\dots,k-1}$ for $1 \leq k \leq l+1$, then the restriction of $a'(\sigma,\sigma')$ to each of the sets $|\sigma_{k-1,\dots,\widehat{j},\dots,l}|$ with $j \geq k$ is already determined by the induction hypothesis applied to $\tau = \sigma_{0,\dots,\widehat{j},\dots,l}$ and $\sigma'$. To summarise, the restriction of $a'(\sigma,\sigma')$ to the set

$$\partial\,|\sigma_{k-1,\dots,l}| \setminus \mathring{\sigma}_{k,\dots,l} \tag{1.49}$$

is already defined. Similarly, the restriction of $a'(\sigma,\emptyset)$ to $\partial\,|\sigma|$ is already known.

We will construct $a'(\sigma,\sigma')$ by induction on $l-k$. First, we have $a'(\sigma,\sigma) = 0$ by (2), and by (3), (4) reduces to the trivial equation

$$a'(\sigma,\sigma_{0,\dots,l-1}) - a(\sigma) = a(\sigma) - a(\sigma) = 0\,.$$

Moreover, (5) follows from Corollary 1.43 (3).

Now fix $\sigma' = \sigma_{0,\dots,k-1}$ for $1 \leq k \leq l$, and suppose that $a'(\sigma,\sigma_{0\dots k})$ has already been constructed on $|\sigma_{k\dots l}| \subset |\sigma_{k-1,\dots,l}| = |\sigma \setminus \sigma'|$. Writing $a'(\,\cdot\,)$ instead of $a'(\sigma,\,\cdot\,)$ as in Remark 1.48, we define

$$a''(\sigma') \in \Omega^*\bigl(|\sigma_{k\dots l}|\,;\operatorname{End} V\bigr)$$

by rewriting (4) as

$$\begin{aligned}(-1)^{k+1}\,a''(\sigma') = \sum_{j=0}^{k-1}(-1)^j\,a'(\sigma_{0\dots\widehat{j}\dots k}) + (-1)^{k+1}\,a(\sigma_{0\dots k})&\\ + \sum_{j=1}^{k}(-1)^{(k+1)(j-1)}\,a(\sigma_{0\dots j})\,a'(\sigma_{j\dots k})&\\ + (-1)^{k+1}\bigl(\nabla^V + a(\sigma_0)\bigr)\,a'(\sigma_{0\dots k}) + a'(\sigma_{0\dots k})\,\bigl(\nabla^V + a'(\emptyset)\bigr)&\end{aligned} \tag{1.50}$$

Then $a''(\sigma')$ clearly satisfies (5). For $j > k$, the induction hypothesis applied to the face $\sigma_{0\dots\widehat{j}\dots l}$ shows that $a'(\sigma')$ satisfies the same equation on $|\sigma_{k\dots\widehat{j}\dots l}|$, thus $a''(\sigma')$ coincides with $a'(\sigma')$ on

$$\partial\,|\sigma_{k\dots l}| \setminus \mathring{\sigma}_{k+1,\dots,l}\ , \tag{1.51}$$

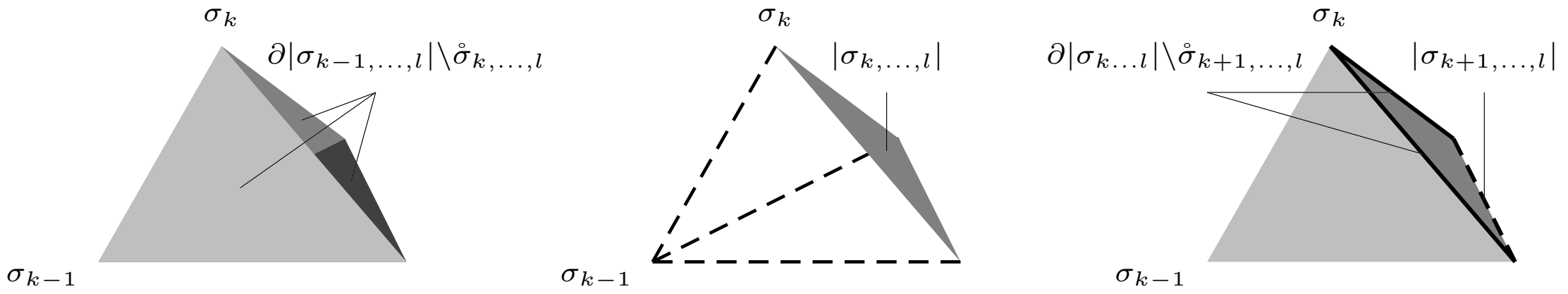

Figure 1.52. The domains of definition of $a'(\sigma')$ and $a''(\sigma')$.

see Figure 1.52.

Because the analogue of (1.50) holds for $a'(\sigma_{0\dots k})$ on $|\sigma_{k+1,\dots,l}|$, and because the superconnections $\nabla^{\overline{V}} + a(\sigma_0)$ and $\nabla^{\overline{V}} + a'(\emptyset)$ are flat, we can rewrite (1.50) for $a''(\sigma')$ as

$$\begin{aligned}
&(-1)^{k+1}a''(\sigma')\\
&\quad=\sum_{j=0}^{k-1}(-1)^j a'(\sigma_{0\dots\widehat{j}\dots k}) + (-1)^{k+1}\,a(\sigma_{0\dots k}) + \sum_{j=1}^{k}(-1)^{(k+1)(j-1)}\,a(\sigma_{0\dots j})\,a'(\sigma_{j\dots k})\\
&\qquad-\left(\nabla^V + a(\sigma_0)\right)\Biggl(\sum_{j=0}^{k}(-1)^j\,a'(\sigma_{0\dots\widehat{j}\dots k+1}) + (-1)^k\,a(\sigma_{0\dots k+1})\\
&\qquad\qquad+\sum_{j=1}^{k+1}(-1)^{k(j-1)}\,a(\sigma_{0\dots j})\,a'(\sigma_{j\dots k+1})\Biggr)\\
&\qquad+(-1)^k\Biggl(\sum_{j=0}^{k}(-1)^j\,a'(\sigma_{0\dots\widehat{j}\dots k+1}) + (-1)^k\,a(\sigma_{0\dots k+1})\\
&\qquad\qquad+\sum_{j=1}^{k+1}(-1)^{k(j-1)}\,a(\sigma_{0\dots j})\,a'(\sigma_{j\dots k+1})\Biggr)\left(\nabla^V + a'(\emptyset)\right)\,.
\end{aligned}$$

We now apply (4) inductively as well as Corollary 1.43 (2) to the last two terms on the right hand side, and find

$$\begin{aligned}
&(-1)^{k+1}a'(i_0\dots i_{k-1},\,\cdot\,)\\
&\quad=\sum_{j=0}^{k-1}(-1)^j a'(\sigma_{0\dots\widehat{j}\dots k}) + (-1)^{k+1}\,a(\sigma_{0\dots k}) + \sum_{j=1}^{k}(-1)^{(k+1)(j-1)}\,a(\sigma_{0\dots j})\,a'(\sigma_{j\dots k})\\
&\qquad-\sum_{j=0}^{k}(-1)^{j+k}\Biggl(\sum_{j'=0}^{j-1}(-1)^{j'}\,a'(\sigma_{0\dots\widehat{j'}\dots\widehat{j}\dots k+1})\\
&\qquad\qquad-\sum_{j'=j+1}^{k+1}(-1)^{j'}\,a'(\sigma_{0\dots\widehat{j}\dots\widehat{j'}\dots k+1}) + (-1)^{k+1}\,a(\sigma_{0\dots\widehat{j}\dots k+1})\\
&\qquad\qquad+\sum_{j'=1}^{j-1}(-1)^{(k+1)(j'-1)}\,a(\sigma_{0\dots j'})\,a'(\sigma_{j'\dots\widehat{j}\dots k+1})\\
&\qquad\qquad+\sum_{j'=j+1}^{k+1}(-1)^{(k+1)j'}\,a(\sigma_{0\dots\widehat{j}\dots j'})\,a'(\sigma_{j'\dots k+1})\Biggr)
\end{aligned}$$

$$
\begin{aligned}
&- a(\sigma_0)\, a'(\sigma_{1\dots k+1}) + a(\sigma_{0\dots k+1})\, a'(\emptyset) \\
&- \sum_{j'=0}^{k+1} (-1)^{j'}\, a(\sigma_{0\dots\widehat{j'}\dots k+1}) - \sum_{j'=1}^{k+1} (-1)^{(k+1)(j'-1)}\, a(\sigma_{0\dots j'})\, a(\sigma_{j'\dots k+1}) \\
&- \sum_{j=1}^{k+1} (-1)^{(k+1)(j-1)} \Bigg( \sum_{j'=0}^{j} (-1)^{j'}\, a(\sigma_{0\dots\widehat{j'}\dots j}) \\
&\qquad\qquad + \sum_{j'=1}^{j} (-1)^{j(j'-1)}\, a(\sigma_{0\dots j'})\, a(\sigma_{j'\dots j}) \Bigg)\, a'(\sigma_{j\dots k+1}) \\
&- \sum_{j=1}^{k+1} (-1)^{kj}\, a(\sigma_{0\dots j}) \Bigg( \sum_{j'=j}^{k+1} (-1)^{j'-j}\, a'_{(}\sigma_{j,\dots\widehat{j'}\dots k+1}) + (-1)^{k-j}\, a(\sigma_{j\dots k+1}) \\
&\qquad\qquad + \sum_{j'=j}^{k+1} (-1)^{(k-j)(j'-j-1)}\, a(\sigma_{j\dots j'})\, a'(\sigma_{j'\dots k+1}) \Bigg) \\
= {}& (-1)^{k+1}\, a'(\sigma_{0\dots k-1}) = (-1)^{k+1}\, a'(\sigma')\ .
\end{aligned}
$$

Note that all terms in the intermediate step except $(-1)^{k+1}\, a'(\sigma_{0\dots k-1})$ cancel by their signs. In particular, we can now extend $a'(\sigma')$ to $\partial\,|\sigma_{0,\dots,k-1}|$ by $a''(\sigma')$. We choose an arbitrary smooth continuation $a'(\sigma,\sigma')$ of $a'(\sigma')|_{\partial|\sigma_{0,\dots,k-1}|}$ to $|\sigma_{0,\dots,k-1}|$ that is of total degree $1-k$, cf. [I2], Lemma 2.2.3. By Remark 1.16 (2), the relation $\prec_\sigma$ refines on the corresponding relations for the faces of $\sigma$, so we can choose $a'(\sigma,\sigma')$ such that (5) still holds.

In the last step, we define $a''(\emptyset)$ on $|\sigma|$ by

$$
\nabla^V + a''(\emptyset) = \bigl(\mathrm{id} + a'(\sigma_0)\bigr)^{-1} \Bigl(\nabla^V + a(\sigma_0)\Bigr) \bigl(\mathrm{id} + a'(\sigma_0)\bigr)\ . \tag{1.53}
$$

As above, we automatically get $a''(\emptyset) = a'(\emptyset)$ on $\partial\,|\sigma| \setminus \mathring{\sigma}_{1\dots l}$, and a similar calculation as above shows $a''(\emptyset) = a'(\emptyset)$ on the whole boundary of $|\sigma|$, so $a'(\emptyset)$ naturally extends to $|\sigma|$. Moreover, since $a(\sigma_0)$ is a parallel endomorphism with $a(\sigma_0)^2 = 0$, the superconnection $\nabla^V + a(\sigma_0)$ on $|\sigma|$ is flat, and thus the superconnection $\nabla^{\overline{V}} + a'(\emptyset)$ on $|\sigma|$ is also flat, which means that (4) is satisfied also for $a'(\sigma,\emptyset)$. $\square$

Given the coefficients $a'(\sigma,\sigma')$ of a "mixed" superconnection, we now have to construct the coefficients $I'(\sigma,\sigma')$ of a "mixed" chain homomorphism in such a way that the analogue of Corollary 1.43 (1) holds. Let $I(\sigma)$ be the integration maps of Corollary 1.43.

**1.54. Proposition.** *For each $0 \le k \le l+1$ and each (possibly empty) face $\sigma' \in S_{k-1}$ of $\sigma \in S_l$, there exists a map*

$$
I'(\sigma,\sigma')\colon\ \Omega^*(M;F) \longrightarrow \Omega^*\bigl(|\sigma \setminus \sigma'|\,;V\bigr)
$$

*of total degree $-k$, such that*

(1) *if $\sigma$ is a face of $\tau$, then*

$$
I'(\sigma,\sigma')\,\omega = \bigl(I'(\tau,\sigma')\,\omega\bigr)|_{|\sigma\setminus\sigma'|} \qquad \text{for all } \omega \in \Omega^*(M;F);
$$

(2) *the map $I'(\sigma,\sigma)$ vanishes;*

(3) *if $\sigma \in S_l$, then*

$$
\bigl(I'(\sigma,\sigma_{0,\dots,l-1})\,\omega\bigr)|_{|\sigma_l|} = \bigl(I(\sigma)\,\omega\bigr)|_{|\sigma_l|} \qquad \text{for all } \omega \in \Omega^*(M;F);
$$

(4) *we have* $I'(\sigma,\emptyset)\circ\nabla^F=\nabla^V\circ I'(\sigma,\emptyset)$, *and for* $\sigma'\in S_k$ *with* $k\geq 1$,

$$\begin{aligned}I'(\sigma,\sigma')\circ\nabla^F&=(-1)^k\nabla^V\circ I'(\sigma,\sigma')\\&\quad+\sum_{j=0}^{k-1}(-1)^j\,I'(\sigma,\sigma'_{0,\dots,\widehat{j},\dots,k-1})+(-1)^k\,I(\sigma')\\&\quad+\sum_{j=0}^{k-1}(-1)^{k(j-1)}\,a(\sigma_{0,\dots,j})\,I'(\sigma,\sigma_{j,\dots,k-1})+a'(\sigma,\sigma_{0,\dots,k-1})\,I'(\sigma,\emptyset)\;;\end{aligned}$$

(5) *for each face* $\sigma''\in S_j$ *of* $\sigma$ *and* $\alpha$, $\beta\in P(\sigma)$, *there exist differential forms*

$$b_{\alpha\beta}(\sigma,\sigma',\sigma'')\in\Omega^{\operatorname{ind}_\beta-\operatorname{ind}_\alpha+j-k}(|\sigma\setminus\sigma'|)\;,$$

*with* $b_{\alpha\beta}(\sigma,\sigma',\sigma'')=0$ *unless* $\beta\prec_\sigma\alpha$ *or* $\beta=\alpha$, *such that*

$$I'_\alpha(\sigma,\sigma')=\sum_\beta\sum_{\sigma''}b_{\alpha\beta}(\sigma,\sigma',\sigma'')\,I_\beta(\sigma'')$$

*for the component* $I'_\alpha(\sigma,\sigma')$ *of* $I(\sigma,\sigma')$ *with values in* $\Omega^*\big(\overline{B}(\sigma,\sigma');\overline{V}_\alpha\big)$;

(6) *if* $\omega\in\Omega^*(M,F)$ *has* $\operatorname{supp}(\omega)\subset h^{-1}(h_\alpha-\varepsilon^2,\infty)$, *then on any* $|\sigma\setminus\sigma'|$, *we have*

$$I'_\alpha(\sigma,\sigma')\,\omega=\begin{cases}\displaystyle\int_{D_{2\varepsilon}T^uX|_{C_\alpha(|\sigma\setminus\sigma'|)}/|\sigma\setminus\sigma'|}\omega & \text{if } \sigma'=\emptyset\text{, and}\\ 0 & \text{otherwise.}\end{cases}$$

Note that the analogue of Remark 1.48 applies to the maps $I'(\sigma,\sigma')$ as well: We may regard $I'(\sigma,\sigma')$ as a restriction of a map $I'(\sigma')$ from $\Omega^*(M|_{U(\sigma')};F)$ to the $\overline{V}$-valued forms on the union of all $\overline{B}(\sigma,\sigma')$. Moreover, if we introduce the notation $I(\sigma'_{0,\dots,k-1,\,\cdot\,})=I'(\sigma')$ for $\sigma'\in S_{k-1}$, then (4) above is related to Corollary 1.43 (1) just as Proposition 1.47 (4) is related to Corollary 1.43 (2). Finally, properties (5) and (6) will be crucial when we prove Theorem 0.1 in Section 4.

*Proof.* We basically repeat the proof of Proposition 1.47. For all $\sigma\in S_0$, we set $I'(\sigma,\sigma)=0$ and $I'(\sigma,\emptyset)=I(\sigma)$. By Corollary 1.43, this agrees with (2)–(6) above.

Fix $\sigma\in S_l$ and assume that all $I'(\tau,\sigma')$ have already been defined for all faces $\tau$ of $\sigma$ and $\sigma'$ of $\tau$. Then by the same argument as in the previous proof, for any face $\sigma'=\sigma_{i_0,\dots,i_{k-1}}$ of $\sigma$ with $i_{k-1}>k-1$, the map $I'(\sigma,\sigma')$ is already inductively defined by (1), and (4)–(6) hold for $I'(\sigma,\sigma')$ on $|\sigma\setminus\sigma'|$. Note also that $I'(\sigma,\sigma)=0$ by (2), which is compatible with (3)–(6).

Now fix $\sigma'=\sigma_{0,\dots,k-1}$ for $0<k\leq l$ and assume that $I'(\sigma,\sigma_{0,\dots,k})$ has already been defined. Then $I'(\sigma,\sigma')$ is already known on $\partial\,|\sigma_{k-1\dots l}|\setminus\mathring{\sigma}_{k\dots l}$, cf. (1.49). Writing $I'(\,\cdot\,)$ shorthand for $I'(\sigma,\,\cdot\,)$, we define $I''(\sigma')\colon\Omega^*(M;F)\to\Omega^*(|\sigma_{k\dots l}|\,;V)$ by

$$\begin{aligned}(-1)^{k+1}\,I''(\sigma')&=\sum_{j=0}^{k-1}(-1)^j\,I'(\sigma_{0\dots\widehat{j}\dots k})+(-1)^{k+1}\,I(\sigma_{0\dots k})\\&\quad+\sum_{j=1}^{k}(-1)^{(k+1)(j-1)}\,a(\sigma_{0\dots j})\,I'(\sigma_{j\dots k})\\&\quad+(-1)^{k+1}\big(\nabla^V+a(\sigma_0)\big)\,I'(\sigma_{0\dots k})+I'(\sigma_{0\dots k})\,\nabla^F\;,\end{aligned}\tag{1.55}$$

cf. (1.50).

Note that since (5) holds for $I'(\sigma_{0\dots k})$ by the construction in the previous step, we have

$$(1.56)\quad \begin{aligned}(-1)^{k+1}\,\nabla^V\,I'_\alpha(\sigma_{0\dots k}) + I'_\alpha(\sigma_{0\dots k})\,\nabla^F\\ =\sum_\beta \sum_{\substack{\sigma''\in S\\ \sigma''\text{ is a face of }\sigma}} \biggl((-1)^{k+1}\,db_{\alpha\beta}(\sigma,\sigma_{0\dots k},\sigma'')\,I_\beta(\sigma'')\\ + b_{\alpha\beta}(\sigma,\sigma_{0\dots k},\sigma'')\Bigl((-1)^{j+\operatorname{ind}_\beta-\operatorname{ind}_\alpha}\,\nabla^V\,I_\beta(\sigma'')+I_\beta(\sigma'')\,\nabla^F\Bigr)\biggr)\;.\end{aligned}$$

By Corollary 1.43, (1.55) and (1.56), we see that $I''(\sigma')$ satisfies (5), where the degrees of the coefficient forms can be calculated using the fact that $I''(\sigma')$ is of total degree $-k$.

As in the proof of Proposition 1.47, $I'(\sigma')$ and $I''(\sigma')$ coincide on $\partial\,|\sigma_{k\dots l}|\setminus\mathring{\sigma}_{k+1\dots l}$, cf. (1.51). In order to check that $I'(\sigma')=I''(\sigma')$ on $|\sigma_{k+1\dots l}|$, we proceed by computations similar to those in the previous proof. In fact, comparing (1.50) with its analogue (1.55), we find that we can re-use the computations in the previous proof after making the following changes.

(1) First, replace each occurrence of $\nabla^V + a'(\emptyset)$ on the right hand side of a summand by $\nabla^F$.
(2) Then replace the rightmost occurrence of $a'(\,\cdot\,)$ or $a(\,\cdot\,)$ by $I'(\,\cdot\,)$ or $I(\,\cdot\,)$, respectively.

Then all calculations remain valid and show that we can use $I''(\sigma')$ to extend $I'(\sigma')$ to $\partial\,|\sigma_{k-1\dots l}|$. We now choose smooth continuations of the coefficients $b_{\alpha\beta}(\sigma,\sigma',\sigma'')$ in (5) to $\partial\,|\sigma_{k-1\dots l}|$, thus extending $I'(\sigma')$ to $|\sigma\setminus\sigma'|$.

In order to determine $I'(\emptyset)$, we may rewrite (4) as

$$(1.57)\qquad I''(\emptyset)=\bigl(\operatorname{id}+a'(\sigma_0)\bigr)^{-1}\Bigl(\bigl(\nabla^V+a(\sigma_0)\bigr)\circ I'(\sigma_0)+I'(\sigma_0)\circ\nabla^F+I(\sigma_0)\Bigr)\;.$$

As before, we can show that $I''(\emptyset)=I'(\emptyset)$ on $\partial\,|\sigma|$, so we can set $I'(\emptyset)=I''(\emptyset)$ on $|\sigma|$. By the argument following (1.56), property (5) then also holds for $I'(\emptyset)$. Moreover, using Corollary 1.43, (1.53) and (4), we finally find that

$$\begin{aligned}I'(\emptyset)\circ\nabla^F &= \bigl(\operatorname{id}+a'(\sigma_0)\bigr)^{-1}\Bigl(\bigl(\nabla^V+a(\sigma_0)\bigr)\circ I'(\sigma_0)+I(\sigma_0)\Bigr)\circ\nabla^F\\ &= \bigl(\operatorname{id}+a'(\sigma_0)\bigr)^{-1}\bigl(\nabla^V+a(\sigma_0)\bigr)\\ &\qquad\Bigl(\bigl(\operatorname{id}+a'(\sigma_0)\bigr)\,I'(\emptyset)-I(\sigma_0)-\bigl(\nabla^V+a(\sigma_0)\bigr)\,I'(\sigma_0)+I(\sigma_0)\Bigr)\\ &= \bigl(\nabla^V+a'(\emptyset)\bigr)\,I'(\emptyset)\;.\end{aligned}$$

We have now proved Proposition 1.54 (1)–(5). Claim (6) follows from (5) by Corollary 1.43 (4), where we use that for $\sigma'\in S_{k-1}$ and $\sigma''\in S_j$, $b_{\alpha\alpha}(\sigma,\sigma',\sigma'')=0$ unless $j\ge k$. □

**1.f. A smoothing construction.** We note that by Proposition 1.47 (4) and Proposition 1.54 (4), we obtain a superconnection $\nabla^V + a'(\sigma,\emptyset)$ on $\Omega^*(|\sigma|\,;V)$ and a cochain-map from $\Omega^*(M;F)$ to $\Omega^*(|\sigma|\,;V)$. By Remark 1.48, these objects do not patch together to yield a smooth flat superconnection and a smooth cochain-map over all of $B$. There are several possible ways out.

For example, one can define an algebra of "$S$-piece-wise smooth" differential forms, given as families $\alpha=(\alpha_\sigma)_{\sigma\in S}$ with $\alpha_\sigma\in\Omega^*(|\sigma|)$, such that $\alpha_\sigma|_\tau=\alpha_\tau$ whenever $\tau$ is a face of $\sigma$. This means that near a simplex $\sigma\in S$, the components of $\alpha$ tangential to $|\sigma|$ are continuous, while the remaining components are only bounded with all their derivatives. One can show that it is possible to do analysis in such a setup. E.g., the equation $(\nabla^V + a'(\,\cdot\,,\emptyset))^2 = 0$ is still meaningful. It is possible to define an $S$-piece-wise smooth torsion form on $B$.

We will instead give a smoothing construction that leads to a smooth superconnection $A'$ on $V$. Choose a partition of unity $\{\,\varphi_\sigma\colon M \to [0,1] | \sigma \in S_0\,\}$ subordinate to the open covering $\mathcal{U}(S)$ of (1.36), i.e., with $\operatorname{supp}(\varphi_\sigma) \subset U_\sigma$ for all $\sigma \in S_0$. If $\sigma \in S_k$, then all these functions vanish on $|\sigma|$ except $\varphi_{\sigma_0}, \dots, \varphi_{\sigma_k}$, and of course $\varphi_{\sigma_0} + \dots + \varphi_{\sigma_k} = 1$ on $|\sigma|$. If we interpret $\varphi_{\sigma_0}, \dots, \varphi_{\sigma_k}$ as barycentric coordinates on $\Delta_k$, we obtain a smooth map $\varphi\colon B \to B$ by setting

$$\bar\varphi(b) = \sigma\bigl(\varphi_{\sigma_0}(b), \varphi_{\sigma_1}(b), \dots\bigr) \in |\sigma| \qquad \text{for } b \in \mathring{\sigma}. \tag{1.58}$$

Note that $\bar\varphi$ maps a small neighbourhood

$$U'(\sigma) = \bigl\{\, b \in B \bigm| b \notin \operatorname{supp}(\varphi_{\sigma'}) \text{ unless } \sigma' \text{ is a vertex of } \sigma \,\bigr\}$$

of each simplex $\sigma \in S$ to $|\sigma|$.

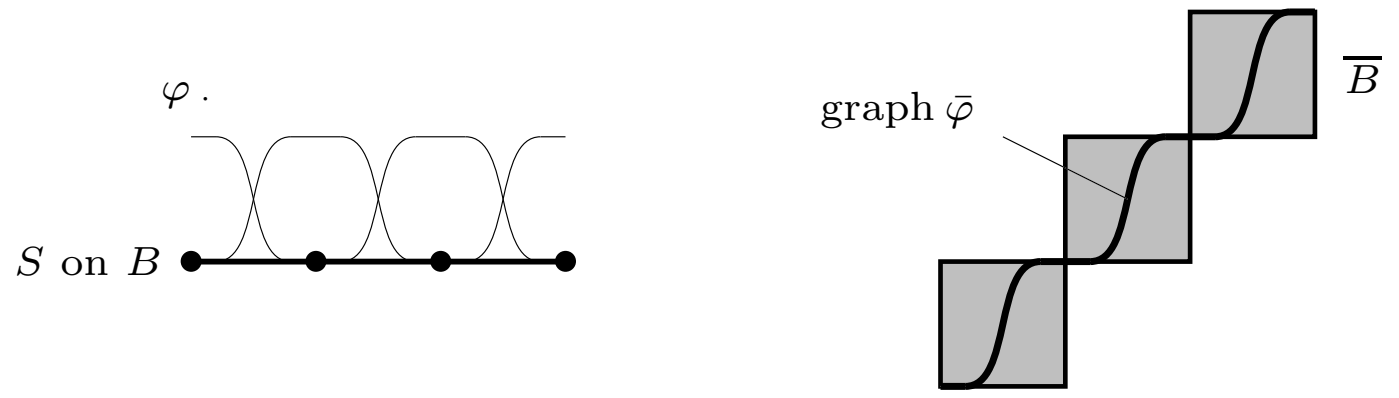

Figure 1.59. A partition of unity and the corresponding section of $\pi_1^*$.

**1.60. Definition.** Define a superconnection $A' = \nabla^V + a'$ on $V \to B$ with $a' \in \Omega^*(B, \operatorname{End}_+^* V)$, and a map $I\colon \Omega^*(M;F) \to \Omega^*(B;V)$, such that restricted to $|\sigma|$,

$$A' = \nabla^V + \bar\varphi^* a'(\sigma, \emptyset)\ ,$$
$$\text{and} \qquad I = \sum_{\sigma'' \text{ is a face of } \sigma} \bar\varphi^* b(\sigma, \emptyset, \sigma'') \cdot I(\sigma'')\ ,$$

where the coefficients $b(\sigma, \emptyset, \sigma'')$ are given in Proposition 1.54 (5).

Let $h^V$ denote the endomorphism of $V \to B$ induced by $h|_C$, then we will call $(V, A', h^V)$ a *family Thom-Smale complex* for $p\colon M \to B$ and $F \to M$. The complexes $\Omega^*(M;F)$ and $\Omega^*(B;V)$ are filtered by $F^*\Omega^*(M;F)$ and $F^*\Omega^*(B;V)$ according to the horizontal degree, see Definition 1.4. Recall that $(E_k^{*,*}(M;F), d_k)$ and $(E_k^{*,*}(B;V), d_k)$ denote the accociated spectral sequences, and that $E_*^{*,*}(M;F)$ is the Leray spectral sequence in de Rham theory with coefficients in $F \to M$.

**1.61. Theorem.** *The superconnection $A' = \nabla^V + a'$ is flat and of total degree one, and $I$ is a cochain map of total degree zero. The endomorphism $a'$ decomposes locally over $|\sigma|$ into components $a'_{\alpha\beta} \in \Omega^*(B, \operatorname{Hom}(V_\beta, V_\alpha))$ with*

$$a'_{\alpha\beta} = 0 \qquad \textit{unless } \beta \prec_\sigma \alpha. \tag{1}$$

*Locally on $|\sigma|$, the map $I$ is an $\Omega^*(B;\operatorname{End} V)$-linear combination of the maps $I_\alpha(\sigma')$ of (1.42), with $\sigma'$ a face of $\sigma$. If $\omega \in \Omega^*(M;F)$ has support in $h^{-1}(h_\alpha - \varepsilon^2, \infty)$ locally on $\overline{U(\sigma)} \subset B$, then the component $I_\alpha\omega$ of $I\,\omega$ in $\Omega^*(\overline{U(\sigma)};V_\alpha)$ is given by*

$$I_\alpha\omega = \int_{D_{2\varepsilon}T^uX|_{C_\alpha(\sigma)}/\overline{U(\sigma)}} \omega \quad \in \Omega^*\bigl(\overline{U(\sigma)}; V_\alpha\bigr)\ . \tag{2}$$

*The map $I$ commutes with the action of $\Omega^*(B)$, maps $F^p\Omega^*(M;F)$ to $F^p\Omega^*(B;V)$ for all integers $p \geq 0$, and induces cochain homomorphisms*

$$I_k^{p,q}\colon \big(E_k^{*,*}(M;F), d_k\big) \longrightarrow \big(E_k^{*,*}(B;V), d_k\big) \tag{3}$$

*for $k$, $p$, $q \geq 0$ that are isomorphisms for $k \geq 1$. In particular, $I\colon (\Omega^*(M;F), \nabla^F) \to (\Omega^*(B;V), A')$ is a quasi-isomorphism.*

*Proof.* We know that near a simplex $\sigma \in S$, the components of $a'(\,\cdot\,,\emptyset)$ and $I'(\,\cdot\,,\emptyset)$ that are tangential to $|\sigma|$ are continuous, together with all their higher derivatives in tangential directions. On the other hand, since $\bar\varphi$ is smooth and maps an open neighbourhood of $|\sigma|$ to $|\sigma|$, pulling back by $\bar\varphi$ kills all components and all higher derivatives of $a'(\,\cdot\,,\emptyset)$ and $I'(\,\cdot\,,\emptyset)$ that are not purely tangential to $|\sigma|$. This implies that $A'$ and $I$ are smooth.

It is immediate from Proposition 1.47 that $A'$ is a flat superconnection of total degree 1, and that (1) holds.

Furthermore, (1.57) and Definition 1.60 imply that for all $\sigma$, there exists a map

$$I'_\sigma = \sum_{\sigma'' \text{ is a face of } \sigma} \bar\varphi^* b(\sigma, \sigma_0, \sigma'')\, I(\sigma'')\;,$$

such that

$$I = \big(\mathrm{id} + \bar\varphi^* a'(\sigma_0)\big)^{-1} \Big(\big(\nabla^V + a(\sigma_0)\big) \circ I'_\sigma + I'_\sigma \circ \nabla^F + I(\sigma_0)\Big)\;, \tag{1.62}$$

$$\text{and} \qquad I \circ \nabla^F = \big(\nabla^V + \bar\varphi^* a'(\emptyset)\big) \circ I = A' \circ I\;.$$

From the construction and Proposition 1.54 (5) and (6), we conclude that $I|_{|\sigma|}$ is indeed an $\Omega^*(B;\operatorname{End} V)$-linear combination of the maps $I_\alpha(\sigma')$ of (1.42), and that (2) holds. Since the maps $I_\alpha(\sigma')$ are defined by integration over the fibres of certain smooth families of submanifolds with conical singularities, these maps clearly commute with $\Omega^*(B)$, and so does $I$.

As an $\Omega^*(B)$-homomorphism, $I$ maps $F^p\Omega^*(M;F)$ to $F^p\Omega^*(B;V)$ for all $p$, so it induces cochain homomorphisms of spectral sequences as claimed. It suffices to show that $I_1^{*,*}$ is an isomorphism, then all $I_k^{*,*}$ with $k \geq 1$ are also isomorphisms, and $I$ is thus a quasi-isomorphism. By (1.62), for each $\sigma \in S_k$, we have a diagramme of cochain maps

$$\begin{array}{ccc} \big(\Omega^*(M;F), \nabla^F\big) & \xrightarrow{\;I\;} & \big(\Omega^*(|\sigma|;V), A'\big) \\ & \searrow^{I(\sigma_0)} & \Big\downarrow{\scriptstyle \wr\|}\; {\scriptstyle \mathrm{id} + \bar\varphi^* a'(\sigma,\sigma_0)} \\ & & \big(\Omega^*(|\sigma|;V), \nabla^V + a(\sigma,\sigma_0)\big) \end{array}$$

that commutes up to the cochain homotopy $I'_\sigma$. Note that by classical Morse-theory, the diagonal arrow induces an isomorphism from $H^*(X;F)$ to $H^*(V, a(\sigma_0))$ for each fibre $X$ over $|\sigma|$.

Passing to $E_1^{p,q}$, the map

$$\big(\nabla^V + \bar\varphi^* a(\sigma_j)\big)\, I'_\sigma + I'_\sigma\, \nabla^F\colon\; \Omega^*(M;F) \longrightarrow \Omega^*\big(|\sigma|\,;V\big)$$

induces 0, because in the grading induced by the filtrations in Definition 1.4, it splits into a fibre-wise homotopy and a part that increases the horizontal degree. In particular, the diagramme

$$\begin{array}{ccc} \Omega^p\big(B; H^q(M/B;F)\big) & \xrightarrow{\;I_1^{p,q}\;} & \Omega^p\big(|\sigma|\,; H^q(V, a'_0)\big) \\ & \searrow^{\sim} & \Big\downarrow{\scriptstyle \wr} \\ & & \Omega^p\big(|\sigma|\,; H^q\big(V, a(\sigma,\sigma_0)^{[0]}\big)\big) \end{array}$$

commutes, where $a_0'$ and $a(\sigma,\sigma_0)^{[0]}$ denote the components of horizontal degree 0. It follows that $I_1^{p,q}$ is an isomorphism for all $p$, $q$, which finishes the proof of Theorem 1.61. □

*1.63. Remark.* As we will see in Section 4, the Witten deformation will give another proof of the fact that $I_1^{*,*}$ is a quasi-isomorphism.

*1.64. Remark.* Let $\sigma' \in S_k$ be a face of $\sigma \in S_l$. It is clear from the proof of Proposition 1.47 that we can always choose

$$a'(\sigma,\sigma') \in \bigoplus_{j=0}^{l-k-1} \Omega^j\left(|\sigma \setminus \sigma'|\,;\operatorname{End}_+^{1-j-k} V\right) .$$

However, this means that $A'$ may have components of horizontal degree up to $\dim X + 1$. If we assume that the dimension of the fibres is higher than the dimension of the base—which we are forced to assume when we apply the main theorem of [I1] in order to construct $h$—then $A'$ can have components in any possible horizontal degree.

The construction so far has involved many choices. We had to pick a diffeomorphism $F$ of the bundles $\pi_1^*M$ and $\pi_2^*M$ over a neighbourhood $U$ of $\Delta B$ in $B \times B$, a sufficiently fine smooth ordered simplical complex $S$ on $B$, a continuous vector field $W$ that over each simpex of $S$ is smooth and satisfies the Smale transversality condition, and finally, we had to choose smooth continuations of the coefficients of $a'(\sigma')$ and $I'(\sigma')$ from the boundary of certain simplices to their interiors. Still, the superconnections $A'$ and integration maps $I$ of Definition 1.60 coming from different choices are not completely unrelated.

By Theorem 1.61, any two pairs of complexes and integration maps

$$I_i\colon \left(\Omega^*(M;F),\nabla^F\right) \longrightarrow \left(\Omega^*(B;V_i),A_i'\right)$$

for $i = 0$, 1 are related by a quasi-isomorphism. One possible right inverse for $I_i$ is the inverse of the restriction of $I_i$ to the space of very small eigenvalues under the Witten deformation, cf. Section 4.b. However, we can make the relation of $I_0$, $A_0'$ and $I_1$, $A_1'$ much more concrete.

**1.65. Proposition.** *Assume that $A_0'$ and $A_1'$ are two superconnections and $I_0$, $I_1$ are integration maps for these superconnections, all constructed using the methods of this chapter. Then the same methods allow us to construct a superconnection $\overline{A}'$ on the pull-back $V \times [0,1]$ of $V$ to $B \times [0,1]$, together with an integration map $\overline{I}\colon \Omega^*(M \times [0,1]; F \times [0,1]) \to \Omega^*(B \times [0,1]; V \times [0,1])$, such that Theorem 1.61 holds for $\overline{A}'$ and $\overline{I}$, and $\overline{A}'|_{B\times\{s\}} = A_s'$ and $\overline{I}|_{B\times\{s\}} = I_s$ for $s = 0$, 1.*

*Proof.* One has to check that all choices made so far on $B \times \{0,1\}$ can be extended to the interior of $\overline{B} = B \times [0,1]$, still satisfying all required properties. For the diffeomorphism $F$ of Remark 1.33, this is clear because it is just a perturbation of $\mathrm{id}_{\pi_1^*M}|_{\Delta B}$. Suppose that $\overline{F}$ extends $F$ on a neighbourhood $U_{\Delta\overline{B}}$ of the diagonal in $\overline{B} \times \overline{B}$, then we can construct a smooth ordered simplicial complex $\overline{S}$ as fine as (1.34) that restricts to the given complexes on $B \times \{0\}$ and $B \times \{1\}$.

In the next step, we extend the vertical gradient field $-\nabla^{TX} h$ and the cut-off function $\varphi_M$ over $B \times \{0,1\}$ that were used in (1.20) smoothly to $B \times [0,1]$, thus defining a preliminary vector field $W$ on $M \times [0,1]$. Now the constructions of Proposition 1.27, Proposition 1.47, and Proposition 1.54 are all by induction over the degree of simplices, so we can keep every choice we made over $B \times \{0,1\}$. For Definition 1.60, it suffices to find a partition of unity on $\overline{B}$ that restricts to the given partitions on $B \times \{0\}$ and $B \times \{1\}$. This finishes the proof of our proposition. □

We will need this proposition in [G], to show that a certain "classifying" map from $B$ to the Whitehead space is well-defined up to homotopy. Proposition 1.65 can also be used to compare the torsion forms $T(A_s', g^V, h^V)$ of Section 2 for $s = 0$, 1 coming from superconnections $A_s'$ related as above. Of course, a formula relating $T(A_0', g^V, h^V)$ and $T(A_1', g^V, h^V)$ will ultimately follow from Theorem 0.1.

## 2. Torsion forms for family Thom-Smale complexes

Inspired by the previous section, we define the notion of a family Thom-Smale complex and construct torsion forms for it. These torsion forms generalise the torsion forms for families of finite-dimensional complexes of [BL], Section 2.f.

Before we recall the definitions of odd Chern classes and transgression forms for flat superconnections from [BL] in Sections 2.b, and 2.c, we introduce some sign conventions in Section 2.a. In Section 2.d, we define family Thom-Smale complexes $(V, A', h^V)$ in Definition 2.35. Their torsion forms $T(A', g^V, h^V)$ are constructed in Definition 2.46 using a finite-dimensional version of the Witten deformation. Theorem 2.48 states the most important properties of $T(A', g^V, h^V)$, which are similar to those of $T(A', g^V)$ in [BL]. We show that these forms are rigid under deformation of the underlying flat superconnection in Section 2.e, cf. [BG1]. In Section 2.f, we use $T(A', g^V, h^V)$ to define the higher torsion classes $T(M/B; F, h)$ when $F$ and the fibrewise cohomology $H$ carry parallel metrics, see Definition 2.64.

In Section 2.g, we collect some facts about generalised metrics (cf. [BG1]) that will be needed in Section 4. We recall the construction of analytic torsion forms from [BL] and the Chern normalisation of the odd classes and torsion forms from [BG1] in Sections 2.h and 2.i. We also define a higher analytic torsion class $\mathcal{T}(M/B; F)$ if $F$ and $H$ are unitarily flat in Definition 2.88.

**2.a. Notes on Sign Conventions.** We recall the sign conventions of [BL] and [BG1] regarding superconnections and their adjoints. In particular, we hope to disarm a few well-concealed traps here and in Section 2.g.

For most of this paper, we strictly follow the Koszul sign convention, by which swapping to terms of degrees $p$ and $q$ modulo 2 introduces a sign factor of $(-1)^{pq}$. However, there are a few exceptions in connection with metrics and their generalisations. This is best motivated by the fact that a Hermitian sesquilinear form in the $\mathbb{Z}_2$-graded sense would have to be Hermitian on the even part and skew-Hermitian on the odd degree part of a vector space, whereas for analytic reasons, we need positive definite symmetric forms both in even and in odd degree.

Let $V = V^+ \oplus V^- \to M$ be $\mathbb{Z}_2$-graded complex vector bundle, and let $\overline{V}^*$ be its antidual. We define $\Omega^*(M; V) = \Gamma\big(\Lambda^* M \mathbin{\widehat{\otimes}} V\big)$ to be the space of sections of a graded tensor product. The natural $\mathbb{C}$-sesquilinear pairing $\overline{V}^* \otimes V \to \mathbb{C}$ extends to a pairing

$$\Omega^*\big(M; \overline{V}^*\big) \otimes \Omega^*(M; V) \longrightarrow \Omega^*(M) \tag{2.1}$$

$$\text{with} \qquad \big(\alpha \mathbin{\widehat{\otimes}} f\big) \otimes \big(\beta \mathbin{\widehat{\otimes}} v\big) \longmapsto \big(\alpha \mathbin{\widehat{\otimes}} f, \beta \mathbin{\widehat{\otimes}} v\big) = (-1)^{\left[\frac{\deg \beta}{2}\right]} \alpha\overline{\beta}\, f(v)$$

for all $\alpha$, $\beta \in \Omega^*(M)$, $f \in \Gamma(\overline{V}^*)$ and $v \in \Gamma(V)$. The factor $(-1)^{\left[\frac{\deg \beta}{2}\right]}$ comes from reversing the one-form components of $\overline{\beta}$. Note the absence of a sign factor $(-1)^{\deg V \cdot \deg \beta}$; in particular, the pairing $(\,\cdot\,,\,\cdot\,)$ violates the sign convention.

A metric $g^V = g^{V^+} \oplus g^{V^-}$ on $V$ which respect the $\mathbb{Z}_2$-grading can be regarded as an even complex linear isomorphism $g^V \colon V \to \overline{V}^*$. Together with (2.1), this induces a pairing

$$g^V\colon\ \Omega^*(M, V) \otimes \Omega^*(M; V) \longrightarrow \Omega^*(M) \tag{2.2}$$

$$\text{with} \qquad \big(\alpha \mathbin{\widehat{\otimes}} v\big) \otimes \big(\beta \mathbin{\widehat{\otimes}} w\big) \longmapsto \langle \alpha \mathbin{\widehat{\otimes}} v, \beta \mathbin{\widehat{\otimes}} w\rangle = (-1)^{\left[\frac{\deg \beta}{2}\right]} \alpha\overline{\beta}\, \langle v, w\rangle\ .$$

The fact that $g^V$ is Hermitian is equivalent to the relation

$$\langle \beta \mathbin{\widehat{\otimes}} w, \alpha \mathbin{\widehat{\otimes}} v\rangle = (-1)^{\left[\frac{\deg(\alpha\beta)}{2}\right]} \overline{\langle \alpha \mathbin{\widehat{\otimes}} v, \beta \mathbin{\widehat{\otimes}} w\rangle}\ , \tag{2.3}$$

where the sign factor is once more equivalent to a reversal of the order of the one-form factors in each decomposable summand.

As in [BL] and [BG1], we extend complex conjugation to $\Omega^*(M; \mathbb{C})$ and $\Omega^*(M; \operatorname{End} V)$.

**2.4. Definition.** Define a $\mathbb{C}$-anti-linear operator on $\Omega^*(M;\mathbb{C})$ by

$$\overline{\alpha}^* = (-1)^{\left[\frac{N^B+1}{2}\right]}\overline{\alpha}\ . \tag{2.5}$$

*2.6. Remark.* This operator has the following properties:

$$\begin{aligned} \overline{f}^* &= \overline{f} && \text{for all } f\in\Omega^0(B;\mathbb{C}),\\ \overline{\alpha}^* &= -\overline{\alpha} && \text{for all } \alpha\in\Omega^1(B;\mathbb{C}),\\ \text{and}\qquad \overline{(\alpha\wedge\beta)}^* &= \overline{\beta}^*\wedge\overline{\alpha}^*\ . && \text{for all } \alpha,\ \beta\in\Omega^*(B;\mathbb{C}). \end{aligned}$$

Note that this operator differs from the operator $(-1)^{\left[\frac{\deg(\,\cdot\,)}{2}\right]}\overline{\,\cdot\,}$ encountered in (2.1)–(2.3) by a factor of $(-1)^{\deg(\,\cdot\,)}$.

We let $E^*$ denote the classical (ungraded) adjoint of an endomorphism $E\in\operatorname{End}V$ with respect to $g^V$, and we let $\nabla^{V,*}$ denote the adjoint of a connection $\nabla^V$.

**2.7. Definition.** The adjoint of a differential form $\alpha\mathbin{\widehat{\otimes}}E\in\Omega^*(B;\operatorname{End}V)$ is defined as

$$(\alpha\mathbin{\widehat{\otimes}}E)^* = (-1)^{\deg\alpha\cdot\deg E}\,\overline{\alpha}^*\,E^*\ .$$

The adjoint of a superconnection $A=\nabla^V+a$ with $a\in\Omega^*(M;\operatorname{End}V)$ is given by

$$A^* = \nabla^{V,*}+a^*\ .$$

We call $E$ or $A$ *self-adjoint* (*skew-adjoint*) iff $E^*=E$ or $A^*=A$ ($E^*=-E$ or $A^*=-A$) respectively.

*2.8. Remark.* The definition of $A^*$ is independent of the splitting $A=\nabla^V+a$. For odd elements $E\in\Omega^*(M;\operatorname{End}V)^{\mathrm{odd}}$ and for superconnections $A$ on $V$, Definition 2.7 is compatible with the sign convention of (2.2) in the sense that

$$\begin{aligned} 0 &= \langle E^*v,w\rangle + (-1)^{\deg E\,(\deg v+\deg w)}\,\langle v,Ew\rangle\\ \text{and}\qquad d\langle v,w\rangle &= \langle A^*v,w\rangle + (-1)^{\deg E\,(\deg v+\deg w)}\,\langle v,Aw\rangle \end{aligned}$$

for all $v,\ w\in\Omega^*(M;V)$.

Symmetry of $g^V$ together with Definition 2.7 implies that $(E^*)^*=E$ and $(A^*)^*=A$ for all superconnections $A$ and all $E\in\Omega^*(M;\operatorname{End}V)$.

*2.9. Remark.* Sometimes it is useful to regard $g^V$ as an element of $\Gamma\big(\operatorname{Hom}(V,\overline{V}^*)\big)$. In this setting, the adjoint of a superconnection $A$ is given as

$$A^* = \big(g^V\big)^{-1}\,\overline{A}^*\,g^V\ , \tag{2.10}$$

where $\overline{A}^*$ is the superconnection that is naturally induced on $\overline{V}^*$ by $A$, cf. [BL], Sections 1(c), 1(d). Note that once more $\overline{A}^*$ is compatible with (2.1) in the sense that

$$d(f,v) = (\overline{A}^*f,v) + (-1)^{\deg f+\deg v}\,(f,Av)$$

for all $f\in\Omega^*(M;\overline{V}^*)$ and all $v\in\Omega^*(M,V)$.

**2.b. Odd characteristic classes for flat superconnections.** We recall the construction of Kamber-Tondeur classes in de Rham theory, following [BL] and [BG1] (cf. [KT], [Du]). Let $B$ be a smooth manifold, and let $N^B$ denote the number operator on $\Omega^*(B)$.

We will fix for the rest of the paper a generating function

$$f(z) = z\, e^{z^2} \; , \tag{2.11}$$

although in a few explicitly marked steps, we need a different choice of $f$. Then $f$ is an odd holomorphic function, which is real, i.e., $f(\mathbb{R}) \subset \mathbb{R}$, and there exists $c > 0$ such that

$$\left| \frac{\partial^k f(z)}{\partial z^k} \right| = O\bigl((1 + |\operatorname{Re} z|)^{-l}\bigr) \qquad \text{uniformly on } \bigl\{\, z \in \mathbb{C} \bigm| |\operatorname{Re}(z)| < c \,\bigr\} \tag{2.12}$$

for all pairs of integers $k$, $l \geq 0$.

Now let $A'$ be a flat superconnection on a $\mathbb{Z}_2$-graded vector bundle $V \to B$. Fix a Hermitian metric $g^V$ on $V$, and let $A'' = (A')^*$ be the adjoint superconnection in the sense of Definition 2.7. Define a self-adjoint superconnection $A$ and a skew-adjoint operator $X$ by

$$A = \frac{1}{2}\bigl(A'' + A'\bigr) \qquad \text{and} \qquad X = \frac{1}{2}\bigl(A'' - A'\bigr) \quad \in \quad \Omega^*\bigl(B; \operatorname{End} V\bigr) \; , \tag{2.13}$$

then clearly

$$[A, X] = 0 \quad \in \quad \Omega^*\bigl(B; \operatorname{End} V\bigr)$$

because $A'$ and $A''$ are flat.

**2.14. Definition.** Let $f(A', g^V) \in \Omega^*(B)$ be defined by

$$f\bigl(A', g^V\bigr) = (2\pi i)^{\frac{1-N^B}{2}} \operatorname{str}_V\bigl(f(X)\bigr) \; .$$

Because supertraces vanish on supercommutators, we have $\operatorname{str}_V(X^{2k}) = \frac{1}{2} \operatorname{str}_V([X, X^{2k-1}]) = 0$ for $k \geq 1$, so we have lost no information by specifying $f$ to be an odd function. Since $f$ and $X$ are odd and $\operatorname{str}_V$ vanishes on $\operatorname{End}^{\mathrm{odd}} V$, only the part of $f(X)$ in $\Omega^{\mathrm{odd}}(B; \operatorname{End}^{\mathrm{even}} V)$ contributes, so $f(A', g^V)$ is an odd form.

Because $[A, X] = 0$, we have

$$d^B \operatorname{str}_V\bigl(f(X)\bigr) = \operatorname{str}_V\bigl([A, f(X)]\bigr) = 0 \; ,$$

which means that $f(A', g^V)$ is closed.

Finally, since $f$ is odd and real, it follows that $f(X)$ is skew-adjoint in the sense of Definition 2.7. Using also that $\operatorname{str}_V(f(X))$ is odd, we find

$$\operatorname{str}_V\bigl(f(X)\bigr) = -\operatorname{str}_V\bigl(\overline{f(X)}^*\bigr) = (-1)^{\frac{N^B-1}{2}} \overline{\operatorname{str}_V\bigl(f(X)\bigr)}$$

by (2.5), so the form $f(A', g^V)$ is real.

**2.15. Theorem** ([BL], Theorem 1.8). *The form $f(A', g^V)$ is odd, closed and real. Its cohomology class does not depend on $g^V$.*

Let $g_0^V$ and $g_1^V$ be two metrics on $V$. Then there exists a path $(g_s^V)_{s \in [0,1]}$ of metrics connecting $g_0^V$ and $g_1^V$. Let $A_s''$ denote the adjoint of $A'$ with respect to $g_s^V$, and define $A_s$ and $X_s$ as in (2.13). Let $\bar{g}^V$ be the metric on the pullback of $V$ to $B \times [0,1]$ that restricts to $g_s^V$ on $B \times \{s\}$. We write $\overline{A}'$ for the pull-back of $A'$ to $B \times [0,1]$. Then we have

$$f\bigl(\overline{A}', \bar{g}^V\bigr)_{p,s} = f\bigl(A', g_s^V\bigr)_p + (2\pi i)^{-\frac{N^B}{2}} \operatorname{str}_V\left( \bigl(g_s^V\bigr)^{-1} \frac{\partial g_s^V}{\partial s} f'\bigl(X_s\bigr) \right)_p ds \; .$$

**2.16. Definition.** Let $\tilde f(A',g_0^V,g_1^V)$ be defined by

$$\tilde f\big(A',g_0^V,g_1^V\big)=(2\pi i)^{-\frac{N^B}{2}}\int_0^1 \operatorname{str}_V\bigg(\frac12\,\big(g_s^V\big)^{-1}\,\frac{\partial g_s^V}{\partial s}\,f'\big(X_s\big)\bigg)\,ds\qquad\in\Omega^*(B)/d^B\Omega^*(B)\;.$$

Then by a standard argument, Theorem 2.15 immediately implies [BL], Theorem 1.11.

**2.17. Corollary.** *The class $\tilde f(A',g_0^V,g_1^V)\in\Omega^*(B)/d^B\Omega^*(B)$ is independent of the path $g_s^V$ chosen in its definition. It is even and real, and satisfies*

$$d^B\tilde f\big(A',g_0^V,g_1^V\big)=f\big(A',g_1^V\big)-f\big(A',g_0^V\big)\;.$$

We recall the behaviour of $f(A',g^V)$ under smooth changes of $A'$, following [BG1], Chapter 2. Let $(A'_s)_{s\in[0,1]}$ be a family of flat superconnections on $V\to B$, let $(g_s^V)_{s\in[0,1]}$ be a family of metrics that respect the $\mathbb Z_2$-grading, and define $A_s$, $X_s$ as in (2.13). Define an odd holomorphic function $k$ by

$$k(z)=\frac{f'(z)-f'(0)}{2z}\;.$$

On $\overline B=B\times[0,1]$, we calculate

(2.18)
$$\begin{aligned}d^{\overline B}\bigg((2\pi i)^{-\frac{N^B}{2}}\operatorname{str}_V\bigg(k(X_s)\,\frac{\partial A_s}{\partial s}\bigg)\,ds\bigg)&=(2\pi i)^{\frac{1-N^B}{2}}\operatorname{str}_V\bigg(\Big[A_s,k(X_s)\,\frac{\partial A_s}{\partial s}\Big]\bigg)\,ds\\&=-(2\pi i)^{\frac{1-N^B}{2}}\operatorname{str}_V\bigg(k(X_s)\,\frac{\partial A_s^2}{\partial s}\bigg)\,ds\\&=\frac{\partial}{\partial s}\Big(f\big(A'_s,g_s^V\big)-\operatorname{str}_V\big(X_s\big)\,f'(0)\Big)\,ds\;.\end{aligned}$$

This motivates the following definition.

**2.19. Definition.** Set

$$L_k\big(A'_s,g_s^V\big)=(2\pi i)^{\frac{1-N^B}{2}}\int_0^1\operatorname{str}_V\bigg(k(X_s)\,\frac{\partial A_s}{\partial s}\bigg)\,ds\qquad\in\Omega^*(B;\mathbb C)\;.$$

These classes resemble the Bott-Chern classes in complex geometry. We now state Theorem 2.5 of [BG1] without proof.

**2.20. Theorem.** *The form $L_k(A'_s,g_0^V,g_1^V)$ is even and real. If $\deg k\ge3$, then $L_k(A'_s,g_s^V)$ depends modulo $d^B\Omega^*(B)$ only on $(A'_s,g_s^V)$ at $s=0,1$ and on the homotopy class of the path $(A'_s)_{s\in[0,1]}$ in the space of flat superconnections on $V$. For $k(z)=\frac{f'(z)-f'(0)}{2z}$,*

$$d^BL_k\big(A'_s,g_0^V,g_1^V\big)=f\big(A'_1,g_1^V\big)^{[\ge3]}-f\big(A'_0,g_0^V\big)^{[\ge3]}\;.$$

**2.c. Superconnections of total degree 1 and the form $S_f(A',g^V)$.** We recall the definition of a superconnection of total degree 1. A flat superconnection of degree 1 on a Hermitian vector bundle $V\to B$ gives rise to a fibre-wise cochain complex with a metric $g^H$ and a connection $\nabla^H$ on its cohomology bundle $H\to B$. Then we define a natural form $S_f(A',g^V)\in\Omega^*(B)$ relating the form $f(\nabla^H,g^H)$ to $f(A',g^V)$.

**2.21. Definition.** A superconnection $A$ on a $\mathbb{Z}$-graded vector bundle $V = \bigoplus_j V^j$ is *of total degree* 1 iff

$$A\big(\Omega^j(V^k)\big) \subset \bigoplus_{l\geq 0} \Omega^{j+l}\big(V^{k+1-l}\big) ,$$

and $A$ is *of total degree* $-1$ iff

$$A\big(\Omega^j(V^k)\big) \subset \bigoplus_{l\geq 0} \Omega^{j+l}\big(V^{k-1+l}\big) .$$

*2.22. Remark.* Let $V$ be a $\mathbb{Z}$-graded vector bundle with a Hermition metric $g^V$. We say that $V$ is a *$\mathbb{Z}$-graded Hermitian vector bundle* iff $g^V$ respects the grading, i.e., iff $V^j \perp V^k$ for all $j \neq k$. In this cases, a superconnection is of total degree 1 iff its adjoint is of total degree $-1$.

Now suppose that $V \to B$ is a $\mathbb{Z}$-graded Hermition vector bundle, and let $\nabla^V$ be any connection on $V$ with respect to which all subbundles $V^j$ are parallel. Let $A'$ be a flat superconnection of total degree 1, then

$$A' = \nabla^V + \sum_j a'_j \qquad \text{with} \qquad a'_j \in \bigoplus_k \Omega^j\big(B; \mathrm{Hom}(V^k, V^{k+1-j})\big) , \tag{2.23}$$

and if $\nabla^{V,*}$ denotes the adjoint of $\nabla^V$ and $A''$ the adjoint of $A'$ with respect to $g^V$, then

$$A'' = \nabla^{V,*} + \sum_j a''_j \qquad \text{with} \qquad a''_j = (a'_j)^* \in \bigoplus_k \Omega^j\big(B; \mathrm{Hom}(V^k, V^{k+1-j})\big) . \tag{2.24}$$

Flatness of $A'$ implies in particular

$$(a'_0)^2 = 0 , \qquad \big[\nabla^V + a'_1, a'_0\big] = 0 , \qquad \text{and} \qquad \big(\nabla^V + a'_1\big)^2 + [a'_0, a'_2] = 0 . \tag{2.25}$$

This says first of all that $(V^*, a'_0)$ is a bundle of cochain complexes. Next, the coboundary operator $a'_0$ is parallel with respect to $\nabla^V + a'_1$, so the fibre-wise cohomology $H^*(V^*, a'_0)$ forms a bundle over $B$ denoted $H \to B$, and $\nabla^V + a'_1$ induces a connection $\nabla^H$ on $H$. Finally, the term $a'_2$ provides a cochain homotopy between the curvature of $\nabla^V + a'_1$ and 0, so the induced connection $\nabla^H$ is flat. The connection $\nabla^H$ is called the *Gauß-Manin connection* on $H$.

If $(V, g^V) \to B$ is a $\mathbb{Z}$-graded Hermitian vector bundle as in Remark 2.22, let $a''_0$ denote the adjoint of $a'_0$. By finite dimensional Hodge theory, we may identify $H_b$ with $\ker(a''_0 + a'_0)_b = \ker(a''_0 - a'_0)_b$ for $b \in B$, so we get an identification of vector bundles

$$H = \ker\big(a''_0 + a'_0\big) = \ker\big(a''_0 - a'_0\big) \subset V \longrightarrow B .$$

Let $P^0 \in \Gamma(\mathrm{Hom}(V, H))$ denote the orthogonal projection onto $H$, and let $g^H$ denote the restriction of $g^V$ to $H$. If we take adjoints with respect to $g^V$ and $g^H$, then by [BL], Proposition 2.6, we have

$$\nabla^H = P^0 \left(\nabla^V + a'_1\right) P^0 \qquad \text{and} \qquad \nabla^{H,*} = P^0 \left(\nabla^{V,*} + a''_1\right) P^0 . \tag{2.26}$$

We replace $B$ by $\overline{B} = B \times \mathbb{R}$, and let $\overline{V}$ be the pullback of $V$ to $\overline{B}$, with induced superconnection

$$\overline{A}' = A' + \frac{\partial}{\partial t}\, dt$$

as in (2.23). Let $N^V \in \Gamma(\operatorname{End} V)$ be the number operator that acts by multiplication with $i$ on $V^i$, and consider the metric $g^{\overline{V}}$ on $\overline{V}$ with

$$g^{\overline{V}}(\,\cdot\,,\,\cdot\,)_{(b,t)} = g^V\bigl(t^{N^V}\,\cdot\,,\,\cdot\,\bigr)_b\;.$$

If $A''$ is the adjoint of $A'$ with respect to $g^V$ as in (2.24), then the adjoint of $\overline{A}'$ with respect to $g^{\overline{V}}$ is given by

$$\overline{A}'' = t^{-N^V}\left(A'' + \frac{\partial}{\partial t}\,dt\right) = \nabla^{V,*} + \sum_j t^{1-j} a_j + \left(\frac{\partial}{\partial t} + \frac{N^V}{t}\right) dt\;.$$

We have an operator

$$\overline{X} = \frac{1}{2}\bigl(\overline{A}'' - \overline{A}'\bigr) \qquad \in \Omega^*\bigl(\overline{B};\operatorname{End}\overline{V}\bigr)\;,$$

and we denote its restriction to $B \times \{t\}$ by $X_t$.

If we conjugate by $t^{\frac{N^V}{2}}$, we obtain superconnections

$$\widetilde{A}' = t^{\frac{N^V}{2}}\,\overline{A}'\,t^{-\frac{N^V}{2}} = \nabla^V + \sum_j t^{\frac{1-j}{2}} a'_j + \left(\frac{\partial}{\partial t} - \frac{N^V}{2t}\right) dt\;,$$

$$\text{and} \qquad \widetilde{A}'' = t^{\frac{N^V}{2}}\,\overline{A}''\,t^{-\frac{N^V}{2}} = \nabla^{V,*} + \sum_j t^{\frac{1-j}{2}} a''_j + \left(\frac{\partial}{\partial t} + \frac{N^V}{2t}\right) dt\;,$$

and $\widetilde{A}''$ and $\widetilde{A}'$ are adjoint to each other with respect to the pull-back of $g^V$ with to $\overline{V} \to \overline{B}$ that we still denote by $g^V$. Let

$$\omega^V = \nabla^{V,*} - \nabla^V = \bigl(g^V\bigr)^{-1}\bigl[\nabla^V, g^V\bigr] \qquad \in \Omega^1\bigl(B;\operatorname{End} V\bigr)\;, \tag{2.27}$$

and as in (2.13), set

$$\widetilde{X} = \frac{1}{2}\bigl(\widetilde{A}'' - \widetilde{A}'\bigr) = \frac{\omega^V}{2} + \sum_j t^{\frac{1-j}{2}}\,\frac{a''_j - a'_j}{2} + N^V\,\frac{dt}{2t}\;,$$

and we write $\widetilde{X}_t$ for the restriction of $\widetilde{X}$ to the submanifold $B \times \{t\}$.

We define $\hat{f}$ to be twice the form with the same name in [BL],

$$\hat{f}\bigl(A', g^V\bigr) = (2\pi i)^{-\frac{N^B}{2}}\,\operatorname{str}_{\overline{V}}\Bigl(N^V f'\bigl(\widetilde{X}\bigr)\Bigr)\;. \tag{2.28}$$

If $g^V_t = g^V(t^{N^V}\,\cdot\,,\,\cdot\,)$ denotes the restriction of $g^{\overline{V}}$ to $B \times \{t\}$, then clearly

$$f\bigl(\overline{A}', g^{\overline{V}}\bigr) = f\bigl(\widetilde{A}', g^V\bigr) = f\bigl(\widetilde{A}'_t, g^V\bigr) + \hat{f}\bigl(\widetilde{A}'_t, g^V\bigr)\,\frac{dt}{2t} = f\bigl(A', g^V_t\bigr) + \hat{f}\bigl(A', g^V_t\bigr)\,\frac{dt}{2t}\;, \tag{2.29}$$

and this form is closed on $B \times \mathbb{R}$.

Again differing from [BL], we define

$$\chi(V) = \operatorname{str}_V(1) = \sum_j (-1)^j\,\operatorname{rk}(V^i) \qquad \text{and} \qquad \chi'(V) = \operatorname{str}_V\bigl(N^V\bigr) = \sum_j (-1)^j\,j\,\operatorname{rk}(V^i)\;,$$

so $\chi(V) = \chi(H)$. Then by [BL], Theorem 2.9, as $t \to \infty$,

$$f\bigl(A', g^V_t\bigr) = f\bigl(\nabla^H, g^H\bigr) + O\Bigl(t^{-\frac{1}{2}}\Bigr)\;, \tag{2.30}$$

$$\text{and} \qquad \hat{f}\bigl(A', g^V_t\bigr) = f'(0)\,\chi'(H) + O\Bigl(t^{-\frac{1}{2}}\Bigr)\;.$$

We define $S_f(A', g^V)$ as in [BL], [BG1], where convergence is granted by (2.30).

**2.31. Definition.** Set

$$S_f(A',g^V) = -\int_1^\infty \Big(\hat f(A',g_t^V) - f'(0)\,\chi'(H)\Big)\,\frac{dt}{2t} \quad \in \Omega^*(B)\;.$$

Then we have the following important result of [BL], Theorems 2.16, and [BG1], Theorem 1.25.

**2.32. Theorem.** *The form $S_f(A',g^V)$ is even and real, and satisfies*

$$d^B\,S_f(A',g^V) = f(A',g^V) - f(\nabla^H,g^H)\;.$$

*Proof.* This follows from (2.29) and (2.30). □

Let now $g^V$, $g'^V$ be two metrics on $V\to B$, and let $g^H$, $g'^H$ be the corresponding metrics on $H\to B$. From Theorem 2.32 and Definition 2.16, we easily deduce that

$$S_f(A',g'^V) - S_f(A',g^V) = \tilde f(A',g^V,g'^V) - \tilde f(A',g^H,g'^H)\;. \tag{2.33}$$

In order to control the behaviour of $S_f(A',g'^V)$ under variation of $A'$, we assume for a moment that $f$ is an odd homolomorphic function of degree $\geq 3$ that satisfies (2.12). Then $k(z) = \frac{f'(z)}{2z}$ is also holomorphic. Let $(A'_s,g_s^V)_{s\in[0,1]}$ be a family of flat superconnections and $\mathbb{Z}$-graded metrics on $V$. We assume that the dimension of the cohomology bundles $H_s$ is independent of $s$, then these bundle combine to form a family over $B\times[0,1]$, and we get a family $(\nabla^{H,s},g_s^H)_{s\in[0,1]}$ of flat connections and metrics. We recall Theorem 2.15 of [BG1], again without proof.

**2.34. Theorem.** *Modulo exact forms on $B$, we have*

$$S_f(A'_1,g_1^V) - S_f(A'_0,g_0^V) = L_k(A'_s,g_0^V,g_1^V) - L_k(\nabla^{H,s},g_0^H,g_1^H)\;.$$

**2.d. Family Thom-Smale complexes and their torsion forms.** We define a family Thom-Smale complex to be a $\mathbb{Z}$-graded complex vector bundle $V\to B$ with a flat superconnection $A' = \nabla^V + a'$ of total degree 1 such that the connection $\nabla^V$ is diagonal and the coefficients $a'$ are upper triangular with respect to the eigenspace decomposition of an endomorphism $h^V$, in a sense to be made precise. For such a vector bundle, together with an adapted metric $g^V$, we define a torsion form $T_f(A',g^V,h^V)$ relating $f(\nabla^H,g^H)$ to $f(\nabla^V,g^V)$. In the special case that $a'_j = 0$ for $j\neq 0$, we may take $h^V = N^V$, and our torsion forms equal those of [BL], Chapter 2, that appeared in Theorem 0.1 of [BG1].

**2.35. Definition.** Let $V$ be a $\mathbb{Z}$-graded complex vector bundle as in Remark 2.22, let $A' = \nabla^V + a'$ be a flat superconnection of total degree 1, and let $h^V \in \Gamma(\operatorname{End}^0 V)$ be an endomorphism of $V$ of degree 0. Then $(V,A',h^V)$ is called a *family Thom-Smale complex* iff for each compact subset $K\subset B$ there exists $\varepsilon>0$, and for each point $b\in K$ there exists a neighbourhood $U$ of $b$, such that

1. for each $j$, there is a splitting $V^j|_U = \bigoplus_\alpha V_\alpha^j$, which is preserved by $\nabla^V$;
2. the endomorphism $h^V$ acts on $V_\alpha^j$ by a function $h_\alpha^j\colon U\to\mathbb{R}$;
3. for each $j$ and each $\alpha$,

$$a'(V_\alpha^j) \subset \bigoplus_{\substack{k,\beta\\ h_\alpha^j+\varepsilon\leq h_\beta^k}} \Omega^*(B;V_\beta^k)\;.$$

If $F\colon M\to B$ is a smooth map, then the pullback $F^*(V,A',h^V) = (F^*V,F^*A',F^*h^V)$ is defined in the obvious way.

It seems that we should include $\nabla^V$ in the notation for a family Thom-Smale complex. However, the following statement shows that this is not necessary.

**2.36. Proposition.** *The connection $\nabla^V$ in Definition 2.35 is uniquely determined by $A'$ and $h^V$. Moreover, it is flat.*

*Proof.* Fix a compact subset $K \subset B$, a point $b \in K$ and a small open neighbourhood $U$ of $b$ in $B$ as above. Let $V|_U = \bigoplus_j \bigoplus_\alpha V_\alpha^j$ be the local splitting of (1) above, and let $P_\alpha^j \colon V|_U \to V_\alpha^j$ be the corresponding product projections. By (3), we clearly have

$$P_\alpha^j \, A' \, P_\alpha^j = P_\alpha^j \, \nabla^V \, P_\alpha^j = \nabla^{V_\alpha^j}$$

over $U$, and since $\nabla^V$ preserves $V_\alpha^j$, the connection $\nabla^{V_\alpha^j}$ is simply the restriction of $\nabla^V$ over $U$ to $V_\alpha^j$. In particular, $\nabla^V|_U = \bigoplus_j \bigoplus_\alpha \nabla^{V_\alpha^j}$ is uniquely determined by $A'$.

Similarly, we can calculate the curvature of $\nabla^{V_\alpha^j}$ by

$$(\nabla^{V_\alpha^j})^2 = P_\alpha^j \, (\nabla^V)^2 \, P_\alpha^j = P_\alpha^j \, (A')^2 \, P_\alpha^j = 0 \; ,$$

again using (3) to show that the derivatives of $a'$ do not contribute to $(\nabla^{V_\alpha^j})^2$. □

For example, if $A' = \nabla^V + a_0'$ is a flat family of complexes as in [BL], section 2.f, that is, $a_0' \in \Omega^0(B; \operatorname{End}^1 V)$, we may take $h^V = N^V$ in order to obtain a family Thom-Smale complex. If we constructed $V$, $\nabla^V$ and $A'$ as in Section 1, then we may choose $h^V$ to act on $V_\alpha \subset V$ by $h_\alpha$.

**2.37. Definition.** Let $(V, A', h^V)$ be a family Thom-Smale complex. An *adapted metric* on $V$ is a metric

$$g^V = \bigoplus_i g^{V^i} \qquad \text{on } V$$

such that the splitting of Definition 2.35 (1) is $g^V$-orthogonal.

*2.38. Remark.* If $(V, A', h^V)$ is a family Thom-Smale complex and $g^V$ is an adapted metric, then the endomorphisms $h^V$ and $N^V$ are self-adjoint with respect to $g^V$. Because $\nabla^V$ and $g^V$ respect the splitting of $V$ in Definition 2.35 (1), so do $\nabla^{V,*}$ and $\omega^V$. If we set

$$\nabla^{V,\mathrm{u}} = \frac{1}{2}\left(\nabla^{V,*} + \nabla^V\right) ,$$

then $\nabla^{V,\mathrm{u}}$ also preserves that splitting. If we write

$$dh^V = \left[\nabla^V, h^V\right] = \left[\nabla^{V,*}, h^V\right] = \left[\nabla^{V,\mathrm{u}}, h^V\right] ,$$

then $dh^V$ acts on $V_\alpha^j$ by $dh_\alpha^j$. We have the obvious relations

$$\left[\nabla^{V,\mathrm{u}}, dh^V\right] = \left[\nabla^{V,\mathrm{u}}, \omega^V\right] = \left[h^V, \omega^V\right] = \left[dh^V, \omega^V\right] = 0 \; .$$

Let $\overline{A}' = A' + \frac{\partial}{\partial t} \, dt + \frac{\partial}{\partial T} \, dT$ denote the trivial extension of $A'$ to the pullback $\overline{V}$ of $V$ to $B \times \mathbb{R}_> \times \mathbb{R}$. Define a metric $g^{\overline{V}}$ on $\overline{V}$ by

$$g^{\overline{V}} = g^V\left(t^{N^V} \, T^{h^V} \, \cdot\, , \cdot\,\right) . \tag{2.39}$$

If we take the adjoint of $A'$ with respect to $g^{\overline{V}}$, we obtain a flat superconnection

$$\overline{A}'' = \nabla^{V,*} + \sum_j t^{1-j} \, T^{-h^V} \, a_j'' \, T^{h^V} + \log T \, dh^V + \left(\frac{\partial}{\partial t} + \frac{N^V}{t}\right) dt + \left(\frac{\partial}{\partial T} + \frac{h^V}{T}\right) dT$$

and an operator

$$\overline{X} = \frac{1}{2}\left(\overline{A}'' - \overline{A}'\right) .$$

We denote the restrictions of $g^{\overline{V}}$, $\overline{A}'$, $\overline{A}''$, $\overline{X}$ to $B \times \{(t,T)\}$ by $g^V_{t,T}$, $A'_{t,T}$, $A''_{t,T}$ and $X_{t,T}$.

If we conjugate everything with $t^{\frac{N^V}{2}}\,T^{\frac{h^V}{2}}$, we obtain two superconnections

(2.40)

$$\widetilde{A}' = \nabla^V + \sum_j t^{\frac{1-j}{2}}\,T^{\frac{h^V}{2}}\,a'_j\,T^{-\frac{h^V}{2}} - \frac{\log T}{2}\,dh^V + \left(\frac{\partial}{\partial t} - \frac{N^V}{2t}\right) dt + \left(\frac{\partial}{\partial T} - \frac{h^V}{2T}\right) dT$$

and $$\widetilde{A}'' = \nabla^{V,*} + \sum_j t^{\frac{1-j}{2}}\,T^{-\frac{h^V}{2}}\,a''_j\,T^{\frac{h^V}{2}} + \frac{\log T}{2}\,dh^V + \left(\frac{\partial}{\partial t} + \frac{N^V}{2t}\right) dt + \left(\frac{\partial}{\partial T} + \frac{h^V}{2T}\right) dT$$

that are adjoints of each other with respect to $g^V$, and a $g^V$-skew-adjoint operator

$$\widetilde{X} = \frac{1}{2}\left(\widetilde{A}'' - \widetilde{A}'\right) . \tag{2.41}$$

Fix a compact subset $K$ of $B$ and let $\varepsilon$ be the constant of Definition 2.35. From now on, we keep $t$ fixed. Then

$$\begin{aligned} \widetilde{X}_{t,T} &= \frac{\omega^V}{2} + \frac{1}{2}\sum_j t^{\frac{1-j}{2}}\left(T^{-\frac{h^V}{2}}\,a''_j\,T^{\frac{h^V}{2}} - T^{\frac{h^V}{2}}\,a'_j\,T^{-\frac{h^V}{2}}\right) + \frac{\log T}{2}\,dh^V \\ &= \frac{\omega^V}{2} + \frac{\log T}{2}\,dh^V + O\!\left(T^{\frac{\varepsilon}{2}}\right) \end{aligned} \tag{2.42}$$

uniformly on $K$ as $T \to 0$.

Define $\nabla^{V,\mathrm{u}}$ as in Remark 2.38, then since $f$ and $\omega^V + \log T\,dh^V$ are odd with $[\omega^V, dh^V] = (dh^V)^2 = 0$, and $\omega^V + \log T\,dh^V$ commutes with $h^V$, we have

$$\operatorname{str}_V\left(h^V f'\Bigl(\frac{\omega^V + \log T\,dh^V}{2}\Bigr)\right) = \operatorname{str}_V\bigl(h^V\bigr)\,f'(0)\,,$$

and $$\operatorname{str}_V\left(f\Bigl(\frac{\omega^V + \log T\,dh^V}{2}\Bigr)\right) = \operatorname{str}_V\left(f\Bigl(\frac{\omega^V}{2}\Bigr)\right) + \frac{\log T}{2}\operatorname{str}_V\left(dh^V f'\Bigl(\frac{\omega^V}{2}\Bigr)\right),$$

with $$\begin{aligned}\operatorname{str}_V\left(dh^V f'\Bigl(\frac{\omega^V}{2}\Bigr)\right) &= \operatorname{str}_V\left(\Bigl[\nabla^{V,\mathrm{u}}, h^V f'\Bigl(\frac{\omega^V}{2}\Bigr)\Bigr]\right) \\ &= d^B \operatorname{str}_V\left(h^V f'\Bigl(\frac{\omega^V}{2}\Bigr)\right) = \operatorname{str}_V\bigl(dh^V\bigr)\,f'(0)\end{aligned}$$

by Remark 2.38, so

$$\begin{aligned} f\Bigl(\nabla^V + \frac{\partial}{\partial T}\,dT, g^{\overline{V}}\Bigr) &= f\left(\nabla^V - \frac{\log T}{2}\,dh^V + \left(\frac{\partial}{\partial T} - \frac{h^V}{2T}\right) dT, g^V\right) \\ &= (2\pi i)^{\frac{1-N^B}{2}}\operatorname{str}_V\left(f\Bigl(\frac{\omega^V + \log T\,dh^V}{2}\Bigr)\right) \\ &\qquad + (2\pi i)^{-\frac{N^B}{2}}\operatorname{str}_V\left(h^V f'\Bigl(\frac{\omega^V + \log T\,dh^V}{2}\Bigr)\right)\frac{dT}{2T} \\ &= f\bigl(\nabla^V, g^V\bigr) + d^{\overline{B}}\left(\frac{\log T}{2}\operatorname{str}_V\bigl(h^V\bigr)\right) f'(0)\,. \end{aligned}$$

Define an odd closed and real form

$$\alpha = f\Big(\overline{A}', g^{\overline{V}}\Big) - f\Big(\nabla^V + \frac{\partial}{\partial t}\,dt + \frac{\partial}{\partial T}\,dT, g^{\overline{V}}\Big) \quad \in \Omega^*\big(B \times (0,\infty) \times \mathbb{R}\big)\;, \tag{2.43}$$

then for fixed $t$,

$$\begin{aligned} \alpha_t &= f\big(\widetilde{A}'_{t,T}, g^V\big) - f\big(\nabla^V, g^V\big) \\ &\qquad + (2\pi i)^{-\frac{N^B}{2}}\,\mathrm{str}_V\Big(h^V\, f'\big(\widetilde{X}_{t,T}\big)\Big)\,\frac{dT}{2T} - d^{\overline{B}}\Big(\frac{\log T}{2}\,\mathrm{str}_V\big(h^V\big)\Big)\,f'(0)\;. \end{aligned} \tag{2.44}$$

For a fixed compact $K \subset B$, let $\varepsilon$ be as in Definition 2.35, and choose $0 < \varepsilon' < \frac{\varepsilon}{2}$. Using (2.42), we compute

$$\begin{aligned} \alpha_t &= (2\pi i)^{\frac{1-N^B}{2}}\,\mathrm{str}_V\bigg(f\Big(\frac{\omega^V + \log T\,dh^V}{2} + O\big(T^{\frac{\varepsilon}{2}}\big)\Big)\bigg) - f\big(\nabla^V, g^V\big) \\ &\qquad + (2\pi i)^{-\frac{N^B}{2}}\,\mathrm{str}_V\bigg(h^V\,f'\Big(\frac{\omega^V + \log T\,dh^V}{2} + O\big(T^{\frac{\varepsilon}{2}}\big)\Big)\bigg)\,\frac{dT}{2T} \\ &\qquad - d^{\overline{B}}\Big(\frac{\log T}{2}\,\mathrm{str}_V\big(h^V\big)\Big)\,f'(0) \\ &= O\big(T^{\varepsilon'}\big) + O\big(T^{\varepsilon'}\big)\,\frac{dT}{2T} \end{aligned} \tag{2.45}$$

as $T \to 0$, where the terms $O(T^{\varepsilon'})$ do not involve the exterior variable $dT$. Using this formula and (2.30), we find that the integrals in the following definition converge.

**2.46. Definition.** Let $(V, A', h^V)$ be a family Thom-Smale complex with an adapted metric $g^V$, then we define its *torsion form* in $\Omega^*(B)/d^B\Omega^*(B)$ by

$$\begin{aligned} T_f\big(A', g^V, h^V\big) = -(2\pi i)^{-\frac{N^B}{2}} \int_1^\infty &\Big(\mathrm{str}_V\big(N^V f'(X_{t,1})\big) - \chi'(H) f'(0) \\ &\qquad - \big(\chi'(V) - \chi'(H)\big)\,f'\big(\sqrt{-t}/2\big)\Big)\,\frac{dt}{2t} \\ - (2\pi i)^{-\frac{N^B}{2}} \int_0^1 &\Big(\mathrm{str}_V\big(h^V f'(X_{1,T})\big) - \mathrm{str}_V\big(h^V\big)\,f'(0) \\ &\qquad - \big(\chi'(V) - \chi'(H)\big)\,\big(f'(\sqrt{-T}/2) - f'(0)\big)\Big)\,\frac{dT}{2T}\;. \end{aligned}$$

*2.47. Remark.* Suppose that $A' = \nabla^V + a_0$ is a flat superconnection on a $\mathbb{Z}$-graded Hermitian vector bundle $(V, g^V)$, then the form $T_f(A', g^V, N^V)$ clearly equals the form

$$\begin{aligned} T_f\big(A', g^V\big) = -(2\pi i)^{-\frac{N^B}{2}} \int_0^\infty &\Big(\mathrm{str}_V\big(N^V f'(X_{t,1})\big) \\ &\qquad - \chi'(H) f'(0) - \big(\chi'(V) - \chi'(H)\big)\,f'\big(\sqrt{-t}/2\big)\Big)\,\frac{dt}{2t} \end{aligned} \tag{(3.30)}$$

of [BL], [BG1]. This motivates some of the correction terms in Definition 2.46. In particular, the component of degree 0 can be written as

$$T_f\big(A', g^V, N^V\big)^{[0]} = \frac{1}{2} \sum_j \mathrm{str}_{H^\perp}\big(N^V\ \log(a_0' + a_0'')^2\big)\;.$$

**2.48. Theorem.** *The form $T_f(A', g^V, h^V) \in \Omega^*(B)$ is even and real, and if $F\colon M \to B$ is a smooth map, then*

$$T_f\big(F^*A', F^*g^V, F^*h^V\big) = F^*T_f\big(A', g^V, h^V\big)\,.$$

*We have*

$$d^B T_f\big(A', g^V, h^V\big) = f\big(\nabla, g^V\big) - f\big(\nabla^H, g^H\big)\,.$$

*The component of $T_f(A', g^V, h^V)$ in degree zero is given by*

$$T_f\big(A', g^V, h^V\big)^{[0]} = \frac{1}{2}\,\mathrm{str}\big(N^V\,\log(a_0' + a_0'')^2|_{\ker(a_0'+a_0'')^\perp}\big)\,.$$

*Proof.* The form $T_f(A', g^V, h^V)$ is even and real because $f\big(\overline{A}', g^{\overline{V}}\big)$ is odd on $B \times \mathbb{R}_> \times \mathbb{R}$ and real by Theorem 2.15. The first equation of Theorem 2.48 is clear from Definition 2.46.

It follows from Definition 2.31 and Theorem 2.32 that

$$\begin{aligned}(2.49)\qquad -d^B\Bigg((2\pi i)^{-\frac{N^B}{2}}\int_1^\infty\Big(&\mathrm{str}_V\big(N^V f'(X_{t,1})\big) - \chi'(H)f'(0)\\ &- \big(\chi'(V)-\chi'(H)\big)\,f'\big(\sqrt{-t}/2\big)\Big)\,\frac{dt}{2t}\Bigg)\\ &= d^B\,S_f\big(A', g^V\big) = f\big(A', g^V\big) - f\big(\nabla^H, g^H\big)\,.\end{aligned}$$

Fix $K \subset B$ compact and $\varepsilon' > 0$ as above, then equations (2.44) and (2.45) imply that similarly,

$$\begin{aligned}(2.50)\qquad -d^B\Bigg((2\pi i)^{-\frac{N^B}{2}}\int_0^1\Big(&\mathrm{str}_V\big(h^V f'(X_{1,T})\big) - \mathrm{str}_V\big(h^V\big)\,f'(0)\\ &- \big(\chi'(V)-\chi'(H)\big)\,\big(f'(\sqrt{-T}/2) - f'(0)\big)\Big)\,\frac{dT}{2T}\Bigg)\\ &= -d^B\int_{B\times[0,1]/B}\alpha_1 = -\alpha_{1,1} = f\big(\nabla^V, g^V\big) - f\big(A', g^V\big)\,.\end{aligned}$$

The second equation in Theorem 2.48 follows from (2.49) and (2.50).

Now fix $b \in B$ and consider the $\mathbb{Z}$-graded Hermitian vector bundle $(V_b \times [0,1], g^{V_b}) \to [0,1]$ with the trivial flat connection $\nabla^0$. This bundle admits a flat superconnection $A_b' = a_0' + \nabla^0$. Let $u$ denote the coordinate of $[0,1]$, then $(V_b \times [0,1], A_b', u\,h^V + (1-u)\,N^V)$ is a family Thom-Smale complex with adapted metric $g^{V_b}$. We have

$$dT_f\big(A_b', g^{V_b}, u\,h^V + (1-u)\,N^V\big) = f\big(\nabla^0, g^{V_b}\big) - f\big(\nabla^H, g^H\big) = 0 \quad\in\quad \Omega^*([0,1])\,.$$

On the other hand, Theorem 2.25 in [BL] shows that at $u = 0$, we obtain the classical Reidemeister torsion of the complex $(V, a_0')$, which proves the last assertion. $\square$

*2.51. Remark.* Suppose that $\gamma \in \Gamma(\mathrm{End}\,V)$ acts by fibre-wise isometries on $V$, leaving $A'$ and $h^V$ invariant. Then there is an obvious equivariant analogue $T_{f,\gamma}(A', g^V, h^V)$ satisfying an equivariant version of Theorem 2.48.

**2.e. Rigidity of the forms $T_f(A', g^V, h^V)$.** We analyse the behaviour of $T_f(A', g^V, h^V)$ under smooth variation of $A'$, $g^V$ and $h^V$.

**2.52. Definition.** Let $V \to B$ be a $\mathbb{Z}$-graded complex vector bundle, and let $(A'_s, g^V_s, h^V_s)_{s\in[0,1]}$ be family of flat superconnections of total degree 1 with $A'_s = \nabla^{V,s} + a'_s$, $\mathbb{Z}$-graded Hermitian metrics, and $\mathbb{Z}$-graded endomorphisms on $V$. We call $(A'_s, g^V_s, h^V_s)_{s\in[0,1]}$ a *family of family Thom-Smale complexes* iff for each compactum $K \subset B \times [0,1]$ there exists $\varepsilon > 0$, and for each point $b \in K$ there exists a neighbourhood $U$ of $b$ in $B \times [0,1]$, such that Definition 2.35 (1)–(3) hold over $U$.

Note that both in Definition 2.35 and here, the splittings of $V^j|_U = \bigoplus_\alpha V^j_\alpha$ are allowed to depend on $U$. In particular, it may happen that several eigenvalues of $h^V$ run together over a certain subset of $B \times [0,1]$.

Let us assume for the moment that $f$ is an odd holomorphic function of degree $\geq 3$ that satisfies (2.12), cf. Theorem 2.34, and set $k(z) = \frac{f'(z)}{2z}$ as before. We also assume that the dimension of the cohomology bundles $H_s$ is independent of $s \in [0,1]$, so that they combine to form a bundle on $B \times [0,1]$. Let $(\nabla^{H,s}, g^H_s)_{s\in[0,1]}$ be the induced family of flat connections and Hermitian metrics.

**2.53. Proposition.** *Modulo exact forms on $B$, we have*

$$T_f\bigl(A'_1, g^V_1, h^V_1\bigr) - T_f\bigl(A'_0, g^V_0, h^V_0\bigr) = L_k\bigl(\nabla^{V,s}, g^V_0, g^V_1\bigr) - L_k\bigl(\nabla^{H,s}, g^H_0, g^H_1\bigr) .$$

*Proof.* Because we already know Theorem 2.34, we only have to care about the integral over $T$ in Definition 2.46. For fixed $t$, define a form

$$\ell_t = (2\pi i)^{-\frac{N^B}{2}} \left(\operatorname{str}_V\left(k\bigl(\overline{X}_s\bigr)\,\frac{\partial\overline{A}_s}{\partial s}\right) - \operatorname{str}_V\left(k\left(\frac{\omega^V + \log T\, dh^V_s}{2} + \frac{h^V_s}{2T}\,dT\right)\frac{\partial\nabla^{V,s,\mathrm{u}}}{\partial s}\right)\right) ds$$

on $B \times [0,1] \times \{t\} \times \mathbb{R}$. As in (2.18), we have

$$d^{B\times[0,1]\times\{t\}\times\mathbb{R}}\ell_t = f\left(\overline{A}'_s, g^V_s\right) - f\left(\nabla^{V,s} + \frac{\partial}{\partial T}\,dT, g^V_s\right) = \alpha_{s,t}\,ds ,$$

where $\alpha_{s,t}$ is the form defined in (2.43). In particular, the form $\ell_t + \alpha_t$ is closed on $B\times[0,1]\times\{t\}\times\mathbb{R}$.

As $T \to 0$, we have

$$\begin{aligned} \frac{\partial\widetilde{A}_{s,t}}{\partial s} &= \nabla^{V,s,\mathrm{u}} + O\Bigl(T^{\frac{\varepsilon}{2}}\Bigr) \\ \text{and}\qquad \widetilde{X}_{s,t} &= \frac{\omega^{V,s} + \log T\, dh^V_s}{2} + \frac{h^V_s}{2T}\,dT + O\Bigl(T^{\frac{\varepsilon}{2}}\Bigr) . \end{aligned}$$

This implies that

$$\ell_t = O\Bigl(T^{\frac{\varepsilon}{2}}\Bigr) + O\Bigl(T^{\frac{\varepsilon}{2}}\Bigr)\frac{dT}{2T} \tag{2.54}$$

uniformly in $s$ as $T \to 0$, where the forms $O(T^{\frac{\varepsilon}{2}})$ do not contain the exterior variable $dT$.

By (2.45), (2.54) and Stokes' theorem,

$$d^B \int_{B\times[0,1]\times[0,1]/B} (\ell_1 + \alpha_1) = \int_{B\times[0,1]/B} (\alpha_{1,1} - \alpha_{0,1}) - \int_{B\times[0,1]/B} \ell_{1,1}$$

for $t = 1$ fixed, where the first integration on the right hand side is in the direction of $T$, and the second in the direction of $s$. Together with Theorem 2.34, this proves our proposition. ☐

The next step towards the general case is an analogue of Theorem 2.18 in [BG1]. Let $f$ be an odd holomorphic function satisfying (2.12). Let $a > 0$, let $R$ be a polynomial, and set

$$f_a(z) = R\Bigl(\frac{\partial}{\partial a}\Bigr)\Big|_{a=1} \frac{h(\sqrt{a}\,z)}{\sqrt{a}} . \tag{2.55}$$

Then $f_a$ also satisfies (2.12).

**2.56. Proposition.** *The torsion with respect to the generating function $f_a$ is up to exact forms on $B$ given by*

$$T_{f_a}\big(A', g^V, h^V\big) = R\bigg(\frac{\partial}{\partial a}\bigg)\Big|_{a=1}\Big(a^{\frac{N^B}{2}}\, T_f\big(A', g^V, h^V\big)\Big)\;.$$

*Proof.* Let us regard the form $\alpha$ of (2.43) for variable $t$, but $T=1$. Because

$$f\bigg(\nabla^V + \frac{\partial}{\partial t}\,dt + \frac{\partial}{\partial T}\,dT, g^{\overline{V}}\bigg)\Big|_{T=1} = f\big(\nabla^V, g^V_t\big) + \chi'(V)\,\frac{dt}{2t}\;,$$

we have

$$\alpha|_{T=1} = f\big(\overline{A}', g^{\overline{V}}\big) - f\big(\nabla^V, g^V_t\big) - \chi'(V)\, f'(0)\,\frac{dt}{2t}\;.$$

Setting

$$\beta = \big(\chi'(V) - \chi'(H)\big)\,\Big(f'\big(\sqrt{-t}/2\big) - f'(0)\Big)\,\bigg(\frac{dt}{2t} + \frac{dT}{2T}\bigg)\;,$$

we may therefore rewrite Definition 2.46 as

$$T_f\big(A', g^V, h^V\big) = -\int_{B\times\big(\{1\}\times[0,1]\,\cup\,(1,\infty)\times\{1\}\big)/B} (\alpha - \beta)\;.$$

Rescaling $t$ to $at$ gives

$$\begin{aligned}\widetilde{X}_{at} &= \frac{\omega^V}{2} + \sum_j (at)^{\frac{1-j}{2}}\Big(T^{-\frac{h^V}{2}}\, a''_j\, T^{\frac{h^V}{2}} - T^{\frac{h^V}{2}}\, a'_j\, T^{-\frac{h^V}{2}}\Big) + \frac{\log T}{2}\, dh^V + \frac{h^V}{2T}\, dT\\ &= a^{\frac{1-N^B}{2}}\, \widetilde{X}_t\, a^{\frac{N^B}{2}}\;.\end{aligned}$$

This implies that

$$\bigg(\frac{f\big(\sqrt{a}\,\cdot\,\big)}{\sqrt{a}}\bigg)'(X_{t,T}) = f'\big(\sqrt{a}\,X_{t,T}\big) = a^{\frac{N^B}{2}}\, f'(X_{at,T})\, a^{-\frac{N^B}{2}}\;.$$

We may thus express the torsion form with respect to $\frac{f(\sqrt{a}\,\cdot\,)}{\sqrt{a}}$ as

$$\begin{aligned} &T_{\frac{f(\sqrt{a}\,\cdot\,)}{\sqrt{a}}}\big(A', g^V, h^V\big)\\ &\quad= -a^{\frac{N^B}{2}}\int_{B\times(1,\infty)\times\{1\}/B}\Big(\operatorname{str}_V\big(f'(X_{at,T})\big) - \chi'(H)\, f'(0)\\ &\qquad\qquad\qquad - \big(\chi'(V)-\chi'(H)\big)\, f'\big(\sqrt{-at}/2\big)\Big)\,\frac{dt}{2t}\\ &\qquad - a^{\frac{N^B}{2}}\int_{B\times\{a\}\times[0,1]/B}\bigg(\alpha - \big(\chi'(V)-\chi'(H)\big)\,\Big(f'\big(\sqrt{-aT}/2\big) - f'(0)\Big)\,\frac{dT}{2T}\bigg)\\ &\quad= -a^{\frac{N^B}{2}}\int_{B\times\big([1,a]\times\{0\}\,\cup\,\{a\}\times[0,1]\,\cup\,(a,\infty)\times\{1\}\big)/B}(\alpha-\beta)\;.\end{aligned}$$

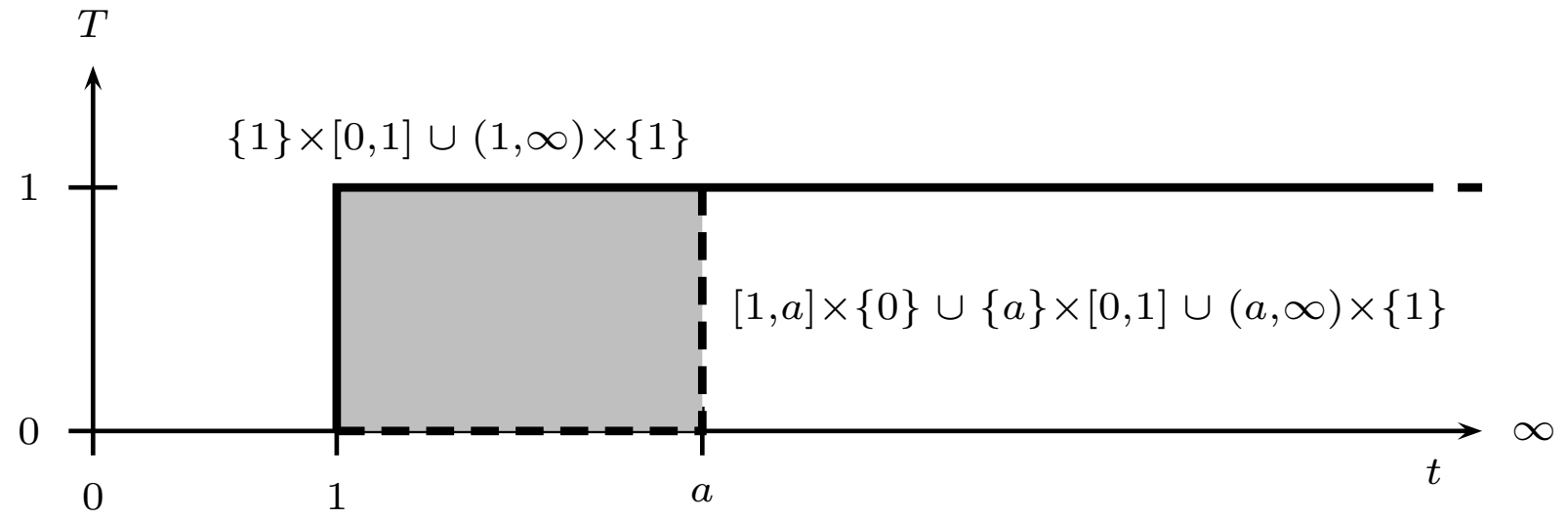

Figure 2.57. Paths in the construction of $T_f$ and $T_{f_\alpha}$.

Here we have used that $\alpha$ vanishes on $B \times [1,a] \times \{0\}$ by (2.45), and that $\beta$ is symmetric in $t$ and $T$. Since the integrand above is closed, we may integrate along any path from $(\infty,1)$ to $(1,0)$, changing the result only by an exact form, so

$$T_{\frac{f(\sqrt{a}\,\cdot)}{\sqrt{a}}}\bigl(A',g^V,h^V\bigr) = a^{\frac{N^B}{2}}\,T_f\bigl(A',g^V,h^V\bigr) + a^{\frac{N^B}{2}}\,d^B \int_{B\times[1,a]\times[0,1]/B} (\alpha-\beta)\;.$$

Our theorem now follows, since

$$\begin{aligned} T_{f_a}\bigl(A',g^V,h^V\bigr) &= R\biggl(\frac{\partial}{\partial a}\biggr)\Big|_{a=1} T_{\frac{f(\sqrt{a}\,\cdot)}{\sqrt{a}}}\bigl(A',g^V,h^V\bigr)\\ &= R\biggl(\frac{\partial}{\partial a}\biggr)\Big|_{a=1}\Bigl(a^{\frac{N^B}{2}}\,T_f\bigl(A',g^V,h^V\bigr)\Bigr)\;. \qquad\square \end{aligned}$$

We can now finally prove the analogue of [BG1], Theorem 2.20. We return to the convention (2.11) that $f(z) = z\,e^{z^2}$, and we set $k(z) = \frac{f'(z)-f'(0)}{2z}$.

**2.58. Theorem.** *Modulo exact forms on $B$, we have*

$$\Bigl(T_f\bigl(A_1',g_1^V,h_1^V\bigr) - T_f\bigl(A_0',g_0^V,h_0^V\bigr)\Bigr)^{[\geq 2]} = L_k\bigl(\nabla^{V,s},g_0^V,g_1^V\bigr) - L_k\bigl(\nabla^{H,s},g_0^H,g_1^H\bigr)\;.$$

*2.59. Remark.* Note that the right hand side is independend modulo $d\Omega^*(B)$ of the path connecting $\nabla^{V,0}$ with $\nabla^{V,1}$ and $\nabla^{H,0}$ with $\nabla^{H,1}$, even though both $L_k$-terms may depend on the homotopy class of these paths in degree $\geq 4$, and on the paths themselves in degree 2.

*Proof.* Choose $f(z) = z\,e^{z^2}$ as in (2.11) and set $R(X) = X$ in (2.55), so that $f_a(z) = z^3\,e^{z^2}$ is now of degree 3, and let $k_a(z) = \frac{f_a'(z)}{2z}$. By Proposition 2.53 and Proposition 2.56,

$$\begin{aligned} &L_{k_a}\bigl(\nabla^{V,s},g_0^V,g_1^V\bigr) - L_{k_a}\bigl(\nabla^{H,s},g_0^H,g_1^H\bigr)\\ &\qquad = \frac{\partial}{\partial a}\Big|_{a=1}\biggl(a^{\frac{N^B}{2}}\Bigl(T_f\bigl(A_1',g_1^V,h_1^V\bigr) - T_f\bigl(A_0',g_0^V,h_0^V\bigr)\Bigr)\biggr)\\ &\qquad = \frac{N^B}{2}\Bigl(T_f\bigl(A_1',g_1^V,h_1^V\bigr) - T_f\bigl(A_0',g_0^V,h_0^V\bigr)\Bigr)\;. \end{aligned}$$

On the other hand, we have

$$k_a(z) = \frac{3z+z^3}{2}\,e^{z^2} = \frac{\partial}{\partial a}\Big|_{a=1}\biggl(\frac{(1+2az^2)\,e^{az^2}-1}{2z}\biggr) = \frac{\partial}{\partial a}\Big|_{a=1}\Bigl(\sqrt{a}\,k\bigl(\sqrt{a}\,z\bigr)\Bigr)\;.$$

It is now easy to check that

$$\begin{aligned} L_{k_a}\big(\nabla^{V,s},g_0^V,g_1^V\big) &= \frac{\partial}{\partial a}\Big|_{a=1}\int_0^1 \operatorname{str}_V\left(k\Big(\sqrt{a}\,\frac{\omega^{V,s}}{2}\Big)\sqrt{a}\,\frac{\partial\nabla^{V,s,\mathrm{u}}}{\partial s}\right)ds \\ &= \frac{\partial}{\partial a}\Big|_{a=1}\int_0^1 a^{\frac{N^B}{2}}\operatorname{str}_V\left(k\Big(\frac{\omega^{V,s}}{2}\Big)\frac{\partial\nabla^{V,s,\mathrm{u}}}{\partial s}\right)ds \\ &= \frac{N^B}{2}\,L_k\big(\nabla^{V,s},g_0^V,g_1^V\big)\;. \end{aligned}$$

Together, these equations give

$$\frac{N^B}{2}\Big(T_f\big(A_1',g_1^V,h_1^V\big)-T_f\big(A_0',g_0^V,h_0^V\big)\Big)=\frac{N^B}{2}\Big(L_k\big(\nabla^{V,s},g_0^V,g_1^V\big)-L_k\big(\nabla^{H,s},g_0^H,g_1^H\big)\Big)\;,$$

from which our theorem finally follows. □

**2.f. Fibre-wise Morse functions and torsion classes.** We show how to associate torsion classes to families of compact manifolds with fibre-wise Morse functions. If we start with a flat bundle $(F,\nabla^F)\to M$ that admits a parallel metric $g^F$, such that the fibre-wise cohomology $(H,\nabla^H)\to B$ also admits a parallel metric, then we obtain a well-defined cohomology class $T_f(\nabla^F,h,g^F,g^H)\in H^{\mathrm{even},\geq 2}(B)$ that is invariant under smooth deformations of $h$ in the space of all fibre-wise Morse-functions on $M$.

Let $p\colon M\to B$ be a bundle with compact fibres, and assume that $(F,\nabla^F)\to M$ is a flat bundle and $h\colon M\to\mathbb{R}$ is a Morse function on each fibre of $p$. If $\hat p\colon C\to B$ denotes the covering of $B$ by the fibre-wise critical points of $h$, we define the bundle $V=\hat p_*(\Lambda^{\max}T^uX\otimes F|_C)$, with the $\mathbb{Z}$-grading by the fibre-wise index of $h$ as in Definition 1.1. Then $h$ acts by multiplication on $\Lambda^{\max}T^uX\otimes F|_C$, and we let $h^V\in\Gamma(\operatorname{End}^0V)$ be the pushdown of this action. Let $A'=\nabla^V+a$ be the superconnection of Definition 1.60, then $a\in\Omega^*(B;\operatorname{End}_+^*V)$ is of total degree 1 by Theorem 1.61.

*2.60. Remark.* The bundle $(V,A',h^V)$ is a family Thom-Smale complex in the sense of Definition 2.35, and any Hermitian metric $g^F$ on $F\to M$ induces by restriction to $F_C$ and push-down to $V\to B$ an adapted metric $g^V$ on $V$. In fact, if we denote by $C_\alpha^j$ the various local leaves of $\hat p^j$ where $h$ has index $j$, we may take the $\nabla^V$-parallel splitting

$$V^j=\bigoplus_\alpha F|_{C_\alpha^j}\;,$$

such that $h^V$ acts on $C_\alpha^j$ by multiplication with $h_\alpha^j=h|_{C_\alpha^j}$, so points (1) and (2) of Definition 2.35 are satisfied. We have $A'=\nabla^V+a'$ with $a\in\Omega^*(B;\operatorname{End}_+^*V)$ by Definition 1.3 and Theorem 1.61, so (3) also holds.

Suppose that $A_0'$ and $A_1'$ are two flat superconnections of total degree one on $V\to B$, both constructed as in Section 1. In particular, both $(V,A_0',h^V)$ and $(V,A_1',h^V)$ are family Thom-Smale complexes. Then we know by Proposition 1.65 that there exists a flat superconnection

$$\overline{A}'=\nabla^{\overline{V}}+\sum_j\bar a_j'$$

of total degree one on the pullback $\overline{V}$ of $V$ to $\overline{B}=B\times[0,1]$, also constructed as in Section 1, such that $\overline{A}'|_{B\times\{s\}}=A_s'$ for $s=0$, 1.

We fix a metric $g^F$ on $F \to M$, let $g^V$ be the induced metric o $V$, and let $g^{\overline{V}}$ be its pull-back to $\overline{V} \to \overline{B}$. Similarly, $h^V$ pulls back to an endomorphism $h^{\overline{V}}$ of $\overline{V}$, and $\big(\overline{V}, \overline{A}', h^{\overline{V}}\big)$ is a family Thom-Smale complex with adapted metric $g^{\overline{V}}$.

Let $\overline{H} = H^*\big(\overline{V}, \bar{a}_0'\big) \to \overline{B}$ denote the fibre-wise cohomology with its Gauß-Manin connection $\nabla^{\overline{H}}$, and let $g^{\overline{H}}$ denote the metric induced on $\overline{H}$ by $g^{\overline{V}}$ using finite-dimensional Hodge theory. Then for $s = 0$, $1$, the restriction of $\big(\overline{H}, \nabla^{\overline{H}}, g^{\overline{H}}\big)$ to $B \times \{s\}$ is naturally isomorphic to the cohomology bundle $H \to B$ with Gauß-Manin connection $\nabla^H$ and metric $g_s^H$ induced by $g^V$. Thus, we have

$$\tilde{f}\big(\nabla^H, g_0^H, g_1^H\big) = \int_0^1 f\big(\nabla^{\overline{H}}, g^{\overline{H}}\big)$$

as in Definition 2.16.

**2.61. Proposition.** *Let $A_0'$ and $A_1'$ be two superconnections on $V$ that are constructed as in Definition 1.60, and let $g^V$ be induced by a metric $g^F$ on $F \to M$. Then*

$$T_f\big(A_1', g^V, h^V\big) - T_f\big(A_0', g^V, h^V\big) = -\tilde{f}\big(\nabla^H, g_0^H, g_1^H\big) \qquad \textit{modulo } d^B\Omega^*(B).$$

*Proof.* By Theorem 2.48 and Stokes' Theorem,

$$\begin{aligned} T_f\big(A_1', g^V, h^V\big) - T_f\big(A_0', g^V, h^V\big) &= d \int_{\overline{B}/B} T_f\big(\overline{A}', g^{\overline{V}}, h^{\overline{V}}\big) \\ &\quad + \int_{\overline{B}/B} \Big( f\big(\nabla^{\overline{V}}, g^{\overline{V}}\big) - f\big(\nabla^{\overline{H}}, g^{\overline{H}}\big)\Big) \\ &= -\tilde{f}\big(\nabla^H, g_0^H, g_1^H\big) \qquad \text{modulo } d^B\Omega^*(B), \end{aligned}$$

because $f\big(\nabla^{\overline{V}}, g^{\overline{V}}\big)$ is just the pullback of $f(\nabla^V, g^V)$ to $\overline{B}$, so $\int_{\overline{B}/B} f\big(\nabla^{\overline{V}}, g^{\overline{V}}\big)$ vanishes. □

In particular, if we fix metrics $g^F$ on $F$ and $g^H$ on $H$, we can define a torsion form for the bundle $M \to B$, equipped with a flat vector bundle $(F, \nabla^F) \to M$ and a fibre-wise Morse function $h$ that is independent of all choices made during the constructions in Section 1, by setting

$$T_f\big(\nabla^F, h, g^F, g^H\big) = T_f\big(A', g^V, h^V\big) + \tilde{f}\big(\nabla^H, g^H, g_V^H\big) \quad \in \Omega^*(B)/d^B\Omega^*(B)\,. \tag{2.62}$$

It follows from Proposition 2.61 that $T_f(\nabla^F, h, g^F, g^H)$ is indeed well-defined.

**2.63. Corollary.** *Assume that the metrics $g^F|_C$ and $g^H$ are parallel, then $T_f(\nabla^F, h, g^F, g^H)$ defines a cohomology class in $H^{\mathrm{even}}_{\mathrm{dR}}(B)$. Moreover, only the homogeneous part of degree $0$ depends on the particular choice of $g^F|_C$ and $g^H$.*

*Proof.* The first assertion is immediate from Theorem 2.48 and Definition 2.64. The second claim follows because any two parallel metrics $g_0^F$, $g_1^F$ on a flat bundle are related by family of flat metrics $g_s^F = (1-s)\, g_0^F + s\, g_1^F$. Let $g^{\overline{F}}$ be the metric induced by this family on the pullback $\overline{F} \to \overline{B} = B \times [0,1]$. Let $\nabla^{\overline{F},*}$ be the adjoint of the pullback $\nabla^{\overline{F}}$ of $\nabla^F$ to $\overline{F}$ with respect to $g^{\overline{F}}$. Then

$$\omega^{\overline{F}} = \nabla^{\overline{F},*} - \nabla^{\overline{F},*} \quad \in \Omega^1\big(\overline{B}, \operatorname{End}\overline{F}\big)$$

vanishes on all vectors in $T\overline{B}$ that are tangent to the factor $B$. In particular, $\big(\omega^{\overline{F}}\big)^k = 0$ for $k > 1$, so

$$\tilde{f}\big(\nabla^F, g_0^F, g_1^F\big) = \tilde{f}\big(\nabla^F, g_0^F, g_1^F\big)^{[0]}\,. \qquad \square$$

**2.64. Definition.** Assuming that the metrics $g^F|_C$ on $F|_C \to C$ and $g^H$ on $H \to B$ are parallel, we define the *higher Morse torsion*

$$T_f(M/B;F,h) = T_f\big(A', g^V, h^V\big)^{[\geq 2]} + \tilde{f}\big(\nabla^H, g^H, g^H_V\big)^{[\geq 2]} \quad \in H^{\mathrm{even},\geq 2}_{\mathrm{dR}}(B)\;.$$

**2.g. Generalised metrics and the classes $S_f(A', \mathbf{g}^V, g^H)$.** We recall the notion of a generalised metric following [BG1], Section 2.9. Our presentation will follow the lines of Section 2.a above. We also recall the definition of $S_f$ for generalised metrics. The content of this subsection will be needed only in Section 4, where generalised metrics appear naturally as restrictions of "classical" metrics to certain subbundles of $\Lambda^*TB \mathbin{\widehat{\otimes}} \Omega(M/B;F)$ formed by "generalised eigenspaces" of the Witten-deformed operator $\mathbb{X}_{t,T}$.

Recall that in Section 2.a, we extended a metric $g^V$ on a vector bundle $V \to B$ to an $\mathbb{R}$-bilinear map $g^V \colon (\Lambda^*TB \mathbin{\widehat{\otimes}} V)^2 \to \Lambda^*TB$ with the properties (2.2) and (2.3). A natural generalisation of metrics is the following.

**2.65. Definition.** A *generalised Hermitian metric* $\mathbf{g}^V$ on a $\mathbb{Z}_2$-graded complex vector bundle $V \to B$ is an $\mathbb{R}$-bilinear map $\mathbf{g}^V \colon (\Lambda^*TB \mathbin{\widehat{\otimes}} V)^2 \to \Lambda^*TB$, such that

(1) the component $(\mathbf{g}^V)^{[0]} \colon (g^V)^2 \to \mathbb{C}$ is a Hermitian metric;
(2) $\mathbf{g}^V$ is even with respect to the total grading of $\Lambda^*TB \mathbin{\widehat{\otimes}} V$;
(3) for $v$, $w \in V$, we have

$$\mathbf{g}^V(w,v) = (-1)^{\left[\frac{\deg(\mathbf{g}^V(w,v))}{2}\right]}\,\overline{\mathbf{g}^V(v,w)}\;;$$

(4) for $\alpha$, $\beta \in \Lambda^*TB$ and $v$, $w \in V$, we have

$$\mathbf{g}^V\big(\alpha \mathbin{\widehat{\otimes}} v, \beta \mathbin{\widehat{\otimes}} w\big) = (-1)^{\left[\frac{\deg(\beta)}{2}\right]}\,\alpha \wedge \mathbf{g}^V\big(v,w\big) \wedge \overline{\beta}\;.$$

Properties (3) and (4) are the natural generalisations of (2.2) and (2.3). Together, they imply that (3) still holds for $v$, $w \in \Lambda^*TB \mathbin{\widehat{\otimes}} V$. Note that as in Remark 2.9, we may regard $\mathbf{g}^V$ as an even element of $\Omega^*\big(B; \mathrm{Hom}(V, \overline{V}^*)\big)$. Clearly, there is an analogous notion of a generalised Euclidean metric. Note that by (2), the part $(\mathbf{g}^V)^{[0]}$ of horizontal degree 0 respects the $\mathbb{Z}^2$-grading of $V$.

*2.66. Remark.* Since any convex combination of generalised metrics is again a generalised metric, the space of generalised metrics on a given vector bundle is contractible. In particular, any generalised metric can be deformed onto its component $(\mathbf{g}^V)^{[0]}$.

When we define adjoints of endomorphisms and superconnections with respect to general metrics, we have to take care of the fact that $\mathbf{g}^V(v,w)$ may be an odd form for $v$, $w \in V$.

**2.67. Definition.** For an endomorphism $E \in \Omega^*(B; \mathrm{End}\, V)$, we define the adjoint endomorphism by

$$E^* = (\mathbf{g}^V)^{-1}\,\overline{E}^*\,\mathbf{g}^V\;,$$

and for a superconnection $A$ on $V$, we set

$$A^* = (\mathbf{g}^V)^{-1}\,\overline{A}^*\,\mathbf{g}^V\;.$$

*2.68. Remark.* For an endomorphism $E \in \mathrm{End}\, V$, Definition 2.67 is equivalent to

$$\mathbf{g}^V(E^*w, v) = (-1)^{|E|\,(|E|+|w|+|v|)}\,\mathbf{g}^V(w, Ev)$$

for all $v$, $w \in V$. Note that $|E| + |w| + |v|$ has the same parity as the exterior degree of the form $\mathbf{g}^V(w, Ev)$.

For a connection $\nabla^V$, we get

$$\nabla^{V,*} = \nabla^V + \omega^V \qquad \text{with} \qquad \omega^V = \big(\mathbf{g}^V\big)^{-1} \big[\nabla^V, \mathbf{g}^V\big] \quad \in \quad \Omega^*(B; \operatorname{End} V)\;.$$

Definition 2.65 (3), (4) and Definition 2.67 imply that $(E^*)^* = E$ and $(A^*)^* = A$ for all superconnections $A$ and all $E \in \Omega^*(M; \operatorname{End} V)$.

For odd elements $\alpha \widehat{\otimes} E \in \Omega^*(M; \operatorname{End} V)^{\mathrm{odd}}$ and for superconnections $A$ on $V$, Definition 2.65 (4) and Definition 2.67 imply that

$$\begin{aligned} 0 &= \langle E^* v, w\rangle + (-1)^{\deg E\,(\deg v + \deg w)}\, \langle v, Ew\rangle \\ \text{and} \qquad d\langle v, w\rangle &= \langle A^* v, w\rangle + (-1)^{\deg E\,(\deg v + \deg w)}\, \langle v, Aw\rangle \end{aligned}$$

for all $v$, $w \in \Omega^*(M; V)$. Thus, most of the properties of metrics, superconnections and their adjoints generalise to the present setting.

Let $f\colon \mathbb{C} \to \mathbb{C}$ be a holomorphic function satisfying (2.12), e.g., $f(z) = z\, e^{z^2}$ as in (2.11). If $\mathbf{g}^V$ is a generalised Hermitian metric and $A'$ is a superconnection, let $A''$ be the adjoint of $A'$ with respect to $\mathbf{g}^V$, and set $X = \frac{1}{2}\,(A'' - A')$ as in (2.13). We have a natural generalisation of Definition 2.14.

**2.69. Definition.** Let $f(A', \mathbf{g}^V) \in \Omega^*(B)$ be defined by

$$f\big(A', \mathbf{g}^V\big) = (2\pi i)^{\frac{1-N^B}{2}}\, \operatorname{str}_V\big(f(X)\big)\;.$$

Theorem 2.15 extends literally to this generalised setting, so in particular, $f(A', \mathbf{g}^V)$ is odd, closed, and real. Similarly, we define $\tilde f(A', \mathbf{g}^V_0, \mathbf{g}^V_1) \in \Omega^*(B)/d^B\Omega^*(B)$ as in Definition 2.16. Again, Corollary 2.17 extends literally to this situation and gives in particular

$$d^B \tilde f\big(A', \mathbf{g}^V_0, \mathbf{g}^V_1\big) = f\big(A', \mathbf{g}^V_1\big) - f\big(A', \mathbf{g}^V_0\big)\;.$$

We also need a generalisation of Definition 2.31 to this setting. Therefore, let $(V, g^V) \to B$ be a $\mathbb{Z}$-graded Hermitian vector bundle, let $(\mathbf{g}^V_t)_{t\in(0,\infty)}$ be a smooth family of generalised metrics on $V$, and let $g' \in \Omega^*(B; \operatorname{End} V)$ be even, such that as $t \to \infty$,

$$t^{\frac{n-N^V}{2}}\, \mathbf{g}^V_t\, t^{-\frac{N^V}{2}} = g^V\Bigg(\bigg(\mathrm{id}_V + \frac{g'}{\sqrt{t}} + O\Big(\frac{1}{t}\Big)\bigg)\,\cdot\,,\,\cdot\Bigg)\;, \tag{2.70}$$

$$\text{and} \qquad t^{\frac{N^V}{2}}\, \Big(\mathbf{g}^V_t\Big)^{-1}\, \frac{\partial \mathbf{g}^V_t}{\partial t}\, t^{-\frac{N^V}{2}} = \Big(N^V - \frac{n}{2}\Big)\,\frac{1}{t} + O\Big(t^{-\frac{3}{2}}\Big)\;,$$

where both $O(\,\cdot\,)$ are even elements of $\Omega^*(B; \operatorname{End} V)$. Let $\bar{\mathbf{g}}^V$ be the generalised metric on $\overline{V} = V \times (0,\infty) \to \overline{B} = B \times (0,\infty)$ induced by $\mathbf{g}^V_t$, in particular, $\bar{\mathbf{g}}^V$ does not contain the exterior variable $dt$.

Let $A'$ be a flat superconnection of total degree 1 as in Definition 2.21, so

$$A' = \nabla^V + \sum_i a_i \qquad \text{with} \qquad a_i \in \Omega^i\big(B; \operatorname{End}^{1-i} V\big)\;,$$

and let $\overline{A}' = A' + \frac{\partial}{\partial t}\, dt$ be the trivial extension of $A'$ to $\overline{V}$. Let $\overline{A}''$ be the adjoint of $\overline{A}'$ with respect to $\bar{\mathbf{g}}^V$. Note that $\overline{A}''$ is not necessarily of total degree $-1$, since $\bar{\mathbf{g}}^V$ can shuffle degrees. Nevertheless, using that $\overline{A}'$ is of total degree 1 together with (2.70), we establish the following properties of $A'$ and $A''_t$.

**2.71. Proposition** ([BG1], Proposition 2.46). *As $t \to \infty$,*

$$t^{\frac{N^V}{2}} A' t^{-\frac{N^V}{2}} = \sqrt{t}\, a_0 + \nabla^V + a_1 + O\Big(t^{-\frac{1}{2}}\Big) ,$$
$$\text{and} \qquad t^{\frac{N^V}{2}} A''_t t^{-\frac{N^V}{2}} = \sqrt{t}\, a_0^* + \nabla^{V,*} + a_1^* + [a_0^*, g'] + O\Big(t^{-\frac{1}{2}}\Big) ,$$

*where adjoints are taken with respect to $g^V$.* □

As in (2.13), we define

$$\overline{X} = \frac{1}{2} \big(\overline{A}'' - \overline{A}'\big) = \frac{1}{2} \big(A''_t - A'\big) + \frac{1}{2} \Big(\mathbf{g}_t^V\Big)^{-1} \frac{\partial \mathbf{g}_t^V}{\partial t}\, dt .$$

We define an odd closed differential form on $\overline{B}$ by

$$(2.72) \quad \beta = f\big(\overline{A}', \bar{\mathbf{g}}^V\big) - \Big(\chi'(H) - \frac{n}{2}\, \chi(V)\Big)\, f'(0)\, \frac{dt}{2t}$$
$$= f\big(A', \mathbf{g}_t^V\big) + (2\pi i)^{-\frac{N^B}{2}} \Bigg( \operatorname{str}_V \bigg( \frac{1}{2} \Big(\mathbf{g}_t^V\Big)^{-1} \frac{\partial \mathbf{g}_t^V}{\partial t}\, f'\big(X_t\big) \bigg) - \Big(\chi'(H) - \frac{n}{2}\, \chi(V)\Big) \frac{f'(0)}{2t} \Bigg)\, dt .$$

Let $g^H$ be the metric induced by $g^V$ on $H = H^*(V, a_1) \to B$, and let $\nabla^H$ be the Gauß-Manin connection on $H$ induced by $A'$. If $P^H\colon V \to H$ denotes the $g^V$-orthogonal projection onto the $(a_0 + a_0^*)$-harmonic elements of $V$, then clearly $P^H\, [a_0^*, g']\, P^H = 0$. Now the proof of [BL], Theorem 2.9, translates to this situation and gives an estimate for $\beta$ as $t \to \infty$, cf. [BG1], Proposition 2.48.

**2.73. Proposition.** *As $t \to \infty$,*

$$\beta = f\big(\nabla^H, g^H\big) + O\Big(t^{-\frac{1}{2}}\Big) + O\Big(t^{-\frac{3}{2}}\Big)\, dt ,$$

*where the $O(\,\cdot\,) \in \Omega^*\big(\overline{B}; \operatorname{End} \overline{V}\big)$ do not contain $dt$.* □

Let $\mathbf{g}_{s,t}^V$ be two generalised metrics on $V$, together with $g_s^V$ and $g'_s$, such that (2.70) holds for $s = 0$, 1. Then there exist families $\mathbf{g}_{s,t}^V$, $g_s^V$ and $g'_s$ for $s \in [0,1]$ that restrict to the given metrics for $s = 0$, 1, such that (2.70) holds uniformly in $s$. We define a form $\overline{\beta}$ on $B \times [0,1] \times (0, \infty)$, and we let $\beta_s$ denote its restriction to $B \times \{s\} \times (0, \infty)$. Let $g_s^H$ denote the metric induced on $H$ by $g_s^V$, then the obvious generalisation of (2.33) is

$$(2.74) \qquad \int_{(B \times (1,\infty))/B} \big(\beta_1 - \beta_0\big) = \tilde{f}\big(A', \mathbf{g}_{0,1}^V, \mathbf{g}_{1,1}^V\big) - \tilde{f}\big(\nabla^H, g_0^H, g_1^H\big) \quad \in \Omega^*(B)/d^B\Omega^*(B) .$$

Writing $\mathbf{g}^V = \mathbf{g}_1^V$, me can now define a generalisation of the class $S_f(A', g^V)$.

**2.75. Definition.** Set

$$S_f\big(A', \mathbf{g}^V, g^H\big) = \int_{(B \times (1,\infty))/B} \beta \quad \in \Omega^*(B)/d^B\Omega^*(B) .$$

Note that the class $S_f(A', \mathbf{g}^V, g^H)$ was called $U_f(A', \mathbf{g}_t^V)$ in [BG1].

**2.76. Theorem** ([BG1], Theorem 2.50). *The form $S_f(A', \mathbf{g}^V, g^H)$ is even and satisfies*

$$d^B S_f\big(A', \mathbf{g}^V, g^H\big) = f\big(A', \mathbf{g}^V\big) - f\big(\nabla^H, g^H\big) \,.$$

*If $\mathbf{g}'^V$ is another generalised metric on $V$ and $g'^H$ is another metric on $H$, them*

$$S_f\big(A', \mathbf{g}'^V, g'^H\big) - S_f\big(A', \mathbf{g}^V, g^H\big) = \tilde{f}\big(A', \mathbf{g}^V, \mathbf{g}'^V\big) - \tilde{f}\big(\nabla^H, g^H, g'^H\big) \,.$$

Let $g^V$ be a fixed metric on $V$ that respects the $\mathbb{Z}$-grading. Note that the family $g^V\, t^{N^V}$ does not quite satisfy (2.70). Me may however set

$$g_t^V = g^V\, t^{N^V - \frac{n}{2}} \,,$$

and since $t^{-\frac{n}{2}}$ commutes with all operators involved, we have

$$A_t'' = t^{-N^V + \frac{n}{2}}\, A''\, t^{N^V - \frac{n}{2}} = t^{-N^V}\, A''\, t^{N^V} \,.$$

Because of this, our form $\beta$ of (2.72) becomes

$$\begin{aligned}\beta = f\big(A', g^V\, t^{N^V}\big) + (2\pi i)^{-\frac{N^B}{2}} \bigg(&\mathrm{str}_V\Big(\Big(N^V - \frac{n}{2}\Big)\, f'(X_t)\Big) \\ &- \Big(\chi'(H) - \frac{n}{2}\,\chi(V)\Big)\, f'(0)\bigg)\, \frac{dt}{2t} \\ = f\big(A', g^V\, t^{N^V}\big) + (2\pi i)^{-\frac{N^B}{2}} &\Big(\mathrm{str}_V\big(N^V\, f'(X_t)\big) - \chi'(H)\, f'(0)\Big)\, \frac{dt}{2t} \,.\end{aligned}$$

Let again $g^H$ be the metric induced on $H$ by $g^V$. A comparison with our definition of $S_f(A', g^V)$ in Definition 2.31 shows that

$$S_f\big(A', g^V, g^H\big) = S_f\big(A', g^V\big) \quad \in \Omega^*(B)/d^B\Omega^*(B) \,, \tag{2.77}$$

which justifies our notation.

**2.h. Analytic torsion forms.** We recall the construction of analytic torsion forms in [BL].

Let $p\colon M \to B$ be a smooth fibre bundle, and let $(F, \nabla^F) \to M$ be a flat complex vector bundle. Let $\Omega^*(M; F)$ denote the space of $F$-valued differential forms on $M$. If we fix a horizontal subbundle $T^HM \subset TM$ with $T^HM \oplus TX = TM$, then $T^HM \cong p^*TB$, and we may split $\Lambda^*TM = p^*\Lambda^*TB \mathbin{\widehat{\otimes}} \Lambda^*TX$. We denote the space of vertical forms by $\Omega^*(M/B; F) = p_*(\Lambda^*TX \otimes F) \to B$, then we have

$$\Omega^*(M; F) = \Omega^*\big(B; \Omega^*(M/B; F)\big) \,. \tag{2.78}$$

Let $\Omega \in \Gamma(p^*\Lambda^2TB \otimes TX)$ denote the fibre-bundle curvature of $T^HM$, given by

$$\Omega(v, w) = -\big[\bar{v}, \bar{w}\big]^\perp \in TX \,,$$

where $\bar{v}$, $\bar{w} \in \Gamma(T^HM)$ are horizontal lifts of vectorfields $v$, $w$ on $B$ and $(\,\cdot\,)^\perp$ denote the projection $TM \to TX$ along $T^HM$. With respect to the splitting (2.78), we may interpret the exterior differential $d^M$ induced by $\nabla^F$ on $\Omega^*(M;F)$ as a flat superconnection $\mathbb{A}'$ of total degree 1 on the infinite-dimensional vector bundle $\Omega^*(M/B;F)\to B$. Classically, this superconnection decomposes as

$$\mathbb{A}'=d^M=d^X+\nabla^{\Omega^*(M/B;F)}+\iota_\Omega \tag{2.79}$$

with respect to horizontal degree, cf. [BGV].

Now let us fix a vertical metric $g^{TX}$, then $g^{TX}$ induces a metric on $\Lambda^*TX$ that we still denote by $g^{TX}$. We obtain a natural $L_2$-metric on $\Omega^*(M/B;F)$, which we write as

$$g^{TX,F}(\omega,\omega')_b=\int_{p^{-1}(b)}\big(g^{TX}\otimes g^F\big)(\omega,\omega')_x\,d\operatorname{Vol}(x)\;, \tag{2.80}$$

where $d\operatorname{Vol}$ is the fibre-wise volume element induced by $g^{TX}$. Let $\mathbb{A}''$ be the adjoint of $\mathbb{A}'$ with respect to $g^{TX,F}$ in the sense of Definition 2.7. Then

$$\mathbb{X}=\frac12\big(\mathbb{A}''-\mathbb{A}'\big)=\frac12\Big(\big(d^{X,*}-d^X\big)+\omega^{\Omega^*(M/B;F)}-\varepsilon_\Omega-\iota_\Omega\Big)\in\Omega^*\big(B;\operatorname{End}\Omega^*(M/B;F)\big) \tag{2.81}$$

is a fibre-wise elliptic differential operator on $\Omega^*(M;F)$ with coefficients in $\Omega^*(B)$. Moreover, $\mathbb{X}$ is skew-adjoint in the sense of Definition 2.7.

By Hodge theory, we can identify the fibre-wise cohomology $H=H^*_{\mathrm{dR}}(M/B;F)$ with the bundle of $(d^{X,*}-d^X)$-harmonic forms in $\Omega^*(M/B;F)$. Then $\nabla^{\Omega^*(M/B;F)}$ induces the Gauß-Manin connection $\nabla^H$ on $H$. We let $g^H_{L_2}$ denote the restriction of $g^{TX,F}$ to $H$.

Because $-\mathbb{X}^2$ has the same spectrum as the Hodge Laplacian $-(\mathbb{X}^{[0]})^2=\frac14\,(d^{X,*}+d^X)^2$, the odd heat kernel $\mathbb{X}\,e^{\mathbb{X}^2}$ is of trace-class. In particular, as in Definition 2.14, we obtain an odd, closed and real form

$$f\big(\mathbb{A}',g^{TX,F}\big)=\operatorname{str}\big(\mathbb{X}\,e^{\mathbb{X}^2}\big)\;. \tag{2.82}$$

Let $N^X$ denote the fibre-wise number operator on $\Omega^*(M;F)$ and $\Omega^*(M/B;F)$. We set $g^{TX}_t=g^{TX}(t^{-N^X}\,\cdot\,,\,\cdot\,)$. Let $g^{TX,F}_t$ denote the induced $L_2$-metric on $\Omega^*(M/B;F)$, then

$$g^{TX,F}_t=g^{TX,F}\big(t^{N^X-\frac n2}\,\cdot\,,\,\cdot\,\big)\;, \tag{2.83}$$

where the additional $t^{\frac n2}$ is due to the change of the fibre-wise volume element. Let $\bar g^{TX,F}$ denote the metric on $\Omega^*(M/B;F)\times(0,\infty)$ such that $\bar g^{TX,F}|_{B\times\{t\}}=g^{TX,F}_t$. Let $\overline{\mathbb{A}}''$ be the adjoint of the obvious extension $\overline{\mathbb{A}}'$ of the superconnection $\mathbb{A}'$ to the bundle $\Omega^*(M/B;F)\times(0,\infty)\to B\times(0,\infty)$, and set $\overline{\mathbb{X}}=\frac12\big(\overline{\mathbb{A}}''-\overline{\mathbb{A}}'\big)$. As before, restriction to $B\times\{t\}$ will be written $(\,\cdot\,)_t$.

The form

$$f\big(\overline{\mathbb{A}}',\bar g^{TX,F}\big)-\chi'(H)\,f'(0)\,\frac{dt}{2t}=f\big(\mathbb{A}'_t,g^{TX,F}\big)+(2\pi i)^{-\frac{N^B}{2}}\Big(\operatorname{str}\big(N^X\,f'(\mathbb{X}_t)\big)-\chi'(H)\,f'(0)\Big)\,\frac{dt}{2t}$$

is closed and satisfies

$$f\big(\overline{\mathbb{A}}',\bar g^{TX,F}\big)-\chi'(H)\,f'(0)\,\frac{dt}{2t}=f\big(\nabla^H,g^H_{L_2}\big)+O\Big(t^{-\frac12}\Big)++O\Big(t^{-\frac32}\Big)\,dt\;,$$

as $t\to\infty$ by [BL], Theorem 3.21, where the $O(\,\cdot\,)$-terms do not involve $dt$.

**2.84. Definition.** Set

$$S_f\big(T^HM,g^{TX},\nabla^F,g^F\big)=-(2\pi i)^{-\frac{N^B}{2}}\int_1^\infty\Big(\operatorname{str}\big(N^X\,f'(\mathbb{X}_t)\big)-\chi'(H)\,f'(0)\Big)\,\frac{dt}{2t}\;.$$

Then clearly, $S_f(T^HM,g^{TX},\nabla^F,g^F)$ is even and real, and satisfies

$$d^BS_f\big(T^HM,g^{TX},\nabla^F,g^F\big)=\operatorname{ch}^o\big(\mathbb{A}',g^{TX,F}\big)-\operatorname{ch}^o\big(\nabla^H,g^H\big)\;.$$

Let $\nabla^{TX}$ be the natural connection on $TX\to M$ that is induced by $T^HM$ and $g^{TX}$, and let $e(\nabla^{TX})$ be its Euler form. As $t\to 0$, again by [BL], Theorem 3.21,

$$\begin{aligned}f\big(\overline{\mathbb{A}}',\bar g^{TX,F}\big)-\chi'(H)\,\frac{dt}{2t}&\\&=e\big(\nabla^{TX}\big)\,f\big(\nabla^F,g^F\big)+\Big(\frac n2\,\chi(H)-\chi'(H)\Big)\,f'(0)\,\frac{dt}{2t}+O\big(t^{\frac12}\big)+O\big(t^{-\frac12}\big)\,dt\;,\end{aligned}$$

where the $O(\,\cdot\,)$-terms do not contain $dt$.

**2.85. Definition.** Set

$$\begin{aligned}\mathcal{T}_f\big(T^HM,g^{TX},\nabla^F,g^F\big)=-(2\pi i)^{-\frac{N^B}{2}}\int_0^\infty\Bigg(&\operatorname{str}\big(N^X\,f'(\mathbb{X}_t)\big)-\chi'(H)\,f'(0)\\&-\Big(\frac n2\,\chi(H)-\chi'(H)\Big)\,f'\Big(\frac{\sqrt{-t}}{2}\Big)\Bigg)\,\frac{dt}{2t}\;.\end{aligned}$$

Let $T^{H,s}M$, $g^{TX}_s$, $\nabla^{F,s}$ and $g^F_s$ denote families of horizontal subbundles and vertical metrics on $TM$, and of flat connections and Hermitian metrics on $F$ for $s\in[0,1]$, where we assume that the rank of $H_s=H^*_{\mathrm{dR}}(M/B;(F,\nabla^{F,s}))$ does not depend on $s$. Let $\nabla^{H,s}$ and $g^{H,s}_{L_2}$ denote the families of induced flat connections and $L_2$-metrics. Let $\nabla^{TX,s}$ denote the induced connections on the vertical tangent bundle, and let $\tilde e(\nabla^{TX,0},\nabla^{TX,1})$ be the Chern-Simons class associated to the Euler class. Let $k(z)=\frac{f'(z)-f'(0)}{2z}$. We now summarise [BL], Theorems 3.23, 3.24, 3.29 and [BG1], Theorem 3.45.

**2.86. Theorem.** *The form $\mathcal{T}_f(T^HM,g^{TX},\nabla^F,g^F)$ is even and real, and satisfies*

$$d^B\mathcal{T}_f\big(T^HM,g^{TX},\nabla^F,g^F\big)=\int_{M/B}e\big(\nabla^{TX}\big)\,f\big(\nabla^F,g^F\big)-f\big(\nabla^H,g^H_{L_2}\big)\;.\tag{1}$$

*The component $\mathcal{T}_f(T^HM,g^{TX},\nabla^F,g^F)^{[0]}$ of degree $0$ is equal to the Ray-Singer analytic torsion $\mathcal{T}_{\mathrm{RS}}(g^{TX},\nabla^F,g^F)$ of the fibre. Moreover, modulo exact forms on $B$,*

$$\begin{aligned}&\mathcal{T}_f\big(T^{H,1}M,g^{TX}_1,\nabla^F,g^F_1\big)-\mathcal{T}_f\big(T^{H,0}M,g^{TX}_0,\nabla^F,g^F_0\big)\\&\qquad=\int_{M/B}\tilde e(\nabla^{TX,0},\nabla^{TX,1})\,f\big(\nabla^F,g^V_0\big)\\&\qquad\quad+\int_{M/B}e(\nabla^{TX,1})\,\tilde f\big(\nabla^F,g^V_0,g^V_1\big)-\tilde f\big(\nabla^H,g^{H,0}_{L_2},g^{H,1}_{L_2}\big)\;,\end{aligned}\tag{2}$$

$$(3) \qquad \textit{and} \qquad \mathcal{T}_f\big(T^HM, g^{TX}, \nabla^{F,1}, g^F_1\big)^{[\geq 2]} - \mathcal{T}_f\big(T^HM, g^{TX}, \nabla^{F,0}, g^F_0\big)^{[\geq 2]} \\ = \int_{M/B} e(\nabla^{TX}_0)\, L_k\big(\nabla^{F,s}, g^V_0, g^V_1\big) - L_k\big(\nabla^{H,s}, g^{H,0}_{L_2}, g^{H,1}_{L_2}\big)\ .$$

*2.87. Remark.* Note again that the right hand side of (3) is independend modulo $d\Omega^*(B)$ of the paths connecting $\nabla^{F,0}$ with $\nabla^{F,1}$ and $\nabla^{H,0}$ with $\nabla^{H,1}$, even though both $L_k$-terms may depend on the homotopy class of these paths in degree $\geq 4$, and on the paths themselves in degree 2.

Suppose now that both the bundles $F \to M$ and $H \to B$ carry parallel metrics $g^F$ and $g^H$. Then as in Corollary 2.63, the form

$$\mathcal{T}_f\big(T^HM, g^{TX}, \nabla^F, g^F\big) + \tilde f\big(\nabla^H, g^H, g^H_{L_2}\big)$$

is closed, and its cohomology class is independend of the choice of the parallel metrics $g^F$ and $g^H$ in degree $\geq 2$.

**2.88. Definition.** If the metrics $g^F$ and $g^H$ are parallel, define the *higher analytic torsion*

$$\mathcal{T}_f\big(M/B; F\big) = \mathcal{T}_f\big(T^HM, g^{TX}, \nabla^F, g^F\big)^{[\geq 2]} + \tilde f\big(\nabla^H, g^H, g^H_{L_2}\big)^{[\geq 2]} \quad \in H^{\mathrm{even},\geq 2}_{\mathrm{dR}}(B)\ .$$

**2.i. The Chern normalisation.** We recall the Chern normalisation of the odd characteristic classes of Definition 2.14 and the torsion forms of Definition 2.46 and Definition 2.85 following [BG1]. Note that "normalisation" means the choice of a particular generating function $f$. Theorem 0.1 does not depend on $f$, and the construction of torsion forms using the function $f$ of (2.11) is very convenient for proving it. Nevertheless, results like Theorem 0.1 of [BG2] depend on the choice of normalisation. This suggests that the most natural normalisation of torsion forms is the one we are introducing now.

Let $(F, \nabla^F) \to M$ be a flat vector bundle. A metric $g^F$ on $F$ determines an adjoint connection $\nabla^{F,*}$ which is also flat, so

$$\mathrm{ch}\big(\nabla^F\big) = \mathrm{tr}\Big(e^{-\frac{(\nabla^F)^2}{2\pi i}}\Big) = \mathrm{ch}\big(\nabla^{F,*}\big) = \mathrm{rk}\, F$$

as differential forms. Let $\widetilde{\mathrm{ch}}(\nabla^F, \nabla^{F,*})$ be the natural Chern-Simons class in $\Omega^*(M)/d^M\Omega^*(M)$, then

$$d^M\widetilde{\mathrm{ch}}\big(\nabla^F, \nabla^{F,*}\big) = \mathrm{ch}\big(\nabla^{F,*}\big) - \mathrm{ch}\big(\nabla^F\big) = 0\ ,$$

so $\widetilde{\mathrm{ch}}(\nabla^F, \nabla^{F,*})$ is in fact a cohomology class.

In order to construct a natural representative for this class, we choose a connection

$$\nabla^{\overline{F}} = (1-s)\,\nabla^F + s\,\nabla^{F,*} + \frac{\partial}{\partial s}\,ds$$

on $\overline{F} = F \times [0,1] \to \overline{B} = B \times [0,1]$. Set $\omega^F = \nabla^{F,*} - \nabla^F$ as before. Because $\nabla^F$ and $\nabla^{F,*}$ are flat, the curvature of $\nabla^{\overline{F}}$ is given by

$$\big(\nabla^{\overline{F}}\big)^2 = s(1-s)\,\big[\nabla^F, \nabla^{F,*}\big] - \big(\nabla^{F^*} - \nabla^F\big)\,ds = -s(1-s)\,(\omega^F)^2 - \omega^F\,ds\ .$$

**2.89. Definition.** Set

$$\mathrm{ch}^{\mathrm{o}}\big(\nabla^F,g^F\big)=\int_0^1\mathrm{tr}\bigg(\frac{\omega^V}{2}\,e^{\frac{s(1-s)\,(\omega^V)^2}{2\pi i}}\bigg)\,ds\;.$$

**2.90. Theorem.** *The form* $\mathrm{ch}^{\mathrm{o}}(\nabla^F,g^F)$ *is real, odd and closed and represents the Chern-Simons class* $\pi i\,\widetilde{\mathrm{ch}}(\nabla^F,\nabla^{F,*})$ *in* $H^*_{\mathrm{dR}}(M)$*. Moreover, for* $f(z)=z\,e^{z^2}$,

$$\mathrm{ch}^{\mathrm{o}}\big(\nabla^F,g^F\big)=\int_0^1\big(4s(1-s)\big)^{\frac{N^B-1}{2}}\,f\big(\nabla^F,g^F\big)\,ds\;.$$

*Proof.* It is clear that $\mathrm{ch}^{\mathrm{o}}(\nabla^F,g^F)$ represents $\pi i\,\widetilde{\mathrm{ch}}(\nabla^F,\nabla^{F,*})$ by the definition of Chern-Simons classes. Moreover,

$$\begin{aligned}\mathrm{ch}^{\mathrm{o}}\big(\nabla^F,g^F\big)&=\int_0^1\mathrm{tr}\bigg(\frac{\omega^V}{2}\,e^{\frac{s(1-s)\,(\omega^V)^2}{2\pi i}}\bigg)\,ds\\&=\int_0^1\big(4s(1-s)\big)^{\frac{N^B-1}{2}}\,\mathrm{tr}\bigg(\frac{\omega^V}{2}\,e^{\frac{(\omega^V/2)^2}{2\pi i}}\bigg)\,ds\\&=\int_0^1\big(4s(1-s)\big)^{\frac{N^B-1}{2}}\,f\big(\nabla^F,g^F\big)\,ds\;.\end{aligned}$$

Finally, $\mathrm{ch}^{\mathrm{o}}(\nabla^F,g^F)$ is real, odd and closed because $f(\nabla^F,g^F)$ is. □

This motivates the following definitions, where we make the same assumptions as before. In particular, let $(V,A',h^V)\to B$ be a family Thom-Smale complex with an adapted metric $g^V$. Let $p\colon M\to B$ be a smooth bundle with compact fibres, with horizontal subbundle $T^HM$ and vertical metric $g^V$. Let $(F,\nabla^F)\to M$ be a flat complex vector bundle, and let $g^F$, $g^F_0$ and $g^F_1$ be metrics on $F$. Let $(V,A'_s,h^V_s)$ be a family of family Thom-Smale complexes with adapted metrics $g^V_s$ for $s\in[0,1]$, and let $\nabla^{F,s}$ be a family of flat connections on $F\to M$, such that the fibre-wise cohomology $H_s\to B$ has constant rank, with $\nabla^{H,s}$, $g^{H,s}_V$ and $g^{H,s}_{L_2}$ in both cases denoting the respective induced family of metrics and connections on $H_s$, respectively.

**2.91. Definition.** Let $\tilde f(\nabla^F,g^F_0,g^F_1)$, $L_k(A'_s,g^V_0,g^V_1)$, $T_f(A',g^V,h^V)$, and $\mathcal{T}_f(T^HM,g^{TX},\nabla^F,g^F)$ be defined as in Definition 2.16, Definition 2.19, Definition 2.46 and Definition 2.85, and set

$$\widetilde{\mathrm{ch}}^{\mathrm{o}}\big(\nabla^F,g^F_0,g^F_1\big)=\int_0^1\big(4s(1-s)\big)^{\frac{N^B}{2}}\,\tilde f\big(\nabla^F,g^F_0,g^F_1\big)\,ds\qquad\in\Omega^*(M)\;,$$

$$L\big(A'_s,g^V_0,g^V_1\big)=\int_0^1\big(4s(1-s)\big)^{\frac{N^B}{2}}\,L_k\big(A'_s,g^V_0,g^V_1\big)\,ds\;,$$

$$T\big(A',g^V,h^V\big)=\int_0^1\big(4s(1-s)\big)^{\frac{N^B}{2}}\,T_f\big(A',g^V,h^V\big)\,ds\;,$$

and

$$\mathcal{T}\big(T^HM,g^{TX},\nabla^F,g^F\big)=\int_0^1\big(4s(1-s)\big)^{\frac{N^B}{2}}\,\mathcal{T}_f\big(T^HM,g^{TX},\nabla^F,g^F\big)\,ds\qquad\in\Omega^*(B)\;.$$

Let us assume that $(H,\nabla^H)$ is the fibre-wise cohomology of $(V,A')$ (or $(\Omega^*(M/B;F),d^M)$ respectively), with induced metric $g^H_V$ (or $g^H_{L_2}$). Then we have the obvious analogues of Corollary 2.17, Theorem 2.48, Theorem 2.58 and Theorem 2.86.

**2.92. Corollary.** *The forms* $\widetilde{\mathrm{ch}}^{\mathrm{o}}(\nabla^F, g^F, g^{F\prime})$, $L(A'_s, g^V_0, g^V_1)$, $T(A', g^V, h^V)$ *and* $\mathcal{T}(T^H M, g^{TX}, \nabla^F, g^F)$ *are even and real, and satisfy*

$$d^M \widetilde{\mathrm{ch}}^{\mathrm{o}}\big(\nabla^F, g^F_0, g^F_1\big) = \mathrm{ch}^{\mathrm{o}}\big(\nabla^F, g^F_1\big) - \mathrm{ch}^{\mathrm{o}}\big(\nabla^F, g^F_0\big)\;, \tag{1}$$

$$d^B L\big(A'_s, g^V_0, g^V_1\big) = \mathrm{ch}^{\mathrm{o}}\big(A'_1, g^V_1\big)^{[\geq 3]} - \mathrm{ch}^{\mathrm{o}}\big(A'_0, g^V_0\big)^{[\geq 3]}\;, \tag{2}$$

$$d^B T\big(A', g^V, h^V\big) = \mathrm{ch}^{\mathrm{o}}\big(\nabla^V, g^V\big) - \mathrm{ch}^{\mathrm{o}}\big(\nabla^H, g^H_V\big)\;, \tag{3}$$

*and*
$$d^B \mathcal{T}\big(T^H M, g^{TX}, \nabla^F, g^F\big) = \int_{M/B} e\big(\nabla^{TX}\big)\; \mathrm{ch}^{\mathrm{o}}\big(\nabla^F, g^F\big) - \mathrm{ch}^{\mathrm{o}}\big(\nabla^H, g^H_{L_2}\big)\;. \tag{4}$$

*Modulo exact forms on* $B$,

$$\begin{aligned} &T\big(A', g^V_1, h^V\big) - T\big(A', g^V_0, h^V\big) \\ &\qquad = \widetilde{\mathrm{ch}}^{\mathrm{o}}\big(\nabla^V, g^V_0, g^V_1\big) - \widetilde{\mathrm{ch}}^{\mathrm{o}}\big(\nabla^H, g^{H,0}_V, g^{H,1}_V\big)\;, \end{aligned} \tag{5}$$

$$\begin{aligned} &\mathcal{T}\big(T^H_1 M, g^{TX}_1, \nabla^F, g^F_1\big) - \mathcal{T}\big(T^H_0 M, g^{TX}_0, \nabla^F, g^F_0\big) \\ &\qquad = \int_{M/B} \tilde{e}\big(\nabla^{TX,0}, \nabla^{TX,1}\big)\; \mathrm{ch}^{\mathrm{o}}\big(\nabla^F, g^F_1\big) \\ &\qquad\quad + \int_{M/B} e\big(\nabla^{TX,0}\big)\, \widetilde{\mathrm{ch}}^{\mathrm{o}}\big(\nabla^F, g^F_0, g^F_1\big) - \widetilde{\mathrm{ch}}^{\mathrm{o}}\big(\nabla^H, g^{H,0}_{L_2}, g^{H,1}_{L_2}\big)\;, \end{aligned} \tag{6}$$

$$\begin{aligned} &T\big(A'_1, g^V_1, h^V_1\big)^{[\geq 2]} - T\big(A'_0, g^V_0, h^V_0\big)^{[\geq 2]} \\ &\qquad = L\big(\nabla^{V,s}, g^V_0, g^V_1\big) - L\big(\nabla^{H,s}, g^{H,0}_V, g^{H,1}_V\big)\;, \end{aligned} \tag{7}$$

*and*
$$\begin{aligned} &\mathcal{T}\big(T^H M, g^{TX}, \nabla^{F,1}, g^F_1\big)^{[\geq 2]} - \mathcal{T}\big(T^H M, g^{TX}, \nabla^{F,0}, g^F_0\big)^{[\geq 2]} \\ &\qquad = \int_{M/B} e\big(\nabla^{TX}\big) L\big(\nabla^{F,s}, g^F_0, g^F_1\big) - L\big(\nabla^{H,s}, g^{H,0}_{L_2}, g^{H,1}_{L_2}\big)\;. \end{aligned} \tag{8}$$

Theorem 2.48 and Corollary 2.92 (7) imply that $T(A', g^V, h^V)$ is independent of $h^V$ modulo exact forms in all degrees. In [G], we will give a slightly more complicated definition of $T(A', g^V, h^V)$ that does not depend on $h^V$ explicitly, but uses only the splitting $A' = \nabla^V + a'$ induced by $h^V$.

Assume that the metrics $g^F$ on $F \to B$ and $g^H$ on $H \to B$ are parallel. Then Corollary 2.92 implies in particular that we get higher torsion classes

$$\begin{aligned} T(M/B; F, h) &= T\big(A', g^V, h^V\big)^{[\geq 2]} + \widetilde{\mathrm{ch}}^{\mathrm{o}}\big(\nabla^H, g^H, g^H_V\big)^{[\geq 2]} \\ \mathcal{T}(M/B; F) &= \mathcal{T}\big(T^H M, g^{TX}, \nabla^F, g^F\big)^{[\geq 2]} + \widetilde{\mathrm{ch}}^{\mathrm{o}}\big(\nabla^H, g^H, g^H_{L_2}\big)^{[\geq 2]} \end{aligned} \tag{2.93}$$

in $H^{\mathrm{even},\geq 2}_{\mathrm{dR}}(B)$ as in Definition 2.64 and Definition 2.88. Moreover, the classes $T(M/B; F, h)$ and $\mathcal{T}(M/B; F)$ are independent of the particular parallel metrics on $F$ and $H$ and of all other secondary data needed to define the respective analytic torsion forms.

## 3. Higher analytic torsion and Morse functions

We recall our main result, Theorem 0.1, and show that it is compatible with other known results on higher analytic torsion. We also name a few conclusions in [BG1] that do not hold without Smale's transversality condition. Lott used the higher analytic torsion to define a pushdown in a secondary $K$-theory in [L]. In Section 3.c, we apply Theorem 0.1 to simplify this pushdown in the presence of a fibre-wise Morse function. In Section 3.d, we adapt the proof of Theorem 0.1 in [BG1] to the current situation. Section 4 contains a new proof for the only intermediate result of [BG1] that is affected by dropping the Smale condition.

**3.a. Statement of the main result.** Let $p\colon M\to B$ be a smooth fibre bundle with compact fibres, let $(F,\nabla^F)\to M$ be a flat complex vector bundle, and assume that there exists a smooth function $h\colon M\to\mathbb{R}$ such that $h$ is a Morse function on each fibre of $p$.

Let $g^F$ be a Hermitian metric on $F$. Then we have an odd characteristic form $\mathrm{ch}^{\mathrm{o}}(\nabla^F,g^F)$ of Definition 2.89 that vanishes if $g^F$ is parallel with respect to $\nabla^F$.

Let $T^HM\subset TM$ be a horizontal subbundle complementary to the vertical tangent bundle $TX$, and let $g^{TX}$ be vertical metric on $TX$. Let $\nabla^{TX}$ be the natural connection on $TX$ induced by these data.

Let $H=H^*(M/B;F)\to B$ denote the bundle of the fibre-wise cohomology with values in $F$, then $H$ carries a natural Gauß-Manin connection $\nabla^H$, such that $(H,\nabla^H)$ is the $E_1$-term in the Leray spectral sequence in de Rham theory for $p$ with coefficients in $F$. Moreover, fibre-wise Hodge theory gives an $L_2$-metric $g^H$ on $H$ as in [BG1]. Let $\mathrm{ch}^{\mathrm{o}}(\nabla^H,g^H)$ be the corresponding odd characteristic form.

Let $\mathcal{T}(T^HM,g^{TX},\nabla^F,g^F)\in\Omega^*(B)/d\Omega^*(B)$ be the analytic torsion form of Definitions 2.85 and 2.91 in the Chern normalization of [BG1], such that

$$d^B\mathcal{T}\big(T^HM,g^{TX},\nabla^F,g^F\big)=\int_{M_g/B}e\big(\nabla^{TX}\big)\,\mathrm{ch}^{\mathrm{o}}\big(\nabla^F,g^F\big)-\mathrm{ch}^{\mathrm{o}}\big(\nabla^H,g^H\big)\;, \tag{3.1}$$

as explained in Theorem 2.86 and Corollary 2.92.

Let $h\colon M\to\mathbb{R}$ be a smooth function that is Morse on each fibre of $p$, and let $\nabla^{TX}h$ be the fibre-wise gradient of $h$ with respect to $g^{TX}$. Let $C\subset M$ denote the set of fibre-wise critical points of $h$, then $\hat p=p|_C\colon C\to B$ is a finite covering of $B$. We denote the stable and unstable part of the vertical tangent bundle $TX|_C$ of $C$ by $T^sX$, $T^uX\to C$.

Let $V=(p|_C)_*o(T^uX)\otimes F|_C\to B$ be the bundle constructed in Definition 1.1, equipped with a flat superconnection $A'$ constructed as in Theorem 1.61. The component $a_0$ of $A'$ of horizontal degree 0 in $\Omega^*(B)$ defines a cochain complex on $V$. Let $h^V$ be the endomorphism of $V$ induced by $h|_C$. Let $g^V$ denote the Hermitian metric on $V$ induced by $g^F$, and let $g^H_V$ denote the $L_2$-metric on $H$ that is constructed as the restriction of $g^V$ to the $(a_0+a_0^*)$-harmonic elements of $V$. Then we have the Chern-Simons class $\widetilde{\mathrm{ch}}^{\mathrm{o}}(\nabla^H,g^H_{L_2},g^H_V)\in\Omega^*(B)/d^B\Omega^*(B)$ of Definitions 2.16 and 2.91, which satisfies

$$d^B\widetilde{\mathrm{ch}}^{\mathrm{o}}\big(\nabla^H,g^H_{L_2},g^H_V\big)=\mathrm{ch}^{\mathrm{o}}\big(\nabla^H,g^H_V\big)-\mathrm{ch}^{\mathrm{o}}\big(\nabla^H,g^H_{L_2}\big)\;. \tag{3.2}$$

Let $T(A',g^V,h^V)\in\Omega^*(B)/d^B\Omega^*(B)$ be the torsion form of Definitions 2.46 and 2.91, such that

$$d^BT\big(A'g^V,h^V\big)=\mathrm{ch}^{\mathrm{o}}\big(\nabla^V,g^V\big)-\mathrm{ch}^{\mathrm{o}}\big(\nabla^H,g^H_V\big)\quad\in\Omega^*(B)_h\;. \tag{3.3}$$

Let $\delta_0$ be the $\pi^*o(TX)$-valued current on the total space of the vertical tangent bundle $\pi\colon TX\to M$, given by integration over the zero section $M$ of $TX$. Let $\Omega^*(TX)_0$ denote currents that are smooth away from the zero section of $TX$, and whose wave front set is contained in $T^*X|_M\subset T^*TX$. Let $\psi(\nabla^{TX},g^{TX})$ denote the $\pi^*o(TX)$-valued Mathai-Quillen current on $T^*X$ of [MQ], see also [BZ1], such that

$$d^{TX}\psi\big(\nabla^{TX},g^{TX}\big)=\pi^*e\big(\nabla^{TX}\big)-\delta_0\quad\in\Omega^*(TX)_0\;. \tag{3.4}$$

Let $\delta_C$ be the $o(TX)$-valued current on $M$, given by integration over $C$. Then the wave front sets of $\delta_C$ and $(\nabla^{TX}h)^*\psi(\nabla^{TX},g^{TX})$ are contained in $T^*X|_C$, and we have

$$d^M\Big(\big(\nabla^{TX}h\big)^*\psi\big(\nabla^{TX},g^{TX}\big)\Big)=e\big(\nabla^{TX}\big)-(-1)^{\mathrm{ind}_h}\,\delta_C\;. \tag{3.5}$$

Finally, let ${}^0J$ denote the additive genus introduced in [BG1]. Let $\zeta$ denote the Riemann $\zeta$-function, and for a real vector bundle $E \to M$, let $\operatorname{ch}(E)^{[4k]}$ denote the component of the Chern character of $E$ in $H^{4k}(M)$. Then

$$ {}^0J(E) = \frac{1}{2} \sum_{k=1}^{\infty} \zeta'(-2k) \operatorname{ch}(E)^{[4k]} . \tag{3.6} $$

For convenience, we recall that Theorem 0.1 states that

$$ \begin{aligned} &\mathcal{T}\big(T^H M, g^{TX}, \nabla^F, g^F\big) - T\big(A', g^V, h^V\big) - \widetilde{\operatorname{ch}}^{o}\big(\nabla^H, g^H_{L_2}, g^H_V\big) \\ &\quad = -\int_{M/B} \operatorname{ch}^{o}\big(\nabla^F, g^F\big)\, \big(\nabla^{TX} h\big)^* \psi\big(\nabla^{TX}, g^{TX}\big) + \hat p_*\big((-1)^{\operatorname{ind}_h}\, {}^0J(T^sX - T^uX)\big) \operatorname{rk} F . \end{aligned} \tag{3.7} $$

modulo exact forms on $B$. Theorem 0.1 will be proved at the end of this section.

*3.8. Remark.* The formulation of Theorem 0.1 is the same as that of Theorem 0.1 in [BG1]. The only difference is in the definition of the "combinatorial torsion form" $T$.

We discuss an important special case of Theorem 0.1. Assume therefore that there are parallel metrics $g^F$ on $F \to M$ and $g^H$ on $H \to B$, and recall the definition of the higher torsion classes $T(M/B; F, h)$ and $\mathcal{T}(M/B; F)$ in (2.93), cf. Definition 2.64 and Definition 2.88. Theorem 0.1 immediately implies

**3.9. Corollary.** *Assume that the vector bundles $F \to M$ and $H \to B$ carry parallel metrics. Then*

$$ \mathcal{T}(M/B; F) - T(M/B; F, h) = \hat p_*\big((-1)^{\operatorname{ind}_h}\, {}^0J(T^sX - T^uX)\big) \operatorname{rk} F . $$

We will give an application of Corollary 3.9 in Section 5.

**3.b. Compatibility with known results on analytic torsion forms.** We check that Theorem 0.1 is compatible with the transgression and variation formulas for higher analytic torsion forms of [BL] and [BG1], and with Poincaré duality, cf. [BG1], sections 7.4–7.6. On the other hand, we remark that most of the results of sections 7.8–7.10 in [BG1] fail if $h$ does not admit a fibre-wise Morse-Smale gradient field.

*3.10. Remark.* From equations (3.1)–(3.5) it is clear that applying the exterior differential $d^B$ on $B$ to both sides of (3.7) gives the trivial identity

$$ \int_{M/B} e\big(\nabla^{TX}\big) \operatorname{ch}^{o}\big(\nabla^F, g^F\big) - \operatorname{ch}^{o}\big(\nabla^V, g^V\big) = \int_{M/B} \operatorname{ch}^{o}\big(\nabla^F, g^F\big)\, e\big(\nabla^{TX}\big) - \hat p_* \operatorname{ch}^{o}\big(\nabla^F, g^F\big) . $$

Thus, Theorem 0.1 is compatible with the transgression formulas satisfied by $\mathcal{T}$, $T$, $\widetilde{\operatorname{ch}}^{o}$ and $\psi$. This automatically implies that Theorem 0.1 is also compatible with the variation formulas for $\mathcal{T}$, $T$, $\operatorname{ch}^{o}$, $\widetilde{\operatorname{ch}}^{o}$ and $\psi$ under smooth variations of $T^HM$, $g^{TX}$, and $g^F$. In particular, we may choose these data such that the simplifying assumptions of Proposition 1.10 are satisfied when we prove Theorem 0.1.

Let $\nabla^{F,s}$ be a smooth one-parameter family of flat connections on $F$, and let $g^{F,s}$ be a smooth one-parameter family of Hermitian metrics, with $s \in [0,1]$. Let $(H_s, \nabla^{H,s}, g^{H,s}_{L_2})$ be the fibre-wise cohomology with values in $(F, \nabla^{F,s})$ with its induced $L_2$-metric. We assume that $\dim H_s$ is constant, so that $(H_s, \nabla^{H,s})$ becomes a smooth vector bundle on $B \times [0,1]$.

For $k \geq 0$, let $(\,\cdot\,)^{[\geq k]}$ denote the sum of the homogeneous parts of degree $\geq k$ of a differential form. In [BG1], Definition 2.4, we defined the so-called Bott-Chern classes $L(\nabla^{F,s}, g^{F,s}) \in \Omega^{\geq 2}(M)/d\Omega^*(M)$ and $L(\nabla^{H,s}, g^{H,s}_{L_2}) \in \Omega^{\geq 2}(B)/d\Omega^*(B)$, see Definitions 2.19 and 2.91, such that

$$d^M L\big(\nabla^{F,s}, g^{F,s}\big) = \mathrm{ch}^{\mathrm{o}}\big(\nabla^{F,1}, g^{F,1}\big)^{[\geq 3]} - \mathrm{ch}^{\mathrm{o}}\big(\nabla^{F,0}, g^{F,0}\big)^{[\geq 3]}$$
$$\text{and} \qquad d^B L\big(\nabla^{H,s}, g^{H,s}\big) = \mathrm{ch}^{\mathrm{o}}\big(\nabla^{H,1}, g^{H,1}\big)^{[\geq 3]} - \mathrm{ch}^{\mathrm{o}}\big(\nabla^{H,0}, g^{H,0}\big)^{[\geq 3]} , \tag{3.11}$$

see Theorem 2.20 and Corollary 2.92 (2). In [BG1], Theorem 3.45, it was shown that

$$\mathcal{T}\big(T^H M, g^{TX}, \nabla^{F,1}, g^{F,1}\big)^{[\geq 2]} - \mathcal{T}\big(T^H M, g^{TX}, \nabla^{F,0}, g^{F,0}\big)^{[\geq 2]}$$
$$= \int_{M/B} e\big(\nabla^{TX}\big)\, L\big(\nabla^{F,s}, g^{F,s}\big) - L\big(\nabla^{H,s}, g^{H,s}\big) , \tag{3.12}$$

see Theorem 2.86 (3) and Corollary 2.92 (8). Similarly, we have shown in Theorem 2.58 and Corollary 2.92 (7) that

$$T\big(A'^1, g^{V,1}, h^V\big)^{[\geq 2]} - T\big(A'^0, g^{V,0}, h^V\big)^{[\geq 2]} = L\big(\nabla^{V,s}, g^{V,s}\big) - L\big(\nabla^{H,s}, g^{H,s}_V\big) . \tag{3.13}$$

*3.14. Remark.* By equations (3.11)–(3.13), subtracting (3.7) for $\overline{F}_0$ from (3.7) for $\overline{F}_1$ in degree $\geq 2$ after a straightforward partial integration gives the trivial identity

$$\int_{M/B} e\big(\nabla^{TX}\big)\, L\big(\nabla^{F,s}, g^{F,s}\big) - L\big(\nabla^{V,s}, g^{V,s}\big)$$
$$= \int_{M/B} L\big(\nabla^{F,s}, g^{F,s}\big)\, e\big(\nabla^{TX}\big) - \hat{p}_* L\big(\nabla^{F,s}, g^{F,s}\big) .$$

In particular, Theorem 0.1 is compatible with the known variation formulas for analytic torsion forms in dependence of the flat connection on $F$.

Let $o(TX) \to M$ denote the orientation line bundle of $TX$, and let $n$ denote the dimension of the fibres of $p\colon M \to B$. Fibre-wise Poincaré duality implies that there exists a canonical isomorphism

$$\overline{H^k(M/B;F)}^* \cong H^{n-k}\big(M/B; \overline{F}^* \otimes o(TX)\big) .$$

As in [BG1], one verifies that the analytic torsion form satisfies

$$\mathcal{T}\Big(T^H M, g^{TX}, \nabla^{\overline{F}^* \otimes o(TX)}, g^{\overline{F}^* \otimes o(TX)}\Big)$$
$$= (-1)^{n+1}\, \mathcal{T}\big(T^H M, g^{TX}, \nabla^F, g^F\big) \quad \in \quad \Omega^*(B)/d\Omega^*(B) . \tag{3.15}$$

If we replace the fibrewise Morse function $h$ by $-h$ and $F$ by $\overline{F}^* \otimes o(TX)$, then the bundles $V^k \to B$ get replaced by $\overline{V^{n-k}}^* \to B$. The superconnections $A'$ and $A''$ on $V$ get replaced by $\overline{A'}^*$ and $\overline{A''}^*$, which are conjugates of $A''$ and $A'$ by the metric $g^V$. As in [BG1], Definition 2.46 implies that

$$T\Big(\overline{A'}^*, g^{\overline{V}^*}, -h^{\overline{V}^*}\Big) = (-1)^{n+1}\, T\big(A', g^V, h^V\big) \quad \in \quad \Omega^*(B)/d\Omega^*(B) . \tag{3.16}$$

*3.17. Remark.* It has been checked in [BG1] that the remaining terms in (3.7) satisfy equations analogous to (3.15) and (3.16) if one replaces $h$ by $-h$ and $F$ by $\overline{F}^* \otimes o(TX)$. Thus, Theorem 0.1 is compatible with Poincaré duality.

*3.18. Remark.* In [BG1], sections 7.8 and 7.9, we drew some conclusions in the special case that the metric $g^F|_C$ is parallel. We found that then the combinatorial torsion vanishes in higher degrees on $B$, and that the metric $g_V^H$ induced by $g^V$ on $H \to B$ is parallel; both due to the fact that $A'$ takes the special form $a_0 + \nabla^V$ with $a_0$ of degree 0 on $B$ and $\nabla^V$ unitary with respect to $g^V$.

If $\nabla^{TX} h$ does not satisfy Smale's transversality condition, then $A'$ is more complicated, and we will give an example in [G] where $T(A', g^V, h^V)$ carries cohomological information even though $\nabla^V$ is unitarily flat and $H = 0$.

Moreover in our situation, $\nabla^H$ is induced by $\nabla^V + a_1$, which is not necessarily unitary, so we cannot conclude any longer that $g^H$ is a parallel metric. Thus, the arguments in [BG1], sections 7.8 and 7.9 fail in our context.

Note also that in [BG1], Thm 7.11, we deduced that if the fibres are even-dimensional, then the bundle

$$\hat{p}_*\big((-1)^{\operatorname{ind}_h} TX|_C\big)$$

becomes trivial in $KR_0(B) \otimes \mathbb{Q}$. By Theorem 0.3 that will be proved in [G], the above bundle gives the difference between the higher analytic torsion and the higher Franz-Reidemeister torsion

$$\int_{M/B} e(TX)\ {}^0J(TX)\ \operatorname{rk} F = \hat{p}_*\big((-1)^{\operatorname{ind}_h}\ {}^0J(TX|_C)\big)\ \operatorname{rk} F\ ,$$

which may be nonzero.

*3.19. Remark.* In [BG1], section 7.10, we used Smale's transversality condition to show that for every subcovering $C' \subset C \to B$ that leads to an acyclic direct summand $V' \subset V$, we have

$$\big(\hat{p}|_{C'}\big)_*\big((-1)^{\operatorname{ind}_h}\ {}^0J(T^sX - T^uX)\big) = 0\ .$$

We will see in Remark 5.14 that this is not the case in Hatcher's example, which shows that the Morse function in Hatcher's example does not admit a global Smale gradient field, and that Theorem 7.12 in [BG1] becomes false without the transversality condition.

**3.c. Morse functions and Lott's pushdown.** In this subsection, we give a short review of Lott's $\overline{K}_R^0$-groups and describe the pushdown map $p_!\colon \overline{K}_R^0(M) \to \overline{K}_R^0(B)$ associated to a proper submersion $p\colon M \to B$ in the special case that $M$ carries a fibre-wise Morse function. In Section 5.f, we will apply the $\overline{K}_R^0$-functor to Hatcher's example and discuss a conjectural relation of $p_!$ to a Becker-Gottlieb transfer.

Let $R$ be a right-regular, right-Noetherian ring, e.g., $R = \mathbb{Z}$, $R = \mathbb{R}$ or $R = \mathbb{C}$. Let $\rho\colon R \to \operatorname{End}(\mathbb{C}^n)$ be a complex representation of $R$ such that $\mathbb{C}^n$ becomes a flat $R$-module. If $M$ is a manifold and $F \to M$ is a local system of right-$R$-modules, let $F_{\mathbb{C}} = F \otimes_\rho \mathbb{C}^n$ denote the complex flat vector bundle associated to $F$ with flat connection $\nabla^F$.

Lott first defines a Grothendieck group $\widehat{K}_R^0(M)$ with generators $(F, g^F, \eta)$, where $F$ is a local system of right-$R$-modules, $g^F$ is a Hermitian metric on $F_{\mathbb{C}}$, and $\eta \in \Omega^{\text{even}}(M)/d\Omega^{\text{odd}}(M)$ is a class of even differential forms modulo exact forms on $M$.

To define the relations in $\widehat{K}_R^0(M)$, let

$$0 \longrightarrow F_1 \longrightarrow F_2 \longrightarrow F_3 \longrightarrow 0 \tag{3.20}$$

be a short exact sequence, choose metrics $g^{F_i}$ on $F_{i,\mathbb{C}}$ and forms $\eta_i \in \Omega^{\text{even}}(M)/d\Omega^{\text{odd}}(M)$, and let $T$ denote the Bismut-Lott analytic torsion of the finite-dimensional flat family of acyclic complexes

$$0 \longrightarrow (F_{1,\mathbb{C}}, g^{F_1}) \longrightarrow (F_{2,\mathbb{C}}, g^{F_2}) \longrightarrow (F_{3,\mathbb{C}}, g^{F_3}) \longrightarrow 0 \tag{3.21}$$

over $M$. Then in $\widehat{K}^0_R(M)$,

$$[F_2, g^{F_2}, \eta_2] = [F_1, g^{F_1}, \eta_1] \oplus [F_3, g^{F_3}, \eta_3] \qquad \text{if} \qquad \eta_2 = \eta_1 + \eta_3 + T\ . \tag{3.22}$$

There is a well-defined map $\mathrm{ch}^o\colon \widehat{K}^0_R(M) \to \Omega^{\mathrm{odd}}(M)$ given on generators by

$$\mathrm{ch}^o[F, g^F, \eta] = \mathrm{ch}^o(\nabla^F, g^F) - d\eta\ ,$$

and $\overline{K}^0_R(M) \subset \widehat{K}^0_R(M)$ is defined as the kernel of $\mathrm{ch}^o$. Then $\overline{K}^0_R(M)$ fits into an exact sequence

$$H^{\mathrm{even}}(M, \mathbb{R}) \to \overline{K}^0_R(M) \to K^0_R(M) \to H^{\mathrm{odd}}(M, \mathbb{R})\ , \tag{3.23}$$

where $K^0_R(M)$ is simply the Grothendieck group of local systems of right-$R$-modules on $M$, and the rightmost arrow is the characteristic class $\mathrm{ch}^o$. For any smooth map $f\colon M \to N$ there is a natural pullback $f^*\colon \widehat{K}^0_R(N) \to \widehat{K}^0_R(M)$ that restricts to a pullback $f^*\colon \overline{K}^0_R(N) \to \overline{K}^0_R(M)$ and is compatible with (3.23).

Now let $p\colon M \to B$ be a smooth proper submersion, and fix a horizontal subbundle $T^HM$ and a vertical metric $g^{TX}$. Let $H = H^*(M/B; F)$ denote the fibre-wise cohomology of $F$, which is again a local system of right-$R$-modules over $B$, and let $g^H_{L_2}$ denote the $L_2$-metric on $H_{\mathbb{C}} \to B$. Then a pushdown in $\widehat{K}^0_R(M)$ is given on generators by

$$p_![F, g^F, \eta] = \Bigl[H, g^H_{L_2}, \int_{M/B} e(TX, \nabla^{TX})\,\eta - \mathcal{T}(T^HM, g^{TX}, \nabla^F, g^F)\Bigr]\ . \tag{3.24}$$

It is proved in [L] that $p_!$ is compatible with (3.22). Moreover, the pushdown restricts to a map $p_!\colon \overline{K}^0_R(M) \to \overline{K}^0_R(B)$ that is independent of $T^HM$ and $g^{TX}$.

If $h\colon M \to B$ is a fibre-wise Morse-function, we define $C$, $(V, \nabla^V, A', g^V) \to B$ and $g^H_V$ as above. We define another pushdown on $\widehat{K}^0_R(M)$ by

$$p_{!h}[F, g^F, \eta] = \Bigl[H, g^H_V, \hat{p}_*\Bigl((-1)^{\mathrm{ind}_h}\bigl({}^0J(T^sX - T^uX)\,\mathrm{rk}\,F_{\mathbb{C}} + \eta\bigr)\Bigr) - T(A', g^V, h^V)\Bigr]\ . \tag{3.25}$$

As in [L], one can show that this pushdown is compatible with (3.22).

**3.26. Theorem.** *Suppose that $h\colon M \to \mathbb{R}$ is a fibre-wise Morse function. Then $p_{!h}$ restricts to the pushdown $p_!\colon \overline{K}^0_R(M) \to \overline{K}^0_R(B)$.*

*Proof.* The relations (3.22)imply that

$$[F, g'^F, \eta'] = [F, g^F, \eta' - \widetilde{\mathrm{ch}}^o(\nabla^F, g^F, g'^F)] \quad \in \quad \widehat{K}^0_R(M)\ . \tag{3.27}$$

Also, a straightforward partial integration shows that

$$\int_{M/B} e(TX, \nabla^{TX})\,\eta - \hat{p}_*\eta = \int_{M/B} (\nabla^{TX}h)^*\psi(\nabla^{TX}, g^{TX})\,d\eta$$

modulo exact forms on $B$. Using (3.7) and the two observations above, we may rewrite (3.24) as

$$\begin{aligned} p_![F, g^F, \eta] &= \Bigl[H, g^H_V, \hat{p}_*\eta + \int_{M/B} (\nabla^{TX}h)^*\psi(\nabla^{TX}, g^{TX})\,d\eta \\ &\qquad\qquad - \mathcal{T}(T^HM, g^{TX}, \nabla^F, g^F) - \widetilde{\mathrm{ch}}^o(\nabla^H, g^H_{L_2}, g^H_V)\Bigr] \\ &= p_{!h}[F, g^F, \eta] + \Bigl[0, 0, \int_{M/B} (\nabla^{TX}h)^*\psi(\nabla^{TX}, g^{TX})\,\Bigl(d\eta - \mathrm{ch}^o(\nabla^F, g^F)\Bigr)\Bigr]\ . \end{aligned}$$

The theorem follows because in $\overline{K}^0_R(M)$, we have

$$d\eta - \mathrm{ch}^o(\nabla^F, g^F) = -\,\mathrm{ch}^o[F, g^F, \eta] = 0\ . \qquad \square$$

**3.d. A proof of Theorem 0.1.** In order to prove Theorem 0.1, we refer to the proof of Theorem 0.1 of [BG1]. However, we have to replace all arguments that rely on the Smale transversality condition for the negative gradient flow, or on its consequences. We start by summarizing the places in [BG1] that depend on the Smale condition, then we indicate how to extend [BG1] to the more general case of fibre-wise Morse functions that do not necessarily admit a fibre-wise Morse-Smale gradient field, and finally, we show how to fit the contents of this paper into the general framework of [BG1], thereby proving Theorem 0.1.

In Chapter 5 of [BG1], we used the transversality condition to prove a families' version of the result of Laudenbach in the appendix of [BZ1]. In particular, we constructed a vector bundle $V = C^*(W^u; F) \to B$ with a flat superconnection

$$A' = \nabla^V + v \qquad \text{with} \qquad v \in \Omega^0\big(B; \operatorname{End}^1 V\big) \tag{3.28}$$

of total degree one, and an $\Omega^*(B)$-linear quasi-isomorphism

$$\mathbf{P}^\infty\colon\ \Omega^*(M; F) \longrightarrow \Omega^*\big(B; V\big)\ , \tag{3.29}$$

given as integration over the fibres of the unstable cells.

In Chapters 1 and 2 of [BG1], we recalled and extended the construction of torsion forms from Section 2 (f) in [BL]. The torsion forms

$$T_f\big(A', g^V\big) \in \Omega^*(B)/d^B\Omega^*(B) \tag{3.30}$$

are only defined for flat superconnections of the type (3.28), in particular, they converge only if there is no component in $\Omega^{\geq 2}(B; \operatorname{End} V)$. Moreover, the identity

$$dT_f\big(A', g^V\big) = f\big(\nabla^V, g^V\big) - f\big(\nabla^{H^*(V,v)}, g^{H^*(V,v)}\big) \tag{3.31}$$

only holds if the connection part of $A'$ equals $\nabla^V$.

In Theorem 9.8 of [BG1], we proved that up to a few divergences and correction terms, the analytic torsion form of [BL] converges to the form $T_f(A', g^V)$ of (3.30). We gave two proofs of this theorem in Chapters 10 and 11 of [BG1]. While the proof in Chapter 11 needs the transversality condition in its full strength, the proof in Chapter 10 uses only that $A'$ is of the form (3.28), and this assumption is needed only in very few places.

We now indicate the changes that lead to Theorem 0.1. In Section 1 of this paper, we constructed a flat superconnection

$$A' = \nabla^V + \sum_j a_j \qquad \text{with} \qquad a_j \in \Omega^j\big(B; \operatorname{End}^{1-j} V\big)$$

of total degree one on $V$, which replaces the one of (3.28). If we view $h|_C$ as the "height" of a critical point, then the coefficients $a_j$ "move" elements of $V$ "up" by at least $\varepsilon^2 > 0$. We also constructed an $\Omega^*(B)$-linear quasi-isomorphism

$$I\colon\ \Omega^*(M; F) \longrightarrow \Omega^*(B; V)$$

that takes the role of $\mathbf{P}^\infty$ of (3.29). The main properties of $I$ and $A'$ are summarised in Theorem 1.61.

In Section 2, we generalised the torsion forms of [BL], [BG1] to general flat superconnections of total degree one that are compatible with an endomorphism $h^V$ given by the Morse function $h|_C$.

The necessary condition is that $(V, A', h^V)$ is a family Thom-Smale complex in the sense of Definition 2.35 with an adapted metric $g^V$, and the torsion forms $T_f(A', g^V, h^V)$ are constructed in Definition 2.46.

In Section 4, we summarise Chapter 10 of [BG1] and prove our Theorem 4.31, which replaces Theorem 9.8 of [BG1]. We still use the Witten deformation as in [BZ2] and [BG1]. In fact, our proof is a bit more direct because our definition of $T_f(A', g^V, h^V)$ is modelled itself on a "Witten deformation" of the complex $(V, A')$ with metric $g^V$.

Let us recall Theorem 9.8 of [BG1]. The definition of $T_f(A', g^V)$ in [BL], Chapter 2 and [BG1], Definition 1.29 was recalled in (3.30). It implies in particular that

$$\begin{aligned}
T_f\bigl(A', g^V\bigr) - \bigl(\chi'(V) - \chi'(H)\bigr) \int_1^\infty f'\bigl(\sqrt{-t}/2\bigr)\, \frac{dt}{2t}& \\
- \bigl(\chi'(V) - \chi'(H)\bigr) \int_0^1 \Bigl(f'\bigl(\sqrt{-t}/2\bigr) - f'(0)\Bigr)\, \frac{dt}{2t}& \\
= -(2\pi i)^{-\frac{N^B}{2}} \int_1^\infty \Bigl(\operatorname{str}_V\bigl(N^V \gamma\, f'(X_{t,1})\bigr) - \chi'(H) f'(0)\Bigr)\, \frac{dt}{2t}& \\
- (2\pi i)^{-\frac{N^B}{2}} \int_0^1 \Bigl(\operatorname{str}_V\bigl(N^V \gamma\, f'(X_{t,1})\bigr) - \chi'(V) f'(0)\Bigr)\, \frac{dt}{2t}& \;,
\end{aligned}$$

where $f(z) = z\, e^{z^2}$, so in particular, $f'(0) = 1$.

We must also recall the definition of $\tilde\chi'^{\pm}(F)$ in [BG1]. Starting from (9.6) in [BG1], an elementary calculation shows that

$$\begin{aligned}
\tilde\chi'^{+}(F) &= \sum_{c \in C/B} (-1)^{\operatorname{ind}_h(c)} \bigl(n - \operatorname{ind}_h(c)\bigr)\, \operatorname{rk} F \\
&= \operatorname{str}_V\bigl((n - N^V)\bigr) = n\,\chi(V) - \chi'(V)\;, \\
\text{and} \qquad \tilde\chi'^{-}(F) &= \sum_{c \in C/B} (-1)^{\operatorname{ind}_h(c)} \operatorname{ind}_h(c)\, \operatorname{rk} F \\
&= \operatorname{str}_V\bigl(N^V\bigr) = \chi'(V)\;.
\end{aligned}$$

With this preliminaries, we can now reformulate Theorem 9.8 of [BG1] in our notation.

**3.32. Theorem.** *In the limit $T \to \infty$, we have*

$$\begin{aligned}
\tilde f\bigl(\nabla^H, g^H, g^H_T\bigr) + S_f\bigl(\mathbb{A}', g^{TX,F}_T\bigr) - \biggl(\frac{\chi'(V)}{2} - \frac{n\,\chi(V)}{4}\biggr) \log T + T\, \operatorname{str}_V\bigl(h^V \gamma\bigr)& \\
\longrightarrow\; \tilde f\bigl(\nabla^H, g^H, g^H_V\bigr) - \biggl(\frac{\chi'(V)}{2} - \frac{n\,\chi(V)}{4}\biggr) \log \pi + T_f\bigl(A', g^V\bigr)& \\
- \bigl(\chi'(V) - \chi'(H)\bigr) \int_1^\infty f'\bigl(\sqrt{-t}/2\bigr)\, \frac{dt}{2t}& \\
- \bigl(\chi'(V) - \chi'(H)\bigr) \int_0^1 \Bigl(f'\bigl(\sqrt{-t}/2\bigr) - 1\Bigr)\, \frac{dt}{2t}&
\end{aligned}$$

*modulo exact forms on $B$.* □

In Section 4, we prove Theorem 4.31, which equals the theorem above, except that we do no longer require the existence of a Smale gradient field (which is used implicitly in [BG1], Theorem 9.8), and that $T_f(A', g^V)$ gets replaced by $T_f(A', g^V, h^V)$.

*Proof of Theorem 0.1.* By Remark 3.10, we may assume that the simplifying assumptions of Proposition 1.10 are satisfied. By the definitions and arguments of Section 2.i, we may replace the Chern normalization in Theorem 0.1 by the normalisation $f(z) = z\, e^{z^2}$ of (2.11). If we define

$$ {}^0I_f(E) = \frac{1}{2} \sum_{k=1}^{\infty} \frac{(4k+1)!}{2^{4k}\,(2k)!^2}\, \zeta'(-2k)\ \mathrm{ch}(E)^{[4k]} \tag{3.33} $$

in analogy with (3.6), then an elementary calculation shows that

$$ {}^0J(E) = \int_0^1 \big(4s(1-s)\big)^{\frac{N^B}{2}}\ {}^0I_f(E)\, ds\ . $$

It remains to establish

$$ \begin{aligned} &\mathcal{T}_f\big(T^HM, g^{TX}, \nabla^F, g^F\big) - T_f\big(A', g^V, h^V\big) - \tilde f\big(\nabla^H, g^H_{L_2}, g^H_V\big) \\ &\quad = -\int_{M/B} f\big(\nabla^F, g^F\big)\, \big(\nabla^{TX} h\big)^* \psi\big(\nabla^{TX}, g^{TX}\big) + \hat p_*\big((-1)^{\mathrm{ind}_h}\ {}^0I_f(T^sX - T^uX)\big)\ \mathrm{rk}\, F\ . \end{aligned} \tag{3.34} $$

As indicated above, this equation now follows from the proof given in Chapter 9 of [BG1], where Theorem 4.31 replaces Theorem 9.8. □

## 4. The Witten deformation

We consider the Witten deformation of the fibre-wise de Rham complex with parameter $T > 0$. Proceeding as in [BG1], we show that for $T$ sufficiently large, the de Rham complex splits naturally into two sub-complexes. One of these sub-complexes represents all eigenvalues that grow at least linearly in $T$ (this will be made precise), so that it does not contribute to the torsion integral in the limit. The other one can be identified with the sections of a finite-dimensional vector bundle $V \to B$, generalising the families' Thom-Smale complex of [BG1]. While in [BG1], the identification was provided by the integration over the unstable bundles of the fibre-wise gradient flow, we now use the integration map $I\colon \Omega^*(M;F) \to \Omega^*(B;V)$ constructed in Section 1.

In Sections 4.a–4.c, we recall several results of [BG1] related to the model operator near the fibre-wise critical points, the space spanned by the generalised eigenforms with small eigenvalues, and the behaviour of the torsion integrand for large values of $t$, always using $I$ instead of the de Rham map employed in [BG1]. In Section 4.d we show how the definition of $T(A', g^V, h^V)$ via finite-dimensional Witten deformation on $\Omega^*(B;V)$ fits into the picture of the infinite-dimensional Witten deformation on $\Omega^*(M;F)$, thereby establishing Theorem 4.31, which was used in Section 3.d to prove Theorem 0.1.

**4.a. The model operator near the critical points.** We recall constructions from [BG1], sections 4 and 10 related to a model operator $\widehat{\mathbb{X}}_T$ acting on forms on the total space of a $\mathbb{Z}_2$-graded Euclidean vector bundle $E = E^+ \oplus E^-$ over a smooth manifold $C$, with values in a flat vector bundle $F \to E$ that is pulled back from $C$. This will be applied as follows. Assume that $p\colon M \to B$ carries a fibre-wise Morse function $h\colon M \to \mathbb{R}$ with fibre-wise critical set $C$, and there exists a vertical metric $g^{TX}$ on the vertical tangent bundle $TX \to M$, such that $h$ and $g^{TX}$ satisfy the simplifying assumptions of Proposition 1.10. Then the restriction of the vertical tangent bundle to $C$ splits orthogonally as $TX|_C = T^sX \oplus T^uX$. The Witten deformed operator $\widetilde{\mathbb{X}}_T$ acting on $\Omega^*(M;F)$ coincides with $\widehat{\mathbb{X}}_T$ on $E = TX|_C$ in a tubular neighbourhood of $C$. It has been shown in [BG1] that for $T \to \infty$, the operator $e^{\widehat{\mathbb{X}}_T^2}$ provides us with a good approximation for the heat kernel $e^{\widetilde{\mathbb{X}}_T^2}$, in a sense to be made precise later.

Let $\pi\colon (E, g^E) = (E^+, g^{E^+}) \oplus (E^-, g^{E^-}) \to C$ be a $\mathbb{Z}_2$-graded Euclidean vector bundle of dimension $n = n^+ + n^-$. Let $\nabla^E = \nabla^{E^+} \oplus \nabla^{E^-}$ be a Euclidean connection that respects the grading, with curvature $F^E$. If we identify $\pi^* E \to E$ with the vertical tangent bundle of $E \to C$, we can define a vertical radial vector field on the total space $E$ by

$$\mathcal{R} = \mathcal{R}^+ + \mathcal{R}^- \qquad \text{with} \qquad \mathcal{R}^\pm(p, e^+, e^-) = e^\pm \in \pi^* E \subset TE ,$$

for all $p \in C$ and all $(e^+, e^-) \in E_p$.

The connection $\nabla^E$ induces a bundle connection on $TE$, i.e., a projection $TE \to \pi^* E$ with kernel $T^H E$, and $g^E$ becomes a vertical metric. These data induce a connection on the vertical tangent bundle as in [B], Theorem 1.1, which coincides with the pull back $\nabla^{\pi^* E}$ of the connection $\nabla^E$. The fibre bundle curvature of $T^H E \subset TE$ is given by $\Omega^E = F^E \mathcal{R}$, see (1.9).

Let $\Omega^*(E/C) \to C$ denote the bundle of differential forms on the fibres of $E$, equipped with the $L_2$-metric

$$g^{\Omega^*(E/C)}(s, s')_p = \int_{E_p} g^E\bigl(s(e), s'(e)\bigr)\, d\operatorname{Vol}_E(e) ,$$

where $d\operatorname{Vol}_E$ is the volume element induced by the vertical metric $g^E$. Then $\nabla^E$ induces a natural connection $\nabla^{\Omega^*(E/C)}$ on $\Omega^*(E/C)$, and $\nabla^{\Omega^*(E/C)}$ is compatible with the $L_2$-metric $g^{\Omega^*(E/C)}$.

If we identify $E$ with $E^*$, then $E$ acts on $\Omega^*(E/C)$ by exterior multiplication $\varepsilon_e$, interior multiplication $\iota_e$, and two types of Clifford multiplication

$$c_e = \varepsilon_e - \iota_e \qquad \text{and} \qquad \hat{c}_e = \varepsilon_e + \iota_e ,$$

satisfying

$$[c_e, c_{e'}] = -g^E(e, e') , \qquad [c_e, \hat{c}_{e'}] = 0 , \qquad \text{and} \qquad [\hat{c}_e, \hat{c}_{e'}] = g^E(e, e') ,$$

where we use $[A, B] = AB - (-1)^{\deg A \cdot \deg B} BA$ to denote the supercommutator of $A$ and $B$. Fixing a local orthonormal basis $e_1, \dots, e_n$ of $E$, we abbreviate $\varepsilon_j = \varepsilon_{e_j}$, $\iota_j = \iota_{e_j}$ etc.

Finally, let $(F, \nabla^F) \to C$ be a flat complex vector bundle with a Hermitian metric $g^F$ that is not necessarily parallel. Let $\nabla^{F,*}$ be the adjoint of $\nabla^F$ with respect to $g^F$, and set

$$\omega^F = \nabla^{F,*} - \nabla^F = \bigl(g^F\bigr)^{-1}\bigl[\nabla^F, g^F\bigr] , \qquad \text{and} \qquad \nabla^{F,\mathrm{u}} = \frac{1}{2}\bigl(\nabla^{F,*} + \nabla^F\bigr) = \nabla^F + \frac{1}{2}\omega^F .$$

Then $\nabla^{F,\mathrm{u}}$ is unitary with respect to $g^F$. By an abuse of notation, we write $(F, \nabla^F) \to E$ and $g^F$ for the pull-backs $(\pi^* F, \pi^* \nabla^F)$ and $\pi^* g^F$. Then $g^E$ and $g^F$ induce an $L_2$-metric

$$g^{\Omega^*(E/C;F)}(s, s')_p = \int_{E_p} \bigl(g^E \otimes g^F\bigr)\bigl(s(e), s'(e)\bigr)\, d\operatorname{Vol}_E(e)$$

on the bundle $\Omega^*(E/C;F) = \Omega^*(E/C) \otimes F \to C$. The connections $\nabla^F$, $\nabla^{F,*}$ and $\nabla^{F,\mathrm{u}}$ induce connections on $\Omega^*(E/C;F)$ denoted by

$$\nabla^{\Omega^*(E/C;F)} , \qquad \nabla^{\Omega^*(E/C;F),*} = \nabla^{\Omega^*(E/C;F)} + \omega^F ,$$
$$\text{and} \qquad \nabla^{\Omega^*(E/C;F),\mathrm{u}} = \frac{1}{2}\bigl(\nabla^{\Omega^*(E/C;F),*} + \nabla^{\Omega^*(E/C;F)}\bigr) .$$

We are now ready to construct two superconnections on the bundle $\Omega^*(E/C;F) \to C$ by

$$\widehat{\mathbb{A}}'_0 = d^E = \varepsilon_i\,\nabla^F_{e_i} + \nabla^{\Omega^*(E/C;F)} + \iota_{F^E\mathcal{R}}\ , \tag{4.1}$$

$$\text{and}\qquad \widehat{\mathbb{A}}''_0 = -\iota_i\,\nabla^F_{e_i} + \nabla^{\Omega^*(E/C;F),*} - \varepsilon_{F^E\mathcal{R}}\ .$$

Then $\widehat{\mathbb{A}}''$ is the adjoint of $\widehat{\mathbb{A}}'$ with respect to $g^{\Omega^*(E/C;F)}$.

Now suppose that $h^C\colon C \to \mathbb{R}$ is any smooth function, then the function $h^E\colon E \to \mathbb{R}$ given by

$$h^E = h^C \circ \pi + \frac{1}{2}\left(\left|\mathcal{R}^+\right|^2 - \left|\mathcal{R}^-\right|^2\right)$$

is a fibre-wise Morse function. The critical set is the zero section $C \subset E$ of $E \to C$, and its index is $n^- = \dim E^-$. We extend the superconnections $\widehat{\mathbb{A}}'_0$ and $\widehat{\mathbb{A}}''_0$ to superconnection on the pullback $\Omega^*(\widehat{E}/\widehat{C};F)$ of $\Omega^*(E/C;F)$ to $\widehat{C} = C \times (0,\infty)$ by setting

$$\widehat{\mathbb{A}}' = e^{-Th^E}\left(\widehat{\mathbb{A}}'_0 + \frac{\partial}{\partial T}\,dT\right)e^{Th^E}\ ,\qquad \text{and}\qquad \widehat{\mathbb{A}}'' = e^{Th^E}\left(\widehat{\mathbb{A}}''_0 + \frac{\partial}{\partial T}\,dT\right)e^{-Th^E}\ ,$$

where $T$ is the standard coordinate on $(0,\infty)$. Then $\widehat{\mathbb{A}}''$ is the adjoint of $\widehat{\mathbb{A}}''$ with respect to the pullback of the $L_2$-metric $g^{\Omega^*(E/C;F)}$. We denote the restriction of any of these connections to the submanifold $C \times \{T\}$ by adding a subscript $(\,\cdot\,)_T$.

Define an operator $\widehat{\mathbb{X}}$ on $\Omega^*(\widehat{E}/\widehat{C};F)$ by

$$\begin{aligned}\widehat{\mathbb{X}} = \frac{1}{2}\left(\widehat{\mathbb{A}}'' - \widehat{\mathbb{A}}'\right) &= -\frac{1}{2}\big(\hat c_j\,\partial_j + T\,c_{\mathcal{R}^+-\mathcal{R}^-}\big) + \frac{1}{2}\,\omega^F - \frac{1}{2}\,\hat c_{F^E\mathcal{R}}\\ &\qquad - T\,dh^E - \frac{\left|\mathcal{R}^s\right|^2 - \left|\mathcal{R}^u\right|^2}{2}\,dT - h^C\,dT = \widehat{\mathbb{X}}_T - h^E\,dT\ .\end{aligned}$$

The square of this operator is given by

$$\begin{aligned}\widehat{\mathbb{X}}^2 &= -\frac{1}{4}\Big(\Delta^E + T^2\left|\mathcal{R}\right|^2 - \big\langle F^E\mathcal{R}, F^E\mathcal{R}\big\rangle\Big) + \frac{1}{2}\nabla_{F^E\mathcal{R}}\\ &\qquad - \frac{T}{2}\left(N^+ - N^- + \frac{n^- - n^+}{2}\right) + \frac{1}{4}\langle F^E e_j, e_k\rangle \hat c_j \hat c_k\\ &\qquad + \frac{1}{4}\big(\omega^F\big)^2 + \frac{1}{2}\,\hat c_{\mathcal{R}^+-\mathcal{R}^-}\,dT\\ &= \widehat{\mathbb{X}}_T^2 + \frac{1}{2}\,\hat c_{\mathcal{R}^+-\mathcal{R}^-}\,dT\ ,\end{aligned}\tag{4.2}$$

where $\Delta^E$ denotes the fibre-wise Laplacian on $E$, and $N^+$, $N^-$ are number operators that act on $\Lambda^pE^+ \otimes \Lambda^qE^- \subset \Lambda^*E$ by $p$ and $q$ respectively.

The model operator $\widehat{\mathbb{X}}$ is invariant under an $\mathbb{R}^*_+$-action that we now define. Let $o(E^-) = \lambda^{\max}\widehat{E}^- \to \widehat{C}$ be the orientation bundle of $\widehat{E}^- \to \widehat{C}$. Let $N^{dT}$ denote the number operator on $\Lambda^*T(0,\infty) = \mathbb{R} \oplus \mathbb{R}\,dT$.

**4.3. Definition.** For $\lambda > 0$, set

(1) $\quad r_\lambda(\hat p) = \big(p, T/\lambda\big)$ $\quad$ for $\hat p = (p,T) \in \widehat{C}$,

(2) $\quad r_\lambda^{o(E^-)}(o) = o \quad \in o(E^-)_{r_\lambda \hat p}$ $\quad$ for $o \in o(E^-)_{\hat p}$,

(3) $\quad r_\lambda^F(f) = f \quad \in F_{r_\lambda \hat p}$ $\quad$ for $f \in F_{\hat p}$,

(4) $\quad r_\lambda^{\Lambda^*T\widehat{C}}(\alpha) = \lambda^{\frac{N^{\widehat{C}}}{2}}\,\alpha \quad \in \Lambda^*T_{r_\lambda \hat p}\widehat{C}$ $\quad$ for $\alpha \in \Lambda^*T_{\hat p}\widehat{C}$,

(5) $\quad \big(r_\lambda^{\Omega^*(\widehat{E}/\widehat{C})}(\omega)\big)\big(r_\lambda \hat p, e\big) = \lambda^{-\frac{n}{4}}\,\omega\big(\hat p, e/\sqrt{\lambda}\big)$ $\quad$ for $\omega \in \Omega^*\big(\widehat{E}/\widehat{C}\big)_{\hat p}$,

where we agree that the rescaling of the coordinates $T$ and $e$ affects $dT$, but not the vertical forms. For $\omega \in \Omega^*\bigl(\widehat{E};F\bigr) = \Omega^*\bigl(\widehat{C};\Omega^*(\widehat{E}/\widehat{C};F)\bigr)$, set

$$(6)\qquad (r_\lambda^{\Omega^*(\widehat{E};F)}\omega)(\hat{p},e) = \Bigl(r_\lambda^{\Lambda^*T\widehat{C}\otimes\Omega^*(\widehat{E}/\widehat{C};F)}\,\omega_{r_\lambda^{-1}\hat{p}}\Bigr)(e) = \lambda^{\frac{N^{\widehat{C}}}{2}-\frac{n}{4}}\,\omega\Bigl(p,\lambda T,e/\sqrt{\lambda}\Bigr)\;.$$

Note that $r_\lambda^{o(E^-)}$, $r_\lambda^F$, and $r_\lambda^{\Omega^*(\widehat{E}/\widehat{C})}$ are isometries with respect to the pull-backs to $\widehat{C}$ of the bundle metrics $g^{o(E^-)}$, $g^F$ and $g^{\Omega^*(E/C;F)}$. When no confusion can arise, we will drop the superscript and write simply $r_\lambda$.

**4.4. Proposition.** *The operator $\widehat{\mathbb{X}}^2$ satisfies*

$$\widehat{\mathbb{X}}^2 = \lambda\, r_\lambda\, \widehat{\mathbb{X}}^2\, r_\lambda^{-1}\;.$$

*Proof.* This is clear from (4.2) and Definition 4.3. Note that $dT$ is multiplied by $\lambda^{\frac{3}{2}}$ by Definition 4.3 (1) and (4). □

Let

$$\widehat{\mathbb{X}}_1^{[0]} = -\frac{1}{2}\left(\hat{c}_j\nabla^F_{e_j} + c_{\mathcal{R}^+-\mathcal{R}^-}\right)$$

be the homogeneous part of horizontal degree 0 with respect to $\Omega^*\bigl(\widehat{C}\bigr)$. Then

$$-\bigl(\widehat{\mathbb{X}}_1^{[0]}\bigr)^2 = \frac{1}{4}\Bigl(\Delta^E + |\mathcal{R}|^2\Bigr) + \frac{N^+-N}{2} - \frac{n^+-n^-}{4}$$

is a fibre-wise harmonic oscillator with spectrum $\bigl\{0,\frac{1}{2},1,\dots\bigr\}$. Let $o(E^-) = \Lambda^{\max}E^-$ denote the orientation bundle of $E^-$. Then the map

$$(4.5)\qquad \widehat{J}_1^{[0]}:\ o(E^-)\otimes F \longrightarrow \ker\Bigl(\widehat{\mathbb{X}}_1^{[0]}\Bigr)^2 \subset \Omega^*\bigl(E/C;F\bigr)$$

$$\text{with}\qquad \widehat{J}_1^{[0]}(v) = \pi^{-\frac{n}{4}}\,e^{-\frac{|\cdot|^2}{2}}\,\pi^*v$$

is an isometry with respect to $g^F$ and $g^{\Omega^*(E/C;F)}$.

Because $\widehat{\mathbb{X}}_1^2 - (\widehat{\mathbb{X}}_1^{[0]})^2$ is nilpotent, the operators $\widehat{\mathbb{X}}_1^{[0]}$ and $\widehat{\mathbb{X}}_1$ have the same spectrum. By Proposition 4.4, the operator $\widehat{\mathbb{X}}^2$ has spectrum $\operatorname{Spec}(\widehat{\mathbb{X}}^2)_{C\times\{T\}} = \bigl\{0,\frac{T}{2},T,\dots\bigr\}$. Let $\delta\subset\mathbb{C}$ be circle of radius 1 around 0 and consider the map

$$\begin{aligned}\widehat{\mathbb{P}} &= \frac{1}{2\pi i}\int_{\frac{T}{4}\delta}\frac{d\lambda}{\lambda-\widehat{\mathbb{X}}^2}\\ &= \widehat{\mathbb{P}}_T + \widehat{\mathbb{Q}}_T\,dT \quad \in \Omega^*\Bigl(\widehat{C};\operatorname{End}\Omega^*\bigl(\widehat{E}/\widehat{C};F\bigr)\Bigr)\;.\end{aligned}$$

We define another map $\widehat{J}\colon \Lambda^*T\widehat{C}\otimes o(E^-)\otimes F\to\Lambda^*T\widehat{C}\otimes\Omega^*\bigl(\widehat{E}/\widehat{C};F\bigr)$ by

$$(4.6)\qquad \widehat{J} = \widehat{\mathbb{P}}\circ\biggl(\Bigl(r_T^{\Lambda^*T\widehat{C}\otimes\Omega^*(\widehat{E}/\widehat{C};F)}\Bigr)^{-1}\circ J_1^{[0]}\circ r_T^{\Lambda^*T\widehat{C}\otimes o(E^-)\otimes F}\biggr)\;.$$

Then we may rewrite Theorem 10.53 in [BG1] as follows, where for any vector bundle $V\to\widehat{C}$, we set

$$\Omega^*\bigl(\widehat{C}/(0,\infty);V\bigr) = \bigl\{\,\omega\in\Omega^*\bigl(\widehat{C};V\bigr)\bigm|\iota_{\frac{\partial}{\partial T}}\omega = 0\,\bigr\}\;.$$

**4.7. Theorem.** *The map $\widehat{\mathbb{P}}$ is a spectral projection that commutes with $\widehat{\mathbb{A}}'$ and $\widehat{\mathbb{A}}''$, is orthogonal with respect to $g^{\Omega^*(E/C;F)}$, and we have*

$$\widehat{\mathbb{P}} \in \Omega^{\mathrm{even}}\Big(\widehat{C}/(0,\infty); \mathrm{End}^{\mathrm{even}}\,\Omega^*\big(\widehat{E}/\widehat{C}\big) \oplus \mathrm{End}^{\mathrm{odd}}\,\Omega^*\big(\widehat{E}/\widehat{C}\big)\,dT\Big)\;. \tag{1}$$

*For any positive integer $k$,*

$$\widehat{\mathbb{P}} = \frac{1}{2\pi i}\int_{(\frac{T}{4})^{\frac{k}{2}}\delta} \frac{d\lambda}{\lambda - \widehat{\mathbb{X}}^k}\;. \tag{2}$$

*Let $\widehat{\mathbb{F}}$ be the range of $\widehat{\mathbb{P}}$, then $\widehat{\mathbb{F}}$ is a finite dimensional bundle over $\widehat{C}$ which is preserved by the action of $r_\lambda$, and $\widehat{\mathbb{F}} = \ker\big(\widehat{\mathbb{X}}^k\big)$ for $k$ sufficiently large.*

*The map $\widehat{J}$ maps $\Lambda^*T\widehat{C} \otimes o(E^-) \otimes F$ isomorphically to $\widehat{\mathbb{F}}$, and $\widehat{J}^{[0]}\colon o(E^-) \otimes F \to \widehat{\mathbb{F}}^{[0]} = \widehat{\mathbb{F}} \cap \Omega^*\big(\widehat{E}/\widehat{C}; F\big)$ is an isometry. We may write*

$$\widehat{J}(\alpha \otimes o \otimes f) = \alpha \otimes \left(\left(\frac{T}{\pi}\right)^{\frac{n}{4}} e^{-\frac{T|\mathcal{R}|^2}{2}}\,\big(1 + R_{\widehat{E}}\big)\big(\pi^* o\big)\right) \otimes \pi^* f \tag{3}$$

*for $\alpha \in \Lambda^{\mathrm{even}}T\widehat{C}$, $o \in o(E^-)$, and $f \in F$, where*

$$R_{\widehat{E}} \in \Omega^{\mathrm{even}}\Big(\widehat{C}/(0,\infty); \mathrm{End}^{\mathrm{even}}\,\Omega^*\big(\widehat{E}/\widehat{C}\big) \oplus \mathrm{End}^{\mathrm{odd}}\,\Omega^*\big(\widehat{E}/\widehat{C}\big)\,dT\Big) \tag{4}$$

*is an $r_\lambda$-invariant polynomial in $\mathcal{R}^\pm$ with coefficients depending only on $F^{E^\pm}$ and $T$.* □

*Proof.* It follows from Proposition 4.4 that $\widehat{\mathbb{P}}$ and $\widehat{\mathbb{F}}$ are $r_\lambda$-invariant. The other assertions on $\widehat{\mathbb{P}}$ and $\widehat{\mathbb{F}}$ follow as in Theorem 10.53 in [BG1]. The fact that the only term in (4.2) containing $dT$ is

$$\frac{1}{2}\,\hat{c}_{\mathcal{R}^+ - \mathcal{R}^-}\,dT \in \mathrm{End}^{\mathrm{odd}}\,\Omega^*(\widehat{E}/\widehat{C})\,dT$$

accounts for the form of (1) compared to the corresponding statement in Theorem 10.53 in [BG1]. From [BZ2], we know that $\widehat{J}^{[0]}$ takes the form

$$\widehat{J}^{[0]}(v) = \left(\frac{T}{\pi}\right)^{\frac{n}{4}} e^{-\frac{T|\mathcal{R}|^2}{2}}\,\pi^*(v)\;.$$

Then (3) follows immediately from (1), and $R_{\widehat{E}}$ is a polynomial by (10.191) in [BG1], which is $r_\lambda$-invariant since $\widehat{\mathbb{F}}$ and $\widehat{\mathbb{P}}$ are $r_\lambda$-invariant, and of the form (4) by (1). □

Let $\widehat{J}^*$ be the adjoint of $\widehat{J}$ with respect to $g^{o(E^-)\otimes F}$ and $g^{\Omega^*(\widehat{E}/\widehat{C};F)}$ in the sense of Definition 2.7. Integration over the fibres of $\widehat{E}^- \to \widehat{C}$ defines a map

$$\widehat{I}\colon\ \Omega^*\big(\widehat{E}; F\big) \longrightarrow \Omega^*\big(\widehat{C}; o(E^-) \otimes F\big) \qquad \text{with} \qquad \widehat{I}(\alpha) = \int_{\widehat{E}^-/\widehat{C}} \alpha\;.$$

Then we have the following version of [BG1], Proposition 10.58 and Theorem 10.59 for our model situation.

**4.8. Proposition.** *There exist $r_\lambda$-invariant characteristic differential forms $R_{J^*J}$ and $R_{IJ} \in \Omega^{\mathrm{even},\geq 2}\big(\widehat{C}/(0,\infty)\big)$ of $R^{E^\pm}$, such that*

$$\widehat{J}^* \circ \widehat{J} = 1 + R_{J^*J} \tag{1}$$

$$\text{and} \qquad \widehat{I} \circ e^{Th^E} \circ \widehat{J} = e^{Th^C} \left(\frac{T}{\pi}\right)^{\frac{n^+-n^-}{4}} \big(1 + R_{IJ}\big) \; . \tag{2}$$

*In particular, $R_{J^*J}$, $R_{IJ}$ are polynomials in $\frac{1}{T}$ without constant term.*

Note in particular that the right hand sides of both formulas are closed differential forms that do not contain any multiple of $dT$.

*Proof.* These claims follow by elementary calculations using $r_\lambda$-invariance, see also [BG1], Section 10.15. Since $R_{\widehat{E}}$ is $r_\lambda$-invariant,

$$R_{\widehat{E}}|_{(p,T,e)} = T^{-\frac{N^{\widehat{C}}}{2}} R_{\widehat{E}}|_{(p,1,\sqrt{T}\,e)} \; .$$

We will use the symbol $|\cdot|^2$ in the sense of (2.2). Let $v \in (o(E^-) \otimes F)_{p,T}$, then

$$\begin{aligned}\big(\widehat{J}^* \circ \widehat{J}\big)(v) &= \int_{E_p} \left(\frac{T}{\pi}\right)^{\frac{n}{2}} e^{-T|e|^2} \,\big|\big(1 + R_{\widehat{E}}|_{(p,T,e)}\big)(\pi^* v)\big|^2 \, d\,\mathrm{Vol}^E(e) \\ &= T^{-\frac{N^{\widehat{C}}}{2}} \int_{E_p} \pi^{-\frac{n}{2}} e^{-|e|^2} \,\big|\big(1 + R_{\widehat{E}}|_{(p,1,e)}\big)(\pi^* v)\big|^2 \, d\,\mathrm{Vol}^E(e) \; ,\end{aligned}$$

which proves (1). Note that the part of $R_{\widehat{E}}$ that contains $dT$ does not contribute because it is odd as an endomorphism of $\Omega^*\big(\widehat{E}/\widehat{C}\big)$.

Similarly, we get (2) because

$$\begin{aligned}\big(\widehat{I} \circ e^{Th^E} \circ \widehat{J}\big) &= e^{Th^C} \left(\frac{T}{\pi}\right)^{\frac{n}{4}} \\ &\qquad \int_{\widehat{E}^-/\widehat{C}|_p} e^{T(h^E - h^C)(e^-)} \, e^{-\frac{T|e^-|^2}{2}} \big(1 + R_{\widehat{E}}|_{(p,T,e^-)}\big)\big(\pi^* v\big) \\ &= e^{Th^C} \left(\frac{T}{\pi}\right)^{\frac{n^+-n^-}{4}} \pi^{-\frac{n}{2}} \\ &\qquad T^{-\frac{N^{\widehat{C}}}{2}} \int_{\widehat{E}^-/\widehat{C}|_p} e^{-|e^-|^2} \big(1 + R_{\widehat{E}}|_{(p,1,e^-)}\big)\big(\pi^* v\big) \; .\end{aligned}$$

Note that if $v = o \otimes f \in o(E^-) \otimes F$, then $\pi^* o$ is a volume form on the fibres of $\widehat{E}^-$, so only the component of $R_{\widehat{E}}$ in $\Omega^*\big(\widehat{C}; \big\{\mathrm{id}_{\Omega^*(\widehat{E}/\widehat{C};F)}\big\}\big)$ contributes, and this component again does not contain $dT$ since it is even.

Then both $R_{J^*J}$ and $R_{IJ}$ are $r_\lambda$-invariant, even forms of degree $\geq 2$, so they are polynomials in $\frac{1}{T}$ without constant terms. $\square$

**4.b. Small eigenvalues and spectral projections.** Let

$$\overline{M} = M \times (0,\infty) \times \mathbb{R} \;\longrightarrow\; \overline{B} = B \times (0,\infty) \times \mathbb{R}$$

be the pull-back of $M \to B$ to $\overline{B}$. By an abuse of notation, let $TX \to \overline{M}$ and $F \to \overline{M}$ denote the pull-backs of the vector bundles $TX \to M$ and $F \to M$ to $\overline{M}$. Let $g^{TX}$ and $g^F$ denote the pull-back metrics on these bundles, then these metrics induce a fibre-wise $L_2$-metric $g^{TX,F}$ on the infinite dimensional vector bundle $\Omega^*(M/B;F) \to B$. For $t \in (0,\infty)$, $T \in \mathbb{R}$ consider the metrics

$$g_t^{TX} = \frac{1}{t}\, g^{TX} \qquad \text{and} \qquad g_T^F = e^{-2Th}\, g^F \;.$$

These metrics induce a fibre-wise $L_2$-metric $\bar{g}^{TX,F}$ on $\Omega^*\bigl(\overline{M}/\overline{B};F\bigr)$ such that

$$\bar{g}^{TX,F} = g^{TX,F}\left(t^{N^X - \frac{n}{2}}\, e^{-2Th}\,\cdot\,,\,\cdot\right) .$$

Let $\mathbb{A}' = d^M$ denote the total derivative on $\Omega^*(M;F)$. We extend $\mathbb{A}'$ to a superconnection $\overline{\mathbb{A}}'$ on $\Omega^*\bigl(\overline{M};F\bigr)$ by

$$\overline{\mathbb{A}}' = \mathbb{A}' + \frac{\partial}{\partial t}\,dt + \frac{\partial}{\partial T}\,dT \;.$$

Let $\overline{\mathbb{A}}''$ be the adjoint superconnection of $\overline{\mathbb{A}}'$ with respect to $\bar{g}^{TX,F}$. Then we have

$$\mathbb{A}'' = t^{-N^X}\, e^{2Th}\left(\mathbb{A}'' + \frac{\partial}{\partial t}\,dt + \frac{\partial}{\partial T}\,dT\right) e^{-2Th}\, t^{N^X} \tag{4.9}$$

Then we define a self-adjoint superconnection $\overline{\mathbb{A}}$ on $\Omega^*(M/B;F)$ and a skew-adjoint operator $\overline{\mathbb{X}} \in \Omega^*(B;\operatorname{End}\Omega(M/B);F)$—both with respect to $\bar{g}^{TX,F}$—by

$$\overline{\mathbb{A}} = \frac{1}{2}\left(\overline{\mathbb{A}}'' + \overline{\mathbb{A}}'\right) \qquad \text{and} \qquad \overline{\mathbb{X}} = \frac{1}{2}\left(\overline{\mathbb{A}}'' - \overline{\mathbb{A}}'\right) .$$

We will denote the restrictions of $\overline{\mathbb{A}}'$, $\overline{\mathbb{A}}''$ and $\bar{g}^{TX,F}$ to $M \times \{t\} \times \{T\}$ by $\mathbb{A}'_{t,T}$, $\mathbb{A}''_{t,T}$ and $g_{t,T}^{TX,F}$, respectively. We have

$$\overline{\mathbb{X}} = \mathbb{X}_{t,T} + N^X\,\frac{dt}{2t} - h\,dT \;.$$

If $\mathbb{A}_{t,T}^{[0]}$ and $\mathbb{X}_{t,T}^{[0]} \in \Omega^0(B;\operatorname{End}\Omega^*(M/B;F))$ denote the homogeneous part of horizontal degree 0, then $\overline{\mathbb{A}}^2|_{t,T} - (\mathbb{A}_{t,T}^{[0]})^2$ and $\overline{\mathbb{X}}^2|_{t,T} - (\mathbb{A}_{t,T}^{[0]})^2$ are elements of $\Omega^{>0}(B;\operatorname{End}\Omega^*(M/B;F))$ and thus nilpotent in $\Omega^*(B;\operatorname{End}\Omega^*(M/B;F))$. This implies that

$$\begin{aligned}\operatorname{Spec}\overline{\mathbb{A}}^2|_{t,T} &= \operatorname{Spec}\mathbb{A}_{t,T}^2 = \operatorname{Spec}\bigl(\mathbb{A}_{t,T}^{[0]}\bigr)^2\\ &= -\operatorname{Spec}\overline{\mathbb{X}}^2|_{t,T} = -\operatorname{Spec}\mathbb{X}_{t,T}^2 = -\operatorname{Spec}\bigl(\mathbb{X}_{t,T}^{[0]}\bigr)^2 \;.\end{aligned}$$

Recall that a vector bundle $V \to B$ was defined in Definition 1.1 by

$$V^j = \hat{p}_*^j F|_{C^j} \otimes o(N^u) \;,$$

$$\text{and} \qquad V = \hat{p}_* F|_C \otimes o(N^u) = \bigoplus_j V^j \;,$$

and that we have constructed a map $I\colon \Omega^*(M;F) \to \Omega^*(B;V)$ in Definition 1.60, which is an $\Omega^*(B)$-linear combination of maps that arise by integration over the fibres of certain subbundles of $M \to B$. Let $I^{[0]}\colon \Omega^*(M/B;F) \to V$ denote the homogeneous part of horizontal degree 0. Then we have the following result of [BZ1], Theorems 7.8 and 7.9, and [BZ2], Theorem 6.11, cf. [BG1], Theorems 10.2 and 10.59.

**4.10. Theorem.** *There exists $T_0 \geq 0$ such that for all $T \geq T_0$,*

$$\operatorname{Spec}\bigl(\mathbb{A}_{t,T}^2\bigr) = -\operatorname{Spec}\bigl(\mathbb{X}_{t,T}^2\bigr) \subset \left[0, \frac{t}{4}\right] \cup (4t, \infty) \;,$$

*and the direct sums of eigenspaces associated to eigenvalues in $\bigl[0, \frac{t}{4}\bigr]$ form a vector bundle $W \to B \times (0,\infty) \times [T_0, \infty) \subset \overline{B}$, such that*

$$I^{[0]}|_W \colon W \longrightarrow V$$

*is a linear isomorphism.*

*Proof.* We can copy the proof of [BZ2], Theorem 6.11. However, since $I^{[0]}$ is not simply integration over the unstable cell of the gradient flow, we have to use Theorem 1.61. ☐

Let $S^1(t) \subset \mathbb{C}$ denote a circle around the origin of radius $t$, and note that if $\frac{t}{4} < t_0$, $t_1 < 4t$, then

$$\frac{1}{2\pi i} \int_{S^1(t_0)} \frac{d\lambda}{\lambda - \overline{\mathbb{A}}^2} = \frac{1}{2\pi i} \int_{S^1(t_1)} \frac{d\lambda}{\lambda - \overline{\mathbb{A}}^2} \;.$$

This allows us to patch spectral projectors that are defined on $B \times \bigl(\frac{t_0}{4}, 4t_0\bigr) \times [T_0, \infty)$ together as in the following definition, cf. Definition 10.33 in [BG1].

**4.11. Definition.** Define an endomorphism $\overline{\mathbb{P}}$ of the bundle $\Lambda^* T(\overline{B}) \otimes \Omega^*\bigl(\overline{M}/\overline{B}; F\bigr)$ for $T \geq T_0$ in a neighbourhood of $B \times \{(t_0, T)\}$ in $B \times (0, \infty) \times [T_0, \infty)$ by

$$\overline{\mathbb{P}} = \frac{1}{2\pi i} \int_{S^1(t_0)} \frac{d\lambda}{\lambda - \overline{\mathbb{A}}^2} \;.$$

We recall parts of Theorem 10.9 from [BG1].

**4.12. Theorem.** *For $T \geq T_0$, the operator $\overline{\mathbb{P}} \in \Omega^*\bigl(\overline{B}; \operatorname{End}\Omega^*(M/B; F)\bigr)$ is an even self-adjoint projector with range a finite dimensional vector bundle $\overline{\mathbb{W}} \to B \times (0,\infty) \times [T_0, \infty)$, which commutes with the action of $\Lambda^* T(\overline{B})$ and with $\overline{\mathbb{A}}'$ and $\overline{\mathbb{A}}''$. For any integer $k > 0$ and $t \in \bigl(\frac{t_0}{4}, 4t_0\bigr)$, we have*

$$\overline{\mathbb{P}} = \frac{1}{2\pi i} \int_{S^1\left(t_0^{\frac{k}{2}}\right)} \frac{d\lambda}{\lambda - \overline{\mathbb{X}}^k} \;.$$

Let $\mathbf{g}_T^{\mathbb{W}}$ denote the restriction of the $L_2$-metric $g_{1,T}^{TX,F}$ on $\Omega^*\bigl(\overline{M}/\overline{B}; F\bigr)$ to $\mathbb{W}_T = \overline{\mathbb{W}}|_{B \times \{(1,T)\}}$, which is a generalized metric on $\mathbb{W}$ in the sense of Definition 2.65. Similarly, let $g_T^H$ denote the $L_2$-metric on $H$ induced by $g_{1,T}^{TX,T}$.

**4.13. Theorem.** *There exists a constant $\delta \in (0,1)$ such that for $T \geq T_0$,*

$$\int_1^\infty \Bigl( \operatorname{str}\Bigl(N^X f'\bigl(\mathbb{X}_{t,T}\bigr)\Bigr) - \operatorname{str}\Bigl(N^X f'\bigl(\mathbb{P}_{t,T}\, \mathbb{X}_{t,T}\, \mathbb{P}_{t,T}\bigr)\Bigr)\Bigr) \frac{dt}{2t} = O\bigl(T^{-\delta}\bigr) \;, \tag{1}$$

*and for all $T_1 \geq T_0$,*

$$\begin{aligned} &\tilde f\bigl(\nabla^H, g^H, g_{T_1}^H\bigr) - (2\pi i)^{-\frac{N^B}{2}} \int_1^\infty \Bigl( \operatorname{str}\Bigl(N^X f'\bigl(\mathbb{X}_{t,T_1}\bigr)\Bigr) - \chi'(H)\, f'(0)\Bigr) \frac{dt}{2t} \\ &\quad = \tilde f\bigl(\nabla^H, g^H, g_{T_0}^H\bigr) - (2\pi i)^{-\frac{N^B}{2}} \int_1^\infty \Bigl( \operatorname{str}\Bigl(N^X f'\bigl(\mathbb{P}_{t,T_0}\, \mathbb{X}_{t,T_0}\, \mathbb{P}_{t,T_0}\bigr)\Bigr) - \chi'(H)\, f'(0)\Bigr) \frac{dt}{2t} \\ &\qquad + \tilde f\bigl(\mathbb{A}'_{1,\,\cdot}, \mathbf{g}_{T_0}^{\mathbb{W}}, \mathbf{g}_{T_1}^{\mathbb{W}}\bigr) + O\bigl(T_1^{-\delta}\bigr) \end{aligned} \tag{2}$$

*modulo an exact form on $B$.*

*Proof.* Equation (1) is a consequence of Theorem 4.10, cf. [BG1], Theorem 10.5.

We recall that

$$\overline{\mathbb{X}}_T = \mathbb{X}_{t,T} + \big(g_{t,T}^{TX,F}\big)^{-1} \frac{\partial g_{t,T}^{TX,F}}{\partial t}\, dt = \mathbb{X}_{t,T} + \frac{N^X}{2t}\, dt\,.$$

Let $\overline{\mathbb{P}}_T = \mathbb{P}_{t,T} + \mathbb{P}'_{t,T}\, dt$ and note that $\overline{\mathbb{P}}_T^2 = \overline{\mathbb{P}}_T$ implies

$$\mathbb{P}'_{t,T} = \mathbb{P}_{t,T}\, \mathbb{P}'_{t,T} + \mathbb{P}'_{t,T}\, \mathbb{P}_{t,T}\ ,$$

so $\mathbb{P}'_{t,T}$ interchanges the images of $\mathbb{P}_{t,T}$ and $\mathrm{id} - \mathbb{P}_{t,T}$. This implies that

$$\begin{aligned}(4.14)\qquad \operatorname{str}\Big(f\big(\overline{\mathbb{P}}_T\, \overline{\mathbb{X}}_T\, \overline{\mathbb{P}}_T\big)\Big) &= \operatorname{str}\Big(f\big(\mathbb{P}_{t,T}\, \mathbb{X}_{t,T}\, \mathbb{P}_{t,T}\big)\Big)\\ &\quad + \operatorname{str}\Bigg(\Big(-\mathbb{P}'_{t,T}\, \mathbb{X}_{t,T}\, \mathbb{P}_{t,T} + \frac{N^X}{2t} + \mathbb{P}_{t,T}\, \mathbb{X}_{t,T}\, \mathbb{P}'_{t,T}\Big)\\ &\qquad\qquad f'\big(\mathbb{P}_{t,T}\, \mathbb{X}_{t,T}\, \mathbb{P}_{t,T}\big)\Bigg)\, dt\\ &= \operatorname{str}\Big(f\big(\mathbb{P}_{t,T}\, \mathbb{X}_{t,T}\, \mathbb{P}_{t,T}\big)\Big) + \operatorname{str}\Big(N^X\, f'\big(\mathbb{P}_{t,T}\, \mathbb{X}_{t,T}\, \mathbb{P}_{t,T}\big)\Big)\, \frac{dt}{2t}\,.\end{aligned}$$

A similar argument with fixed $t$ and variable $T$ shows that

$$(4.15)\qquad \operatorname{str}\Big(f\big(\overline{\mathbb{P}}_t\, \overline{\mathbb{X}}_t\, \overline{\mathbb{P}}_t\big)\Big) = \operatorname{str}\Big(f\big(\mathbb{P}_{t,T}\, \mathbb{X}_{t,T}\, \mathbb{P}_{t,T}\big)\Big) - \operatorname{str}\Big(h\, f'\big(\mathbb{P}_{t,T}\, \mathbb{X}_{t,T}\, \mathbb{P}_{t,T}\big)\Big)\, dT\,.$$

Now consider the closed differential form on $B \times (0,\infty) \times \{T\}$, which is given by

$$\begin{aligned}\operatorname{str}\Big(f\big(\overline{\mathbb{P}}\, \overline{\mathbb{X}}\, \overline{\mathbb{P}}\big)\Big) &- \chi'(H)\, f'(0)\, \frac{dt}{t}\\ &= \operatorname{str}\Big(f\big(\mathbb{P}_{t,T}\, \mathbb{X}_{t,T}\, \mathbb{P}_{t,T}\big)\Big) + \Big(\operatorname{str}\Big(N^X\, f'\big(\mathbb{P}_{t,T}\, \mathbb{X}_{t,T}\, \mathbb{P}_{t,T}\big)\Big) - \chi'(H)\, f'(0)\Big)\, \frac{dt}{2t}\\ &\qquad - \operatorname{str}\Big(h\, f'\big(\mathbb{P}_{t,T}\, \mathbb{X}_{t,T}\, \mathbb{P}_{t,T}\big)\Big)\, dT + R\, dt\, dT\ ,\end{aligned}$$

with $R \in \Omega^*(B \times (0,\infty) \times [T_0,\infty))$.

Let $\bar{g}^H$ be the metric on the pull-back of $H$ to $B \times [T_0,\infty)$, let $\overline{\nabla}^{H,*}$ denote the adjoint of the pull-back $\overline{\nabla}^H$ of $\overline{\nabla}^H$ with respect to $\bar{g}^H$, and let

$$\overline{\omega}^H = \overline{\nabla}^{H,*} - \overline{\nabla}^H\,.$$

Then we know that as $t \to \infty$,

$$\operatorname{str}\Big(f\big(\overline{\mathbb{P}}\, \overline{\mathbb{X}}\, \overline{\mathbb{P}}\big)\Big) - \chi'(H)\, f'(0)\, \frac{dt}{t} = \operatorname{str}\Big(f\big(\overline{\omega}^H/2\big)\Big) + O\big(t^{-\frac{1}{2}}\big) + O\big(t^{-\frac{3}{2}}\big)\, dt\ ,$$

cf. [BL], Theorems 3.16 and 3.21. Now, equation (2) follows easily. $\square$

**4.c. The limit $t \to \infty$ and the class $S_f(A', g^V)$.** We identify the bundle $\mathbb{W} \to B$ with $\Lambda^*TB \otimes V$. This allows us to compare the term

$$-(2\pi i)^{-\frac{N^B}{2}} \int_1^\infty \Bigl(\operatorname{str}\Bigl(N^X\, f'\bigl(\mathbb{P}_{t,T_0}\, \mathbb{X}_{t,T_0}\, \mathbb{P}_{t,T_0}\bigr)\Bigr) - \chi'(H)\, f'(0)\Bigr) \frac{dt}{2t}$$

of Theorem 4.13 (2) with a class $S_f(A', \mathbf{g}^V, g^H)$ of the bundle $V \to B$.

Let $I\colon \Omega^*(M;F) \to \Omega^*(B;V)$ be the integration map of Definition 1.60, then we have an obvious extension

$$\overline{I}\colon\, \Omega^*\bigl(\overline{M};F\bigr) \;\longrightarrow\; \Omega^*\bigl(\overline{B};V\bigr)\;,$$

for which the analogue of Theorem 1.61 still holds. By Theorem 4.10, $I^{[0]}$ provides us with an isomorphism of the bundle $W_{t,T} \to B$ of the small eigenspaces with the bundle $V$ of Definition 1.1. By $\Omega^*(B)$-linearity, this implies that

$$\overline{\mathbb{I}} = \overline{I}|_{\overline{\mathbb{W}}}\colon\; \overline{\mathbb{W}} \;\longrightarrow\; \Lambda^*T\overline{B} \otimes V \tag{4.16}$$

is an isomorphism as well. Define

$$\bar{\mathbf{g}}^V_{t,T} = \bigl(\overline{\mathbb{I}}^{-1}\bigr)^*\bigl(\mathbf{g}^{\overline{\mathbb{W}}}\bigr) = g^{TX,F}\Bigl(\overline{\mathbb{I}}^{-1}\,\cdot\,, \overline{\mathbb{I}}^{-1}\,\cdot\,\Bigr) \tag{4.17}$$

Let $N^M$ denote the total number operator on $\Omega^*(M;F)$. Recall that we used the subscript $(\,\cdot\,)_T$ as an abbreviation for $(\,\cdot\,)_{1,T}$. We collect a few properties of the maps $\overline{\mathbb{P}}$ and $\overline{\mathbb{I}}$ and of $\bar{\mathbf{g}}^V$.

**4.18. Proposition.** *For $T \geq T_0$, the map $\mathbb{P}_{t,T}$ satisfies*

$$\mathbb{P}_{t,T} = t^{-\frac{N^M}{2}}\, \mathbb{P}_T\, t^{\frac{N^M}{2}}\;. \tag{1}$$

*Moreover, $\mathbb{P}_{t,T}$ is even in $\Omega^*(B;\operatorname{End}\Omega^*(M/B;F))$ and polynomial in $\frac{1}{t}$, and $\mathbb{P}^{[0]}_{t,T} = \mathbb{P}^{[0]}_T$ independent of $t$.*

*The inverse of $\mathbb{I}_{t,T}$ is even in $\Omega^*\bigl(B;\operatorname{Hom}(V,\Omega^*(M/B;F))\bigr)$ and polynomial in $\frac{1}{t}$ with $(\mathbb{I}^{-1}_{t,T})^{[0]} = (I^{[0]}|_{W_T})^{-1}$ independent of $t$, and satisfies*

$$\mathbb{I}^{-1}_{t,T} = t^{-\frac{N^M}{2}}\, \mathbb{I}^{-1}_T\, t^{\frac{N^B+N^V}{2}}\;. \tag{2}$$

*The map $\mathbf{g}^V_{t,T}$ is a generalized metric on $V$ that depends polynomially on $\frac{1}{t}$. If we set $g^V_T = (\mathbf{g}^V_T)^{[0]}$, then there exists an even $g'_T \in \Omega^*(B;\operatorname{End}V)$ such that (2.70) holds for each $T \geq T_0$.*

*Let $A''_{t,T}$ denote the adjoint of $A'$ with respect to $\mathbf{g}^V_{t,T}$ and set $X_{t,T} = \frac{1}{2}\,(A''_{t,T} - A')$, then*

$$\mathbb{I}_{t,T}\, \mathbb{X}_{t,T} = X_{t,T}\, \mathbb{I}_{t,T}\;. \tag{3}$$

*Proof.* Since $\mathbb{A}'$ is of total degree 1 and $\mathbb{A}''$ is of total degree $-1$, we have

$$t^{\frac{1-N^M}{2}}\, \mathbb{A}'\, t^{\frac{N^M}{2}} = \mathbb{A}' \qquad \text{and} \qquad t^{\frac{1-N^M}{2}}\, \mathbb{A}''_T\, t^{\frac{N^M}{2}} = \mathbb{A}''_{t,T}\;.$$

By Definition 4.11, this implies (1).

Next, we note that $\mathbb{A}'_{t,T}$ raises the total degree by 1, while $\mathbb{A}''^{[k]}_{t,T}$ raises it by $2k-1$. This implies that

$$\mathbb{A}^2_{t,T} = \mathbb{A}'_{t,T}\mathbb{A}''_{t,T} + \mathbb{A}''_{t,T}\mathbb{A}'_{t,T}$$

does not lower the total degree. Moreover

$$\big(\mathbb{A}^{[0]}_{t,T}\big)^2 = \mathbb{A}'^{[0]}_{t,T}\mathbb{A}''^{[0]}_{t,T} + \mathbb{A}''^{[0]}_{t,T}\mathbb{A}'^{[0]}_{t,T} = t\,\big(\mathbb{A}^{[0]}_{T}\big)^2$$

preserves the total degree, and $\mathbb{A}^2_{t,T} - (\mathbb{A}^{[0]}_{t,T})^2$ is nilpotent. In particular,

$$\begin{aligned}\mathbb{P}_{t,T} &= \frac{1}{2\pi i}\int_{S^1(t)} \frac{1}{\lambda - (\mathbb{A}^{[0]}_{t,T})^2}\sum_{k\geq 0}\bigg(\Big(\mathbb{A}^2_{t,T} - (\mathbb{A}^{[0]}_{t,T})^2\Big)\frac{1}{\lambda - (\mathbb{A}^{[0]}_{t,T})^2}\bigg)^k d\lambda \\ &= \frac{1}{2\pi i}\int_{S^1(1)} \frac{1}{\lambda - (\mathbb{A}^{[0]}_{T})^2}\sum_{k\geq 0}\bigg(\Big(\frac{1}{t}\,\mathbb{A}^2_{t,T} - (\mathbb{A}^{[0]}_{T})^2\Big)\frac{1}{\lambda - (\mathbb{A}^{[0]}_{T})^2}\bigg)^k d\lambda\;,\end{aligned}$$

where the summation index $k$ is bounded by $\dim B$. Since $\frac{1}{t}\,\mathbb{A}^2_{t,T} - (\mathbb{A}^{[0]}_{T})^2$ is even, it is polynomial in $\frac{1}{t}$ by (1), which proves the remaining statements on $\mathbb{P}_{t,T}$.

Taking the inverse of $I^{[0]}|_{W_T}$, and noting that $\mathbb{I}_{t,T} = I\,\mathbb{P}_{t,T}$ and $\big(I\,\mathbb{P}_{t,T}\,(I^{[0]}|_{W_T})^{-1}\big)^{[0]} = \mathrm{id}$, we may write

$$\mathbb{I}^{-1}_{t,T} = \mathbb{P}_{t,T}\,\big(I^{[0]}|_{W_T}\big)^{-1}\Big(I\,\mathbb{P}_{t,T}\,\big(I^{[0]}|_{W_T}\big)^{-1}\Big)^{-1}\;.$$

This explains why $\mathbb{I}_{t,T}$ is even and polynomial in $\frac{1}{t}$. Claim (2) follows from the analogous statement for $\mathbb{P}_{t,T}$, because $I^{[0]}$ preserves the total degree.

Since $(I^{[0]}|_W)^{-1}\colon V \to \Omega(M/B;F)$ is clearly injective and preserves the total degree, we know that $(\mathbf{g}^V_{t,T})^{[0]} = \big((I^{[0]}|_W)^{-1}\big)^* g^{TX,F}_{t,T}$ is a Hermitian metric on $V$ that respects the $\mathbb{Z}$-grading. Moreover, $\mathbf{g}^V_{t,T}$ is even because $\mathbb{I}^{-1}_{t,T}$ is. Finally, Definition 2.65 (3), (4) are clearly satisfied by (2.2), so $\mathbf{g}^V_{t,T}$ is a generalised metric. It is also clear that $\mathbf{g}^V_{t,T}$ is polynomial in $\frac{1}{t}$.

In order to check (2.70), note that

$$\begin{aligned}\mathbf{g}^V_{t,T} &= t^{-\frac{n}{2}}\,g^{TX,F}_T\Big(t^{\frac{N^X}{2}}\,\mathbb{I}^{-1}_{t,T}\,\cdot\,,t^{\frac{N^X}{2}}\,\mathbb{I}^{-1}_{t,T}\,\cdot\Big)\\ &= t^{-N^B-\frac{n}{2}}\,g^{TX,F}_T\Big(\mathbb{I}^{-1}_T\,t^{\frac{N^B+N^V}{2}}\,\cdot\,,\mathbb{I}^{-1}_T\,t^{\frac{N^B+N^V}{2}}\,\cdot\Big)\\ &= t^{-N^B-\frac{n}{2}}\,\mathbf{g}^V_T\Big(t^{\frac{N^B+N^V}{2}}\,\cdot\,,t^{\frac{N^B+N^V}{2}}\,\cdot\Big)\;.\end{aligned}$$

Since $\mathbb{I}^{-1}_T$ is $\Omega^*(B)$-linear, it cannot lower horizontal degrees, which proves that $t^{\frac{n-N^V}{2}}\,\mathbf{g}^V_{t,T}\,t^{-\frac{N^V}{2}}$ is a polynomial in $\frac{1}{\sqrt{t}}$ with leading term $g^V_T$. This implies the first equation of (2.70). It also implies that

$$\begin{aligned}t^{\frac{N^V}{2}}\,\big(\mathbf{g}^V_{t,T}\big)^{-1}\,\frac{\partial\mathbf{g}^V_{t,T}}{\partial t}\,t^{-\frac{N^V}{2}} &= \Big(t^{\frac{n-N^V}{2}}\,\mathbf{g}^V_{t,T}\,t^{-\frac{N^V}{2}}\Big)^{-1}\,\frac{\partial\big(t^{\frac{n-N^V}{2}}\,\mathbf{g}^V_{t,T}\,t^{-\frac{N^V}{2}}\big)}{\partial t}\\ &\quad+\frac{N^V}{2t} + \Big(t^{\frac{n-N^V}{2}}\,\mathbf{g}^V_{t,T}\,t^{-\frac{N^V}{2}}\Big)^{-1}\,\frac{N^V-n}{2t}\,\Big(t^{\frac{n-N^V}{2}}\,\mathbf{g}^V_{t,T}\,t^{-\frac{N^V}{2}}\Big)\end{aligned}$$

is a polynomial in $\frac{1}{\sqrt{t}}$ with leading term $\frac{N^V}{t} - \frac{n}{2t}$, which implies the second equation of (2.70).

Finally, the adjoint of $A'$ with respect to $\mathbf{g}^V_{t,T}$ can be written as

$$\begin{aligned} A''_{t,T} &= \Bigl(\overline{\mathbb{I}^{-1}_{t,T}}^{*}\, g^{TX,F}_{t,T}\, \mathbb{I}^{-1}_{t,T}\Bigr)^{-1} \overline{A'}^{*} \Bigl(\overline{\mathbb{I}^{-1}_{t,T}}^{*}\, g^{TX,F}_{t,T}\, \mathbb{I}^{-1}_{t,T}\Bigr) \\ &= \mathbb{I}_{t,T}\, \bigl(g^{TX,F}_{t,T}\bigr)^{-1}\, \overline{\mathbb{A}'}^{*}\, g^{TX,F}_{t,T}\, \mathbb{I}^{-1}_{t,T} \\ &= \mathbb{I}_{t,T}\, \mathbb{A}''_{t,T}\, \mathbb{I}^{-1}_{t,T}\ , \end{aligned}$$

which implies (3). $\square$

Let $\overline{A}'$, $\overline{A}''_T$ and $\overline{X}_T$ denote the obvious extensions of $A'$, $A''_{t,T}$ and $X_{t,T}$ to the pullback of $V$ to $B \times (0,\infty) \times \{T\}$ . Note that

$$\operatorname{str}\Bigl(f'\bigl(\mathbb{P}_{t,T}\, \mathbb{X}_{t,T}\, \mathbb{P}_{t,T}\bigr)\Bigr) = \chi(V)\, f'(0)$$

because $f'$ is even and $\mathbb{P}_{t,T}\, \mathbb{X}_{t,T}\, \mathbb{P}_{t,T}$ is odd. With this observed, we can now calculate

$$\begin{aligned} &-(2\pi i)^{-\frac{N^B}{2}} \int_1^\infty \Bigl(\operatorname{str}\Bigl(N^X\, f'\bigl(\mathbb{P}_{t,T}\, \mathbb{X}_{t,T}\, \mathbb{P}_{t,T}\bigr)\Bigr) - \chi'(H)\, f'(0)\Bigr) \frac{dt}{2t} \\ &\quad= -(2\pi i)^{-\frac{N^B}{2}} \int_1^\infty \Biggl(\operatorname{str}\Bigl(\Bigl(N^X - \frac{n}{2}\Bigr) f'\bigl(\mathbb{P}_{t,T}\, \mathbb{X}_{t,T}\, \mathbb{P}_{t,T}\bigr)\Bigr) \\ &\qquad\qquad - \Bigl(\chi'(H) - \frac{n}{2}\,\chi(V)\Bigr) f'(0)\Biggr) \frac{dt}{2t} \\ &\quad= -(2\pi i)^{-\frac{N^B}{2}} \int_{(B\times(1,\infty)\times\{T\})/B\times\{T\}} \Bigl(\operatorname{str}\Bigl(f\bigl(\overline{\mathbb{P}}_T\, \overline{\mathbb{X}}_T\, \overline{\mathbb{P}}_T\bigr)\Bigr) \\ &\qquad\qquad - \Bigl(\chi'(H) - \frac{n}{2}\,\chi(V)\Bigr) f'(0)\, \frac{dt}{2t}\Bigr) \\ &\quad= -(2\pi i)^{-\frac{N^B}{2}} \int_{(B\times(1,\infty)\times\{T\})/B\times\{T\}} \Bigl(\operatorname{str}\Bigl(f\bigl(\overline{X}_T\bigr)\Bigr) \\ &\qquad\qquad - \Bigl(\chi'(H) - \frac{n}{2}\,\chi(V)\Bigr) f'(0)\, \frac{dt}{2t}\Bigr) \\ &\quad= -(2\pi i)^{-\frac{N^B}{2}} \int_1^\infty \Biggl(\operatorname{str}_V\Bigl(\bigl(\mathbf{g}^V_{t,T}\bigr)^{-1} \frac{\partial \mathbf{g}^V_{t,T}}{\partial t}\, f'\bigl(X_{t,T}\bigr)\Bigr) \\ &\qquad\qquad - \Bigl(\chi'(H) - \frac{n}{2}\,\chi(V)\Bigr) f'(0)\Biggr) \frac{dt}{2t} \\ &\quad= S_f\bigl(A', \mathbf{g}^V_T, g^H_T\bigr)\ . \end{aligned}$$

Using Theorem 2.76 and (4.15), we can now rewrite the equation in Theorem 4.13 (2) as

$$\begin{aligned} &\tilde f\bigl(\nabla^H, g^H, g^H_T\bigr) + S_f\bigl(T^HM, g^{TX}, \nabla^F, g^F_T\bigr) \\ &\quad= \tilde f\bigl(\nabla^H, g^H, g^H_{T_0}\bigr) + S_f\bigl(A', \mathbf{g}^V_{T_0}, g^H_{T_0}\bigr) + \tilde f\bigl(A', \mathbf{g}^V_{T_0}, \mathbf{g}^V_T\bigr) + O\bigl(T^{-\delta}\bigr) \\ &\quad= \tilde f\bigl(\nabla^H, g^H, g^H_V\bigr) + S_f\bigl(A', g^V\bigr) + \tilde f\bigl(A', g^V, \mathbf{g}^V_T\bigr) + O\bigl(T^{-\delta}\bigr) \qquad \in \Omega^*(B)/d^B\Omega^*(B)\ . \end{aligned} \tag{4.19}$$

**4.d. The limit $T \to \infty$ and the class $T_f(A', g^V, h^V)$.** We state some estimates for the generalised metric $\mathbf{g}_T^V$ as $T \to \infty$ that were established in [BG1], and we use these estimates to identify the right hand side of (4.19) with the torsion form $T_f(A', g^V, h^V)$, up to divergences and correction terms.

Recall that we have defined an $\Omega^*(C)$-linear map $\widehat{J}\colon \Omega^*(C; o(E^-) \otimes F) \to \Omega^*\big(\widehat{E}; F\big)$ in (4.6) that maps $\Lambda^* T\widehat{C} \otimes o(E^-) \otimes F$ isomorphically onto the bundle $\mathbb{F} \to \widehat{C}$ of the fibre-wise generalised eigenspaces of the model operator $\widehat{\mathbb{X}}$ for the generalised eigenvalue 0. We will apply this definition to $E = TX|_C = T^sX \oplus T^uX$. Restriction to $T \geq T_0$ and push-down to $B \times [T_0, \infty)$ gives a map

$$\overline{J}\colon\ \Omega^*\big(B \times [T_0, \infty), V\big) \longrightarrow \Omega^*\big(TX|_C \times [T_0, \infty), \pi^* F\big)\ . \tag{4.20}$$

Recall that by the simplifying assumptions of Proposition 1.10, the fibre-wise exponential map at $C \subset M$ maps the $2\varepsilon$-disc bundle of $TX|_C$ fibre-wise isometrically onto a tubular neighbourhood of $C$ in $M$. Let us fix a $C^\infty$ cut-off function $\varphi\colon \mathbb{R} \to [0,1]$ with

$$\varphi(r) = \begin{cases} 1 & \text{if } r \leq \varepsilon, \text{ and} \\ 0 & \text{if } r \geq 2\varepsilon. \end{cases}$$

Then we may regard $\varphi(|\mathcal{R}|)\, \overline{J}$ as a map

$$\varphi(|\mathcal{R}|)\, \overline{J}\colon\ \Omega^*\big(B \times [T_0, \infty), V\big) \longrightarrow \Omega^*\big(M \times [T_0, \infty), F\big)\ . \tag{4.21}$$

This allows to define

$$\mathbf{e}_T\colon\ \Lambda^* TB \otimes V \longrightarrow \mathbb{W}_T \qquad \text{by} \qquad \mathbf{e}_T = \mathbb{P}_T\left(e^{Th}\, \varphi(|\mathcal{R}|)\, J_T\right)\ . \tag{4.22}$$

We recall a fundamental result from [BG1].

**4.23. Theorem** ([BG1], Theorem 10.56, cf. [BZ2], Theorem 6.7). *There exists a constant $c > 0$ such that as $T \to \infty$ for any $k \geq 0$ and any $v \in \Omega^*(B; V)$,*

$$\big(\mathbf{e}_T - e^{Th}\, \varphi(|\mathcal{R}|)\, J_T\big)(v) = e^{Th}\, O\big(e^{-cT}\ |v|\big)$$

*uniformly in $C^k$-topology over compact subsets of $B$. In particular, $\mathbf{e}_T$ is an $\Omega^*(B)$-linear vector bundle isomorphism for sufficiently large $T$.* □

This theorem was applied to prove an estimate for $\mathbf{g}_T^V$ in [BG1], Section 10.16, which we now recall.

**4.24. Theorem.** *There exists an $r_\lambda$-invariant and $g^V$-self-adjoint $R_T \in \Omega^{\mathrm{even}, \geq 2}(B; \mathrm{End}^{\mathrm{even}} V)$ commuting with $\nabla^V$, $h^V$ and $N^V$, which is a polynomial in $\frac{1}{T}$ without constant term, and another $g^V$-self-adjoint $R'_T \in \Omega^*(B; \mathrm{End}\, V)$ with $|R'_T| = O(e^{-cT})$ uniformly in $C^k$-topology, such that*

$$\mathbf{g}_T^V = \left(\big(\mathrm{id}_V + R_T + R'_T\big)\, e^{-Th^V} \left(\frac{T}{\pi}\right)^{\frac{N^V}{2} - \frac{n}{4}}\right)^{*} g^V\ .$$

*Proof.* By definition of $\mathbf{g}_T^V$,

$$\begin{aligned} \mathbf{g}_T^V &= g_T^{TX,F}\big(\mathbb{I}_T^{-1}\cdot, \mathbb{I}_T^{-1}\cdot\big) \\ &= g_T^{TX,F}\big(\mathbf{e}_T \circ (I \circ \mathbf{e}_T)^{-1}\cdot, \mathbf{e}_T \circ (I \circ \mathbf{e}_T)^{-1}\cdot\big)\ . \end{aligned} \tag{4.25}$$

We will prove Theorem 4.24 by giving suitable approximations for $I \circ \mathbf{e}_T$ and $\mathbf{e}_T^* \, g_T^{TX,F}$ using Theorem 1.61 and Proposition 4.8.

Let $c \in C_\beta$ be a fibre-wise critical point of $h$, and let $v \in F_c \otimes o(T_c^u X)$. Then

(4.26)
$$\begin{aligned}\big(I \circ \mathbf{e}_T\big)_\alpha(v) &= e^{Th_\alpha} \, I\big(e^{-Th_\alpha} \mathbf{e}_T(v)\big)\\ &= e^{Th_\alpha} \, I\big(e^{T(h-h_\alpha)} \, \varphi(|\mathcal{R}|) \, J(v)\big) + e^{Th_\alpha} \, I\Big(e^{T(h-h_\alpha)} \, O\big(e^{-cT} \; |v|\big)\Big) \, .\end{aligned}$$

Clearly,
$$e^{T(h-h_\alpha)}|_{h^{-1}(-\infty,h_\alpha-\varepsilon^2)} = O\big(e^{-\varepsilon^2 T}\big) \, .$$

Together with Theorem 1.61 (2), this allows us to replace $M$ near $C$ by the total space of the vertical tangent bundle $E = TX|_C \to C$, and thus to approximate $e^{-Th_\alpha} \, I_T \circ \mathbf{e}_T \in$ by $\widehat{I} \circ e^{Th^E} \circ \widehat{J}$ considered in Proposition 4.8 (2). Here $h^E$ replaces $h - h_\alpha$, and we set $h^C = 0$. In particular, if $\alpha \neq \beta$, then (4.26) gives
$$\big(I_T \circ \mathbf{e}_T\big)_\alpha(v) = e^{Th_\alpha} \, O\big(e^{-cT} \, |v|\big)$$
uniformly in $C^k$-topology, possibly with a new constant $c$.

If we assume that $\beta = \alpha$, then Proposition 4.8 (2) and (4.26) imply that

$$\begin{aligned}\big(I_T \circ \mathbf{e}_T\big)_\alpha(v) &= e^{Th_\alpha} \, \Big(\widehat{I} \circ e^{T(h-h_\alpha)} \circ \widehat{J}(v)\big) + O\big(e^{-cT} \; |v|\big)\Big)\\ &= e^{Th_\alpha} \left(\left(\frac{T}{\pi}\right)^{\frac{n^+-n^-}{4}} \big(1 + R_{IJ,T}\big) \, v + O\big(e^{-cT} \; |v|\big)\right) ,\end{aligned}$$

where $R_{IJ}$ is an $r_\lambda$-invariant characteristic differential form of $\nabla^{TX}$ on $C$. We can summarize our results so far by

$$I_T \circ \mathbf{e}_T = e^{Th^V} \left(\frac{T}{\pi}\right)^{\frac{n}{4} - \frac{N^V}{2}} \Big(1 + \hat{p}_* R_{IJ,T} + O\big(e^{-cT}\big)\Big) \tag{4.27}$$

with $O(e^{-cT}) \in \Omega^*(B; \operatorname{End} V)$ for some $c > 0$, cf. [BG1], Theorem 10.59, and [BZ2], Theorem 6.11. Here, $\hat{p}_* R_{IJ,T} \in \Omega^{\geq 2}(B, \operatorname{End} V)$ clearly commutes with $\nabla^V$, $h^V$, and $N^V$. By Proposition 4.8, $R_{IJ,T}$ is a polynomial in $\frac{1}{T}$ without constant term.

We regard the generalised metric $\mathbf{e}_T^* \, g_T^{TX,F}$ on $V$. Taking adjoints with respect to $g^{TX,F}$ and $g^V$, Theorem 4.23 implies

$$\begin{aligned}\mathbf{e}_T^* \, g_T^{TX,F} &= g^{TX,F}\big(e^{-Th} \, \mathbf{e}_T \, \cdot \, , e^{-Th} \, \mathbf{e}_T \, \cdot\big)\\ &= g^V\Big( \, \cdot \, , \Big(J_T^* \, \varphi(\mathcal{R}) + O\big(e^{-cT}\big)\Big) \, \Big(\varphi(\mathcal{R}) \, J_T + O\big(e^{-cT}\big)\Big) \, \cdot \, \Big) \, .\end{aligned}$$

The presence of the cut-off function $\varphi(\mathcal{R})$ allows us again to replace $M$ by $\hat{p}_*(TX|_C)$. By Theorem 4.7 and Proposition 4.8, we find

$$\mathbf{e}_T^* \, g_T^{TX,F} = g^V\Big( \, \cdot \, , \Big(\hat{p}_*\big(\widehat{J}_T^* \, \widehat{J}_T\big) + O\big(e^{-cT}\big)\Big) \, \cdot \, \Big) = g^V\Big( \, \cdot \, , \Big(\mathrm{id}_V + \hat{p}_* R_{J^*J,T} + O\big(e^{-cT}\big)\Big) \, \cdot \, \Big) \tag{4.28}$$

uniformly in $C^k$-topology, cf [BG1], Proposition 10.58. Moreover, $\hat p_* R_{J^*J,T} \in \Omega^{\ge 2}(B, \operatorname{End} V)$ again commutes with $\nabla^V$, $h^V$, and $N^V$, and $\hat p_* R_{J^*J,T}$ is a polynomial in $\frac{1}{T}$ without constant term.

Let $\bigl(\mathrm{id}_V + \hat p_* R_{J^*J,T} + O(e^{-cT})\bigr)^{\frac12}$ denote the square-root of $\bigl(\mathrm{id}_V + \hat p_* R_{J^*J,T} + O(e^{-cT})\bigr)$, which is self adjoint in the sense of Definition 2.7. Combining (4.25), (4.27), and (4.28), we finally find that

$$\begin{aligned}\mathbf{g}^V_T &= \Biggl(\Bigl(\mathrm{id}_V + \hat p_* R_{J^*J,T} + O\bigl(e^{-cT}\bigr)\Bigr)^{\frac12} \\ &\qquad \Bigl(1 + \hat p_* R_{IJ,T} + O\bigl(e^{-cT}\bigr)\Bigr)^{-1} \left(\frac{T}{\pi}\right)^{\frac{N^V}{2} - \frac{n}{4}} e^{-Th^V}\Biggr)^* g^V \\ &= \Biggl(\Bigl(\mathrm{id}_V + R_T + O\bigl(e^{-cT}\bigr)\Bigr) \left(\frac{T}{\pi}\right)^{\frac{N^V}{2} - \frac{n}{4}} e^{-Th^V}\Biggr)^* g^V . \qquad \square\end{aligned}$$

We can now compare the generalised metric $\mathbf{g}^V_T$ with the deformed metric $e^{-2Th}\, g^V$.

**4.29. Theorem.** *As $T \to \infty$,*

$$\tilde f\bigl(A', e^{-2Th}\, g^V, \mathbf{g}^V_T\bigr) = \left(\frac{\chi'(V)}{2} - \frac{n\,\chi(V)}{4}\right) (\log T - \log \pi)\, f'(0) + O\!\left(\frac{1}{T}\right) .$$

*Proof.* We set

$$E_{s,T} = \Bigl(\mathrm{id}_V + s\,\bigl(R_T + R'_T\bigr)\Bigr)\, e^{-Th^V} \left(\frac{T}{\pi}\right)^{s\,\left(\frac{N^V}{2} - \frac{n}{4}\right)} \qquad \in \Omega^*\bigl(B; \operatorname{End} V\bigr) ,$$

and construct a family of metrics

$$\tilde{\mathbf{g}}^V_{s,T} = E^*_{s,T} g^V ,$$

so $g^V = \tilde{\mathbf{g}}^V_{0,T}$ and $\mathbf{g}^V_T = \tilde{\mathbf{g}}^V_{1,T}$. Let $\widetilde A''_{s,T}$ denote the adjoint of $A'$ with respect to $\tilde{\mathbf{g}}^V_{s,T}$, and let $\widetilde X_{s,T} = \frac12\bigl(\widetilde A''_{s,T} - A'\bigr)$, then we have

$$\tilde f\bigl(A', e^{-2Th}\, g^V, \mathbf{g}^V_T\bigr) = (2\pi i)^{-\frac{N^B}{2}} \int_0^1 \operatorname{str}_V\left(\frac12 \bigl(\tilde{\mathbf{g}}^V_{s,T}\bigr)^{-1} \frac{\partial \tilde{\mathbf{g}}^V_{s,T}}{\partial s}\, f'\bigl(\widetilde X_{s,T}\bigr)\right) ds . \tag{4.30}$$

Using Theorem 4.24, we can estimate the various terms appearing in (4.30). First of all,

$$\begin{aligned} E_{s,T}\,\bigl(\tilde{\mathbf{g}}^V_{s,T}\bigr)^{-1} \frac{\partial \tilde{\mathbf{g}}^V_{s,T}}{\partial s}\, E^{-1}_{s,T} &= \Bigl(\mathrm{id}_V + s\bigl(R_T + R'_T\bigr)\Bigr)^{-1} \Biggl(\bigl(R_T + R'_T\bigr) \\ &\qquad + \left(\frac{N^V}{2} - \frac{n}{4}\right)(\log T - \log \pi)\Bigl(\mathrm{id}_V + s\bigl(R_T + R'_T\bigr)\Bigr)\Biggr) \\ &\quad + \Biggl(\left(\frac{N^V}{2} - \frac{n}{4}\right)(\log T - \log \pi)\Bigl(\mathrm{id}_V + s\bigl(R_T + R'_T\bigr)\Bigr) \\ &\qquad + \bigl(R_T + R'_T\bigr)\Biggr)\Bigl(\mathrm{id}_V + s\bigl(R_T + R'_T\bigr)\Bigr)^{-1} \\ &= \left(N^V - \frac{n}{2}\right)(\log T - \log \pi) + 2R_T\,\bigl(\mathrm{id} + s\,R_T\bigr)^{-1} + O\bigl(e^{-cT}\bigr) \end{aligned}$$

uniformly in $s$, since $R_T$ and $N^V$ commute.

Now, we investigate the superconnections $A'$ and $\widetilde{A}''_{s,T}$. Let $A''$ be the adjoint of $A'$ with respect to $g^V$. Then we have

$$A' = \nabla^V + \sum_j a_j \qquad \text{and} \qquad A'' = \nabla^{V,*} + \sum_j a_j^* \; ,$$

where adjoints are taken with respect to $g^V$, and we know that $a_j$ raises eigenvalues of $h^V$ by at least $\varepsilon^2$ because of Theorem 1.61, while $a_j^*$ lowers them by at least the same amount. Since the estimates in Theorem 4.23 and Theorem 4.24 include derivatives of $R'_T$, we find

$$\begin{aligned} E_{s,T}\, A'\, E_{s,T}^{-1} &= \Bigl(\mathrm{id}_V + s\bigl(R_T + R'_T\bigr)\Bigr) \Biggl(\nabla^V + T\,\bigl[\nabla^V, h^V\bigr] \\ &\qquad + \sum_j \Bigl(\frac{T}{\pi}\Bigr)^{\frac{s(1-j)}{2}} e^{-Th^V}\, a_j\, e^{Th^V}\Biggr) \Bigl(\mathrm{id}_V + s\bigl(R_T + R'_T\bigr)\Bigr)^{-1} \\ &= \nabla^V + T\,\bigl[\nabla^{V,*}, h^V\bigr] + O\bigl(e^{-cT}\bigr)\; . \end{aligned}$$

Similarly, we obtain

$$E_{s,T}\, \widetilde{A}''_{s,T}\, E_{s,T}^{-1} = E_{s,T}^{-1}\, A''\, E_{s,T} = \nabla^{V,*} - T\,\bigl[\nabla^V, h^V\bigr] + O\bigl(e^{-cT}\bigr)\; .$$

Now we note that $h^V$ acts as a scalar function on each component $V^\alpha = \hat{V}|_{C_\alpha}$ of $V$, and that $R_T$, $\nabla^V$ and $\nabla^{V,*}$ preserve these components. Writing

$$\omega^V = \nabla^{V,*} - \nabla^V \qquad \text{and} \qquad dh^V = \bigl[\nabla^V, h^V\bigr] = \bigl[\nabla^{V,*}, h^V\bigr]\; ,$$

we finally obtain

$$\begin{aligned} &\mathrm{str}_V\left(\frac{1}{2}\,\bigl(\tilde{\mathbf{g}}^V_{s,T}\bigr)^{-1}\,\frac{\partial \tilde{\mathbf{g}}^V_{s,T}}{\partial s}\, f'\bigl(\widetilde{X}_{s,T}\bigr)\right) \\ &\qquad = \mathrm{str}_V\left(\frac{1}{2}\left(\Bigl(N^V + \frac{n}{2}\Bigr)(\log T - \log\pi) + 2R_T\,\bigl(\mathrm{id} + R_T\bigr)^{-1}\right) f'\Bigl(\frac{\omega^V}{2} - T\,dh^V\Bigr)\right) + O\bigl(e^{-cT}\bigr) \\ &\qquad = \left(\frac{\chi'(V)}{2} - \frac{n\,\chi(V)}{4}\right)(\log T - \log\pi)\, f'(0) + O\Bigl(\frac{1}{T}\Bigr)\; , \end{aligned}$$

because $\frac{\omega^V}{2} + T\,dh^V$ is odd and commutes with $N^V + 2R_T\,(\mathrm{id} + R_T)^{-1}$, while $f'$ is an even function, and supertraces of even powers of odd operators vanish. By (4.30), this finishes our proof of Theorem 4.29. □

Using this result, we can summarize the results of this section, where we now make use of the fact that $f(z) = ze^{z^2}$ by (2.11).

**4.31. Theorem.** *As $T \to \infty$, we have*

$$\begin{aligned} &\tilde{f}\bigl(\nabla^H, g^H, g^H_T\bigr) + S_f\bigl(T^H M, g^{TX}, \nabla^F, g^F_T\bigr) - \left(\frac{\chi'(V)}{2} - \frac{n\,\chi(V)}{4}\right)\log T + T\,\mathrm{str}_V\bigl(h^V\bigr) \\ &\qquad\qquad \longrightarrow\; \tilde{f}\bigl(\nabla^H, g^H, g^H_V\bigr) + T_f\bigl(A', g^V, h^V\bigr) \\ &\qquad\qquad\qquad - \left(\frac{\chi'(V)}{2} - \frac{n\,\chi(V)}{4}\right)\log\pi - \bigl(\chi'(V) - \chi'(H)\bigr)\int_1^\infty f'\bigl(\sqrt{-t}/2\bigr)\,\frac{dt}{2t} \\ &\qquad\qquad\qquad - \bigl(\chi'(V) - \chi'(H)\bigr)\int_0^1 \Bigl(f'\bigl(\sqrt{-t}/2\bigr) - 1\Bigr)\,\frac{dt}{2t} \end{aligned}$$

*modulo exact forms on $B$.*

This theorem is our analogue of Theorem 9.8 in [BG1]. Note however, that a different notion of torsion appears in [BG1].

*Proof.* From (4.19), Theorem 4.29 and the fact that $f'(0) = 1$, we get

$$\begin{aligned}(4.32)\quad \tilde f\big(\nabla^H, g^H, g^H_T\big) + S_f\big(T^H M, g^{TX}, \nabla^F, g^F_T\big) - \left(\frac{\chi'(V)}{2} - \frac{n\,\chi(V)}{4}\right)(\log T - \log\pi)\\ = \tilde f\big(\nabla^H, g^H, g^H_V\big) + S_f\big(A', g^V\big) + \tilde f\big(A', g^V, e^{-2Th^V} g^V\big) + O\big(T^{-\delta}\big) \quad \in \Omega^*(B)/d^B\Omega^*(B)\ .\end{aligned}$$

We change back to the notation of Section 2.d. In particular, $A''_{1,\tau}$ now refers to the adjoint of $A'$ with respect to the metric $g^V\big(T^{h^V}\cdot\,,\,\cdot\big)$, and of course $X_{1,T} = \frac12\,(A''_{1,T} - A')$. Then

$$\begin{aligned}(4.33)\qquad \tilde f\big(A', g^V, e^{-2Th^V} g^V\big) &= -(2\pi i)^{-\frac{N^B}{2}} \int_0^T \operatorname{str}_V\Big(h^V\, f'\big(X_{1,e^{-2\tau}}\big)\Big)\, d\tau\\ &= -(2\pi i)^{-\frac{N^B}{2}} \int_{e^{-2T}}^1 \operatorname{str}_V\Big(h^V\, f'\big(X_{1,T'}\big)\Big)\,\frac{dT'}{2T'}\ .\end{aligned}$$

By (4.32) and (4.33),

$$\begin{aligned}\tilde f\big(\nabla^H, g^H, g^H_T\big) + S_f\big(T^H M, g^{TX}, \nabla^F, g^F_T\big) - \left(\frac{\chi'(V)}{2} - \frac{n\,\chi(V)}{4}\right)(\log T - \log\pi) + T\,\operatorname{str}_V\big(h^V\big)\\ = \tilde f\big(\nabla^H, g^H, g^H_V\big) + S_f\big(A', g^V\big) - (2\pi i)^{-\frac{N^B}{2}} \int_{e^{-2T}}^1 \Big(\operatorname{str}_V\big(h^V\, f'(X_{1,T'})\big)\\ - \operatorname{str}_V\big(h^V\big)\Big)\,\frac{dT'}{2T'} + O\big(T^{-\delta}\big) \quad \in \Omega^*(B)/d^B\Omega^*(B)\ .\end{aligned}$$

In the limit $T\to\infty$, Definition 2.46 gives us

$$\begin{aligned}\tilde f\big(\nabla^H, g^H, g^H_T\big) + S_f\big(T^H M, g^{TX}, \nabla^F, g^F_T\big) - \left(\frac{\chi'(V)}{2} - \frac{n\,\chi(V)}{4}\right)(\log T - \log\pi) + T\,\operatorname{str}_V\big(h^V\big)\\ \longrightarrow\ \tilde f\big(\nabla^H, g^H, g^H_V\big) + T_f\big(A', g^V, h^V\big)\\ - \big(\chi'(V) - \chi'(H)\big) \int_1^\infty f'\big(\sqrt{-t}/2\big)\,\frac{dt}{2t}\\ - \big(\chi'(V) - \chi'(H)\big) \int_0^1 \Big(f'\big(\sqrt{-t}/2\big) - 1\Big)\,\frac{dt}{2t} \quad \in \Omega^*(B)/d^B\Omega^*(B)\ ,\end{aligned}$$

which finishes the proof. □

*4.34. Remark.* Suppose that $\gamma$ is a fibre-wise isometry of $M$ that acts on $F$ and preserves $T^H M$, $g^{TX}$, $\nabla^F$ and $g^F$. Suppose also that $\gamma$ acts on the fibres of $V\to B$ preserving $A'$, $g^V$ and $h^V$, and that the map $I\colon \Omega^*(M;F)\to\Omega^*(B;V)$ is $\gamma$-equivariant. Then there are forms $\tilde f_\gamma(\nabla^H, g^H, g^H_T)$, $S_{f,\gamma}(T^H M, g^{TX}, \nabla^F, g^F_T)$, $T_{f,\gamma}(A', g^V, h^V)$ etc., and all considerations in this section still go through and provide us with a $\gamma$-equivariant version of Theorem 4.31, cf. [BG1], Theorem 9.8

## 5. Exotic smooth fibre bundles

We calculate the analytic torsion form of a sphere bundle $p_\gamma\colon M_\gamma \to S^{4k}$ that is fibre-wise homeomorphic to a trivial bundle, but not diffeomorphic. We then use this bundle to construct "exotic" bundle structures on a quite general class of compact fibre bundles.

The construction of the sphere bundle $p_\gamma$ in Section 5.a goes back to Hatcher, and was kindly explained to us by Igusa, see also [I2], Section 6. The computation of $\mathcal{T}(M_\gamma, S^{4k}; \mathbb{C})$ in Section 5.b and the extension to bundles with arbitrary fibres in Section 5.c are still closely related to Igusa's work.

In Section 5.d, we use Hatcher's example to exhibit a large class of fibre bundles $M \to B$ where we can construct infinitely many bundles that are homeomorphic but not diffeomorphic as bundles to $M \to B$, see Corollary 5.22. In Section 5.e we finally prove a refinement of Theorem 0.2, see Theorem 5.29. We will also relate our results to [DWW].

**5.a. Hatcher's construction.** We will construct a bundle $N_\gamma \to S^{4k}$, with fibre

$$D^{m+n} = (S^m \times D^n) \cup_{S^m \times D^{n-1}} (D^{m+1} \times D^{n-1}) , \tag{5.1}$$

see Figure 5.2. Here, $\gamma \in \ker J$, see (5.3) below. The bundle $M_\gamma$ will then be given by the fibre-wise double of $N_\gamma$. Over the "southern" hemisphere $D^{4k}_- \subset S^{4k}$ of the base sphere, we take a fixed embedding of $D^{n-1}$ into a fixed hemisphere of $S^{n-1} = \partial D^n$ in (5.1), thus obtaining a fixed trivialisation of $D^{m+n} \times D^{4k}_-$.

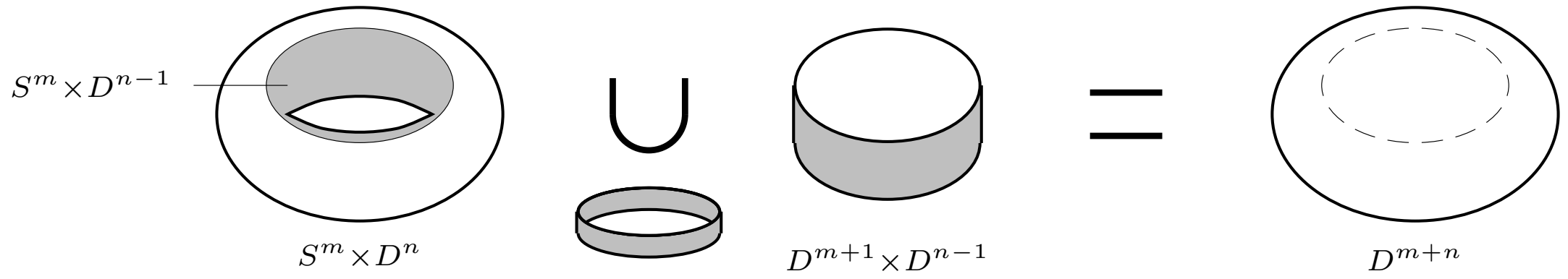

Figure 5.2. A single fibre of $N_\gamma$.

We will now describe the construction over the "northern" hemisphere $D^{4k}_+ \subset S^{4k}$. For fixed $k$ and $m$, we have the $J$-homomorphism, which is given as

$$J\colon\ \pi_{4k-1}\big(O(m)\big) \longrightarrow \pi_{4k-1}\big(C^0((S^m, x_0), (S^m, x_0))\big) \cong \pi_{4k-1+m}(S^m, x_0) , \tag{5.3}$$

where $x_0 \in S^m$ denotes the north pole. The groups $\pi_{4k-1}(O(m))$ and $\pi_{4k-1+m}(S^m, x_0)$ are independent of $m$ for $m$ sufficiently large, and $\pi_{4k-1+m}(S^m, x_0)$ is finite, while $\pi_{4k-1}(O(m))$ is infinite cyclic. Thus, pick $0 \neq \gamma \in \ker J$, which is represented by a map $\gamma\colon S^{4k-1} \to O(m)$ that extends continuously to $\gamma\colon D^{4k}_+ \to C^0((S^m, x_0), (S^m, x_0))$.

For $n$ large enough, using transversality, we can approximate the map

$$D^{4k}_+ \times (S^m, x_0) \xrightarrow{\ \gamma\ } (S^m, x_0) \lhook\joinrel\longrightarrow \big(S^m \times D^{n-1}, (x_0, 0)\big)$$

by a smooth family of embeddings $\tilde\gamma_q\colon (S^m, x_0) \hookrightarrow (S^m \times D^{n-1}, (x_0, 0))$ for $q \in D^{4k}_+$, which equals

$$(qr, x) \longmapsto \big(\gamma_q(x), 0\big) \qquad \text{for } q \in S^{4k-1} = \partial D^{4k}_+,\ r \in \big[\tfrac12, 1\big], \text{ and } x \in S^m .$$

The normal bundle $\mathcal{N}_q \to S^m$ of $\tilde\gamma_q$ is stably trivial for all $q \in D^{4k}_+$, because the tangent bundles of $S^m$ and $S^m \times D^{n-1}$ are stably trivial. Making $n$ even larger if necessary, we may assume that the $\mathcal{N}_q$ are trivial. Then the bundle $\mathcal{N} \to D^{4k}_+ \times S^m$ is of course also trivial. Fixing a trivialisation of $\mathcal{N}$ and using the normal exponential map, we can extend $\widetilde{\gamma}$ to a family of embeddings

$$\widetilde{\gamma}_q\colon\ S^m \times D^{n-1} \longrightarrow S^m \times D^{n-1} \qquad \text{for } q \in D^{4k}_+ \ .$$

We may assume that there exists a map $\gamma'\colon S^{4k-1} \times S^m \to O(n-1)$ such that for all $q \in S^{4k-1} = \partial D^{4k}_+$ and all $r \in \left[\frac{1}{2}, 1\right]$, the embedding $\widetilde{\gamma}_{qr}$ becomes

$$\widetilde{\gamma}_{qr}(x,y) = \bigl(\gamma_q(x), \gamma'_{q,x}(y)\bigr) \ .$$

**5.4. Proposition.** *The family of embeddings $\widetilde{\gamma}_q\colon S^m \times D^{n-1} \to S^m \times D^{n-1}$ for $q \in D^{4k}_+$ is isotopic to a family of embeddings that restricts over $S^{4k-1} = \partial D^{4k}_+$ to*

$$(x,y) \longmapsto \bigl(\gamma_q(x), (\gamma_q^{-1} \oplus \mathrm{id}_{\mathbb{R}^{n-1-m}})(y)\bigr) \qquad \textit{for all } q \in S^{4k-1} \ .$$

With this result, we finally construct $N_\gamma$. We still identify $D^{n-1}$ with a fixed hemisphere of $S^{n-1} = \partial D^n$. For $q \in D^{4k}_+$, we construct the fibre

$$D^{m+n} = (S^m \times D^n) \cup_{S^m \times D^{n-1}} (D^{m+1} \times D^{n-1})$$

as in (5.1), identifying

$$(x,y) \in S^m \times D^{n-1} \subset \partial(D^{m+1} \times D^{n-1}) \qquad \text{with} \qquad \widetilde{\gamma}_q(x,y) \in S^m \times D^{n-1} \subset \partial(S^m \times D^n) \ .$$

Then we can glue $N_\gamma|_{D^{4k}_+}$ with $N_\gamma|_{D^{4k}_-}$ along $S^{4k-1}$ using the identity on $S^m \times D^n$ and identifying

$$(q,x,y) \in N_\gamma|_{D^{4k}_+} \qquad \text{with} \qquad \bigl(q, \gamma_q(x), (\gamma_q^{-1} \oplus \mathrm{id}_{\mathbb{R}^{n-1-m}})(y)\bigr) \in N_\gamma|_{D^{4k}_-} \ .$$

Finally, we construct the sphere bundle

$$p\colon\ M_\gamma = N_\gamma \cup_{\mathrm{id}_{\partial N_\gamma}} N_\gamma \longrightarrow S^{4k} \ . \tag{5.5}$$

*Proof* of Proposition 5.4. Consider the family of rotations $S^{4k-1} \to O(2m)$ given by

$$q \longmapsto \bigl(\gamma_q \oplus \gamma_q^{-1}\bigr)\colon\ S^{4k-1} \longrightarrow O(2m) \ .$$

By a well-known trick from $K$-theory, this family is homotopic to the trivial family $q \mapsto \mathrm{id}_{\mathbb{R}^{2m}}$: apply the Cartan projection $GL(2m,\mathbb{R}) \to O(2m)$ to the homotopy

$$(q,t) \longmapsto \begin{pmatrix} 1 & t\gamma_q \\ 0 & 1 \end{pmatrix} \begin{pmatrix} 1 & 0 \\ -t\gamma_q^{-1} & 1 \end{pmatrix} \begin{pmatrix} 1 & t\gamma_q \\ 0 & 1 \end{pmatrix} \begin{pmatrix} 1-t & -t \\ t & 1-t \end{pmatrix} \tag{5.6}$$

in $GL(2m,\mathbb{R})$, where $t \in [0,1]$. In particular, there exists a map $G\colon D^{4k}_+ \to O(m+n)$ that for all $q \in S^{4k-1}$ and all $r \in \left[\frac{1}{2}, 1\right]$ takes the form

$$G_{rq} = \bigl(\gamma_q \oplus \gamma_q^{-1} \oplus \mathrm{id}_{\mathbb{R}^{n-m}}\bigr) \ .$$

The map $\widetilde{\gamma}_q\colon S^m \times D^{n-1} \to S^m \times D^{n-1}$ is determined by its 1-jet

$$D\widetilde{\gamma}_q|_{S^m \times \{0\}}\colon\ TS^m \times \mathbb{R}^{n-1} \longrightarrow TS^m \times TD^{n-1} \ .$$

We add the trivial line $\mathbb{R}$ on both sides, and identify $TS^m \times \mathbb{R}$ with $\mathbb{R}^{m+1}$ using the standard embedding $S^m \hookrightarrow \mathbb{R}^{m+1}$. Then $D\widetilde{\gamma}_q|_{S^m \times \{0\}}$ acts on $TS^m \times \mathbb{R}$ by the matrix $\gamma_q \oplus \mathrm{id}_{\mathbb{R}} \in O(m+1)$ for $q \in S^{4k-1}$. This implies that $G_q^{-1} \circ D\widetilde{\gamma}_q|_{S^m \times \{0\}}$ acts trivially on $TS^m \times \mathbb{R}$ for $q \in S^{4k-1}$.

Because the Stiefel manifold $O(m+n+1)/O(m+1)$ is $4k$-connected for $n$ sufficiently large, we can deform $G_q^{-1} \circ D\widetilde{\gamma}_q$ on $D^{4k}_+$ relative to $S^{4k-1}$ such that $TS^m \times \mathbb{R} \cong \mathbb{R}^{m+1}$ is trivially mapped to itself. Finally, we can deform the action on the factor $\mathbb{R}^{n-1}$ to the trivial one over $D^{4k}_+$. This shows that the map $G_q^{-1} \circ D\widetilde{\gamma}_q|_{S^m \times \{0\}}$ can be deformed to the identity, relative to $S^{4k-1} \times TS^m \times \mathbb{R}$. Then $D\widetilde{\gamma}_q|_{S^m \times \{0\}}$ can be deformed to $G_q$ relative to $S^{4k-1} \times TS^m \times \mathbb{R}$. In particular, the embedding $\widetilde{\gamma}_q$ is isotopic to an embedding that takes the form given in the proposition. □

In the remainder of this chapter, we will use the example above to construct families of bundles which are pairwise homeomorphic as bundles, but not diffeomorphic in one the following two senses.

**5.7. Definition.** Two fibre bundles $p\colon M \to B$ and $p'\colon M' \to B$ are *homeomorphic (diffeomorphic) as bundles* iff there exists a homeomorphism (diffeomorphism) $\Phi$ such that the left diagramme below commutes.

$$\begin{array}{ccc} M & \xrightarrow[\sim]{\Phi} & M' \\ {\scriptstyle p}\big\downarrow & & \big\downarrow{\scriptstyle p'} \\ B & \xrightarrow{\mathrm{id}} & B \end{array} \qquad\qquad \begin{array}{ccc} M & \xrightarrow[\sim]{\Phi} & M' \\ {\scriptstyle p}\big\downarrow & & \big\downarrow{\scriptstyle p'} \\ B & \xrightarrow[\sim]{\Psi} & B \end{array}$$

We say that $p$ and $p'$ are *homeomorphic (diffeomorphic) as maps* iff there exist homeomorphisms (diffeomorphisms) $\Phi$ and $\Psi$ such that the right diagramme above commutes.

**5.8. Proposition.** *The bundle $p\colon M_\gamma \to S^{4k}$ is homeomorphic as a bundle to a trivial sphere bundle.*

*Proof.* It is well known that the diffeomorphism group of $D^{m+n}$ with the $C^0$-topology retracts to $O(m+n)$, so $M_\gamma \to S^{4k}$ is homeomorphic to some sphere bundle with structure group $O(m+n)$. Our claim then follows because the bundle $p\colon M_\gamma \to S^{4k}$ admits sections $C_1, \dots, C_6$, and the restriction of $TX$ to each of these sections is trivial. $\square$

*5.9. Remark.* In the construction above, we did not care at all about the dimension of the fibre of $p\colon M_\gamma \to S^{4k}$. It should be clear that there are some numbers $p$, $q$ such that one can always keep the dimension of the fibre at most $pk + q$, and it is interesting to know explicite values of $p$ and $q$. This question is adressed in [I2], see Theorem 6.2.2 and Corollary 6.5.8.

**5.b. The higher torsion invariants of Hatcher's sphere bundle.** We calculate Bismut-Lott's analytic torsion form and Igusa's higher Franz-Reidemeister torsion of $p\colon M_\gamma \to S^{4k}$ by specifying a fibre-wise Morse function on $M_\gamma$.

We first construct a fibre-wise Morse function $h\colon N_\gamma \to \mathbb{R}$ such that $dh$ is positive on the outward normal vector to $\partial N_\gamma$. We start with a fixed Morse function on $S^m$, having exactly two critical points: a minimum of $-3$ at the south pole, and a maximum of $-2$ at the north pole. Adding a small positive function on $D^n$ whose only critical point is a minimum 0 at the origin, we obtain a Morse function on the trivial bundle $S^m \times D^n \times S^{4k} \to S^{4k}$.

We then choose a Morse function on $D^{m+1} \times D^{n-1}$ with only one critical point at $(0,0)$ of value $-1$, such that $D^{m+1} \times \{0\}$ is the descending cell, and $\{0\} \times D^{n-1}$ is the ascending cell. Clearly, these fibre-wise Morse functions can be patched together to give a fibre-wise Morse function $h\colon N_\gamma \to \mathbb{R}$ with three leaves of fibre-wise critical points of indices 0, $m$ and $m+1$, see Figure 5.10. We may assume that $h^{-1}(0) = \partial N_\gamma$.

We finally get a Morse function on $M_\gamma = N_\gamma \cup_{\mathrm{id}_{\partial N_\gamma}} N_\gamma$ by taking $h$ on one copy of $N_\gamma$ and $-h$ on the other. Assuming that $n \geq m + 3$, this new Morse function has six fibre-wise critical leaves $C_1$, $\dots$, $C_6$ of indices 0, $m$, $m+1$, $n-1$, $n$ and $m+n$, with values $-3$, $-2$, $-1$, 1, 2 and 3, respectively. In particular, no critical leaf can "cross over" another critical leaf of higher or equal index. For this reason, the "simplicial superconnection" $a$ of Proposition 1.31 has no terms of higher degree, i.e., $a_{\alpha\beta}(\sigma) = 0$ unless $\sigma \in S_0$.

Letting $F = M_\gamma \times \mathbb{C}$ be trivial, the bundle $V = S^{4k} \times \mathbb{C}^6$ is also trivial. We can arrange the superconnection $A'$ of Theorem 1.61 to equal $a_0 + \nabla^V$, with

$$a_0 = \begin{pmatrix} 0 &&&&& \\ & 0 &&&& \\ & 1 & 0 &&& \\ &&& 0 && \\ &&& 1 & 0 & \\ &&&&& 0 \end{pmatrix} \in \Gamma(\operatorname{End} V)$$

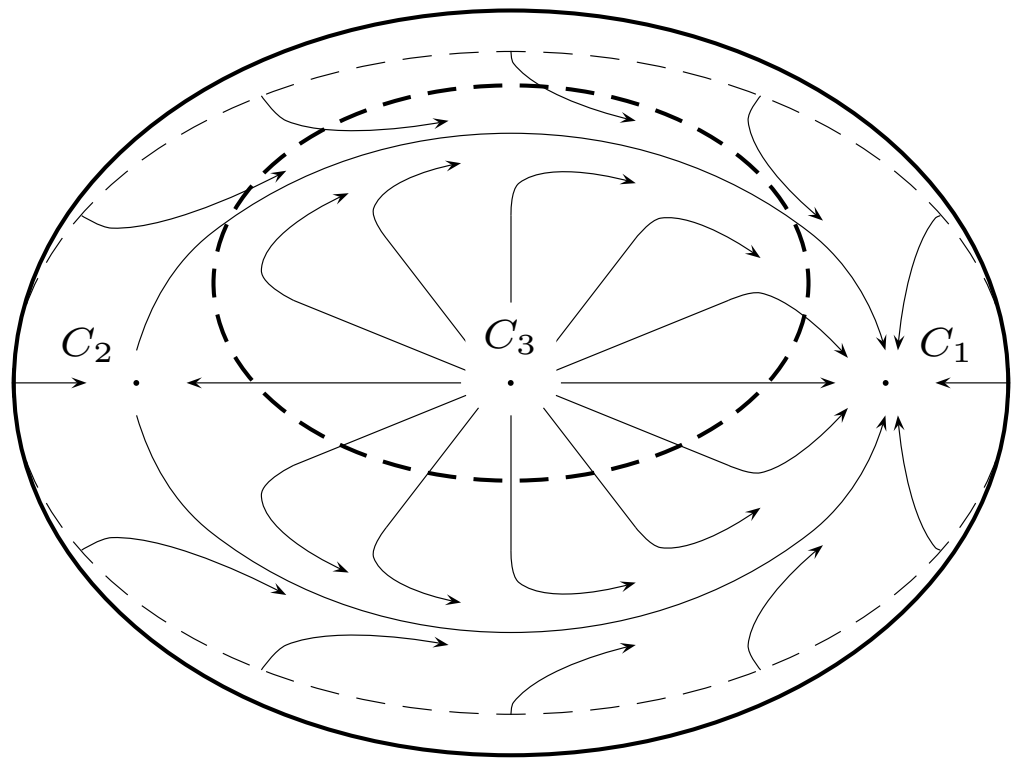

Figure 5.10. The negative gradient flow of the fibre-wise Morse function $h$ on $N_\gamma$.

parallel with respect to the trivial connection $\nabla^V$. Thus

$$\big(\Omega^*(S^{4k};V),A'\big)=\big(\Omega^*(S^{4k}),d\big)\otimes(\mathbb{C}^6,a_0)$$

is a product, and the bundle $H=H^*(M_\gamma/S^{4k};\mathbb{C})\to S^{4k}$ is trivial, too.

We also find that $T(A',g^V,h^V)\in\mathbb{R}\subset\Omega^*(S^{4k})$ is homogeneous of degree 0 (it actually vanishes). However, the term

$$\hat{p}_*{}^0J\big((-1)^{\mathrm{ind}_h}(T^sX-T^uX)\big)$$

in Theorem 0.1 does not vanish. The critical leaves $C_1$, $C_2$, $C_5$ and $C_6$ lie within the trivial bundle

$$\big((S^m\times D^n)\cup_{\mathrm{id}_{S^m\times D^{n-1}}}(S^m\times D^n)\big)\times S^{4k}\;\longrightarrow\;S^{4k}$$

and come with trivial bundles $T^sX$, $T^uX$, so they do not contribute. The bundles $TX|_{C_3}$ and $TX|_{C_4}$ are trivial, because they arise from trivial bundles over $D^{4k}_+$ and $D^{4k}_-$ using the clutching function $q\mapsto(\gamma_q\oplus\gamma_q^{-1}\oplus\mathrm{id})$, which is homotopic to the identity by (5.6). However, the bundles $T^{s,u}X|_{C_3}$ and $T^{s,u}X|_{C_4}$ are not trivial because $\gamma_q$ was chosen non-trivial. The bundle $T^sX|_{C_4}$ is isomorphic to $T^uX|_{C_3}$, and has clutching function $\gamma^{-1}$. If $u$ is a generator of $H^{4k}(S^{4k};\mathbb{Z})\subset H^{4k}_{\mathrm{dR}}(S^{4k})$, if $\gamma$ is the $n_\gamma$-th multiple of a generator of $\ker J$, and if $0\neq b_k\in\mathbb{Q}$ is the coefficient discussed in [I2], Remark 6.4.3, then

$$p_j\big(T^uX|_{C_3}\big)=p_j\big(T^sX|_{C_4}\big)=\begin{cases}n_\gamma\,b_k\,u & \text{if } j=k\text{, and}\\ 0 & \text{otherwise.}\end{cases}\tag{5.11}$$

By triviality of $TX|_{C_3}$ and $TX|_{C_4}$ and multiplicativity of the total Pontrijagin class,

$$p_j\big(T^uX|_{C_4}\big)=p_j\big(T^sX|_{C_3}\big)=-p_j\big(T^uX|_{C_3}\big)=\begin{cases}-n_\gamma\,b_k\,u & \text{if } j=k\text{, and}\\ 0 & \text{otherwise,}\end{cases}\tag{5.12}$$

Because the two flat bundles $F\to M_\gamma$ and $H\to S^{4k}$ are trivial, we have a higher analytic torsion class $\mathcal{T}(M_\gamma/S^{4k};\mathbb{C})\in\Omega^{\geq 2}(S^{4k})$ as in Definition 2.88. By Theorem 0.1 and Corollary 3.9, we get a parallel result to [I2], Theorem 6.4.2. Note that passing from Pontrijagin to Chern classes in the definition of ${}^0J$ gives an extra factor $2^{-k}$.

**5.13. Theorem.** *The analytic torsion class of the bundle $M_\gamma \to S^{4k}$ above is given by*

$$\begin{aligned}\mathcal{T}(M_\gamma/S^{4k};\mathbb{C}) &= (-1)^{m+1}\,{}^0J\big(T^sX|_{C_3}\big) + (-1)^m\,{}^0J\big(T^uX|_{C_3}\big)\\ &\qquad + (-1)^{n-1}\,{}^0J\big(T^sX|_{C_4}\big) + (-1)^n\,{}^0J\big(T^uX|_{C_4}\big)\\ &= \begin{cases} (-1)^m\,2^{2-k}\,b_k\,\zeta'(-2k)\,u \in \Omega^{4k}(S^{4k}) & \text{if } n+m \text{ is odd, and}\\ 0 & \text{if } n+m \text{ is even.}\end{cases}\end{aligned}$$

In particular, neither the sphere bundle $M_\gamma \to S^{4k}$ nor the disc bundle $N_\gamma \to S^{4k}$ are diffeomorphic to a trivial bundle. Thus we recover a result by Bökstedt [Bo], cf. [I2], Chapter 6.

Note on the other hand that using Igusa's framing principle (Theorem 6.1.1 in [I2]), the higher Franz-Reidemeister torsion of [I2] equals

$$\tau(M_\gamma/S^{4k};\mathbb{C}) = (-1)^m\,2\,{}^0J\big(T^uX|_{C_3}\big) + (-1)^n\,2\,{}^0J\big(T^uX|_{C_4}\big) = \mathcal{T}(M_\gamma/S^{4k};\mathbb{C})\;.$$

This is consistent with Theorem 0.3, by which

$$\mathcal{T}(M_\gamma/S^{4k};\mathbb{C}) - \tau(M_\gamma/S^{4k};\mathbb{C}) = \int_{M_\gamma/S^{4k}} e(TX)\,{}^0J(TX) = 0\;,$$

because the bundle $TX|_C$ is trivial by construction.

*5.14. Remark.* Note that the superconnection $A' = \nabla^V + a_0$ on $V \to S^{4k}$ looks like a superconnection of a fibre-wise Morse function with a fibre-wise gradient field satisfying the Smale transversality condition, cf. [BG1], section 5. On the other hand, if there was a fibre-wise Smale gradient field for $h$, we would have

$$\begin{aligned} T^sX|_{C_2} &\cong T^sX|_{C_3} \oplus \mathbb{R}\;, & & T^uX|_{C_2} \oplus \mathbb{R} \cong T^uX|_{C_3}\;,\\ T^sX|_{C_4} &\cong T^sX|_{C_5} \oplus \mathbb{R}\;, &\quad\text{and}\quad & T^uX|_{C_4} \oplus \mathbb{R} \cong T^uX|_{C_5}\end{aligned}$$

by [BG1], section 7.10, contradicting (5.11), (5.12). This shows that the considerations in Section 1 are indeed needed to deal with this example.

*5.15. Remark.* Consider once more the $(m+n)$-disc bundle $N_\gamma \to S^{4k}$ above. If $m+n$ is odd, the boundary $\partial N_\gamma$ becomes a trivial smooth sphere bundle, possibly after replacing $N_\gamma$ by $N_\gamma \times D^{2l} \to S^{4k}$ for some $l \geq 0$. This means that instead of doubling $N_\gamma$, we could glue $N_\gamma$ and a trivial disc bundle to obtain another sphere bundle $\overline{N}_\gamma \to S^{4k}$ with a fibre-wise Morse function with only four fibre-wise critical points, namely $C_1$, $C_2$ and $C_3$ on $N_\gamma$ and a maximum on the trivial disc bundle. Then clearly

$$\mathcal{T}(\overline{N}_\gamma/S^{4k};\mathbb{C}) = \frac{1}{2}\,\mathcal{T}(M_\gamma/S^{4k};\mathbb{C}) = (-1)^m\,2^{1-k}\,b_k\,\zeta'(-2k)\,u \in \Omega^{4k}\big(S^{4k}\big)\;.$$

In order to keep this chapter reasonably self-contained, we will stick to the slightly more complicated bundle $p\colon M_\gamma \to S^{4k}$ constructed above.

**5.c. Exotic bundles over $S^{4k}$ with general fibres.** We extend Hatcher's example to generate examples of bundles with arbitrary fibre. This involves a very simple kind of a vertical gluing formula. A more general kind of vertical gluing will be discussed in [G]. The following construction is closely related to Corollary 6.5.8 in [I2]. We need it as an intermediate step for the general construction in the following subsection.

Fix $k > 0$ and let $X$ be a compact manifold of dimension $\dim X = 2l-1$, such that $l$ is sufficiently large with respect to $k$. Let $p_X \colon X \times S^{4k} \to S^{4k}$ be the trivial bundle. We fix a Morse function on $X$, which gives rise to a trivial fibre-wise Morse function $h_X \colon X \times S^{4k} \to \mathbb{R}$. Any flat vector bundle $F_X \to X$ of rank $r$ with a parallel metric pulls back to a flat bundle $F_X \to X \times S^{4k}$.

Let $p_\gamma \colon M_\gamma \to S^{4k}$ be the bundle with fibre $S^{2l-1}$ constructed in Section 5.a. In particular, there exists a fibre-wise Morse function $h_\gamma \colon M_\gamma \to \mathbb{R}$ that has six critical points on each fibre. Let $S^{4k} \cong C_6 \subset M_\gamma$ be the fibre-wise maximum, then $C_6$ has a trivial neighbourhood $U \cong C_0 \times D^{2l-1}$ in $M_\gamma$.

We can now replace a small neighbourhood of one of the local minima of $h_X$ by $M_\gamma \backslash U$, obtaining a new bundle $p'_\gamma \colon M'_\gamma \to S^{4k}$ with fibre $X$ that is still homeomorphic as a bundle to $p_X \colon X \times S^{4k} \to S^{4k}$. By gluing the bundle $F_X \to X \times S^{4k}$ with the trivial bundle of rank $r$ over $M_\gamma$, we get a flat vector bundle $F'_\gamma \to M'_\gamma$, and $F'_\gamma \to M'_\gamma$ is the pull back of $F$ by the homeomorphism $M'_\gamma \to X \times S^{4k}$ above. In particular, the fibre-wise cohomology of $M'_\gamma \to S^{4k}$ is still a trivial bundle over $S^{4k}$ with fibre $H^*(X; F_X)$.

**5.16. Proposition** (cf. Corollary 6.5.8 in [I2]). *If $F_X \to X$ is trivial, we have*

$$\mathcal{T}(M'_\gamma/S^{4k}; \mathbb{C}) = 2^{2-k}\, j\, b_k\, \zeta'(-2k)\, u \quad \in \quad H^{\text{even},\geq 2}(B; \mathbb{R}) .$$

Thus again, the higher analytic torsion is able to detect the "exotic" differentiable structure of the homeomorphically trivial bundle $p'_\gamma \colon M'_\gamma \to S^{4k}$.

In [G], we will exhibit a spectral sequence argument which gives rise to a gluing formula for higher torsion invariants. In particular, we will find that

$$\mathcal{T}(M'_\gamma/S^{4k}; F'_\gamma) = 2^{2-k}\, j\, b_k\, \zeta'(-2k)\, u \operatorname{rk} F_X \quad \in \quad H^{\text{even},\geq 2}(B; \mathbb{R}) .$$

Note that in the argument below, there are no spectral sequence contributions by (5.18).

*Proof.* After subtracting a suitable constant from $h_\gamma$, gluing with $h_X$, and smoothing along the common boundary, we obtain a new fibre-wise Morse function $h'_\gamma \colon M'_\gamma \to \mathbb{R}$. Let $C'_\gamma \subset M'_\gamma$ denote the critical points of $h'_\gamma$, then

$$C'_\gamma = C'^{+}_\gamma \,\dot\cup\, C'^{-}_\gamma ,$$

where $C'^{-}_\gamma$ are the five remaining critical levels of $h_\gamma$, and $C'^{+}_\gamma$ are the remaining critical points of $h_X$. This induces a splitting

$$V'_\gamma = V'^{+}_\gamma \oplus V'^{-}_\gamma \longrightarrow B . \tag{5.17}$$

Since the descending cells of the critical levels $C'^{+}_\gamma$ can meet the ascending cells of $C'^{-}_\gamma$ in any possible way, we cannot conclude any more that the superconnection $A' = \nabla^{V'_\gamma} + a'$ of Definition 1.60 on the family Thom-Smale complex $\Omega^*(S^{4k}; V)$ still takes the form $A' = \nabla^{V'_\gamma} + a_0$ with $a_0 \in \Omega^0(B, \operatorname{End} V)$. We may however conclude that

$$a' = \begin{pmatrix} b' & d' \\ 0 & c' \end{pmatrix}$$

with respect to (5.17), such that $b' \in \Omega^0(B; \operatorname{End} V'^{+}_\gamma)$ and $c' \in \Omega^0(B; \operatorname{End} V'^{-}_\gamma)$ are both of horizontal degree 0 and parallel with respect to $\nabla^{V'_\gamma}$.

We note that the fibre-wise cohomologies of $(V'^{+}_\gamma, b')$ and $(V'^{-}_\gamma, c')$ are given by

$$H^j\big(V'^{+}_\gamma, b'\big) = \begin{cases} 0 & j = 0 , \\ H^j(X; \mathbb{C}) & j > 0 , \end{cases} \qquad \text{and} \qquad H^j\big(V'^{-}_\gamma, c'\big) = \begin{cases} \mathbb{C} & j = 0 , \\ 0 & j > 0 , \end{cases}$$

because $H^*(V_\gamma'^-,c')$ is the cohomology of a disc. This implies in particular that

$$H^*\bigl(V_\gamma',a_0'\bigr)=H^*(X;\mathbb{C})=H^j\bigl(V_\gamma'^+,b'\bigr)\oplus H^j\bigl(V_\gamma'^-,c'\bigr)\;, \tag{5.18}$$

so the spectral sequence associated with the filtration $V'\gamma\supset V_\gamma'^+\supset 0$ degenerates already at $E_1=E_\infty$.

We now need the analytic analogue of Igusa's splitting principle ([I2], Lemma 5.7.8). Multiplying the (parallel) metric $g^V$ restricted to $V_\gamma'^-$ by a small positive constant $T^{-2}$ has the same effect on $T(A',g^V,h^V)$ as conjugating $a'$ by $\left(\begin{smallmatrix}1&0\\0&T\end{smallmatrix}\right)$. In the limit $T\to\infty$, only terms in $\Omega^0(S^{4k};\operatorname{End}V)$ remain, giving a torsion form in $\Omega^0(S^{4k})$. The class $T(M_\gamma'/S^{4k};\mathbb{C},h)\in H^{\mathrm{even},\geq 2}(S^{4k})$ is unaffected by this deformation for all $T>0$ according to Corollary 2.63. Because the rank of the fibre-wise cohomology does not jump in the limit $T\to\infty$ by (5.18), we find that

$$T(M_\gamma'/S^{4k};\mathbb{C},h)=T\bigl(\nabla^{V_\gamma'^+}+b',g^{V_\gamma'^+}\bigr)^{[\geq 2]}+T\bigl(\nabla^{V_\gamma'^-}+c',g^{V_\gamma'^-}\bigr)^{[\geq 2]}=0\;.$$

Because the bundles $T^sX$ and $T^uX$ are trivial along $C_\gamma'^+$, and the ${}^0J$-classes of $(T^sX-T^uX)|_{C_\gamma'^-}$ are the same as in Theorem 5.13, Theorem 0.1 implies that

$$\mathcal{T}(M_\gamma'/S^{4k};\mathbb{C})=2^{2-k}\,j\,b_k\,\zeta'(-2k)\,u\;. \quad\Box \tag{5.19}$$

**5.d. Exotic structures on general fibre bundles.** We show that we can alter the differentiable structure on a quite general fibre bundle by "horizontally" gluing with the bundle $M_\gamma'$ of the previous subsection. In the following subsection, we use the bundles obtained this way to prove Theorem 0.2.

Let $p=p_0\colon M_0\to B$ be a smooth proper submersion of compact manifolds, and assume that the base is orientable with $\dim B=4k$, and that $\dim X=\dim M_0-\dim B=2l-1$, such that $l$ is sufficiently large with respect to $k$. Let $F\to M_0$ be a unitarily flat vector bundle. If $H=H^*(M_0/B;F)\to B$ admits a parallel metric, the higher torsion $\mathcal{T}(M_0/B;F)\in H_{\mathrm{dR}}^{\mathrm{even},\geq 2}(B)$ is well-defined. Otherwise, we fix an arbitrary metric $g^H$ on $H$ and regard the class

$$\mathcal{T}(M_0/B;F,g^H)=\Bigl(\mathcal{T}\bigl(T^HM_0,g^{TX},\nabla^F,g^F\bigr)+\widetilde{\mathrm{ch}}^{\mathrm{o}}\bigl(\nabla^H,g^H,g^H_{L_2}\bigr)\Bigr)^{[\geq 2]}$$

in $\Omega^{\mathrm{even}}(B)/d\Omega^{\mathrm{odd}}(B)$, which is independent of the choices of $T^HM_0$, $g^{TX}$ and the parallel metric $g^F$, but depends on $g^H$.

Let $X_0=p^{-1}(b_0)\subset M_0$ be a fixed fibre of $p$ for a fixed $b_0\in B$, and let $F_0$ be the restriction of $F$ to $X_0$. Let $M_\gamma'\to S^{4k}$ be constructed as in the previous subsection with fibre $X=X_0$. We may assume that the bundles $M_0\to B$ and $M_\gamma'\to S^{4k}$ are trivialised in small discs $D_0^{4k}$ around $b_0\in B$ and $D_-^{4k}\subset S^{4k}$, and we pick the horizontal subbundles $T^HM_0$ and $T^HM_\gamma'$ according to that trivialisation. Fixing a metric $g^X$ on $X$ and $g_0^H$ on $H(X;F)$, we define $g^{TX}$ and $g_0^H$ over $D_0^{4k}$ and $D_-^{4k}$ using the trivialisation above together with the given $g^X$ and $g_0^H$. Then the data defining the classes $\mathcal{T}(M_0/B;F,g^H)$ and $\mathcal{T}(M_\gamma'/S^{4k};E_i)$ coincide on $\partial D_0^{4k}$ and $\partial D_-^{4k}$, which allows us to glue $M_0|_{B\setminus D_0^{4k}}$ and $M_\gamma'|_{S^{4k}\setminus D_+^{4k}}$ along their common boundary, obtaining a new bundle

$$p_1\colon\, M_1\;\longrightarrow\; B_1\;,$$

where $B_1$ is naturally diffeomorphic to $B$. In other words, we replace the restriction of $M_0$ to $D_0^{4k}$ by the restriction $M_\gamma'$ to $D_+^{4k}$.

Then $p_1\colon M_1 \to B_1 \cong B$ is homeomorphic to $p_0\colon M_0 \to B$ as a bundle by Proposition 5.8. Because the bundle $F|_X$ is diffeomorphic and isometric to $E_j|_{p^{-1}(b_j)}$, we can also identify the flat vector bundles $F$ and $E_j$ along $M|_{\partial D}$, obtaining a new flat vector bundle $F_1 \to M_1$.

We note that all secondary data $T^H M$, $g^{TX}$, $g^F$ and $g^H_0$ has carefully been chosen so as to be compatible with the gluing process. We now restrict our attention to the trivial bundle $F = M_0 \times \mathbb{C}$. Proposition 5.16 implies that

$$(5.20)\quad \mathcal{T}\big(T^H M_1, g^{TX_1}, d, g^{\mathrm{triv}}\big) + \widetilde{\mathrm{ch}}^{\mathrm{o}}\big(\nabla^H, g^H_{0,L_2}, g^H_{1,L_2}\big) = \mathcal{T}\big(T^H M_0, g^{TX_0}, d, g^{\mathrm{triv}}\big) + 2^{2-k}\, b_k\, \zeta'(-2k)\, u \quad \in \quad \Omega^*(B)/d\Omega^*(B)\ ,$$

where $g^H_{0,L_2}$ and $g^H_{1,L_2}$ are the $L_2$-metrics on $H \cong H^*(M_0/B) \cong H^*(M_1/B)$, and where $u$ represents a generator of $H^{4k}(B;\mathbb{Z}) \subset H^{4k}_{\mathrm{dR}}(B)$. Note that both $\widetilde{\mathrm{ch}}^{\mathrm{o}}(\nabla^H, g^H_{0,L_2}, g^H_{1,L_2})$ and $u$ are given by smooth forms supported on $D^{4k}_0 \subset B$, which are well-defined modulo differentials of forms supported on $D^{4k}_0$.

Because the bundles $p_1$ and $p_0$ are naturally homeomorphic, we can identify the fibre-wise cohomology bundles $H(M_1/B;\mathbb{C}) \to B_1 \cong B$ and $H(M_0/B;\mathbb{C}) \to B$ as flat vector bundles, and this identification is canonical on $B \subset D^{4k}_0 = B_1 \subset D^{4k}_+$. In particular, the metric $g^H$ can be regarded as a metric on $H(M_1/B;\mathbb{C})$ as well. Passing to classes, equation (5.20) becomes

$$\mathcal{T}\big(M_1/B;\mathbb{C}, g^H\big) = \mathcal{T}\big(M_0/B;\mathbb{C}, g^H\big) + 2^{2-k}\, b_k\, \zeta'(-2k)\, u \quad \in \quad \Omega^*(B)/d\Omega^*(B)\ .$$

Iterating this process gives bundles $p_j\colon M_j \to B$ for all $j \geq 0$. If we repeat the constructions above with $\gamma$ replaced by $\gamma^{-1}$ in (5.5), we also obtain bundles $p_j\colon M_j \to B$ for all $j < 0$, such that

$$(5.21)\qquad \mathcal{T}\big(M_j/B;\mathbb{C}, g^H\big) = \mathcal{T}\big(M_0/B;\mathbb{C}, g^H\big) + 2^{2-k} j\, b_k\, \zeta'(-2k)\, u \quad \in \quad \Omega^*(B)/d\Omega^*(B)$$

for all $j \in \mathbb{Z}$. Let us summarize our examples so far.

**5.22. Corollary.** *Let $p\colon M \to B$ be a smooth fibre bundle with compact fibre $X$ of dimension $2l-1$ and an orientable, $4k$-dimensional, compact base space $B$, such that $l$ is sufficiently large with respect to $k$, and let $g^H$ be a fixed metric on the fibre-wise cohomology $H(M/B;\mathbb{C}) \to B$. Then there exist infinitely many bundles $p_j\colon M_j \to B$ for $j \in, j > 0$, with fibres diffeomorphic to $X$, that are all homeomorphic to $p = p_0$ as bundles, but with pairwise different values of*

$$\int_B \mathcal{T}\big(M_j/B;\mathbb{C}, g^H\big)\ .$$

*5.23. Remark.* The result above is related to Proposition 6.5.7 in [I2] for trivial bundles over arbitrary compact orientable base spaces. Indeed, if the bundle $H(M/B;\mathbb{C}) \to B$ is trivial, then we could have stated Corollary 5.22 using Igusa's higher Franz-Reidemeister torsion $\tau(M_j/B;\mathbb{C})$. However, Igusa's $\tau(M_j/B;\mathbb{C})$ has so far only been defined if $H(M/B;\mathbb{C}) \to B$ is trivial, so we use $\mathcal{T}(M_j/B;F,g^H)$ to get a more general statement. Since we used Theorem 0.1 to relate the higher analytic torsion of Hatcher's sphere bundles to $\tau(M_\gamma/B,\mathbb{C})$, Corollary 5.22 is in fact an application of the methods developped in this paper.

*5.24. Remark.* At this point, it is natural to ask for those cohomology classes on $B$ that can be realised as differences of higher torsion classes of fibre bundles that are homeomorphic to a given $M \to B$ as bundles. Corollary 5.22 only realises certain multiples of the fundamental class as such differences. However, the gluing techniques used in the proof immediately give a slightly more general result, which applies in particular to the trivial bundle $M \times B \to B$.

Let $p\colon M \to B$ be a smooth fibre bundle with compact fibre $X$ of dimension $2l-1$ and base space $B$. Assume that $B$ contains a properly immersed oriented submanifold $A$ of codimension $4k > 0$ with trivial normal bundle, which has a lift $N \subset M$, such that $l$ is sufficiently large with respect to $k$, and the vertical tangent bundle $TX|_N$ is trivial. Let $\alpha \in H^{4k}(B,\mathbb{R})$ be the Poincaré dual of $A$. Let $g^H$ be a fixed metric on the fibre-wise cohomology $H(M/B;\mathbb{C}) \to B$. Then there exist infinitely many bundles $p_j\colon M_j \to B$ for $j \in$, $j > 0$, with fibres diffeomorphic to $X$, that are all homeomorphic to $p = p_0$ as bundles, but with

$$\mathcal{T}\big(M_1/B;\mathbb{C},g^H\big) = \mathcal{T}\big(M_0/B;\mathbb{C},g^H\big) + 2^{2-k}\, b_k\, \zeta'(-2k)\, \alpha \quad \in \quad \Omega^*(B)/d\Omega^*(B)\ . \tag{5.25}$$

Clearly, equation (5.21) is contained here as a special case where $A$ is a single point.

To proof this claim, we may assume that $N$ is embedded in $M$ by transversality. Then let $U$ denote a tubular neighbourhood of $N$, fix $p \in N$, and let $D_p$ be a fibre of the bundle $U \to N$. We may assume that $U$ fibres over a $\dim B$-dimensional submanifold $V \subset M$ with fibre $D^{2l-1}$, where $V$ is transversal to the fibres of $M \to B$. By our assumptions, the bundle $U \to V$ is trivial, and we have a projection

$$\pi\colon V \longrightarrow D^{4k} \qquad \text{with} \qquad U = \pi^*\big(D^{2l-1} \times D^{4k}\big) \longrightarrow V\ .$$

Thus, we can replace $U$ by the pullback of the bundle $M_\gamma|_{D^{4k}_+} \to D^{4k}$ with a neighbourhood of the critical level $C_1$ removed as in Section 5.c. Then (5.25) follows by the same arguments as (5.19) and (5.21).

**5.e. Homeomorphism classes and diffeomorphism classes of submersions, and a proof of Theorem 0.2.** In the last section, we have exhibited a method to construct an infinite family of pairwise homeomorphic fibre-bundles with pairwise different higher analytic torsion classes, starting from an arbitrary fibre-bundle with compact fibre and compact orientable base of the right dimensions. To prove Theorem 0.2, we have to show that an infinite subfamily has diffeomorphic total spaces and that these bundles are pairwise not diffeomorphic as bundles, which is again a consequence of Corollary 5.22 above. In Theorem 5.29, we give a necessary condition for a diagramme

$$\begin{array}{ccc} M & \xrightarrow[\sim]{\ \Phi\ } & M \\ {\scriptstyle p_i}\big\downarrow & & \big\downarrow{\scriptstyle p_j} \\ B & \xrightarrow[\sim]{\ \Psi\ } & B \end{array}$$

as in Definition 5.7 to commute, where $\Phi$ and $\Psi$ are now diffeomorphisms. This proves Theorem 0.2 and gives the refinement mentioned in the introduction.

Therefore, let $p\colon M \to B$ be a smooth fibre bundle with compact fibre $X$ of dimension $2l-1$ and an orientable, $4k$-dimensional, compact base space $B$, such that $l$ is sufficiently large with respect to $k$. Let $F \to M$ be a unitarily flat vector bundle, and let $p_j\colon M_j \to B$ and $F_j \to M_j$ be constructed as in Corollary 5.22. We do not require that the pair-wise isomorphic flat vector bundles $H_j = H^*(M_j/B;\mathbb{C}) \to B$ carry a parallel metric.

A careful investigation of the construction in Section 5.d shows that $M_j$ arises from $M$ by removing several subsets $D_i$ diffeomorphic to $D^{4k} \times D^{2l-1}$ and gluing them back in a different way, such that before and after the gluing, the fibres of $p$ respectively $p_j$ intersect $D^{4k} \times D^{2l-1}$ either in $\{b\} \times D^{2l-1}$ for some $b \in D^{4k}$ or not at all. After sliding the sets $D_i$ around within $M$, we may assume that all $D_i$ are contained in a larger set $\overline{D} \subset M$ diffeomorphic to $D^{4k} \times D^{2l-1}$, such that each inclusion $D_i \hookrightarrow \overline{D}$ projects to the identity on $D^{4k}$. If we denote the corresponding subset of $M_j$ by $\overline{D}_j$, then $\overline{D}_j$ is still diffeomorphic to $D^{4k} \times D^{2l-1}$, where the fibre of $p_j$ intersect $D^{4k} \times D^{2l-1}$ either in $\{b\} \times D^{2l-1}$ for some $b \in D^{4k}$ or not at all.

We construct $M_j$ for $j < 0$ similarly by glueing in $|j|$ copies of another exotic sphere bundle where we start with the inverse of $\gamma$ in $\pi_{4k-1}(O(m))$. Then (5.21) still holds.

**5.26. Proposition.** *There exists $j_0 \in \mathbb{Z}$ such that there exists a diffeomorphism*

$$\Theta_j \colon M \longrightarrow M_j \qquad \text{with} \qquad \Theta_j|_{M\setminus\overline{D}} = \mathrm{id}_{M\setminus\overline{D}} \colon M \setminus \overline{D} \longrightarrow M_j \setminus \overline{D}_j$$

*for all $j \in j_0\mathbb{Z}$. In particular, $\Theta_j$ lifts to a parallel map*

$$\overline{\Theta}_j \colon F \longrightarrow F_j \qquad \text{with} \qquad \overline{\Theta}_j|_{F|_{M\setminus\overline{D}}} = \mathrm{id}_{F|_{M\setminus\overline{D}}} \colon F|_{M\setminus\overline{D}} \longrightarrow F_j|_{M_j\setminus\overline{D}_j} \,.$$

*Proof.* By the constructions in Subsections 5.c and 5.d, there exists a possibly exotic $(4k+2l-1)$-dimensional sphere $\Sigma$ such that

$$M_j = M \# \underbrace{\Sigma \# \cdots \# \Sigma}_{j \text{ times}}$$

for all $j > 0$, where the connected sum operation takes places inside $\overline{D}$. For $j < 0$, we have to take the connected sum with $|j|$ copies of another possibly exotic sphere $\Sigma'$. Because the group of exotic spheres is finite in each dimension, the first claim now follows.

The second claim is then clear because $F$ and $F_j$ are up to isomorphism the unique extensions of $F|_{M\setminus\overline{D}}$ to $M$ and $M_j$, respectively, thus $\Theta_j^* F_j$ and $F$ are isomorphic as flat bundles and can thus be identified. □

It is an immediate consequence of Corollary 5.22 that the maps $p$ and $p_j \circ \Theta_j$ differ somewhere on $\overline{D}$ for $j \neq 0$. On the other hand, recall that by Proposition 5.8, there always exists a homeomorphism $\Theta_j' \colon M \to M_j$ supported on $\overline{D}$ such that $p = p_j \circ \Theta_j'$. In the remainder of this chapter, we want to generalise and strengthen this observation. But first of all, we need some facts on Chern-Simons classes.

Let $(F, \nabla^F) \to M$ be a flat vector bundle, and let $\mathcal{D}\mathit{iff}(M;F)$ denote the group of fibre preserving, fibre-wise linear, parallel diffeomorphisms $\overline{\Phi} \colon F \to F$. An element $\Phi \in \mathcal{D}\mathit{iff}(M;F)$ is a parallel lift of a diffeomorphism of $M$, which by an abuse of notation will again be named $\Phi$. If $M$ is oriented, let $\mathcal{D}\mathit{iff}^+(M;F) \subset \mathcal{D}\mathit{iff}(M;F)$ denote the subgroup of lifts of orientation preserving diffeomorphisms of $M$.

Let $\overline{F} = F \times [0,1]$ denote the trivial extension of $F$ to $\overline{M} = M \times [0,1]$. For $\Phi \in \mathcal{D}\mathit{iff}(M;F)$, the bundle $F$ naturally extends to a flat vector bundle

$$T_\Phi F = \overline{F}/\sim \quad \longrightarrow \quad T_\Phi M = \overline{M}/\sim$$

over the mapping torus $T_\Phi M$, where "$\sim$" is generated by $(v,1) \sim (\Phi(v),0)$ for all $v \in F$, respectively $(p,1) \sim (\Phi(p),0)$ for all $p \in M$. If $g^F$ is any metric on $F$, let $g_t^F$ be a family of metrics that interpolates between $g_0^F = g^F$ and $g_1^F = \Phi^* g^F$. Let $g^{\overline{F}}$ denote the corresponding metric on $\overline{F}$, then $g^{\overline{F}}$ descends to a metric $g^{T_\Phi F}$ on $T_\Phi F$.

Now we assume in addition that $M$ is oriented and that $\Phi$ preserves the orientation. By the definition of the Chern-Simons class $\widetilde{\mathrm{ch}}^\mathrm{o}$, we have

$$\int_M \widetilde{\mathrm{ch}}^\mathrm{o}\big(\nabla^F, g^F, \Phi^* g^F\big) = \int_{M\times[0,1]} \mathrm{ch}^\mathrm{o}\Big(\nabla^{\overline{F}}, g^{\overline{F}}\Big) = \mathrm{ch}^\mathrm{o}\big(\nabla^{T_\Phi F}\big)[T_\Phi M] \,, \tag{5.27}$$

independent of the choice of $g^F$.

**5.28. Proposition.** *Assume that $M$ is oriented, then the map*

$$\mathcal{D}\mathit{iff}^+(M;F) \longrightarrow \mathbb{R} \qquad \text{with} \qquad \Phi \longmapsto \mathrm{ch}^o\big(\nabla^{T_\Phi F}\big)[T_\Phi M]$$

*is a group homomorphism.*

*Assume in addition that* $\dim M > 0$. *Then* $\mathrm{ch}^o(\nabla^{T_\Phi F})[T_\Phi M]$ *depends only on the homotopy type of the action of* $\Phi$ *on* $M$. *If* $F$ *carries a parallel metric or if* $\Phi^k \sim \mathrm{id}_M$ *for some* $k \neq 0$, *then* $\mathrm{ch}^o(\nabla^{T_\Phi F})[T_\Phi M] = 0$.

To see that the invariant $\mathrm{ch}^o(\nabla^{T_\Phi F})[T_\Phi M]$ is nontrivial, one may consider a surface $M$ of genus $\geq 2$ and a diffeomorphism $\Phi$ such that $T_\Phi M = \mathbb{H}^3/\Gamma$ becomes hyperbolic with fundamental group $\Gamma \subset SL(2,\mathbb{C})$. We then choose $F \to M$ and a lift $\Phi$ to $F$ such that $T_\Phi F$ becomes the flat bundle associated to the adjoint representation of $SL(2,\mathbb{C})$ restricted to $\Gamma$. In this case, $\mathrm{ch}^o(\nabla^{T_\Phi F})[T_\Phi M]$ is a nonzero multiple of the hyperbolic volume of $T_\Phi M$.

*Proof.* The map $\widetilde{\mathrm{ch}}^o(M;F)$ is a homomorphism by (5.27), because for $\Phi$, $\Phi' \in \mathcal{D}\mathit{iff}^+(M;F)$, we have

$$\widetilde{\mathrm{ch}}^o\big(\nabla^F, g^F, (\Phi\circ\Phi')^* g^F\big) = \widetilde{\mathrm{ch}}^o\big(\nabla^F, g^F, \Phi^* g^F\big) + \Phi^*\widetilde{\mathrm{ch}}^o\big(\nabla^F, g^F, (\Phi')^* g^F\big)$$

modulo exact forms.

We now assume that $\dim M > 0$. Let the action of $\Phi$ on $M$ be homotopic to the identity, then $T_\Phi M$ is diffeomorphic to $M \times S^1$. Because flat vector bundles on $M \times S^1$ are given by representations of $\pi_1(M) \times \mathbb{Z}$, the bundle $T_\Phi F$ decomposes into parallel subbundles of the form

$$\pi_1^* F_M \otimes \pi_2^* F_{S^1} \ ,$$

where $F_M \to M$ and $F_{S^1} \to S^1$ are flat vector bundles. By Proposition 1.13 in [BL], we have

$$\mathrm{ch}^o\Big(\nabla^{\pi_1^* F_M \otimes \pi_2^* F_{S^1}}\Big) = \pi_1^*\,\mathrm{ch}^o\big(\nabla^{F_M}\big)\ \mathrm{rk}\,F_{S^1} + \mathrm{rk}\,F_M\ \pi_2^*\,\mathrm{ch}^o\big(\nabla^{F_{S^1}}\big) \ ,$$

thus in particular

$$\mathrm{ch}^o\Big(\nabla^{\pi_1^* F_M \otimes \pi_2^* F_{S^1}}\Big)\big[M \times S^1\big] = 0 \ .$$

If $\Phi^k$ is homotopic to the identity, the homomorphism property implies that

$$k\ \mathrm{ch}^o\big(\nabla^{T_\Phi F}\big)[T_\Phi M] = \mathrm{ch}^o\big(\nabla^{T_{\Phi^k} F}\big)[T_{\Phi^k} M] = 0 \ .$$

Assume that $g^F$ is a parallel metric, then $\Phi^* g^F$ is also parallel, and we may interpolate by a family $g_t^F$ of metrics that are parallel on $M \times \{t\}$ for each $t$. Then $\mathrm{ch}^o(\nabla^{T_\Phi F}, g^{T_\Phi F}) \in \Omega^1(T_\Phi M)$, in particular

$$\mathrm{ch}^o\big(\nabla^{T_\Phi F}\big)[T_\Phi M] = 0 \ . \qquad \square$$

We now want state the main application of Subsections 5.c and 5.d, which allows us to give many examples of manifolds $M$, $B$ that allow infinitely many smooth maps $M \to B$ that are pairwise conjugate by homeomorphisms, but not by diffeomorphisms. Whenever the dimensions match and $B$ is compact and orientable, we can construct infinitely many proper submersions $M \to B$ that are homeomorphic but not diffeomorphic as bundles. Let us write $q_j = p_{jj_0} \circ \Theta_{jj_0}$, where $j_0 > 0$ and $\Theta_{jj_0}$ are given by Proposition 5.26.

**5.29. Theorem.** *Let $p\colon M \to B$ be a proper submersion. Assume that $B$ is compact, orientable, and of dimension $4k > 0$. Let the bundles $(q_j\colon M \to B)_{j\in\mathbb{Z}}$ with $q_0 = p$ be constructed as above. Assume that for $i$, $j \in \mathbb{Z}$ with $i \neq j$, there exists diffeomorphisms $\Phi \in \mathcal{D}\mathit{iff}(M)$ and $\Psi \in \mathcal{D}\mathit{iff}^+(B)$ such that the diagramme*

$$\begin{array}{ccc} M & \xrightarrow{\ \Phi\ } & M \\ {\scriptstyle q_i}\big\downarrow & & \big\downarrow{\scriptstyle q_j} \\ B & \xrightarrow{\ \Psi\ } & B \end{array}$$

*commutes. Then $\Psi \in \mathcal{D}\mathit{iff}^+(B;H)$ for $H = H^*(M/B;\mathbb{C}) \to B$, and*

$$\mathrm{ch}^{\mathrm{o}}\big(\nabla^{T_\Psi H}\big)[T_\Psi B] = 2^{2-k}\,(j-i)\,j_0\,b_k\,\zeta'(-2k)\;.$$

*In particular, this is only possible if $H \to B$ does not carry a parallel metric, and if $\Psi$ does not represent a torsion element in the mapping class group.*

*Proof* of Theorem 0.2. By Proposition 5.8, the maps $q_j$ are pairwise homeomorphic as bundles, in the sense of Definition 5.7. In particular, they are also homeomorphic as maps from $M$ to $B$. Proposition 5.28 and Theorem 5.29 have the following two easy consequences.

(1) If the bundle $H = H^*(M/B;\mathbb{C})$ carries a parallel metric, then the maps $q_j$ are not pairwise conjugate by diffeomorphisms, in particular, they are not pairwise isomorphic as smooth fibre bundles.
(2) If the bundle $H$ does not admit a parallel metric, then the maps $q_j$ are not pairwise isomorphic as smooth fibre bundles, because $\mathrm{id}_B$ clearly is a torsion element in the mapping class group.

This proves Theorem 0.2. □

*Proof* of Theorem 5.29. We know from Proposition 5.8 and the construction in Subsections 5.c and 5.d that the fibre bundles $q_i$, $q_j\colon M \to B$ are homeomorphic to $M \to B$ by a homeomorphism that is supported on a contractible subset of $M$. To avoid confusion, we will write $q_i$ and $q_j$ instead of $M/B$ whenever necessary. In particular, the fibrewise cohomologies $H^*(q_i;F)$ and $H^*(q_j;F) \to B$ can be identified with the original $H \to B$. By functoriality of the fibre-wise cohomology, the diffeomorphism $\Psi$ lifts to $\mathcal{D}\mathit{iff}^+(B;H)$.

If $g_0^H$ is any metric on $H$, the class

$$\mathcal{T}\big(T^HM, g^{TX}, d, g^{\mathrm{triv}}\big) - \widetilde{\mathrm{ch}}^{\mathrm{o}}\big(\nabla^H, g^H_{L_2}, g^H_0\big) \quad \in \quad \Omega^{\mathrm{even}}(B)/d\Omega^{\mathrm{odd}}(B)$$

is independent of the choices of $T^HM$ and $g^{TX}$ by Corollary 2.17 and Theorem 2.48. Since we assume that $\dim B > 0$, we obtain an invariant

$$\mathcal{T}\big(M/B, \mathbb{C}; g_0^H\big) = \int_B \Big(\mathcal{T}\big(T^HM, g^{TX}, d, g^{\mathrm{triv}}\big) - \widetilde{\mathrm{ch}}^{\mathrm{o}}\big(\nabla^H, g^H_{L_2}, g^H_0\big)\Big)\;,$$

which depends on the choice of $g_0^H$ by

$$\mathcal{T}\big(M/B, \mathbb{C}; g_1^H\big) - \mathcal{T}\big(M/B, \mathbb{C}; g_0^H\big) = -\int_B \widetilde{\mathrm{ch}}^{\mathrm{o}}\big(\nabla^H, g^H_0, g^H_1\big)\;. \tag{5.30}$$

By naturality of torsion, we find that

$$\mathcal{T}\big(q_j, \mathbb{C}; g_0^H\big) = \mathcal{T}\big(q_i, \mathbb{C}; \Psi^* g_0^H\big)\;. \tag{5.31}$$

On the other hand, it follows from (5.21) that

$$\mathcal{T}\big(q_j, \mathbb{C}, g_0^H\big) - \mathcal{T}\big(q_i, \mathbb{C}, g_0^H\big) = 2^{2-k}(j-i)\,j_0\,b_k\,\zeta'(-2k)\,u\;, \tag{5.32}$$

where $u$ is the positive generator of $H^{4k}(B;\mathbb{Z})$. The formula in the theorem follows from (5.27) and (5.30)–(5.32). The remaining conclusions are a consequence of Proposition 5.28. □

**5.f. Hatcher's example and Lott's secondary $K$-theory.** We reformulate Corollary 5.22 in the language of Lott's secondary $K$-theory $\overline{K}^0_R$ of [L], cf. Section 3.c. We will see that the different maps $p_j\colon M \to B$ constructed in Subsections 5.d and 5.e give rise to different pushdown maps $p_{j!}\colon \overline{K}^0_R(M) \to \widehat{K}^0_R(B)$, at least if $\pi_1(B)$ is finite. We also discuss this fact in the light of Conjecture 1 in [L], which relates the pushdown $p_!$ to a Becker-Gottlieb transfer.

Thus, let $p\colon M \to B$ be a smooth fibre bundle with compact fibres, and assume that $B$ is compact oriented with $\dim B = 4k$ and the dimension of the fibre is odd and sufficiently large. Let $(q_i\colon M \to B)_{i\in\mathbb{Z}}$ with $q_0 = p$ denote the maps considered in Theorem 5.29. We take $R$ to be one of $\mathbb{Z}$, $\mathbb{R}$, $\mathbb{C}$ for simplicity and let $\rho\colon R \to \mathbb{C}$ be the inclusion. Then by (3.25) and Theorem 3.26,

$$q_{i!}\bigl[M \times R, g^{\mathrm{triv}}, 0\bigr] = q_{0!}\bigl[M \times R, g^{\mathrm{triv}}, 0\bigr] + [0, 0, i\,u] \tag{5.33}$$

for some nonzero $u \in H^{4k}(S^{4k}, \mathbb{R})$. By Proposition 5.35 below, the elements of $\overline{K}^0_R(B)$ given by (5.33) are pairwise different, at least if $\pi_1(B)$ is finite.

In the following, we let $B$ be a smooth manifold with finite fundamental group. Then we claim that there exists a well-defined map

$$\widehat{\mathrm{ch}}^{\mathrm{o}}\colon\ \widehat{K}^0 \longrightarrow \Omega^{\mathrm{even},\geq 2}(B)/d\Omega^{\mathrm{odd}}(B)\,.$$

To define $\widehat{\mathrm{ch}}^{\mathrm{o}}$ on a generator $(F, g^F, \eta)$, we choose a parallel metric $g^F_0$ on $F_{\mathbb{C}}$, which is possible because $B$ is simply connected, and we set

$$\widehat{\mathrm{ch}}^{\mathrm{o}}(F, g^F, \eta) = \widetilde{\mathrm{ch}}^{\mathrm{o}}\bigl(\nabla^F, g^F_0, g^F\bigr)^{[\geq 2]} - \eta^{[\geq 2]}\,. \tag{5.34}$$

This is independent of $g^F_0$ because for two parallel metrics $g^F_0$, $g^F_1$ on $F$, we have $\widetilde{\mathrm{ch}}^{\mathrm{o}}(\nabla^F, g^F_0, g^F_1) \in \mathbb{R} \subset \Omega^0(B)$ modulo exact forms.

It is easy to check that (5.34) complies with (3.27) due to the properties of $\widetilde{\mathrm{ch}}^{\mathrm{o}}$. We check that (5.34) is compatible with the $\widehat{K}^0_R$-relations (3.22). Regard a short exact sequence as in (3.20). By (3.27), we may choose parallel metrics $g^{F_i}$ on the bundles $F_{i,\mathbb{C}}$. Then the associated Bismut-Lott torsion of (3.21) vanishes in degree $\geq 2$, so (5.34) is compatible with (3.22) and the map $\widehat{\mathrm{ch}}^{\mathrm{o}}$ is well-defined.

**5.35. Proposition.** *Assume that $\pi_1(B)$ is finite, then*

$$\overline{K}^0_R(B) \cong \overline{K}^0_R(\mathrm{pt}) \oplus H^{\mathrm{even},\geq 2}(B, \mathbb{R})\,,$$

*and the combined map*

$$H^{\mathrm{even},\geq 2}(B, \mathbb{R}) \longrightarrow \overline{K}^0_R(B) \xrightarrow{\widehat{\mathrm{ch}}^{\mathrm{o}}} H^{\mathrm{even},\geq 2}(B, \mathbb{R})$$

*is multiplication by $-1$.*

*Proof.* Because $B$ is simply connected, the sequence (3.23) becomes

$$H^{\mathrm{even}}(B, \mathbb{R}) \longrightarrow \overline{K}^0_R(B) \longrightarrow K^0_R(\mathrm{pt}) \longrightarrow 0\,.$$

Since any element of $\overline{K}^0_R(B)$ can be represented as a sum of $\widehat{K}^0_R$-generators $(F, g^F, \eta)$ with parallel $g^F$ and $\eta \in H^{\mathrm{even}}(B, \mathbb{R})$, the proposition follows from the fact that $\widehat{\mathrm{ch}}^{\mathrm{o}}$ is well-defined. □

*5.36. Remark.* Lott conjectured in [L] that $p_!$ is related to a Becker-Gottlieb transfer. More precisely, let $\mathcal{F}^*_R$ denote the homotopy fibre of the map

$$K^*(R) \xrightarrow{\rho_*} K^*(\mathbb{C}) \longrightarrow \prod_j K(\mathbb{R}, 2j+1) ,$$

where $K^*$ denotes the spectrum of algebraic $K$-theory, and $K(\mathbb{R}, 2j-1)$ denotes the Eilenberg-McLane spectrum of odd real homology. Here, the left arrow is induced by the representation $\rho: R \to \operatorname{End}(C^n)$, and the right arrow is the Borel regulator map. Then $\mathcal{F}^*_R$ is again a generalised cohomology theory, and the Becker-Gottlieb transfer gives a homomorphism

$$\tau^*_{\mathrm{BG}}\colon \mathcal{F}^0_R(M) \longrightarrow \mathcal{F}^0_R(B) . \tag{5.37}$$

Lott conjectures the existence of a natural map $\overline{K}^0_R \to \mathcal{F}^0_R$ that fits into commutative diagrammes

$$\begin{array}{ccccc} H^{\mathrm{even}}(M;\mathbb{R}) & \longrightarrow & \overline{K}^0_R(M) & \longrightarrow & K^0_R(M) \\ \downarrow & & \downarrow & & \downarrow \\ H^{\mathrm{even}}(M;\mathbb{R}) & \longrightarrow & \mathcal{F}^0_R(M) & \longrightarrow & K^0_{\mathrm{alg},R}(M) , \end{array} \tag{5.38}$$

cf. (3.23), and

$$\begin{array}{ccc} \overline{K}^0_R(M) & \longrightarrow & \mathcal{F}^0_R(M) \\ {\scriptstyle p_!}\downarrow & & \downarrow{\scriptstyle \tau^*_{\mathrm{BG}}} \\ \overline{K}^0_R(B) & \longrightarrow & \mathcal{F}^0_R(B) . \end{array} \tag{5.39}$$

The Becker-Gottlieb transfer itself is a stable map $\tau_{\mathrm{BG}}: B^+ \to M^+$, which can be defined without referring to the smooth structure of $q: M \to B$, see [BeG]. This implies that $\tau^*_{\mathrm{BG}}$ in (5.37) is the same map for all $q_i$, $i \in \mathbb{Z}$, whereas the maps $q_{i!}: \overline{K}^0_R(M) \to \overline{K}^0_R(B)$ are pairwise different by (5.33) and Proposition 5.35 if $\pi_1(B)$ is finite. We conclude that the $\overline{K}^0_R$-pushdown contains more information than the Becker-Gottlieb transfer $\tau^*_{\mathrm{BG}}$ in $\mathcal{F}^0_R$. Moreover, if a natural map $\overline{K}^0_R \to \mathcal{F}^0_R$ satisfying (5.38) and (5.39) above exists, then the element $[0, 0, u] \in \overline{K}^0_R(B)$ of (5.33) lies the kernel of $\overline{K}^0_R(B) \to \mathcal{F}^0_R(B)$ for all closed oriented $B$ of dimension $4k$.

Mathematisches Institut der Uni Tübingen, Auf der Morgenstelle 10, 72076 Tübingen, Germany
*E-mail address*: `sebastian.goette@uni-tuebingen.de`